\documentclass[12pt]{amsart}
\usepackage[utf8]{inputenc}
\usepackage{graphicx}
\usepackage{subcaption}
\graphicspath{ {./figures_bis_v2/} }
\usepackage{tikz}
\usepackage{float}
\allowdisplaybreaks
\usepackage{dsfont}
\usepackage{empheq}
\usepackage{enumerate}
\usepackage{setspace}
\usepackage[T1]{fontenc}
\usepackage[letterpaper,top=2cm,bottom=2cm,left=1.5cm,right=1.5cm,marginparwidth=1.75cm]{geometry}

\usepackage[colorlinks=true, allcolors=blue]{hyperref}

\def\a {{\boldsymbol a}}
\def\x {{\boldsymbol x}}

\newtheorem{theorem}{Theorem}
\newtheorem{lemma}{Lemma}
\newtheorem{remark}{Remark}
\newtheorem{proposition}{Proposition}

\title[Physics-based FE approximations of the PNP equations]
{Physics-based stabilized finite element approximations of the Poisson--Nernst--Planck equations}
\author[J. Bonilla]{Jesús Bonilla$^\ddag$}
\address{$\ddag$ Los Alamos National Laboratory, Los Alamos, NM 87545, USA. E-mail: {\tt \href{mailto:jbonilla@lanl.gov}{jbonilla@lanl.gov}}}
\thanks{$\ddag$ Los Alamos National Laboratory, an affirmative action/equal opportunity employer, is operated by Triad National Security, LLC for the National Nuclear Security Administration of U.S. Department of Energy under contract 89233218CNA000001. Los Alamos National Laboratory strongly supports academic freedom and a researcher's right to publish; as an institution, however, the Laboratory does not endorse the viewpoint of a publication or guarantee its technical correctness. LA-UR-24-31612}

\author[J. V. Gutiérrez-Santacreu]{Juan Vicente Gutiérrez-Santacreu$^\S$}
\address{$\S$Dpto. de Matemática Aplicada I\\
         E. T. S. I. Informática\\
         Universidad de Sevilla\\
         Avda. Reina Mercedes, s/n.\\
         E-41012 Sevilla\\
         Spain\\
         E-mail: {\tt \href{mailto:juanvi@us.es}{juanvi@us.es}}}

\thanks{$^\S$ JVGS was partially supported by the Spanish Grant No. PGC2018-098308-B-I00 from Ministerio de Ciencias e Innovación - Agencia Estatal de Investigación with the participation of FEDER and by the Andalusian Grant No. P20\_01120 from Junta de Andalucía (Consejería de Economía, Conocimiento, Empresas y Universidad)}

\begin{document}

\begin{abstract}

We present and analyze two stabilized finite element methods for solving numerically the Poisson--Nernst--Planck equations. The stabilization we consider is carried out by using a shock detector and a discrete graph Laplacian operator for the ion equations, whereas the discrete equation for the electric potential need not be stabilized.  Discrete solutions stemmed from the first algorithm preserve both maximum and minimum discrete principles. For the second algorithm, its discrete solutions are conceived so that they hold discrete principles and obey an entropy law provided that an acuteness condition is imposed for meshes. Remarkably the latter is found to be unconditionally stable.  We validate our methodology through numerical experiments. 
\end{abstract}

\maketitle

{\bf 2020 Mathematics Subject Classification.}  65M60, 35K61, 65Z05

{\bf Keywords.} Poisson--Nernst--Planck equations; stabilized finite-element approximation; shock detector; maximum and minimum discrete principles;  entropy.

\tableofcontents
\section{Introduction}

The Poisson--Nernst--Planck (PNP) system is a mathematical framework for understanding electrodiffusion, which makes reference to  the dynamics of charged particles under diffusion and electric potential. It combines the Possion equation, which describes the electrostatic potential resulting from charge distributions, with the Nernst--Planck equations, which account for the flux of ions due to both concentration gradients and electric fields. Many applications in physics, engineering, and biology, can be described with the PNP system, some of them are semiconductors, electrochemical devices, and biological membranes.

In this paper we consider the numerical approximation of  the PNP problem for two ion species -- cation and anion -- on a bounded domain $\Omega\subset\mathds{R}^d$, with $d=2$ or $3$, over the time interval $(0, \infty)$. This system reads as 
\begin{equation}\label{PNP}
    \left\{
        \begin{array}{rclcc}
             \partial_t p-\Delta p-\nabla\cdot(p\nabla \phi)&=&0& \mbox{ in } &  \Omega\times(0,\infty),
         \\
             \partial_t n-\Delta n+\nabla\cdot(n\nabla \phi)&=&0& \mbox{ in } & \Omega\times(0,\infty), 
        \\  
            -\Delta \phi&=&p-n&\mbox{ in } & \Omega\times(0,\infty),
    \end{array}
\right.
\end{equation}
complemented with homogeneous Neumann boundary conditions
\begin{equation}\label{BC}
\partial_{\boldsymbol{n}} p = \partial_{\boldsymbol{n}} n = 0\quad\mbox{ on }\quad \partial\Omega\times (0,\infty),
\end{equation}
and initial conditions
\begin{equation}\label{IC}
p(0)=p_0\quad\mbox{ and }\quad n(0)=n_0\quad \mbox{ in }\quad\Omega.
\end{equation}
In the above $p, n: \bar \Omega\times [0,\infty)\to\mathds{R}_+$  are densities corresponding to positive and negative ion charges, respectively, and $\phi: \bar\Omega\times [0,\infty)\to \mathds{R}$ is the electric potential. 

In equations $\eqref{PNP}_{1-2}$, the terms $-\Delta p$ and $- \Delta n$ stand for Fick's law of diffusion, which governs ion fluxes due to diffusion down the concentration gradients. The terms $\nabla\cdot(p\nabla \phi)$ and $- \nabla\cdot(n\nabla \phi)$ represent Ohm's law of drift-diffusion, which causes positive and negative ions to move away from each other, down and up, the  gradient of the electric field, respectively.  Equation $\eqref{PNP}_3$ models the electrical activity in the system and is a consequence of Maxwell's equation, in the absence of magnetic fields or with a quasi-static field, and of Gauss'law. 

For the sake of discussion,  we have normalized the coefficients in \eqref{PNP} concerning the diffusion terms $\nabla\cdot(p\nabla \phi)$ and $- \nabla\cdot(n\nabla \phi)$ and drift-diffusion terms $\nabla\cdot(p\nabla \phi)$ and $- \nabla\cdot(n\nabla \phi)$.   

System \eqref{PNP} is very elusive when it comes to constructing numerical solutions. Standard Galerkin methods based on continuous, piecewise linear finite element spaces can often lead to unsatisfactory approximations. Among the most common reasons for inaccuracy of approximation is failure to satisfy discrete principles, energy and/or entropy laws that continuous solutions enjoy. More particularly, discrete principles are of fundamental importance when dealing with the presence of sharp gradients in the electric potential and/or in the ion densities near charged interfaces or boundaries, as these factors further compromise numerical stability and accuracy. Therefore algorithms must adequately construct numerical solutions that capture these gradients in order to avoid local spurious oscillations and maintain computational reliability. The coupled nature of the ion equations through the electric potential is strongly tied to obtaining discrete principles and the energy and entropy decays when the system is isolated. Moreover, the correct dissipation of energy and entropy is essential to reach steady-state equilibria consistent with physics.   

The numerical resolution of system \eqref{PNP} has been tackle dusing different techniques via finite elements for developing physics-based methods. In this framework, reformulating the continuous or discrete equations seems to be an effective approach. For instance, Prohl and Schmuck \cite{Prohl_Schmuck_2009} used an entropy provider to modify the ion fluxes at the discrete level. Discrete principles were only accomplished for initial conditions limited to $0\le p_0, n_0\le 1$. Huang and Shen \cite{Huang_Shen_2021} developed a SAV-based algorithm, which enforces positivity with the help of a suitable functional transformation, since applying the SAV methodology only supplies energy-  and entropy-decreasing algorithms. In regard to high-order finite elements, Xu and Fu \cite{Fu_Xu_2022} proposed a log-density formulation for reformulating system \eqref{PNP}. In this case, positivity is a byproduct of the reformulation. In the same spirit, Shen and Xu \cite{Shen_Xu_2021} made use of the entropy variable to rewrite the ion fluxes. For interested readers, other numerical algorithms based on the finite difference method exist in the literature 
\cite{Ding_Wang_Zhou_2023, Ding_Wang_Zhou_2020,He_Pan_Yen_2019, Liu_Wang_Wise_Yue_Zhou_2021, Liu_Wang_Wise_Ye_Zhou_2023} aligned with addressing the physical constraints of the PNP equations.  

In recent years there has been a significant progress in algorithmic tools for generating discrete-principle-bounded solutions. In particular, artificial diffusion operators \cite{Badia_Bonilla_2017} have considerably gained some popularity since they only manipulate the algebraic structure resulting from using finite elements. These artificial diffusion operators are designed to mitigate spurious oscillations at local extrema, which can significantly deteriorate the numerical resolution. Nevertheless, the resulting algebraic manipulation methods are inherently nonlinear and not simply adapted to construct energy-and-entropy-based numerical solutions. Therefore, some further elaborations are required. Related examples can be found in the context of the Keller--Segel equations \cite{AS_GG_RG_2023,Badia_Bonilla_GS_2023, Bonilla_GS}, which share physical constraints similar to the PNP equations. 

The physical constraints that will be handled are listed below for the sake of clarity and rightful understanding.

\begin{enumerate}[(a)]
\item Minimum and maximum principle. For all  $(\x,t)\in \Omega\times (0,\infty)$,
$$
\min\{\underset {\x\in\Omega}{\min}\,p^0_h(\x),\underset {\x\in\Omega}{\min }\,n^0_h(\x)\}\le p(\x,t), n(\x,t) \le \max\{\underset {\x\in\Omega}{\max}\,p^0_h(\x),\underset {\x\in\Omega}{\max }\,n^0_h(\x)\}.
$$
\item Mass invariance. For all $t\in(0,\infty)$,
$$
\int_\Omega p(\x,t) d\x=\int_\Omega p_0(\x) d\x
$$
and
$$
\int_\Omega n(\x,t) d\x=\int_\Omega n_0(\x) d\x
$$
\item Energy law. For all $t\in(0,\infty)$,
\begin{equation}\label{Energy_law}
\frac{1}{2}\frac{d}{d{\rm t}} \int_\Omega |\nabla \phi(t)|^2\, d\x+\int_\Omega |(p-n)(t)|^2\,d\x+\int_\Omega (p+n)(t)|\nabla\phi(t)|^2\,d\x=0.
\end{equation}
This energy applies only when the diffusion coefficients are identical. 
\item Entropy law. For all $t\in(0, \infty)$,  
\begin{equation}\label{Entropy_law}
\frac{d}{{\rm d}t}\mathcal{E}(p(t),n(t))=-\mathcal{D}(p(t),n(t)).
\end{equation}
where
$$
\mathcal{E}(p(t),n(t))=\int_{\Omega} \Big(p(t) (\log p(t) -1)+1+n(t) (\log n(t) -1)+1+|\nabla \phi(t)|^2 \Big)\,d\x
$$
and
$$
\mathcal{D}(p(t),n(t))=\int_\Omega \Big(p(t)|\nabla(\log p(t)+\phi(t))|^2+n(t)|\nabla(\log n(t)+\phi(t)|^2\Big) d\x.
$$
\end{enumerate}    

Our objective in this paper is to construct and analyze two algorithms, for approximating solutions of the Poisson--Nernst--Planck equations, based on algebraic manipulation methods via artificial diffusion combined with an Euler time-marching integration. In designing our first algorithm, we seek that its numerical solutions satisfy discrete principles, particularly satisfying both maximum and minimum.  Meanwhile the second algorithm aims to provide numerical solutions that exhibit an entropy-decreasing behavior in accordance with \eqref{Entropy_law} as well.  

After this introductory section, the paper is organized as follows. Section 2 is dedicated to setting out the functional and discrete spaces needed for defining our two algorithms, where  the concept of the shock detector and the corresponding stabilzing terms are introduced. The statements of our two results specifying the physical constraints that each algorithm must satisfy is given in Section 3. Following are Sections 4 and 5 that contain the proof of these physical constraints for the first and second algorithms, respectively. To finish up, Section 6 presents a series of numerical experiments focusing on smooth initial data and ion channel transport.        

\section{Numerical approximation}
\subsection{Notation}
Here and throughout, we adopt the following notation.  If $\Omega$ is a Lebesgue-measurable subset of $\mathds{R}^d$, $L^p(\Omega)$, with $1\le p \le \infty$, consists of all Lebesgue functions $f$ such that $\|f\|_{L^p(\Omega)}$ is finite, where
$$
\|f\|_{L^p(\Omega)}=\left(\int_\Omega |f(\x)|^p {\rm d}\x\right)^{1/p}
$$
or 
$$
\displaystyle
\|f\|_{L^\infty}(\Omega) = \underset{\x\in\Omega}{\rm ess\, sup}\, |f(\x)|.  
$$ 
Further we introduce $L^2_{\int=0}(\Omega)$ to be the space of square-summable functions with zero mean in $\Omega$.

\subsection{Setting}

\begin{figure}
	\centering
	\includegraphics[width=0.25\textwidth]{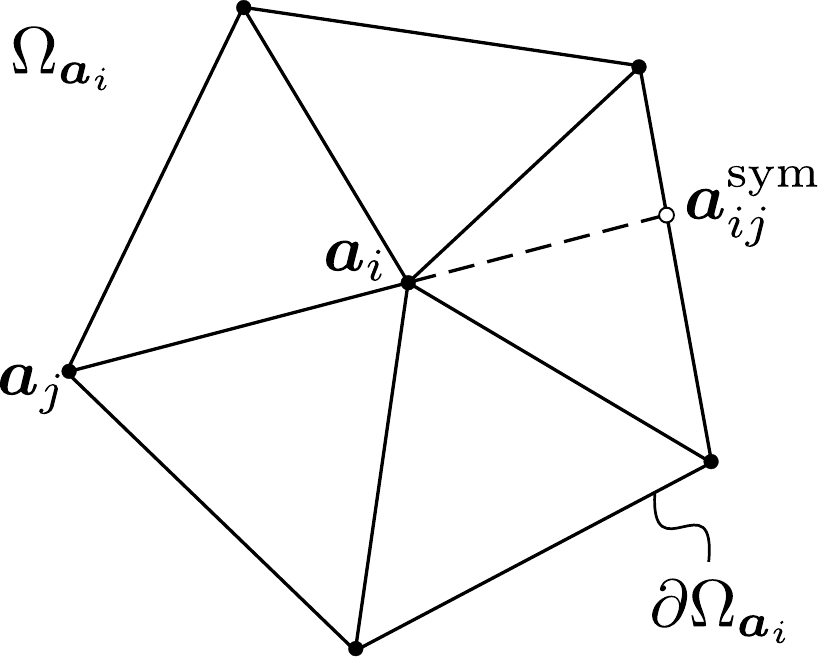}
	\hspace{5em}
	\includegraphics[width=0.21\textwidth]{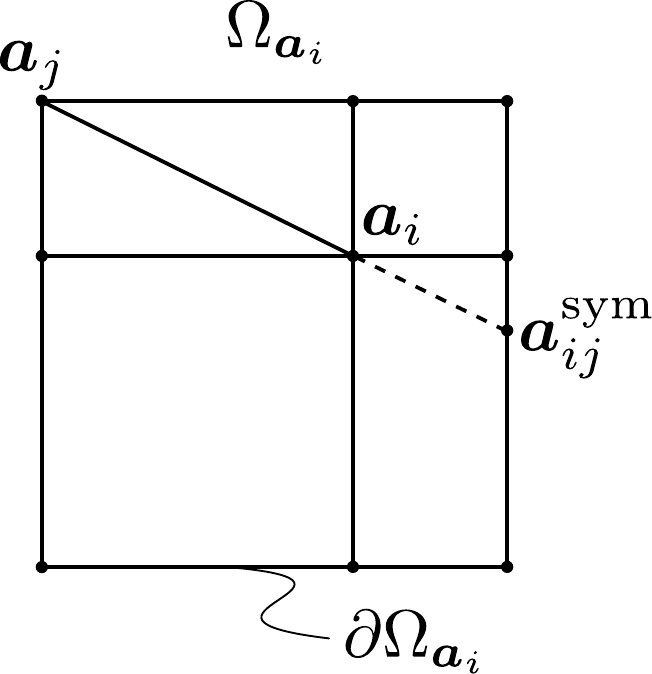}
	\caption{Representation of the symmetric node $\a_{ij}^{\rm sym}$ of $\a_j$ concerning $\a_i$ in a macroelement of triangles (left) and right quadrilaterals (right).}
	\label{fig:usym}
\end{figure}

Henceforth let $\Omega$ be a polygon or Lipschitzian polyhedron. We assume that $\Omega$ is endowed with a regular mesh $\mathcal{E}_h=\{E_h\}$, i.e. $\Omega =\cup_{E_h\in \mathcal{E}_h} $, where $E_h$ is triangle or quadrilateral ($d=2$) and tetrahedron or hexahedron  ($d=3$) and is $h=\max_{E\in\mathcal{E}_h}{h_E}$, denoting the mesh size, with $h_E={\rm diam }\, E$. For such a mesh, consider $X_h$ to be the finite element space of order one, i.e.,
$$
X_h=\{x_h\in C^0(\bar\Omega)\,:\, x_h|_E\in \mathds{K} \quad \forall E\in\mathcal{E}_h\},
$$ 
where $\mathds{K}=\mathds{P}_1$, the set of linear polynomials on $E$, if $\mathds{K}$ is a triangle or tetrahedron, and $\mathds{K}=\mathds{Q}_1$, the set of bilinear, or trilinear, polynomials on $E$, if $E$ is a quadrilateral or hexahedron, respectively. As $X_h$ is finite dimensional, we set $I={\rm dim}\, X_h$.  The nodal basis of $X_h$ is denoted as $\{\varphi_{\a_i}\}_{i=1}^I$, i.e., $\varphi_{\a_i}(\boldsymbol{a}_j)=\delta_{ij}$ for $i,j=1,\cdots, I$, where $\mathcal{N}_h=\{\boldsymbol{a}_i\}_{i=1}^I$ is the set of nodes of $\mathcal{E}_h$.  For each $i\in I$, it is defined $\Omega_{\a_i}={\rm supp }\,\varphi_{\a_i}$  as its macroelement  and $\mathcal{N}_h(\Omega_{\a_i})$ as the set of nodes belonging to $\Omega_{\a_i}$. Additionally, let  $I(\Omega_{\a_i})=\{j\in I\,:\, \a_j\in \Omega_{\a_i}\}$ be the set of indices of $\Omega_{\a_i}$. It will be finally needed to consider the set $\mathcal{N}_h^{\rm sym }(\Omega_{\boldsymbol{a}_i})$, which consists of the symmetric nodes $\boldsymbol{a}_{ij}^{\rm sym}$ concerning  $\a_i$ constructed as the point at the intersection between the line that passes through $\boldsymbol{a}_i$ and $\boldsymbol{a}_j$ and $\partial\Omega_{\boldsymbol{a}_i}$  not being $\boldsymbol{a}_j$; see Figure \ref{fig:usym}.

Two different interpolation operators will be used: the nodal interpolation operator $i_h$ from $C(\bar\Omega)$ into  $X_h$, i.e,  $i_h x(\boldsymbol{a}_i)=x_h(\boldsymbol{a}_i)$ for $i=1,\cdots, I$, and the averaged interpolation operator  $\mathcal{I}_h$ from $L^p(\Omega)$ into $X_h$ such that 
$$
\mathcal{I}_h\psi =\sum_{i\in I} \left(\frac{1}{|E_{\a_i}|}\int_{E_{\a_i}} \psi(\x)\,{\rm d}\x\right)\varphi_{\a_i},
$$  
where $E_{\a_i}\in\mathcal{E}_h$ is  associated with each $\a_i\in\mathcal{N}_h$ once for all. 
A very substantial advantage of $\mathcal{I}_h$ is that of being stable in $L^p(\Omega)$ and preserving discrete principles. More precisely, it is known that there exists $C_{\rm sta}>0$, independent of $h$, such that  
\begin{equation}\label{I-stability}
\|\mathcal{I}_h \psi\|_{L^{p}(\Omega)}\le C_{\rm sta} \| \psi \|_{L^{p}(\Omega)}\quad \mbox{for } 1\le p\le\infty,
\end{equation}
and 
\begin{equation}\label{I-lower-bounds}
\displaystyle
\underset{\x\in\Omega}{\rm ess\, sup}\,\psi(\x) \ge \mathcal{I}_h \psi (\x) \ge\underset {\x\in\Omega}{\rm ess\, inf}\,\psi(\x)\quad \mbox{ for all }\quad\x\in\Omega  \quad\mbox{ if }\quad \psi\in L^\infty(\Omega). 
\end{equation}

It is also introduced the discrete mass-lumping inner product
$$
(x_h,\bar x_h)_h=\int_\Omega i_h(x_h(\x)\bar x_h(\x))\, {\rm d}\x,
$$
which induces the norm $\|\cdot\|_h$, and the nodal values $\{x_h(\boldsymbol{a}_i)\}_{i=1}^I$  are written in a shorthand as  $\{x_i\}_{i=1}^I$.

\subsection{Finite element approximation problem}
Herein we present our two finite element approximations for problem \eqref{PNP}-\eqref{IC}. With the first algorithm we restrict ourselves to devise a stabilized finite method, which consists in modifying algebraically a standard finite element method by adding accordingly stabilizing terms. Such terms will cure the possible breakdown of the discrete minimum and maximum principles. On the contrary, the second algorithm requires manipulating previously the terms responsible for generating the flux induced due to the electric field, aiming a discrete entropy law. This new setting is able to keep minimum and maximum discrete principles as well.  

Before stating our finite element approximation problems, we set $P_h$ and $N_h$ as $X_h$, while $\Phi_h$ is chosen as $X_h\cap L^2_{\int=0}(\Omega)$. Following this we present the fundamental tools to construct the stabilizing terms, which will be added to the finite element equations for discrete ion densities to achieve the above-mentioned properties. Let $\ell(i,j) = 2\delta_{ij}-1$ be the graph-Laplacian operator, where $\delta_{ij}$ is the Kronecker delta. Furthermore, take $q\in \mathds{R}_+$. For each $i\in I$,  let $\alpha_i: X_h\to \mathds{R}$ such that   
\begin{equation}\label{def:alpha_min_max}
\alpha_i(x_h) = \left\{
\begin{array}{cc}  
\left[
\frac{|{\sum_{j\in I(\Omega_{\boldsymbol{a}_i})} [\![\nabla x_h]\!]_{ij}}|}{\sum_{j\in I(\Omega_{\boldsymbol{a}_i})} 2\{\!\!\{|\nabla x_h \cdot \hat{\boldsymbol{r}}_{ij}\}\!\!\}|_{ij}}
\right]^q & \text{if } 
\sum_{j\in I(\Omega_{\boldsymbol{a}_i})} {\!{|\nabla x_h\cdot \hat{\boldsymbol{r}}_{ij}|}\!}_{ij} \neq 0, 
\\
0 & \text{otherwise},
\end{array}
\right. 
\end{equation}
where
$$
[\![\nabla x_h]\!]_{ij} =\frac{x_j - x_i}{|\boldsymbol{r}_{ij}|} + \frac{x_j^{\rm sym} - x_i}{|\boldsymbol{r}_{ij}^{\rm sym}|}, 
$$
and
$$
\{\!\!\{|\nabla x_h\cdot \hat{\boldsymbol{r}}|_{ij}\}\!\!\}_{ij} =\frac{1}{2}
\left(\frac{|x_j-x_i|}{|\boldsymbol{r}_{ij}|}+\frac{|x_j^{\rm sym}-x_i|}{|\boldsymbol{r}_{ij}^{\rm sym}|}\right),
$$
with $\boldsymbol{r}_{ij} = \boldsymbol{a}_j - \boldsymbol{a}_i$ and $\hat{\boldsymbol{r}}_{ij} = \frac{\boldsymbol{r}_{ij}}{|\boldsymbol{r}_{ij}|}$ its normalized vector. Moreover,  $\boldsymbol{r}_{ij}^{\rm sym } = \boldsymbol{a}_{ij}^{\rm sym } - \boldsymbol{a}_i$, and $x_j^{\rm sym } = x_h(\boldsymbol{a}_{ij}^{\rm sym })$. It should be noted that $[\![\nabla x_h]\!]_{ij}$ and $\{\!\!\{|\nabla x_h\cdot \hat{\boldsymbol{r}}|_{ij}\}\!\!\}_{ij}$ can be thought of as  gradient's linear approximation jump or mean at node $\boldsymbol{a}_i$ in the direction of $\boldsymbol{r}_{ij}$. See \cite{Badia_Bonilla_2017}.  

The employment of the shock detector $\alpha_i$ allows us to localize extrema for a given $x_h\in X_h$ at node $\a_i\in\mathcal{N}_h$. This is formalized in the subsequent result.   
\begin{lemma}\label{lm: alpha_i} Let $x_h\in X_h$. Then it follows that $0\leq \alpha_i(x_h)\leq 1$ and that $\alpha_i(x_h)=1$ for any extreme value at $\a_i$.
\end{lemma}
\subsubsection{Algorithm 1}

We first introduce the nonlinear stabilizing operator $\mathcal{B}_1: X_h\times X_h\cap L^2_0(\Omega) \to \mathcal{L}(X_h', X_h)$ defined as  
\begin{equation}\label{B_alg1}
(\mathcal{B}_1^\pm(x_h, \phi_h)  \tilde x_h, \bar x_h)=\sum_{i\in I}\sum_{j\in I(\Omega_{\boldsymbol{a}_i})}\beta^\pm_{ji}(x_h,\phi_h) \tilde x_j \bar x_i \ell(i,j),
\end{equation}
with
\begin{equation}\label{B_alg1:beta_ij}
\beta^\pm_{ij}(x_h,\phi_h)=\max\{\alpha_i(x_h) f^\pm_{ij}(\phi_h), \alpha_j(x_h)f^\pm_{ji}(\phi_h), 0\},
\end{equation}
for $i\not= j$, and
$$
\beta^\pm_{ii}(x_h,\phi_h)=\sum_{j\in I(\Omega_{\boldsymbol{a}_i})\backslash\{i\}}\beta^\pm_{ij}(x_h,\phi_h),
$$
where $f^\pm_{ij}(\phi_h)$ is given by 
$$
f^\pm_{ij}(\phi_h)=k^{-1}(\varphi_{\boldsymbol{a}_j}, \varphi_{\boldsymbol{a}_i})+(\nabla \varphi_{\boldsymbol{a}_j}, \nabla \varphi_{\boldsymbol{a}_i})\pm(\varphi_{\boldsymbol{a}_j}\nabla\phi_h, \nabla \varphi_{\boldsymbol{a}_i}),
$$
respectively.  

From  \eqref{B_alg1}, it is immediate to deduce the following properties: 
\begin{proposition}\label{Prop:B} Given $p_h, n_h\in X_h$ and $\phi_h\in X_h\cap L^2_0(\Omega)$, it follows that 
\begin{itemize}
\item Mass conservation:  
\begin{equation}\label{B_alg1:mass-conservation} 
(\mathcal{B}^\pm_1(x_h, \phi_h)  \tilde x_h, 1)=0\quad\mbox{for all}\quad x_h\in X_h.
\end{equation} 
\item Positive definiteness:
\begin{equation}\label{B_alg1:definiteness}
(\mathcal{B}^\pm_1(x_h, \phi_h)  \tilde x_h, \tilde x_h)\ge 0\quad\mbox{for all}\quad x_h \in X_h.
\end{equation}
\end{itemize}
\end{proposition}

Let $\{t_m\}_{m=0}^M$ be a uniform partitioning of $[0,T]$ with time step $k=\frac{T}{N}$, where $N\in\mathds{N}$.  Starting from $p^0_h=p_{0h}$ and $n^0_h=n_{0h}$, one wants to find $(p^{m+1}_h, n^{m+1}_h, \phi^{m+1}_h)\in P_h\times N_h\times \Phi_h$ such that, for all $(\bar p_h, \bar n_h, \bar \phi_h)\in P_h\times N_h \times \Phi_h$,    
\begin{subequations}\label{Algorithm1}
\begin{empheq}[left=\empheqlbrace]{align}
(\delta_t p^{m+1}_h, \bar p_h)+(\nabla p^{m+1}_h,\nabla\bar p_h)+(p^{m+1}_h\nabla\phi^{m+1}_h,\nabla\bar p_h )+(\mathcal{B}^+_1(p^{m+1}_h, \phi^{m+1}_h) p^{m+1}_h,\bar p_h )&=0,
\label{eq_alg1:ph} 
\\
(\delta_t n^{m+1}_h, \bar n_h)+(\nabla n^{m+1}_h,\nabla \bar n_h)-(n^{m+1}_h\nabla\phi^{m+1}_h,\nabla \bar n_h)+(\mathcal{B}^-_1(n^{m+1}_h, \phi^{m+1}_h)n^{m+1}_h,\bar n_h )&=0,
\label{eq_alg1:nh}
\\
(\nabla \phi^{m+1}_h, \nabla\phi_h)+(n^{m+1}_h-p^{m+1}_h,\bar\phi_h)_h&=0.
\label{eq_alg1:phih}
\end{empheq}
\end{subequations}

\subsubsection{Algorithm 2}  While the algebraic manipulation approach based on artificial diffusion operator provides a powerful tool to design discrete-principle-bounded,  numerical algorithms, it does not have a clear mechanism to preserve entropy laws. An \emph{ad hoc} strategy to enforce numerical solutions to preserve an entropy law is introducing a suitable discretization of ion fluxes, which takes into account the nonlinear structure of the entropy-Lyapunov functional. Additionally, in the construction of the method, and its analysis, the mass-lumping $L^2(\Omega)$-inner product will be of importance. 
 Therefore, for our second algorithm, we want to replace the ion transport term $(x_h\nabla \phi_h, \nabla \bar x_h)$ by $(x_h\nabla\phi_h, \nabla\bar x_h )_*$, which must allow for choosing $\bar x_h = i_h(\log x_h)$, as a test function,  in order to obtain an approximation of $(\nabla \phi_h, \nabla x_h)$. In doing so, it will be used the following \cite{Badia_Bonilla_GS_2023}. Let $\varepsilon>0$ and define
$$
g_\varepsilon(s)=\left\{
\begin{array}{ccl}
s\log s-s+1&\mbox{ if }& s>\varepsilon,
\\
\frac{s^2-\varepsilon^2}{2\varepsilon}+(\log \varepsilon-1)s+1&\mbox{ if }& s\le \varepsilon,
\end{array}
\right.
$$
and hence
$$
g'_\varepsilon(s)=\left\{
\begin{array}{ccl}
\log s&\mbox{ if }& s>\varepsilon,
\\
\frac{s}{\varepsilon}+\log \varepsilon-1&\mbox{ if }& s\le \varepsilon. 
\end{array}
\right.
$$ 
Thus
\begin{equation}\label{Elec_Trans_Term_new}
(x_h\nabla \phi_h, \nabla \bar x_h)_*=\sum_{i<j\in I}\tau_{ji}(x_h) \delta_{ji} \phi_h \delta_{ij} x_h(\nabla \varphi_{\boldsymbol{a}_j},\nabla\varphi_{\boldsymbol{a}_i}),
\end{equation}
with
\begin{equation}\label{def:tau_ij}
\tau_{ji}(x_h)=\left\{
\begin{array}{ccl}
\displaystyle
\frac{\delta_{ji}x_h}{\delta_{ji}g'_\varepsilon(x_h)}&\mbox{ if }& x_j\not=x_i,
\\
\max\{x_i, \varepsilon\}&\mbox{ if }& x_j=x_i,
\end{array}
\right.
\end{equation}
where the functional $[x]_+=\max\{0,x\}$ stands for the positivity part and $\delta_{ji} f(x_h)=f(x_j)-f(x_i)$ for any $f:\mathds{R}\to\mathds{R}$. As $g'_\varepsilon$ is bijective, we are allowed to use $x_j=x_i$ instead  $g'_\varepsilon(x_j)=g'_\varepsilon(x_i)$ when defined $\tau_{ji}(\cdot)$. 

Algorithm $2$ reads as follows. 
Starting from $p^0_h=p_{0h}$ and $n^0_h=n_{0h}$, find $(p^{m+1}_h, n^{m+1}_h, \phi^{m+1}_h)\in P_h\times N_h \times \Phi_h$ such that, for all $(\bar p_h, \bar n_h, \bar \phi_h)\in P_h\times N_h \times \Phi_h$,   
\begin{subequations}\label{Algorithm2}
\begin{empheq}[left=\empheqlbrace]{align}
(\delta_t p^{m+1}_h, \bar p_h)_h+(\nabla p^{m+1}_h,\nabla\bar p_h)+(p^{m+1}_h\nabla\phi^{m+1}_h,\nabla\bar p_h )_*+(\mathcal{B}^+_2(p^{m+1}_h, \phi^{m+1}_h)p^{m+1}_h,\bar p_h )&=0,
\label{eq_alg2:ph} 
\\
(\delta_t n^{m+1}_h, \bar n_h)_h+(\nabla n^{m+1}_h,\nabla \bar n_h)-(n^{m+1}_h\nabla\phi^{m+1}_h,\nabla \bar n_h)_*+(\mathcal{B}^-_2(n^{m+1}_h, \phi^{m+1}_h)n^{m+1}_h,\bar n_h )&=0,
\label{eq_alg2:nh}
\\
(\nabla \phi^{m+1}_h, \nabla\phi_h)+(n^{m+1}_h-p^{m+1}_h,\bar\phi_h)_h&=0.
\label{eq_alg2:phih}
\end{empheq}
\end{subequations}

The stabilizing operators  are now defined as follows \cite{Badia_Bonilla_GS_2023}. Recall $\delta_{ij} x_h =x_j - x_i$ and define $\mathcal{B}_2: X_h \times X_h\cap L^2_0(\Omega) \to \mathcal{L}(X'_h,X_h)$ such that 
\begin{equation}\label{B_alg2}
(\mathcal{B}^\pm_2(x_h, \phi_h) \tilde x_h, \bar x_h)=\sum_{i<j\in I}\tilde\beta^\pm_{ji}(x_h, \phi_h) \delta_{ji} \tilde x_h \delta_{ji}\bar x_h,
\end{equation} 
where 
\begin{equation}\label{beta_alg2:ij}
\tilde \beta^\pm_{ij}(x_h,\phi_h)=\max\{\alpha_i(x_h) f^\pm_{ij}(x_h,\phi_h), \alpha_j(x_h)f^\pm_{ji}(x_h,\phi_h), 0\},
\end{equation}
with
$$
f^\pm_{ij}(x_h,\phi_h)=\left\{
\begin{array}{rcl}
\displaystyle
\left(1\pm\delta_{ji}\phi_h\left[\frac{1}{\delta_{ji} g'_\varepsilon(x_h)} - \frac{\max\{ x_i,\varepsilon\}}{\delta_{ji} x_h}\right]\right)(\nabla \varphi_{\boldsymbol{a}_j}, \nabla \varphi_{\boldsymbol{a}_i})&\mbox{ if }&  x_j\not= x_i,
\\
0&\mbox{ if }& x_j= x_i,
\end{array}
\right.
$$
and
$$
\tilde\beta^\pm_{ii}(x_h, \phi_h)=\sum_{j\in I(\Omega_{\boldsymbol{a}_i})\backslash\{i\}}\tilde \beta_{ij}^\pm(x_h, \phi_h).
$$

 The stabilizing term $\mathcal{B}^\pm_2$ satisfies the same properties as $\mathcal{B}^\pm_1$ does.  
\begin{proposition} It follows that 
\begin{itemize}
\item Mass conservation:
\begin{equation}\label{B_alg2:mass-conservation} 
(\mathcal{B}^\pm_2(x_h, \phi_h)  \tilde x_h, 1)=0.
\end{equation}
\item Positive definiteness: 
\begin{equation}\label{B_alg2:positive-definiteness} 
(\mathcal{B}^\pm_2(x_h,\phi_h)  \tilde x_h, \tilde x_h)\ge 0.
\end{equation}
\end{itemize}
\end{proposition}

\section{Statement of main results}
Herein will be given the statement of our two main results for Algorithms 1 and 2. They will be proved separately in subsequent sections. These proofs will be made clear how the structure of both algorithms interacts with the properties that are required for each.

\begin{theorem}\label{Th:alg1} Let $p_0, n_0\in L^\infty(\Omega)\cap C(\bar \Omega)$ be such that  $p_0, n_0>0$ satisfying the electroneutrality condition $\|p_0\|_{L^1(\Omega)}=\|n_0\|_{L^1(\Omega)}$.  Choose $p_{0h}=\mathcal{I}_{h} p_0$ and $n_{0h}=\mathcal{I}_h n_0$. Moreover, assume that  $k$ is sufficiently small that
\begin{equation}\label{DMP_alg1:smallness_on_k}
1 - k \Big(\max\{\underset {\x\in\Omega}{\max}\,p^0_h(\x),\underset {\x\in\Omega}{\max }\,n^0_h(\x)\} -\min\{\underset {\x\in\Omega}{\min}\,p^0_h(\x),\underset {\x\in\Omega}{\min }\,n^0_h(\x)\}\Big) >0
\end{equation} 
holds.  Then any sequence of discrete solutions $\{(p^m_h, n^m_h, \phi^m_h)\}_{m=1}^{M}\subset P_h\times N_h\times\Phi_h$ constructed through Algorithm $1$ satisfies the following properties:
\begin{enumerate}[$(i)$]
\item Discrete principle. For $s=0, \cdots, M$, 
\begin{equation}\label{DMP_alg1:ph}
\min\{\underset {\x\in\Omega}{\min}\,p^0_h(\x),\underset {\x\in\Omega}{\min }\,n^0_h(\x)\}\le p^s_h(\x) \le \max\{\underset {\x\in\Omega}{\max}\,p^0_h(\x),\underset {\x\in\Omega}{\max }\,n^0_h(\x)\} \quad\mbox{ for all } \x\in \Omega
\end{equation}
and
\begin{equation}\label{DMP_alg1:nh}
\min\{\underset {\x\in\Omega}{\min}\,p^0_h(\x),\underset {\x\in\Omega}{\min }\,n^0_h(\x)\}\le n^s_h(\x) \le  \max\{\underset {\x\in\Omega}{\max}\,p^0_h(\x),\underset {\x\in\Omega}{\max }\,n^0_h(\x)\} \quad\mbox{ for all } \x\in \Omega. 
\end{equation}
\item Mass conservation. For $s=0, \cdots, M$, 
\begin{equation}\label{MC_alg1:ph}
\|p^{m}_h\|_{L^1(\Omega)}=\|p^0_h\|_{L^1(\Omega)}
\end{equation}
and
\begin{equation}\label{MC_alg1:nh}
\|n^{m}_h\|_{L^1(\Omega)}=\|n^0_h\|_{L^1(\Omega)}.
\end{equation}
\end{enumerate}
\end{theorem}

Before stating the theorem for Algorithm 2, we introduce some short-hand notation. Let $g_0(x)= x \log x -x+1$. Then define 
$$
\mathcal{E}_h(p_h, n_h, \phi_h)=(g_0(p_h),1)_h+(g_0(n_h),1)_h+\frac{1}{2}\|\nabla \phi_h\|^2_{L^2(\Omega)},
$$
and 
$$
\begin{array}{rcl}
\mathcal{D}_h(\rho_h, \phi_h)&=&\displaystyle
-\sum_{i<j\in \mathds{I}^*(\rho_h)}\left|\left(\frac{\delta_{ji} g'_0(\rho_h)}{\delta_{ji}\rho_h}\right)^{-\frac{1}{2}}\delta_{ji} \rho_h
-\left(\frac{\delta_{ji}g'_0(\rho_h)}{\delta_{ji} \rho_h}\right)^{\frac{1}{2}}\delta_{ji}\phi_h\right|^2(\nabla\varphi_{\a_j},\nabla\varphi_{\a_i})
\\
&&\displaystyle
-\sum_{i<j\in \mathds{I}^*_c(\rho_h)}\rho_i \delta^2_{ij}\phi_h(\nabla\varphi_{\a_j},\nabla\varphi_{\a_i}),
\end{array}
$$
where $\mathds{I}^*(\rho_h)=\{(i,j)\in I\times I : \rho_j\not = \rho_i \}$ and $\mathds{I}^*_c(\rho_h)$ its complementary. Further $\mathds{I}=I\times I$. 
\begin{theorem}\label{Th:alg2}  Let $p_0, n_0\in L^\infty(\Omega)$ be such that  $p_0, n_0>0$ satisfying the electroneutrality condition $\|p_0\|_{L^1(\Omega)}=\|n_0\|_{L^1(\Omega)}$.  Choose $p_{0h}=\mathcal{I}_{h} p_0$ and $n_{0h}=\mathcal{I}_h n_0$.  Then any sequence of discrete solutions $\{(p^m_h, n^m_h, \phi^m_h)\}_{m=1}^{M}\subset P_h\times N_h\times\Phi_h$ constructed through Algorithm $2$ satisfies the following properties:
\begin{enumerate}[$(i)$]
\item Discrete principles.  For $s=0, \cdots, M$,
\begin{equation}\label{DMP_alg2:ph}
\min\{\underset {\x\in\Omega}{\min }\,p^0_h(\x),\underset {\x\in\Omega}{\min }\,n^0_h(\x)\}\le p^s_h(\x) \le \max\{\underset {\x\in\Omega}{\max }\,p^0_h(\x),\underset {\x\in\Omega}{\max }\,n^0_h(\x)\}
\end{equation} 
and
\begin{equation}\label{DMP_alg2:nh}
\min\{\underset {\x\in\Omega}{\min }\,p^0_h(\x),\underset {\x\in\Omega}{\min }\,n^0_h(\x)\}\le n^s_h(\x) \le \max\{\underset {\x\in\Omega}{\max }\,p^0_h(\x),\underset {\x\in\Omega}{\max }\,n^0_h(\x)\}.
\end{equation} 
\item Mass conservation. For $s=0, \cdots, M$,
\begin{equation}\label{MC_alg2:ph}
\|p^s_h\|_{L^1(\Omega)}=\|p^0_h\|_{L^1(\Omega)}
\end{equation}
and
\begin{equation}\label{MC_alg2:nh}
\|n^s_h\|_{L^1(\Omega)}=\|n^0_h\|_{L^1(\Omega)}.
\end{equation}
\item Discrete entropy law: Suppose that $\mathcal{E}_h$ is strictly acute i.e., there exists $C_{\rm ang}>0$ with
\begin{equation}\label{Alg_2:acuteness}
(\nabla \varphi_{\a_i},\nabla\varphi_{\a_j})\le - C_{\rm ang} \quad\mbox{ for all }\quad 
j\not= i\in I.
\end{equation}
Then,  for $m=0, \cdots, M$, it follows that 
\begin{equation}\label{Entropy:alg2}
\begin{array}{r}
\mathcal{E}_h(p^s_h, n^s_h, \phi^s_h)+k\displaystyle\sum_{m=0}^{s-1}\mathcal{ND}_h(p^{m+1}_h, p^{m}_h, n^{m+1}_h, n^m_h, \phi^{m+1}_h)
\\
+k\displaystyle\sum_{m=0}^{s-1} \Big(\mathcal{D}_h(p^{m+1}_h, \phi^{m+1}_h) + \mathcal{D}_h(n^{m+1}_h, \phi^{m+1}_h) \Big)  \le \mathcal{E}_h(p^{0}_h, n^{0}_h, \phi^0_h),
\end{array}
\end{equation}
where $\alpha_\#(\rho^{m+1}_h)\in\{\alpha_i(\rho^{m+1}_h), \alpha_j(\rho^{m+1}_h), 0\}$ is  such that 
$$\tilde\beta^+_{ji}(\rho^{m+1}_h,\phi^{m+1}_h)=\alpha_\#(\rho^{m+1}_h)\left[1+\delta_{\bar\#\#}\phi^{m+1}_h\left(\frac{1}{\delta_{\bar\#\#}g'_0(\rho^{m+1}_h)}-\frac{\rho^{m+1}_\#}{\delta_{\bar\#\#} \rho^{m+1}_h)}\right)\right](\nabla\varphi_{\a_{\bar \#}},\nabla\varphi_{\a_\#}),$$
with ${\bar \#}=j$ when $\#=i$ and ${\bar \#}=i$ when $\#=j$, and
$$
\begin{array}{rcl}
\mathcal{ND}_h(p^{m+1}_h, p^{m}_h, n^{m+1}_h, n^m_h, \phi^{m+1}_h)&=&\displaystyle\frac{k}{2} (g''_0(p^{m+\theta_p}_h),(\delta_t p^{m+1}_h)^2)_h+\frac{k}{2} (g''_0(n^{m+\theta_n}_h),(\delta_t n^{m+1}_h)^2)_h
\end{array}
$$
with $\theta_\rho\in (0,1)$ such that $ \rho^{m+\theta_\rho}_h=\theta_\rho \rho^{m+1}_h + (1-\theta_\rho) \rho^m_h$ for $\rho=p$ or $n$.
\end{enumerate}
\end{theorem}

\section{Proof of Theorem \ref{Th:alg1}}

$\bullet$ \textbf{Discrete principles } 
Let us consider the following auxiliary algorithm. Find $(p^{m+1}_h, n^{m+1}_h, \phi^{m+1}_h)\in P_h\times N_h\times \Phi_h$ such that, for all $(\bar p_h,\bar n_h, \bar \phi_h)\in P_h\times N_h\times\Phi_h$,  
\begin{subequations}\label{Algorithm1_Aux}
\begin{empheq}[left=\empheqlbrace]{align}
(\delta_t p^{m+1}_h, \bar p_h)+(\nabla p^{m+1}_h,\nabla\bar p_h)+(p^{m+1}_h\nabla\phi^{m+1}_h,\nabla\bar p_h )+(\mathcal{B}^+_1(p^{m+1}_h,\phi^{m+1}_h)p^{m+1}_m,\bar p_h )&=0,
\label{eq_alg2_aux:ph} 
\\
(\delta_t n^{m+1}_h, \bar n_h)+(\nabla n^{m+1}_h,\nabla \bar n_h)-(n^{m+1}_h\nabla\phi^{m+1}_h,\nabla \bar n_h)+(\mathcal{B}^-_1( n^{m+1}_h,\phi^{m+1}_h)n^{m+1}_h,\bar n_h )&=0, 
\label{eq_alg2_aux:nh}
\\
(\nabla \phi^{m+1}_h, \nabla \bar\phi_h)+([n^{m+1}_h]_T-[p^{m+1}_h]_T,\bar\phi_h)_h&=0,
\label{eq_alg2_aux:phih}
\end{empheq}
\end{subequations}
where $[\cdot]_T$ is the nodal truncating operator between $\min\{\underset {\x\in\Omega}{\min}\, p^0_h(\x),  \underset {\x\in\Omega} \min\,n^0_h(\x) \}$ and $\max\{\underset {\x\in\Omega}{\max} \,p^0_h(\x),\break \underset {\x\in\Omega}{\max} \,n^0_h(\x)\}$, respectively. It is forthwith seen that \eqref{Algorithm1_Aux} becomes \eqref{Algorithm1} once a discrete principle is proved with $0<\varepsilon$ satisfying $\varepsilon<\min\{\underset {\x\in\Omega}{\min}\, p^0_h(\x),  \underset {\x\in\Omega} \min\,n^0_h(\x) \}$. 

Assume that \eqref{DMP_alg1:ph} and \eqref{DMP_alg1:nh} are satisfied for $s=m$. Then we wish to show \eqref{DMP_alg1:ph} and \eqref{DMP_alg1:nh} for $s=m+1$ by induction. Let us focus on $p^{m+1}_h$; the case $n^{m+1}_h$ is analogous. 

We begin by proving the discrete minimum principle, i.e, the lower bound in \eqref{DMP_alg1:ph}. Then if this were false, there would exist a node  $\boldsymbol{a_i}\in\mathcal{N}_h$ such that  $p^{m+1}_h$ has a local minimum satisfying  
\begin{equation}\label{DMP_alg1:cont_min}
p^{m+1}_h(\boldsymbol{a}_i):=p^{m+1}_i<\min\{\underset {\x\in\Omega}{\min}\,p^0_h(\x),\underset {\x\in\Omega}{\min }\,n^0_h(\x)\}.
\end{equation}
Selecting $\bar p_h=\varphi_{\boldsymbol{a}_i}$ in \eqref{eq_alg2_aux:ph}  gives 
\begin{equation}\label{DMP_alg1-lab1}
\begin{array}{c}
\displaystyle
\sum_{j \in I(\Omega_{\boldsymbol{a}_i})}p^{m+1}_j\{k^{-1}(\varphi_{\boldsymbol{a}_j}, \varphi_{\boldsymbol{a}_i})+(\nabla\varphi_{\boldsymbol{a}_j},\nabla \varphi_{\boldsymbol{a}_i})+(\varphi_{\boldsymbol{a}_j}\nabla\phi^{m+1}_h,\nabla \varphi_{\boldsymbol{a}_i})
\\
\displaystyle
+(\mathcal{B}^+_1(p^{m+1}_h,\phi^{m+1}_h)\varphi_{\boldsymbol{a}_j}, \varphi_{\boldsymbol{a}_i})\}=k^{-1}\sum_{j \in I(\Omega_{\boldsymbol{a}_i})} p^m_j (\varphi_{\boldsymbol{a}_j}, \varphi_{\boldsymbol{a}_i}).
\end{array} 
\end{equation}
By the very definition of  $\mathcal{B}^+_1$ in \eqref{B_alg1} in conjunction with \eqref{B_alg1:beta_ij}, we obtain, on noting $\alpha_i(p^{m+1}_h)=1$ from Lemma~\ref{lm: alpha_i}, that, for $i\not = j$, 
\begin{equation}\label{DMP_alg1-lab2}
(\varphi_{\boldsymbol{a}_j}, \varphi_{\boldsymbol{a}_i})+(\nabla\varphi_{\boldsymbol{a}_j},\nabla \varphi_{\boldsymbol{a}_i})+(\varphi_{\boldsymbol{a}_j}\nabla\phi^{m+1}_h,\nabla \varphi_{\boldsymbol{a}_i})+(\mathcal{B}^+_1(p^{m+1}_h,\phi^{m+1}_h)\varphi_{\boldsymbol{a}_j}, \varphi_{\boldsymbol{a}_i})\le 0, 
\end{equation}
and hence, from \eqref{DMP_alg1-lab1}, we bound 
$$
\begin{array}{c}
\displaystyle
	\sum_{j \in I(\Omega_{\boldsymbol{a}_i})}p^{m+1}_i\{ k^{-1}(\varphi_{\boldsymbol{a}_j},
	\varphi_{\boldsymbol{a}_i})+(\nabla \varphi_{\boldsymbol{a}_j},\nabla \varphi_{\boldsymbol{a}_i})+	
	(\varphi_{\boldsymbol{a}_j}\nabla\phi^{m+1}_h,\nabla \varphi_{\boldsymbol{a}_i})
\\
\displaystyle
	+(\mathcal{B}^+_1(p^{m+1}_h, n^{m+1}_h)\varphi_{\boldsymbol{a}_j}, \varphi_{\boldsymbol{a}_i})\}\ge k^{-1}\sum_{j\in  I(\Omega_{\boldsymbol{a}_i})} 
	p^m_j (\varphi_{\boldsymbol{a}_j}, \varphi_{\boldsymbol{a}_i}),
\end{array}
$$
since $p^{m+1}_i<p^{m+1}_j$ for all $j\in I(\Omega_{\boldsymbol{a}_i})$. Further use of  
\begin{equation}\label{DMP_alg1-lab3}
\sum_{j \in I(\Omega_{\boldsymbol{a}_i})} (\varphi_{\boldsymbol{a}_j}, \varphi_{\boldsymbol{a}_i})= 
\|\varphi_{\boldsymbol{a}_i}\|_{L^1(\Omega)},
\end{equation}
\begin{equation}\label{DMP_alg1-lab4}
\sum_{j \in I(\Omega_{\boldsymbol{a}_i})} (\nabla \varphi_{\boldsymbol{a}_j}, \nabla\varphi_{\boldsymbol{a}_i})= (\nabla 1,\nabla\varphi_{\boldsymbol{a}_i})=0,
\end{equation}
\begin{equation}\label{DMP_alg1-lab5}
\sum_{j \in I(\Omega_{\boldsymbol{a}_i})} (\varphi_{\boldsymbol{a}_j}\nabla\phi^{m+1}_h, \nabla\varphi_{\boldsymbol{a}_i})= (\nabla\phi^{m+1}_h,\nabla \varphi_{\boldsymbol{a}_i}),
\end{equation}
and, by \eqref{B_alg1:mass-conservation},
\begin{equation}\label{DMP_alg1-lab6}
\sum_{j \in I(\Omega_{\boldsymbol{a}_i})}(\mathcal{B}^+_1(p^{m+1}_h, \phi^{m+1}_h)\varphi_{\boldsymbol{a}_j}, \varphi_{\boldsymbol{a}_i})=(\mathcal{B}^+_1(p^{m+1}_h, \phi^{m+1}_h) 1, \varphi_{\boldsymbol{a}_i})=0 
\end{equation}
is made to get  
\begin{equation}\label{DMP_alg1-lab7}
p^{m+1}_i\left(1+ k \frac{(\nabla \phi^{m+1}_h, \nabla\varphi_{\boldsymbol{a}_i})}{\|
\varphi_{\boldsymbol{a}_i}\|_{L^1(\Omega)}}\right)\ge \sum_{j\in I(\Omega_{\boldsymbol{a}_i})}p^m_j 
\frac{(\varphi_{\boldsymbol{a}_j}, \varphi_{\boldsymbol{a}_i})}{\|\varphi_{\boldsymbol{a}_i}\|_{L^1(\Omega)}}.
\end{equation}
Let us choose  $\bar \phi_h=\varphi_{\boldsymbol{a}_i}$ in \eqref{eq_alg2_aux:phih}, and thereafter insert it into \eqref{DMP_alg1-lab7} and bound to get  
$$
p^{m+1}_i\left(1+ k \frac{([p^{m+1}_h]_T-[n^{m+1}_h]_T,\varphi_{\boldsymbol{a}_i})_h}{\|
\varphi_{\boldsymbol{a}_i}\|_{L^1(\Omega)}}\right)\ge\sum_{j\in I(\Omega_{\boldsymbol{a}_i})}p^m_j 
\frac{(\varphi_{\boldsymbol{a}_j}, \varphi_{\boldsymbol{a}_i})}{\|\varphi_{\boldsymbol{a}_i}\|_{L^1(\Omega)}}\ge \min_{j\in I(\Omega_{\a_i})} p^m_j,
$$
which, in turn, implies, on noting that 
$$
0<1+ k \frac{([p^{m+1}_h]_T-[n^{m+1}_h]_T,\varphi_{\boldsymbol{a}_i})_h}{\|
\varphi_{\boldsymbol{a}_i}\|_{L^1(\Omega)}} <1
$$
holds from \eqref{DMP_alg1:smallness_on_k} and \eqref{DMP_alg1:cont_min}, that
\begin{equation}\label{DMP_alg1-lab9}
\begin{array}{rcl}
\min\{\underset {\x\in\Omega}{\min}\,p^0_h(\x),\underset {\x\in\Omega}{\min }\,n^0_h(\x)\}
&>&\displaystyle p^{m+1}_i\left(1+ k \frac{([p^{m+1}_h]_T-[n^{m+1}_h]_T,\varphi_{\boldsymbol{a}_i})_h}{\|
\varphi_{\boldsymbol{a}_i}\|_{L^1(\Omega)}}\right)
\\
&\ge& \min\{\underset {\x\in\Omega}{\min}\,p^0_h(\x),\underset {\x\in\Omega}{\min }\,n^0_h(\x)\}.
\end{array}
\end{equation}
This contradiction proves that the discrete minimum principle holds for $p^{m+1}_h$  as desired. 

To face the discrete maximum principle for $p^{m+1}_h$, i.e. the upper bound in \eqref{DMP_alg1:ph}, we proceed again by contradiction assuming that $p^{m+1}_h$ attains a local maximum at $\boldsymbol{a}_i\in\mathcal{N}_h$ such that 
\begin{equation}\label{DMP_alg1:cont_max}
p^{m+1}_i>\max\{\underset {\x\in\Omega}{\max}\,p^0_h(\x),\underset {\x\in\Omega}{\max }\,n^0_h(\x)\}.
\end{equation} 
Thus we take $\bar p_h=\varphi_{\boldsymbol{a}_i}$ in \eqref{eq_alg2:ph} to get, on noting \eqref{DMP_alg1-lab3}--\eqref{DMP_alg1-lab6} and using \eqref{eq_alg2_aux:phih} for $\bar\phi_h=\varphi_{\a_i}$, that   
$$
p^{m+1}_i\left(1+ k \frac{([p^{m+1}_h]_T-[n^{m+1}_h]_T,\varphi_{\boldsymbol{a}_i})_h}{\|\varphi_{\boldsymbol{a}_i}\|_{L^1(\Omega)}}\right)\le \sum_{j\in I(\Omega_{\a_i})} p^m_j \frac{(\varphi_{\boldsymbol{a}_j}, \varphi_{\boldsymbol{a}_i})}{\|\varphi_{\boldsymbol{a}_i}\|_{L^1(\Omega)}}.
$$
We next straightforwardly see that  $([p^{m+1}_h]_T-[n^{m+1}_h]_T,\varphi_{\boldsymbol{a}_i})_h\ge0$ from \eqref{DMP_alg1:cont_max} and thereby finding 
$$
\begin{array}{rcl}
\max\{\underset {\x\in\Omega}{\max}\,p^0_h(\x),\underset {\x\in\Omega}{\max }\,n^0_h(\x)\}&<&\displaystyle p^{m+1}_i\left(1+ k \frac{([p^{m+1}_h]_T-[n^{m+1}_h]_T,\varphi_{\boldsymbol{a}_i})_h}{\|\varphi_{\boldsymbol{a}_i}\|_{L^1(\Omega)}}\right)
\\
&\le& \max\{\underset {\x\in\Omega}{\max}\,p^0_h(\x),\underset {\x\in\Omega}{\max }\,n^0_h(\x)\},
\end{array}
$$
a contradiction. 

$\bullet$ \textbf{Mass conservation} The proof is again established by induction. Let \eqref{MC_alg1:ph} and \eqref{MC_alg1:nh} hold for $s=m$.  On substituting $\bar p_h=1$ and $\bar n_h=1$ into \eqref{eq_alg1:ph} and \eqref{eq_alg1:nh}, respectively, it follows, on using \eqref{B_alg1:mass-conservation}, that 
$$
\int_\Omega p^{m+1}_h(\x)\,{\rm d}\x=\int_\Omega p^m_h(\x)\,{\rm d}\x
$$
and
$$
\int_\Omega n^{m+1}_h\,(\x){\rm d}\x=\int_\Omega n^m_h(\x)\,{\rm d}\x.
$$
Now, from \eqref{DMP_alg1:ph} and \eqref{DMP_alg1:nh}, the proof is complete.

\section{Proof of Theorem \ref{Th:alg2}} 
$\bullet $ \textbf{Discrete principles} 
In a similar fashion as in the proof of Theorem \ref{Th:alg1}, we suppose that  there exists a node $\boldsymbol{a}_i\in\mathcal{N}_h$ such that $p^{m+1}_h$ has a local minimum so that  
\begin{equation}\label{DMP_alg2:cont_min}
p^{m+1}_h(\boldsymbol{a}_i):=p^{m+1}_i<\min\{\underset {\x\in\Omega}{\min}\,p^0_h(\x),\underset {\x\in\Omega}{\min }\,n^0_h(\x)\}.
\end{equation}
 
In order to simplify the proof, we may assume  without loss of generality  that $i<j$ for all $ j\in I(\Omega_{\boldsymbol{a}_i})$. Let $I^*(\Omega_{\boldsymbol{a}_i})=\{j\in I(\Omega_{\boldsymbol{a}_i}): p^{m+1}_j=p^{m+1}_i\}$ and let  $I^*_c(\Omega_{\boldsymbol{a}_i})$ be  its complementary. Choosing  $\bar p_h=\varphi_{\boldsymbol{a}_i}$ in \eqref{eq_alg2:ph} yields 
$$
(\delta_t p^{m+1}_h, \varphi_{\boldsymbol{a}_i})_h+(\nabla p^{m+1}_h,\nabla \varphi_{\boldsymbol{a}_i})+(p^{m+1}_h\nabla\phi^{m+1}_h,\nabla \varphi_{\boldsymbol{a}_i})_*+(\mathcal{B}^+_2(p^{m+1}_h,\phi^{m+1}_h)p^{m+1}_h, \varphi_{\boldsymbol{a}_i} )=0, 
$$
which can be rewritten, on noting that 
$$
\begin{array}{rcl}
(\nabla p^{m+1}_h,\nabla \varphi_{\boldsymbol{a}_i})
&=&\displaystyle
\sum_{j\in I^*_c(\Omega_{\boldsymbol{a}_i})} (\nabla\varphi_{\boldsymbol{a}_j},\nabla\varphi_{\boldsymbol{a}_i}) \delta_{ji} p^{m+1}_h,
\end{array}
$$
that, from \eqref{Elec_Trans_Term_new},
\begin{align*}
(p^{m+1}_h\nabla  \phi^{m+1}_h, \nabla \varphi_{\boldsymbol{a}_i})_*=&\displaystyle
\sum_{j\in I^*_c(\Omega_{\boldsymbol{a}_i})} \frac{\delta_{ji}\phi^{m+1}_h}{\delta_{ji}g'_\varepsilon(p^{m+1}_h)}(\nabla \varphi_{\boldsymbol{a}_j},\nabla\varphi_{\boldsymbol{a}_i}) \delta_{ji}p^{m+1} _h 
\\
&\displaystyle
+\sum_{j\in I^*(\Omega_{\boldsymbol{a}_i})} \max\{p^{m+1}_j, \varepsilon\} (\nabla \varphi_{\boldsymbol{a}_j},\nabla\varphi_{\boldsymbol{a}_i}) \delta_{ji}\phi^{m+1}_h 
\\\notag
=&\displaystyle
\sum_{j\in I^*_c(\Omega_{\boldsymbol{a}_i})} \frac{\delta_{ji}\phi^{m+1}_h}{\delta_{ji}g'_\varepsilon(p^{m+1}_h)}(\nabla \varphi_{\boldsymbol{a}_j},\nabla\varphi_{\boldsymbol{a}_i}) \delta_{ji} p^{m+1} _h
\\
&\displaystyle
-\sum_{j\in I^*_c(\Omega_{\boldsymbol{a}_i})} \delta_{ji}\phi^{m+1}_h \frac{\max\{p^{m+1}_i, \varepsilon\}}{\delta_{ji} p^{m+1}_h}  (\nabla \varphi_{\boldsymbol{a}_j},\nabla\varphi_{\boldsymbol{a}_i})\delta_{ji} p^{m+1}_h
\\
&+ \max\{p^{m+1}_i, \varepsilon\} (\nabla \phi^{m+1}_h,\nabla\varphi_{\boldsymbol{a}_i}),
\end{align*}
and that, from \eqref{B_alg2},
$$
\begin{array}{rcl}
(\mathcal{B}_2^+(p^{m+1}_h, \phi^{m+1}_h)p^{m+1}_n, \varphi_{\boldsymbol{a}_i} )
&=&\displaystyle
-\sum_{j\in I^*_c(\Omega_{\boldsymbol{a}_i})} \tilde\beta_{ji}^+(p^{m+1}_h, \phi^{m+1}_h) \delta_{ji} p^{m+1}_h
\end{array}
$$
as 
\begin{equation}\label{DMP_alg2:lab1}
\begin{array}{l}
k^{-1}(1,\varphi_{\boldsymbol{a}_i}) p^{m+1}_i\displaystyle
\displaystyle
+\sum_{j\in I^*_c(\Omega_{\boldsymbol{a}_i})}\left(1+\frac{\delta_{ji}\phi^{m+1}_h}{\delta_{ji}g'_\varepsilon(p^{m+1}_h)}\right)(\nabla\varphi_{\boldsymbol{a}_j},\nabla\varphi_{\boldsymbol{a}_i})  \delta_{ji} p^{m+1}_h 
\\
\displaystyle
-\sum_{j\in I^*_c(\Omega_{\boldsymbol{a}_i})} \delta_{ji}\phi^{m+1}_h \frac{\max\{p^{m+1}_i, \varepsilon\}}{\delta_{ji} p^{m+1}_h}  (\nabla \varphi_{\boldsymbol{a}_j},\nabla\varphi_{\boldsymbol{a}_i})\delta_{ji} p^{m+1}_h
\\
+ \max\{p^{m+1}_i, \varepsilon\} (\nabla \phi^{m+1}_h,\nabla\varphi_{\boldsymbol{a}_i})
\\
\displaystyle
-\sum_{j\in I^*_c(\Omega_{\boldsymbol{a}_i})}
\tilde\beta_{ji}^+(p^{m+1}_h, \phi^{m+1}_h) \delta_{ji} p^{m+1}_h
= k^{-1}  (1,\varphi_{\boldsymbol{a}_i}) p_i^m.
\end{array}
\end{equation}
In accordance with \eqref{B_alg2} -- particularly from \eqref{beta_alg2:ij} --,  we know that, for all $j\in I(\Omega_{\a_i})$,
$$
\left(1+\delta_{ji}\phi^{m+1}_h\left[\frac{1}{\delta_{ji} g'_\varepsilon(p^{m+1}_h)} - \frac{\max\{p^{m+1}_i, \varepsilon\}}{\delta_{ji} p^{m+1}_h}\right]\right)
(\nabla\varphi_{\boldsymbol{a}_j},\nabla\varphi_{\boldsymbol{a}_i})-\tilde\beta_{ji}^+(p^{m+1}_h,\phi^{m+1}_h)\le 0
$$
holds from Lemma \ref{lm: alpha_i}, because of $\alpha_i(p^{m+1}_h)=1$, and hence
\begin{equation}\label{DMP_alg2:lab2}
\begin{array}{l}
\displaystyle
\left[\left(1+\delta_{ji}\phi^{m+1}_h\left[\frac{1}{\delta_{ji} g'_\varepsilon(p^{m+1}_h)} - \frac{\max\{p^{m+1}_i, \varepsilon\}}{\delta_{ji} p^{m+1}_h}\right]\right)
(\nabla\varphi_{\boldsymbol{a}_j},\nabla\varphi_{\boldsymbol{a}_i})\right.
\\
\displaystyle
-\tilde\beta_{ji}^+(p^{m+1}_h,\phi^{m+1}_h)\Big](p^{m+1}_j- p^{m+1}_i)\le 0,
\end{array}
\end{equation}
since $p^{m+1}_h$ possesses a minimum on $\Omega_{\boldsymbol{a}_i}$ by supposition. Combing \eqref{DMP_alg2:lab1} and \eqref{DMP_alg2:lab2} yields
\begin{equation}\label{DMP_alg2:lab3}
p^{m+1}_i +k \frac{(\nabla\phi^{m+1}_h, \nabla\varphi_{\boldsymbol{a}_i})}{\|\varphi_{\boldsymbol{a}_i}\|_{L^1(\Omega)}} \max\{p^{m+1}_i, \varepsilon\}\ge \min_{j\in I(\Omega_{\a_i})} p^{m}_j
\end{equation}
Next on selecting $\bar\phi_h=\varphi_{\boldsymbol{a}_i}$ in \eqref{eq_alg2:phih}, on  truncation, and substituting it into \eqref{DMP_alg2:lab3}, there holds
$$
p^{m+1}_i +k \frac{([p^{m+1}_h]_T-[n^{m+1}_h]_T, \varphi_{\boldsymbol{a}_i})_h}{\|\varphi_{\boldsymbol{a}_i}\|_{L^1(\Omega)}} \max\{p^{m+1}_i, \varepsilon\}\ge \min_{j\in I(\Omega_{\a_i})} p^{m}_j
$$
and therefore 
$$
p^{m+1}_i \ge p^m_i \ge \min\{\underset {\x\in\Omega}{\min }\,p^0_h(\x),\underset {\x\in\Omega}{\min }\,n^0_h(\x)\},
$$
since $([p^{m+1}_h]_T-[n^{m+1}_h]_T, \varphi_{\boldsymbol{a}_i})_h<0$. This is then obviously a contradiction from \eqref{DMP_alg2:cont_min} demonstrating the discrete minimum principle for $p^{m+1}_h$ in \eqref{DMP_alg2:ph}. Likewise, for the discrete maximum principle, one finds 
$$
p^{m+1}_i +k \frac{([p^{m+1}_h]_T-[n^{m+1}_h]_T, \varphi_{\boldsymbol{a}_i})_h}{\|\varphi_{\boldsymbol{a}_i}\|_{L^1(\Omega)}} \max\{p^{m+1}_i, \varepsilon\}\le p^{m}_i,
$$
which implies 
$$
p^{m+1}_i \le p^m_i \le \max\{\underset {\x\in\Omega}{\max }\,p^0_h(\x),\underset {\x\in\Omega}{\max }\,n^0_h(\x)\}.
$$

It is not hard to see that, of course, the discrete principles for $n^{m+1}_h$ in \eqref{DMP_alg2:nh} are obtained with the same procedure. 
\begin{remark}
Upon choice of $\varepsilon$ such that $0<\varepsilon< \min\{\underset {\mathbf{x}\in\Omega}{\min },p^0_h(\mathbf{x}),\underset {\mathbf{x}\in\Omega}{\min },n^0_h(\mathbf{x})\}$, one can substitute $g'_\varepsilon(s)$ with $\log(s)$. This replacement is viable, since $g'_\varepsilon(s)=\frac{s}{\varepsilon}+\log \varepsilon-1$ is not utilized in the subsequent steps.
\end{remark}

$\bullet$ \textbf{Mass convervation} Let us assume that \eqref{MC_alg2:ph} and \eqref{MC_alg2:nh} holds for $s=m$. Picking $\bar p_h=1$ and $\bar n_h=1$ as test functions in \eqref{eq_alg2:ph} and \eqref{eq_alg2:nh}, respectively, leads, on invoking \eqref{B_alg2:mass-conservation}, to 
$$
\int_\Omega p^{m+1}_h(\x){\rm d}\x=\int_\Omega p^m_h(\x){\rm d}\x
$$
and
$$
\int_\Omega n^{m+1}_h(\x){\rm d}\x=\int_\Omega n^m_h(\x){\rm d}\x.
$$
Positivity for $p^{m+1}_h$ and  $n^{m+1}_h$ obtained in \eqref{DMP_alg2:ph} and \eqref{DMP_alg2:nh} concludes the proof for $s=m+1$.     

$\bullet$ \textbf{Discrete entropy law}
Letting $p_h= i_h g_0'(p^{m+1}_h)+\phi^{m+1}_h$ in \eqref{eq_alg2:ph} and $n_h= i_h g'_0(n^{m+1}_h)-\phi^{m+1}_h$ in \eqref{eq_alg2:nh}, there follows that   
\begin{equation}\label{DEtL_alg2:lab1}
\begin{array}{rcl}
(\delta_t p^{m+1}_h, i_h g_0'(p^{m+1}_h)+\phi^{m+1}_h)_h&+&(\nabla p^{m+1}_h, \nabla(i_h g_0'(p^{m+1}_h)+\phi^{m+1}_h))
\\
&&-(p^{m+1}_h\nabla \phi^{m+1}_h, \nabla (i_h g_0'(p^{m+1}_h)+\phi^{m+1}_h))_*
\\
&&+(\mathcal{B}^+_2(p^{m+1}_h, \phi^{m+1}_h) p^{m+1}_h,  i_h g'_0(p^{m+1}_h)+\phi^{m+1}_h)=0
\end{array}
\end{equation}
and
\begin{equation}\label{DEtL_alg2:lab2}
\begin{array}{rcl}
(\delta_t n^{m+1}_h, i_h g'_0(n^{m+1}_h)-\phi^{m+1}_h)_h&-&(\nabla n^{m+1}_h, \nabla(i_h g'_0(n^{m+1}_h)-\phi^{m+1}_h))
\\
&&-(n^{m+1}_h\nabla \phi^{m+1}_h, \nabla (i_h g'_0(n^{m+1}_h)-\phi^{m+1}_h))_*
\\
&&+(\mathcal{B}^-_2(n^{m+1}_h, \phi^{m+1}_h) n^{m+1}_h,  i_h g'_0(n^{m+1}_h)-\phi^{m+1}_h)=0.
\end{array}
\end{equation}
Let us first start with \eqref{DEtL_alg2:lab1}. Recall $\mathds{I}^*(\rho_h)=\{(i,j)\in I\times I : \rho_j\not = \rho_i \}$ and $\mathds{I}^*_c(\rho_h)$ its complementary. Thus 
$$
\begin{array}{rcl}
(\nabla p^{m+1}_h, \nabla(i_h g_0'(p^{m+1}_h)+\phi^{m+1}_h))&=&\displaystyle
-\sum_{i<j\in \mathds{I}^*(p^{m+1}_h)}\frac{\delta_{ji}g'_0(p^{m+1}_h)}{\delta_{ji}p^{m+1}_h} (\delta_{ji}p^{m+1}_h)^2(\nabla \varphi_{\a_j},\nabla\varphi_{\a_i})
\\
&&\displaystyle
-\sum_{i<j\in \mathds{I}^*(p^{n+1}_h)}\delta_{ji}p^{m+1}_h \delta_{ji}\phi^{m+1}_h(\nabla\varphi_{\a_j},\nabla\varphi_{\a_i})
\end{array}
$$
and
$$
\begin{array}{rcl}
(p^{m+1}_h\nabla \phi^{m+1}_h, \nabla (i_h g'_0(p^{m+1}_h)+\phi^{m+1}_h))_*&=&\displaystyle
-\sum_{i<j\in \mathds{I}^*(p^{n+1}_h)} \delta_{ji}p^{n+1}_h \delta_{ji}\phi^{n+1}_h(\nabla\varphi_{\a_j},\nabla\varphi_{\a_i})
\\
&&\displaystyle
-\sum_{i<j\in \mathds{I}^*(p^{m+1}_h)}\frac{\delta_{ji}p^{m+1}_h}{\delta_{ji}g'_0(p^{m+1}_h)}(\delta_{ji}\phi^{m+1}_h)^2(\nabla\varphi_{\a_j},\nabla\varphi_{\a_i})
\\
&&\displaystyle
-\sum_{i<j\in \mathds{I}^*_c(p^{m+1}_h)}p^{m+1}_i(\delta_{ji}\phi^{m+1}_h)^2(\nabla\varphi_{\a_j},\nabla\varphi_{\a_i}),
\end{array}
$$
whereupon 
\begin{align*}
(\nabla p^{m+1}_h,  \nabla(i_h g'_0(p^{m+1}_h)&+\phi^{m+1}_h))+(p^{m+1}_h\nabla \phi^{m+1}_h, \nabla (i_h g_0'(p^{m+1}_h)+\phi^{m+1}_h))_*
\\
=\displaystyle
-\sum_{i<j\in \mathds{I}^*(p^{m+1}_h)}&\left|\left(\frac{\delta_{ji}g'_0(p^{m+1}_h)}{\delta_{ji}p^{m+1}_h}\right)^{\frac{1}{2}}\delta_{ji}p^{m+1}_h
+\left(\frac{\delta_{ji}p^{m+1}_h}{\delta_{ji}g'_0(p^{m+1}_h)}\right)^{\frac{1}{2}}\delta_{ji}\phi^{m+1}_h\right|^2(\nabla\varphi_{\a_j},\nabla\varphi_{\a_i})
\\
\displaystyle
-\sum_{i<j\in \mathds{I}^*_c(p^{n+1}_h)}&p^{m+1}_i \delta^2_{ij}\phi^{m+1}_h(\nabla\varphi_{\a_j},\nabla\varphi_{\a_i})>0,
\end{align*}
on account of \eqref{Alg_2:acuteness}.  Now the stabilizing term in \eqref{DEtL_alg2:lab1} is written out as 
\begin{align}
(B_2(p^{m+1}_h, \phi^{m+1}_h) p^{m+1}_h, &i_h g'_0(p^{m+1}_h)-\phi^{m+1}_h)
\nonumber
\\
=\displaystyle
\sum_{i<j\in I^*}&\tilde\beta^+_{ji}(p^{m+1}_h, \phi^{m+1}_h) \delta_{ji} p^{m+1}_h ( \delta_{ji} g'_0(p^{m+1}_h)+\delta_{ji}\phi^{m+1}_h).
\nonumber
\end{align}
On the one hand, if one assumes that there holds 
$$
\tilde\beta^+_{ji}(p^{m+1}_h, \phi^{m+1}_h)=\alpha_i(p^{m+1}_h)\left(1+\delta_{ji}\phi^{m+1}_h\left[\frac{1}{\delta_{ji} g'_0(p^{m+1}_h)} - \frac{p^{m+1}_i}{\delta_{ji} p^{m+1}_h}\right]\right)(\nabla\varphi_{\a_j},\nabla\varphi_{\a_i}),
$$
we deduce that $\delta_{ji}\phi^{m+1}_h\le 0$, since $\tilde\beta^+_{ji}(p^{m+1}_h, \phi^{m+1}_h)>0$ and $(\nabla\varphi_{\a_j},\nabla\varphi_{\a_i})<0$. 
Next observe (with the help of the mean value theorem) that 
$$
\frac{1}{\delta_{ji} g'_0(p^{m+1}_h)} - \frac{p^{m+1}_i}{\delta_{ji} p^{m+1}_h}=\frac{\theta_{ji}p^{m+1}_j+(1-\theta_{ji}) p_i^{m+1}}{\delta_{ji} p^{m+1}_h}-\frac{p^{m+1}_i}{\delta_{ji} p^{m+1}_h}=\theta_{ji}.
$$
with $\theta_{ji}\in (0,1)$. We now distinguish among two cases:

$\bullet $ \textbf{Case 1}: $\delta_{ji} p^{m+1}_h<0$.  It is clear that 
$$
\tilde\beta^+_{ji}(p^{m+1}_h, \phi^{m+1}_h) \delta_{ji} p^{m+1}_h ( \delta_{ji} g'_0(p^{m+1}_h)+\delta_{ji}\phi^{m+1}_h)>0,
$$ 
since $\delta_{ji} p^{m+1}_h<0$ implies $\delta_{ji} g'_0(p^{m+1}_h)<0$.

$\bullet $ \textbf{Case 2}: $\delta_{ji} p^{m+1}_h>0$. Write out 
\begin{align*}
\left(1+\delta_{ji}\phi^{m+1}_h\left[\frac{1}{\delta_{ji} g'_0(p^{m+1}_h)} - \frac{p^{m+1}_i}{\delta_{ji} p^{m+1}_h}\right]\right)(\delta_{ji}p^{m+1}_h)(\delta_{ji}g'_0(p^{m+1}_h)+\delta_{ij}\phi^{m+1}_h) (\nabla\varphi_{\a_j},\nabla\varphi_{\a_i})
\\
= \left|\left(\frac{\delta_{ji}g'_0(p^{m+1}_h)}{\delta_{ji}p^{m+1}_h}\right)^{\frac{1}{2}}\delta_{ji}p^{m+1}_h
+\left(\frac{\delta_{ji}p^{m+1}_h}{\delta_{ji}g'_0(p^{m+1}_h)}\right)^{\frac{1}{2}}\delta_{ji}\phi^{m+1}_h\right|^2 (\nabla\varphi_{\a_j},\nabla\varphi_{\a_i})
\\
- p^{m+1}_i(\delta_{ij}\phi^{m+1}_h)^2 (\nabla\varphi_{\a_j},\nabla\varphi_{\a_i})
\\
- p^{m+1}_i \delta_{ji}g'_0(p^{m+1}_h)\delta_{ji}\phi^{m+1}_h (\nabla\varphi_{\a_j},\nabla\varphi_{\a_i}).
\end{align*}
As  $\delta_{ji} p^{m+1}_h>0$ and $1<-\theta_{ji}\delta_{ji}\phi^{m+1}_h $, we bound
$$
0\le \delta_{ji}g'_0(p^{m+1}_h)\le \frac{\delta_{ji}p^{m+1}_h}{\theta_{ji}\delta_{ji}p^{m+1}_h+ p^{m+1}_i}\le \frac{1}{\theta_{ji}}\le - \delta_{ij}\phi^{m+1}_h,
$$
which gives
$$
- p^{m+1}_i \delta_{ji}g'_0(p^{m+1}_h)\delta_{ij}\phi^{m+1}_h (\nabla\varphi_{\a_j},\nabla\varphi_{\a_i})\ge p^{m+1}_i (\delta_{ij}\phi^{m+1}_h)^2 (\nabla\varphi_{\a_j},\nabla\varphi_{\a_i}).
$$
Therefore, 
\begin{align*}
\tilde\beta^+_{ji}(p^{m+1}_h, \phi^{m+1}_h) \delta_{ji} p^{m+1}_h (\delta_{ji} g'_0(p^{m+1}_h)+\delta_{ji}\phi^{m+1}_h)
\\
\ge \alpha_{i}(p^{m+1}_h)\left|\left(\frac{\delta_{ji}g'_0(p^{m+1}_h)}{\delta_{ji}p^{m+1}_h}\right)^{\frac{1}{2}}\delta_{ji}p^{m+1}_h
+\left(\frac{\delta_{ji}p^{m+1}_h}{\delta_{ji}g'_0(p^{m+1}_h)}\right)^{\frac{1}{2}}\delta_{ji}\phi^{m+1}_h\right|^2 (\nabla\varphi_{\a_j},\nabla\varphi_{\a_i}). 
\end{align*} 

Alternatively, if we have
$$
\tilde\beta^+_{ji}(p^{m+1}_h, \phi^{m+1}_h)=\alpha_j(p^{m+1}_h)\left(1+\delta_{ij}\phi^{m+1}_h\left[\frac{1}{\delta_{ij} g'_0(p^{m+1}_h)} - \frac{p^{m+1}_j}{\delta_{ij} p^{m+1}_h}\right]\right)(\nabla\varphi_{\a_j},\nabla\varphi_{\a_i}),
$$
the same result holds. Consequently, we are led to 
\begin{equation}\label{DEtL_alg2:lab3}
\begin{array}{c}
(\nabla p^{m+1}_h, \nabla(i_h g'_0(p^{m+1}_h)+\phi^{m+1}_h))+(p^{m+1}_h\nabla \phi^{m+1}_h, \nabla (i_h g'_0(p^{m+1}_h)+\phi^{m+1}_h))_*
\\
+(\mathcal{B}_2(p^{m+1}_h, n^{m+1}_h, \phi^{m+1}_h) p^{m+1}_h, i_h g'_0(p^{m+1}_h)+\phi^{m+1}_h)
\\
\ge\displaystyle
-\sum_{i<j\in \mathds{I}^*(p^{m+1}_h)}\hspace{-0.5cm}(1-\alpha_\#(p^{m+1}_h))\left|\left(\frac{\delta_{ji} g'_0(p^{m+1}_h)}{\delta_{ji}p^{m+1}_h}\right)^{-\frac{1}{2}}\delta_{ji} p^{m+1}_h\right.
\\
\displaystyle
\hspace{9cm}
\left.+\left(\frac{\delta_{ji} p^{m+1}_h}{\delta_{ji}g'_0(p^{m+1}_h)}\right)^{\frac{1}{2}}\delta_{ji}\phi^{m+1}_h\right|^2(\nabla\varphi_{\a_j},\nabla\varphi_{\a_i})
\\
\displaystyle
-\sum_{i<j\in \mathds{I}^*_c(p^{m+1}_h)}p^{m+1}_i \delta^2_{ij}\phi^{m+1}_h(\nabla\varphi_{\a_j},\nabla\varphi_{\a_i}).
\end{array}
\end{equation}

Taylor's theorem of $g_0$ about $p^{m+1}_h$ with the Lagrange form of the remainder yields
\begin{equation}\label{DEtL_alg2:lab4}
g_0(p^m_h)=g_0(p^{m+1}_h)-g'_0(p^{m+1}_h)(p^{m+1}_h-p^m_h)+\frac{g''_\varepsilon(p^{m+\theta_p}_h)}{2}(p^{m+1}_h-p^m_h)^2,
\end{equation}
where $\theta_p\in (0,1)$ is such that $p^{m+\theta_p}_h=\theta_p p^{m+1}_h+(1-\theta_p) p^m_h$. Using \eqref{DEtL_alg2:lab4}, we show  \cite{GS_RG_2021}  that  
$$
(\delta_t p^{m+1}_h, g'_0(p^{m+1}_h))_h=\frac{1}{k}(g_0(p^{m+1}_h),1)_h- \frac{1}{k}(g_0(p^m_h),1)_h+ \frac{k}{2} (g''_0(p^{n+\theta_p}_h),(\delta_t p^{m+1}_h)^2)_h.
$$

On combining \eqref{DEtL_alg2:lab3} and \eqref{DEtL_alg2:lab4} with \eqref{DEtL_alg2:lab1} leads to
\begin{equation}\label{DEtL_alg2:lab5}
\begin{array}{c}
\displaystyle
\frac{1}{k}(g_0(p^{m+1}_h),1)_h- \frac{1}{k}(g_0(p^m_h),1)_h
+\frac{k}{2} (g''_0(p^{m+\theta_p}_h),(\delta_t p^{m+1}_h)^2)_h+(\delta_t p^{m+1}_h, \phi^{m+1}_h)_h
\\
\displaystyle
-\sum_{i<j\in \mathds{I}^*(p^{m+1}_h)}(1-\alpha_\#(p^{m+1}_h))\left|\left(\frac{\delta_{ji} g'_0(p^{m+1}_h)}{\delta_{ji}p^{m+1}_h}\right)^{-\frac{1}{2}}\delta_{ji} p^{m+1}_h\right.
\\
\displaystyle
\hspace{9cm}\left.+\left(\frac{\delta_{ji}g'_0(p^{m+1}_h)}{\delta_{ji} p^{m+1}_h}\right)^{\frac{1}{2}}\delta_{ji}\phi^{m+1}_h\right|^2(\nabla\varphi_{\a_j},\nabla\varphi_{\a_i})
\\
\displaystyle
-\sum_{i<j\in \mathds{I}^*_c(p^{m+1}_h)}p^{m+1}_i \delta^2_{ij}\phi^{m+1}_h(\nabla\varphi_{\a_j},\nabla\varphi_{\a_i}) 
\le 0,
\end{array}
\end{equation}
where $\alpha_\#\in\{\alpha_i, \alpha_j, 0\}$ is  such that 
$$\tilde\beta^+_{ji}(p^{m+1}_h,\phi^{m+1}_h)=\alpha_\#(p^{m+1}_h)\left[1+\delta_{\bar\#\#}\phi^{m+1}_h\left(\frac{1}{\delta_{\bar\#\#}g'_0(p^{m+1}_h)}-\frac{p^{m+1}_\#}{\delta_{\bar\#\#} p^{m+1}_h)}\right)\right](\nabla\varphi_{\a_{\bar \#}},\nabla\varphi_{\a_\#}),$$
with ${\bar \#}=j$ when $\#=i$ and ${\bar \#}=i$ when $\#=j$.

Similarly for  \eqref{DEtL_alg2:lab5}, one finds
\begin{equation}\label{DEtL_alg2:lab6}
\begin{array}{c}
\displaystyle
\frac{1}{k}(g_0(n^{m+1}_h),1)_h- \frac{1}{k}(g_0(n^m_h),1)_h
+\frac{k}{2} (g''_0(n^{m+\theta_n}_h),(\delta_t n^{m+1}_h)^2)_h-(\delta_t n^{m+1}_h, \phi^{m+1}_h)_h
\\
\displaystyle
-\sum_{i<j\in \mathds{I}^*(n^{m+1}_h)}(1-\alpha_\#(n^{m+1}_h))\left|\left(\frac{\delta_{ji} g'_0(n^{m+1}_h)}{\delta_{ji}n^{m+1}_h}\right)^{-\frac{1}{2}}\delta_{ji} n^{m+1}_h\right.
\\
\displaystyle
\hspace{9cm}-\left.\left(\frac{\delta_{ji}g'_0(n^{m+1}_h)}{\delta_{ji} n^{m+1}_h}\right)^{\frac{1}{2}}\delta_{ji}\phi^{m+1}_h\right|^2(\nabla\varphi_{\a_j},\nabla\varphi_{\a_i})
\\
\displaystyle
-\sum_{i<j\in \mathds{I}^*_c(n^{m+1}_h)}n^{m+1}_i \delta^2_{ij}\phi^{m+1}_h(\nabla\varphi_{\a_j},\nabla\varphi_{\a_i}) \le 0,
\end{array}
\end{equation}
where $\alpha_\#\in\{\alpha_i, \alpha_j, 0\}$ is  such that 
$$\tilde\beta^+_{ji}(n^{m+1}_h,\phi^{m+1}_h)=\alpha_\#(n^{m+1}_h)\left[1+\delta_{\bar\#\#}\phi^{m+1}_h\left(\frac{1}{\delta_{\bar\#\#}g'_0(n^{m+1}_h)}-\frac{n^{m+1}_\#}{\delta_{\bar\#\#} n^{m+1}_h)}\right)\right](\nabla\varphi_{\a_{\bar \#}},\nabla\varphi_{\a_\#}),$$
with ${\bar \#}=j$ when $\#=i$ and ${\bar \#}=i$ when $\#=j$.

Addition of  \eqref{DEtL_alg2:lab5} and \eqref{DEtL_alg2:lab6} completes the proof of \eqref{Entropy:alg2}; thereby establishing Theorem 2.

\section{Numerical simulations}

In this section we provide numerical evidence of the good performance of our numerical solution for both algorithms in three different scenarios. To begin with, we use a smooth initial data with critical peaks. As a second example, we deal with three experiments being a schematic representation of an ion channel, namely an initial uniform distribution and two initial distributions concentrated close to the wall for both ions, respectively. The latter leads to two opposite travelling waves.    

As both algorithms are nonlinear, we use Picard iteration as a linealization. More particularly,  
\begin{align}
(\delta_t p^{m+1, i+1}, \bar p_h)_\#+(\nabla p^{m+1,i+1}_h,\nabla\bar p_h)+(p^{m+1,i+\$}_h\nabla\phi^{m+1,i}_h,\nabla\bar p_h)_\#
\nonumber
\\
+(\mathcal{B}_\#(p^{m+1}_h, n^{m+1}_h,\phi^{m+1}_h)p^{m+1}_m,\bar p_h )&=0,
\nonumber
\\
(\delta_t n^{m+1, i+1}_h, \bar n_h)_\#+(\nabla n^{m+1}_h,\nabla \bar n_h)-(n^{m+1,i+\$}_h\nabla\phi^{m+1,i}_h,\nabla \bar n_h)_\#
\nonumber
\\
+(\mathcal{B}_\#(p^{m+1,i}_h, n^{m+1,i}_h,\phi^{m+1,i}_h)n^{m+1,i+1}_h,\bar p_h )&=0,
\nonumber 
\\
(\nabla \phi^{m+1,i+1}_h, \nabla \bar\phi_h)+(n^{m+1,i+1}_h-p^{m+1,i+1}_h,\bar\phi_h)_h&=0,
\label{eq:phi_m+1,i+1}
\end{align}
with $\$=1$ or $0$ and $\#=1$ or $2$, respectively, where 
$$
(\delta_t p^{m+1, i+1}, \bar p_h)_1=(\delta_t p^{m+1, i+1}, \bar p_h)\quad \mbox{ and }\quad (p^{m+1,i}_h\nabla\phi^{m+1,i}_h,\nabla\bar p_h)_1=(p^{m+1,i+1}_h\nabla\phi^{m+1,i}_h,\nabla\bar p_h)
$$
and
$$
(\delta_t p^{m+1, i+1}, \bar p_h)_2=(\delta_t p^{m+1, i+1}, \bar p_h)_h \quad \mbox{ and } \quad (p^{m+1,i}_h\nabla\phi^{m+1,i}_h,\nabla\bar p_h)_2=(p^{m+1,i}_h\nabla\phi^{m+1,i}_h,\nabla\bar p_h)_*
$$
This linealization is conjugated with a  backtracking line search with stopping criteria of $10^{-6}$ for the residual in the $L^\infty(\Omega)$-norm and of $10^{-16}$ for the increment in the $L^2(\Omega)$-norm. In defining the shock detector we take $q=2$ in \eqref{def:alpha_min_max}.

\subsection{Smooth initial data} 
\begin{figure}
\begin{subfigure}[b]{0.25\textwidth}
\centering
\includegraphics[width=1.0\textwidth]{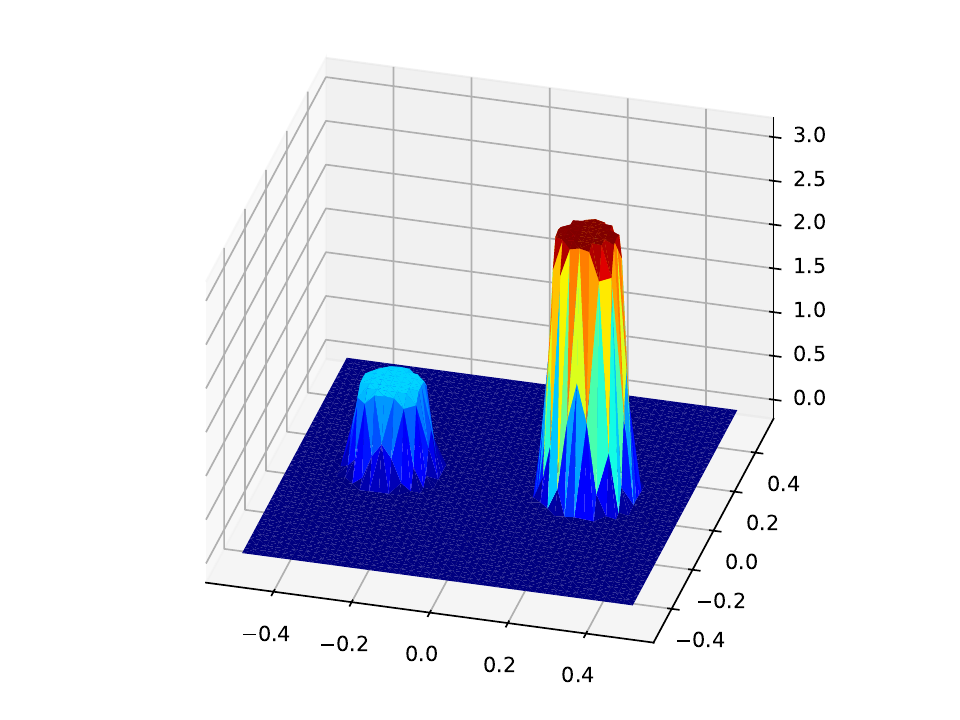}
\subcaption{$n_{0h}$}
\end{subfigure}
\begin{subfigure}[b]{0.25\textwidth}
\centering
\includegraphics[width=1.0\textwidth]{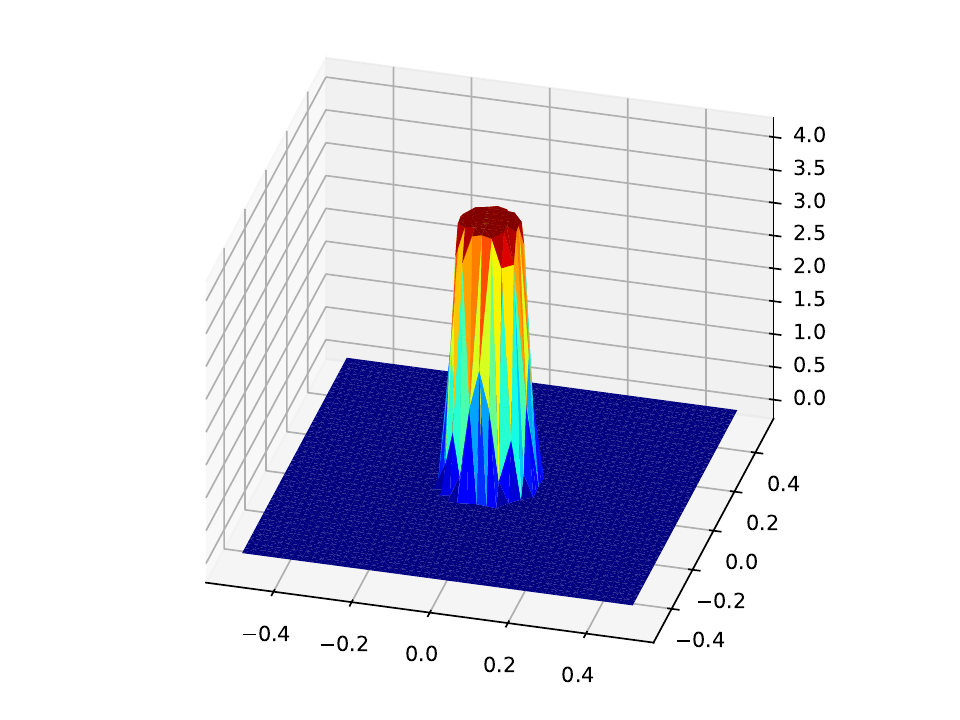}
\subcaption{$p_{0h}$}
\end{subfigure}
\begin{subfigure}[b]{0.25\textwidth}
\centering
\includegraphics[width=1.0\textwidth]{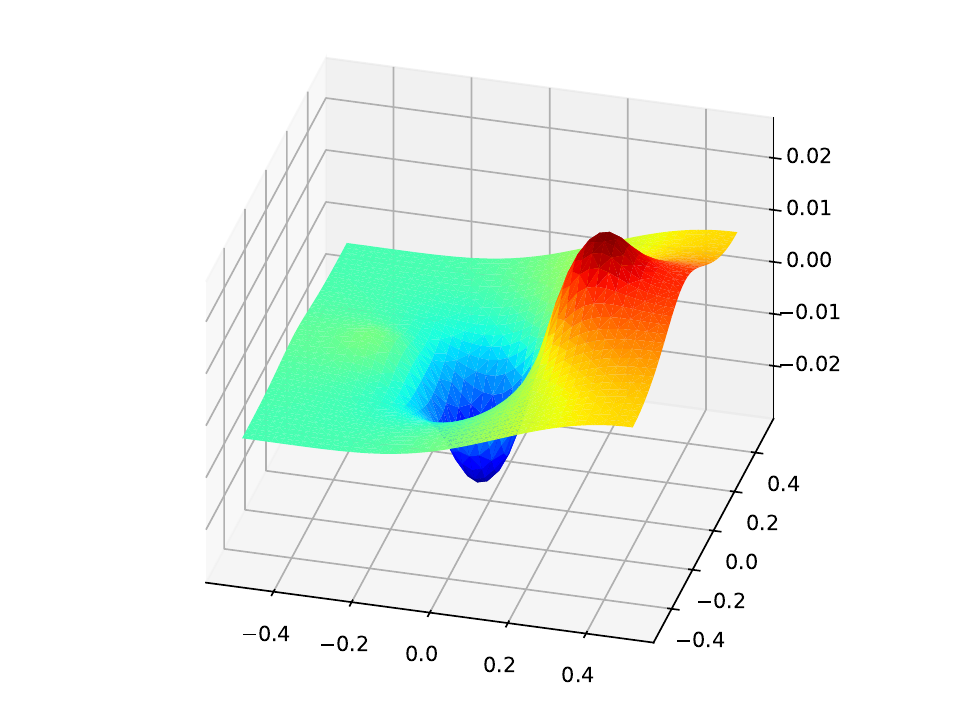}
\subcaption{$\phi_{0h}$}
\end{subfigure}
\caption{Initial data}
\label{smooth:initial_data}
\end{figure}

In this first example, we want to approximate the solution of \eqref{PNP}-\eqref{IC} on $\Omega=(-\frac{1}{2},\frac{1}{2})\times(-\frac{1}{2},\frac{1}{2})$ with  
$$
p_0(x,y)=\frac{1}{2}\tanh\left(\frac{1-10\sqrt{(x+\frac{1}{4})^2+y^2}}{0.1}\right)+\frac{3}{2}\tanh\left(\frac{1-10\sqrt{(x-\frac{1}{4})^2+ y^2}}{0.1}\right)+2
$$ 
and
$$
n_0(x,y)=2\left(\tanh\left(\frac{1-10\sqrt{x^2+ y^2}}{0.1}\right)+1\right).
$$
The initial datum $p_0$ is located around the point $(0,0)$, while $n_0$ has two accumulation points at $(0, \pm \frac{1}{4})$ with different heights $1$ and $3$, respectively. Moreover, the electroneutrality condition $\|p_{0}\|_{L^1(\Omega)}=\|n_{0}\|_{L^1(\Omega)}$ holds. Figure \ref{smooth:initial_data} shows the initial configuration of $p_{0h}$, $n_{0h}$, and $\phi_{h0}$, respectively, where $\phi_{0h}$ is computed by \eqref{eq:phi_m+1,i+1}.     
The expected dynamics involves the convergence toward a steady-state solution according to the decay of the energy and entropy quantities in \eqref{Energy_law} and \eqref{Entropy_law}.  
The mesh $\mathcal{T}_h$ employed is generated by dividing each side of the square domain into $40$ subintervals, resulting in $h\approx0.035$, and the time step is $k=10^{-3}$.

\begin{figure}
    \begin{subfigure}[b]{0.25\textwidth}
        \centering
        \includegraphics[width=1.0\textwidth]{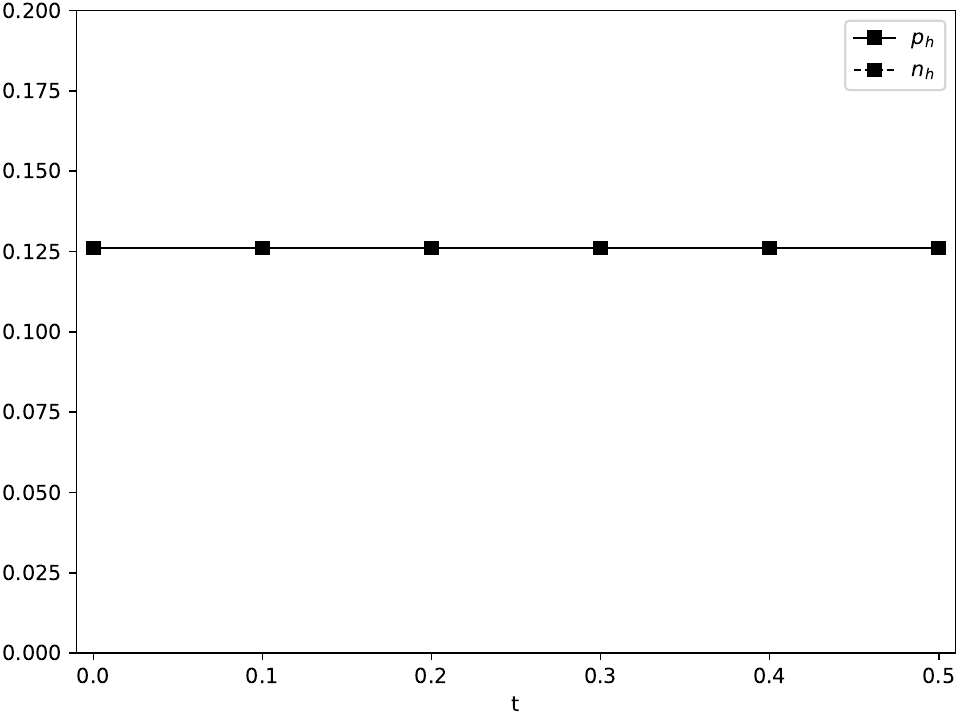}
    \end{subfigure}
    \begin{subfigure}[b]{0.25\textwidth}
        \centering
        \includegraphics[width=1.0\textwidth]{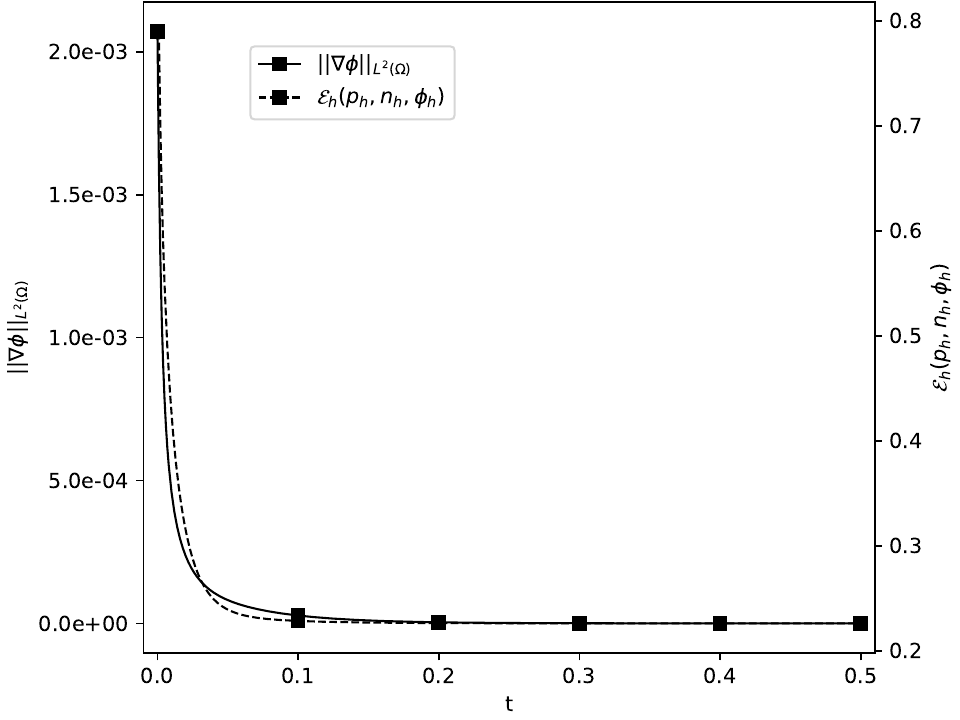}
    \end{subfigure}
    \caption*{Algorithm 1}
    \begin{subfigure}[b]{0.25\textwidth}
        \centering
        \includegraphics[width=1.0\textwidth]{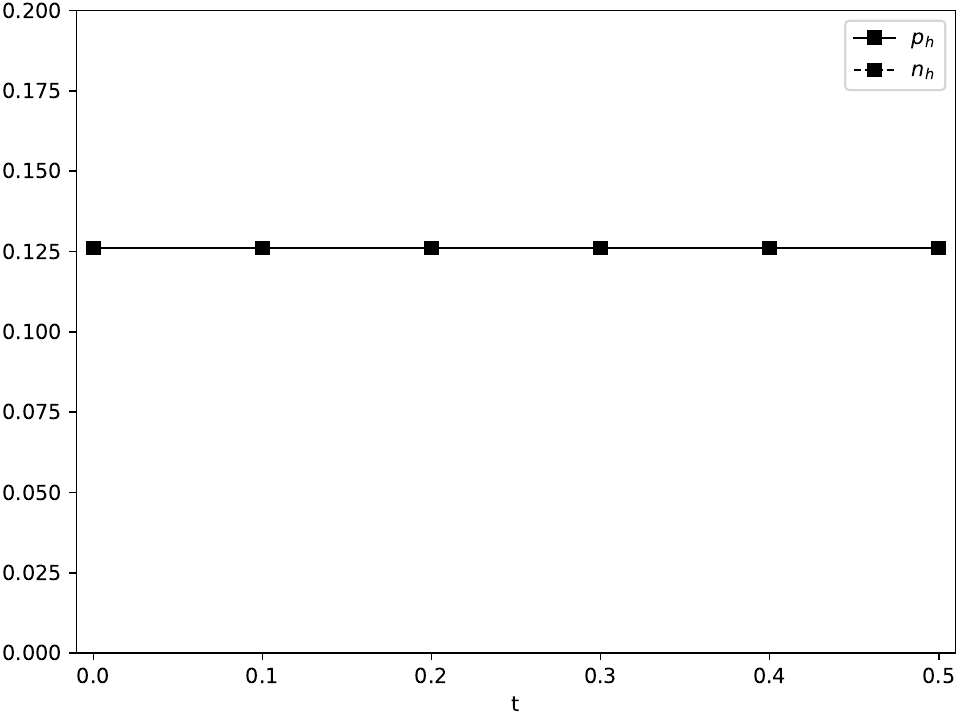}
    \end{subfigure}
    \begin{subfigure}[b]{0.25\textwidth}
        \centering
        \includegraphics[width=1.0\textwidth]{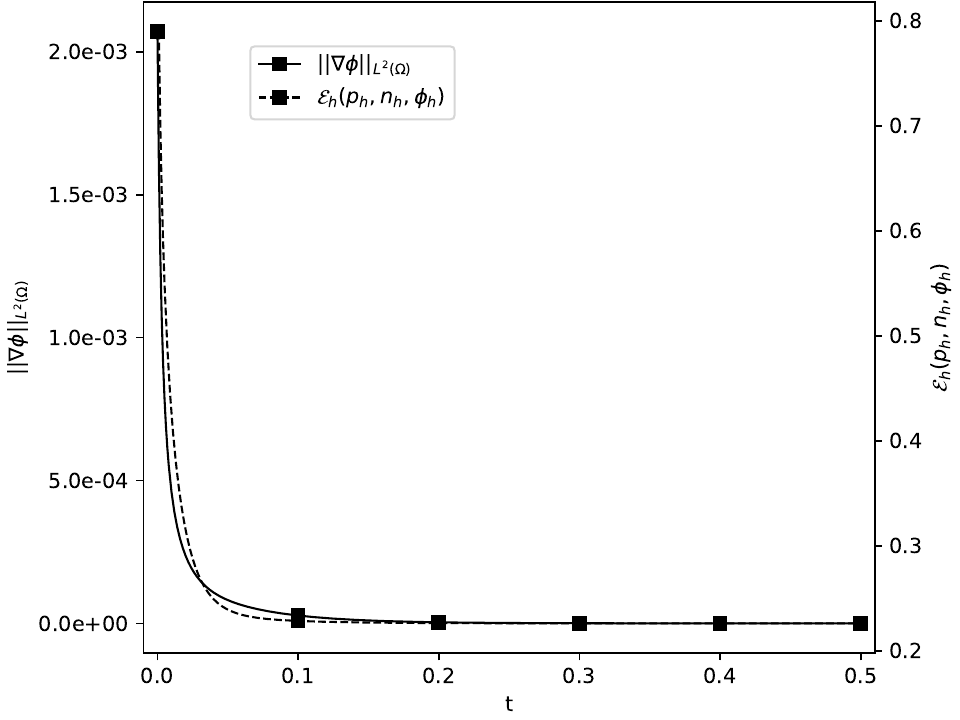}
    \end{subfigure}
    \caption*{Algorithm 2}
    \caption{Mass conservation (left), and energy and entropy evolutions (right)}
    \label{fig_smooth:mass_and_energies}
\end{figure}

\begin{figure}
    \begin{subfigure}[b]{0.25\textwidth}
        \centering
        \includegraphics[width=1.0\textwidth]{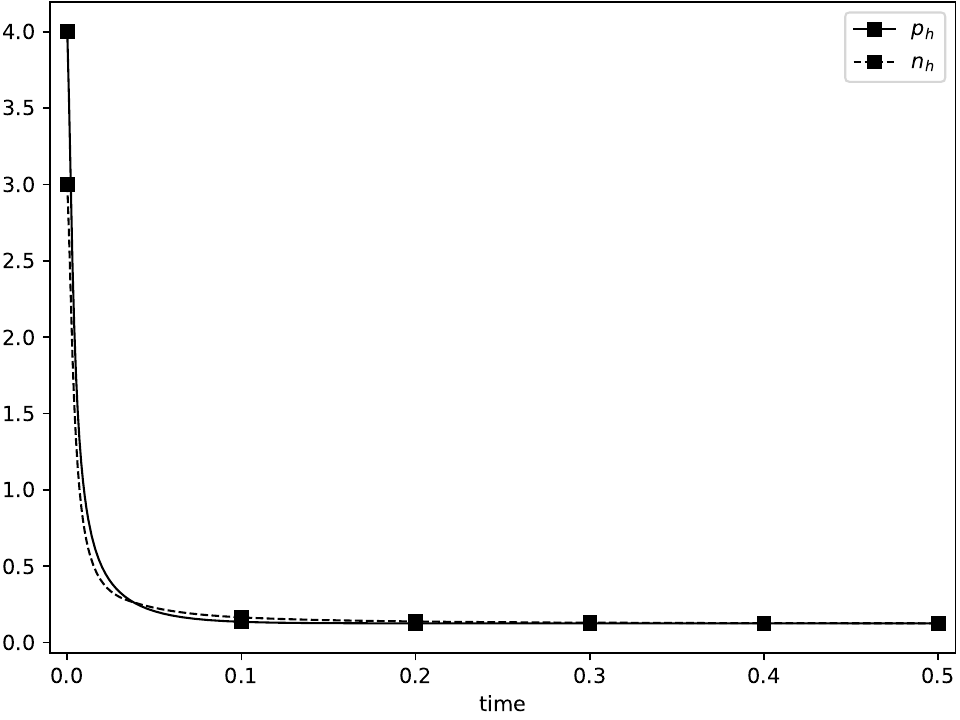}
     \end{subfigure}
    \begin{subfigure}[b]{0.25\textwidth}
        \centering
        \includegraphics[width=1.0\textwidth]{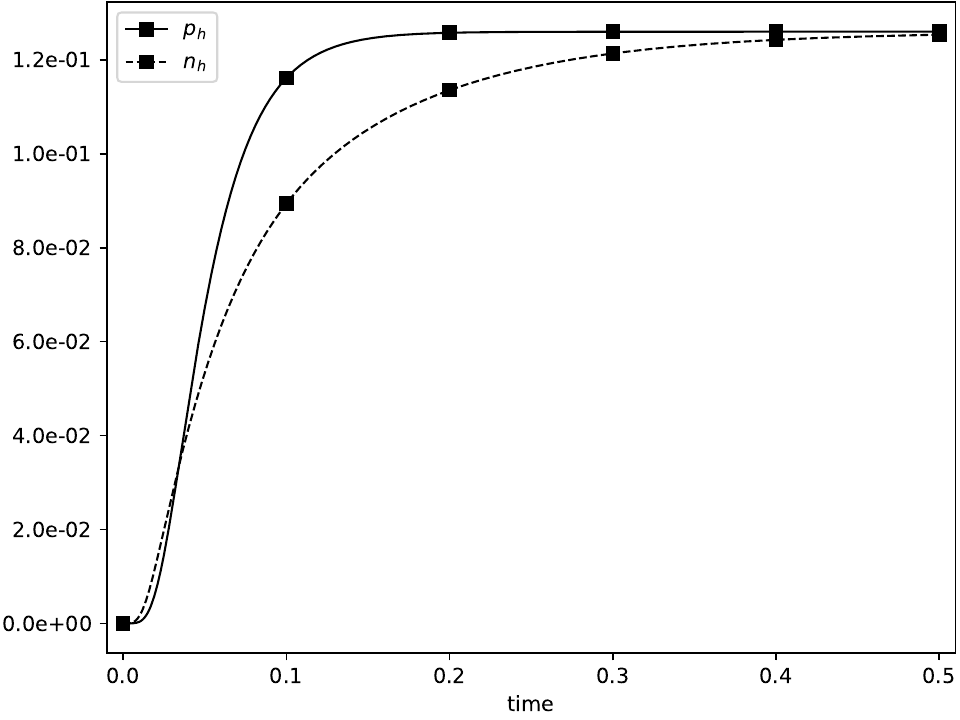}
     \end{subfigure}
  \caption*{Algorithm 1} 
       \begin{subfigure}[b]{0.25\textwidth}
        \centering
        \includegraphics[width=1.0\textwidth]{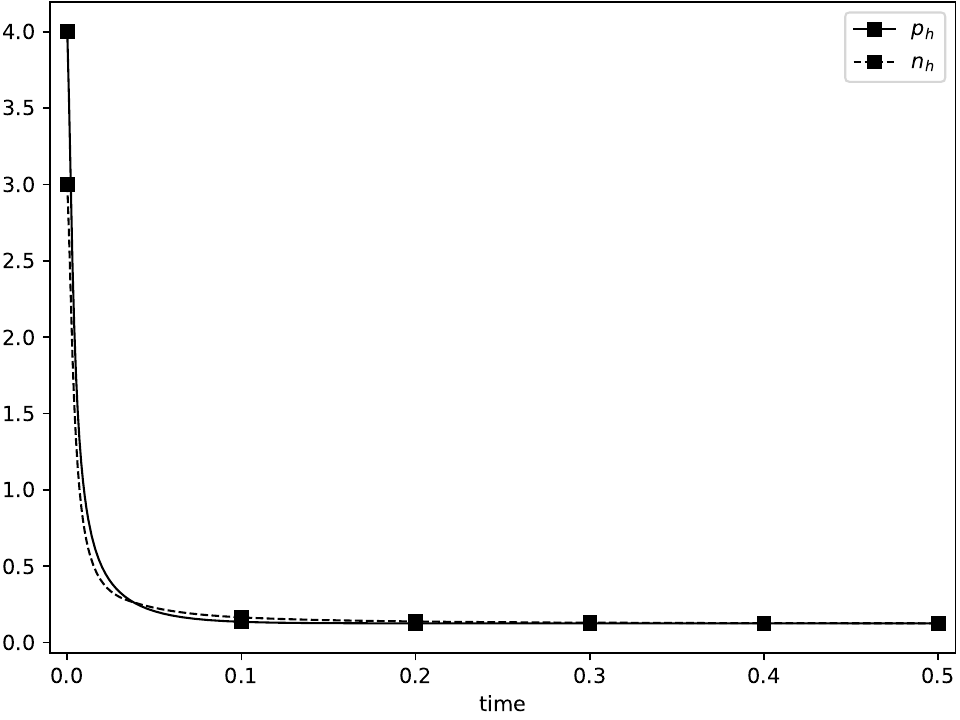}
    \end{subfigure}
    \begin{subfigure}[b]{0.25\textwidth}
        \centering
        \includegraphics[width=1.0\textwidth]{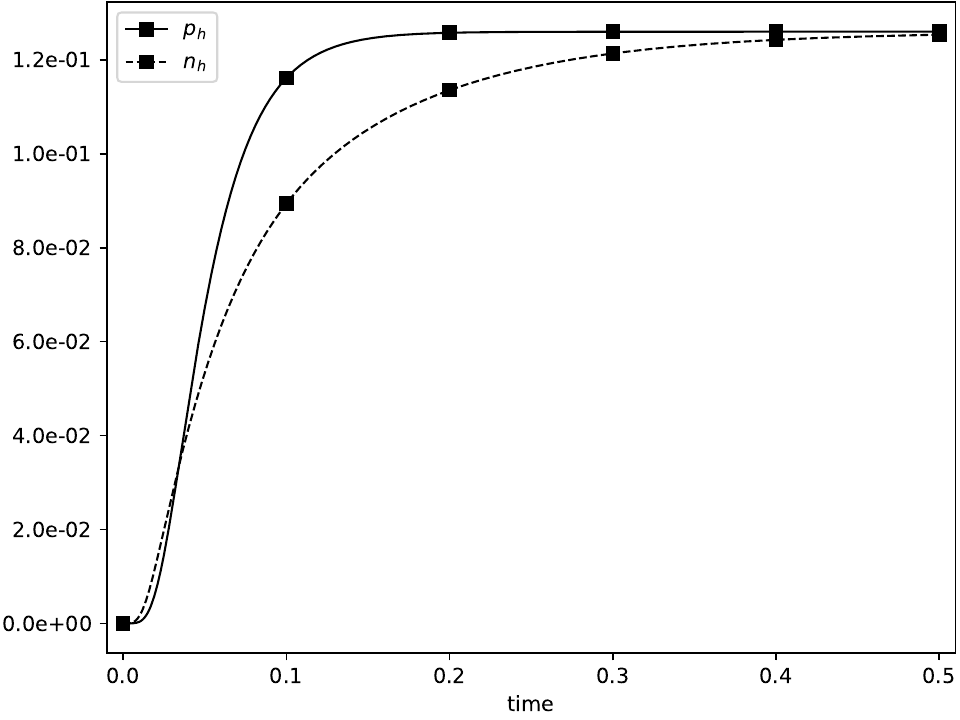}
    \end{subfigure}    
    \caption*{Algorithm 2}
    \caption{Maxima (left) and minima (right)}
    \label{fig_smooth:max_and_min}
\end{figure}

\begin{figure}    
\begin{subfigure}[b]{0.22\textwidth}
        \centering
        \includegraphics[width=1.0\textwidth]{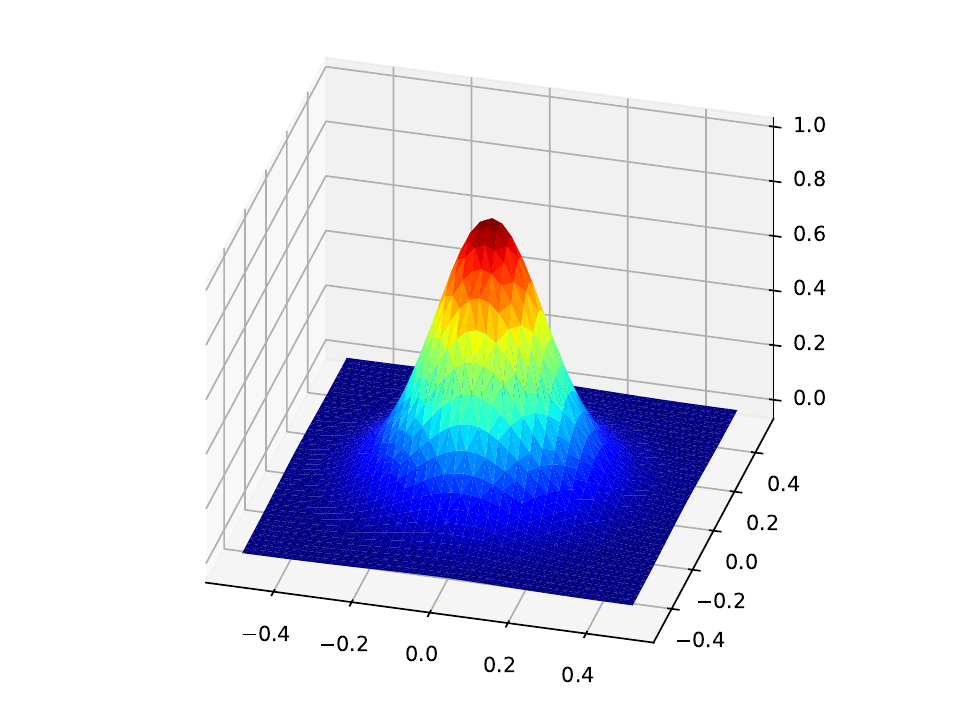}
    \end{subfigure}
    \begin{subfigure}[b]{0.22\textwidth}
        \centering
        \includegraphics[width=1.0\textwidth]{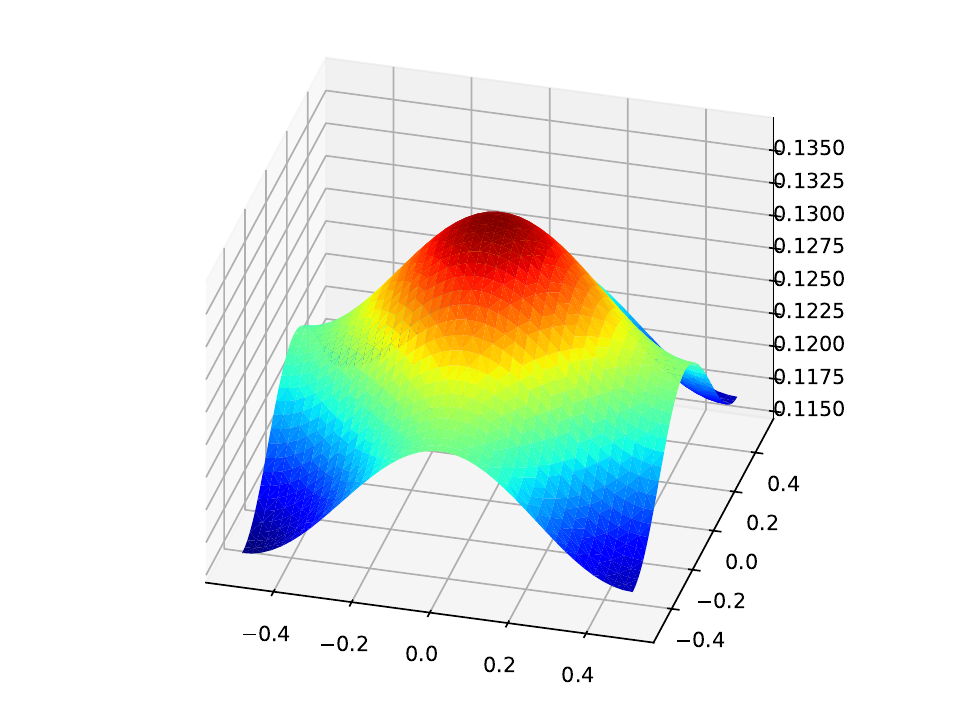}
    \end{subfigure}
    \begin{subfigure}[b]{0.22\textwidth}
        \centering
        \includegraphics[width=1.0\textwidth]{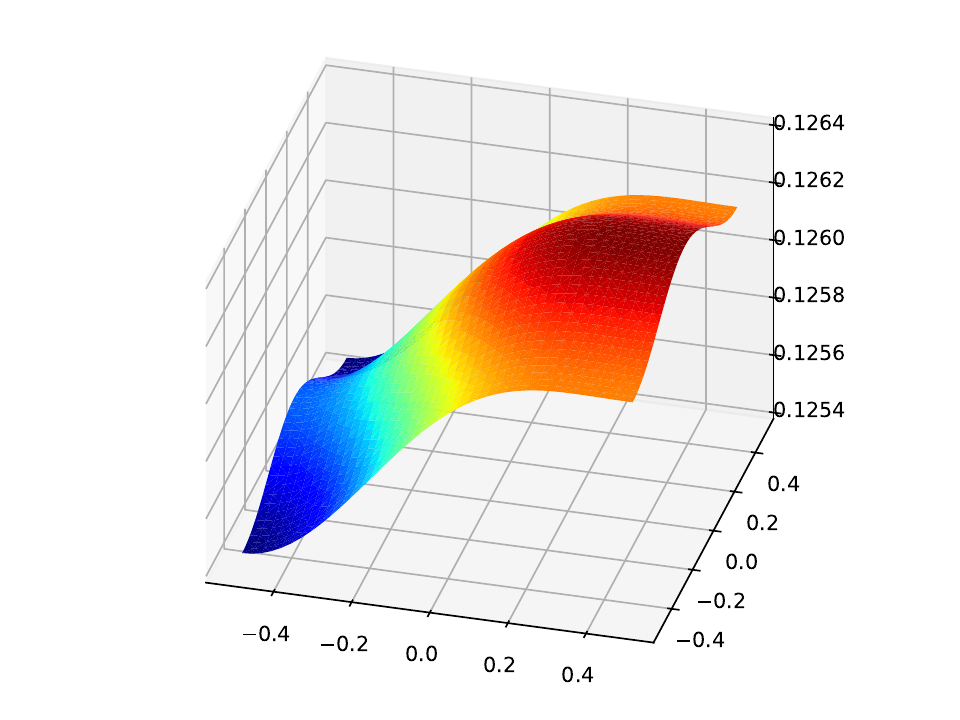}
    \end{subfigure}
    \begin{subfigure}[b]{0.22\textwidth}
        \centering
        \includegraphics[width=1.0\textwidth]{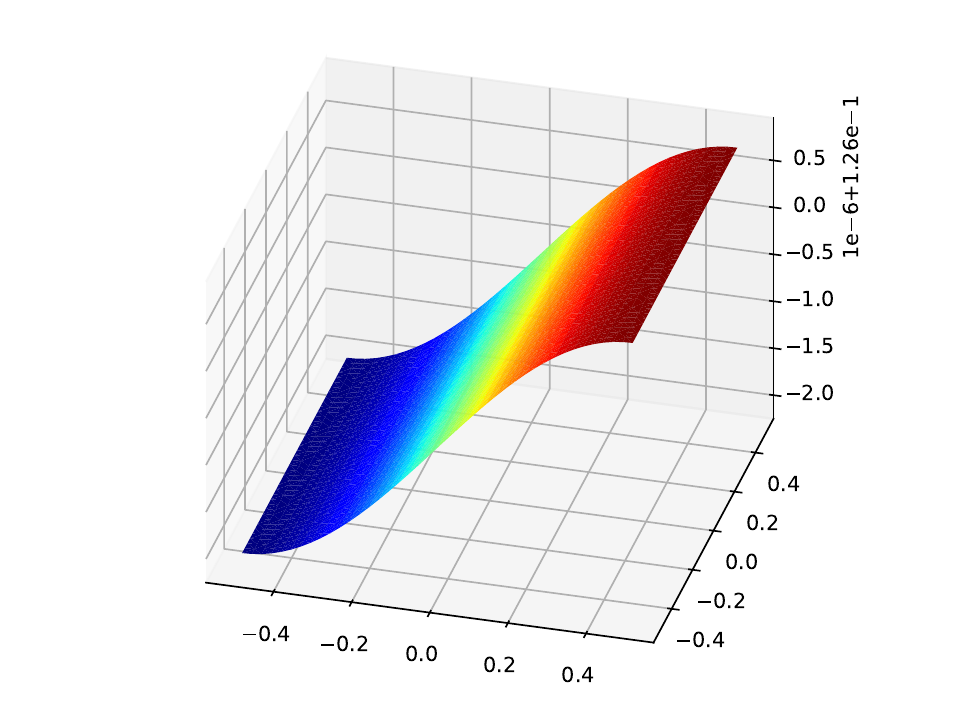}
    \end{subfigure}    
    \caption*{Snapshots of $n_h$ at times $t=0.01$, $0.1$, $0.2$ and $0.5$}
    \begin{subfigure}[b]{0.22\textwidth}
        \centering
        \includegraphics[width=1.0\textwidth]{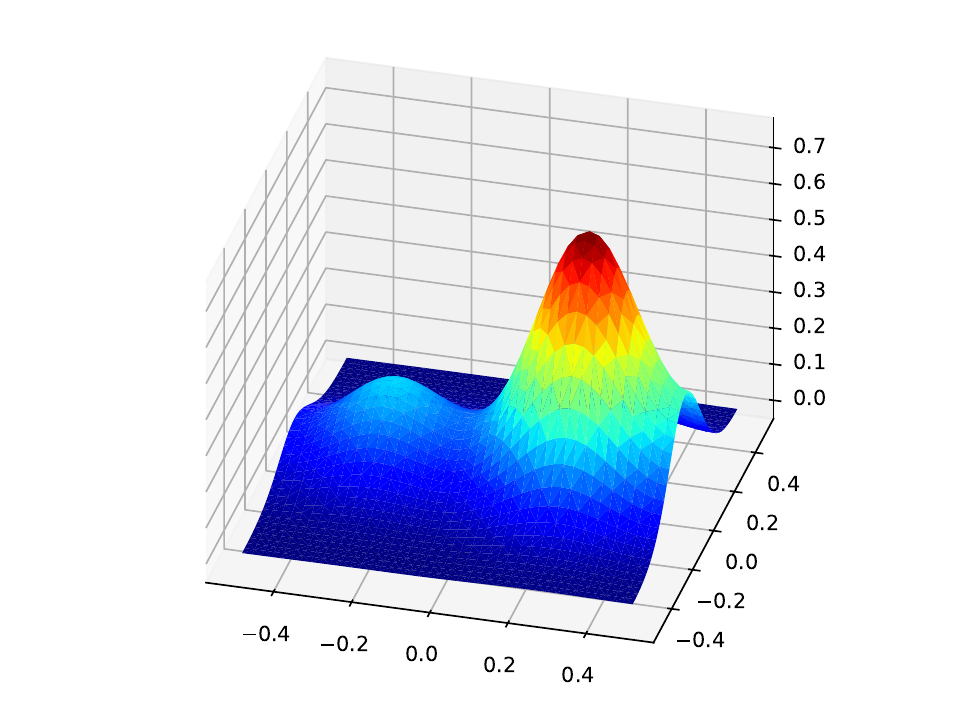}
    \end{subfigure}
    \begin{subfigure}[b]{0.22\textwidth}
        \centering
        \includegraphics[width=1.0\textwidth]{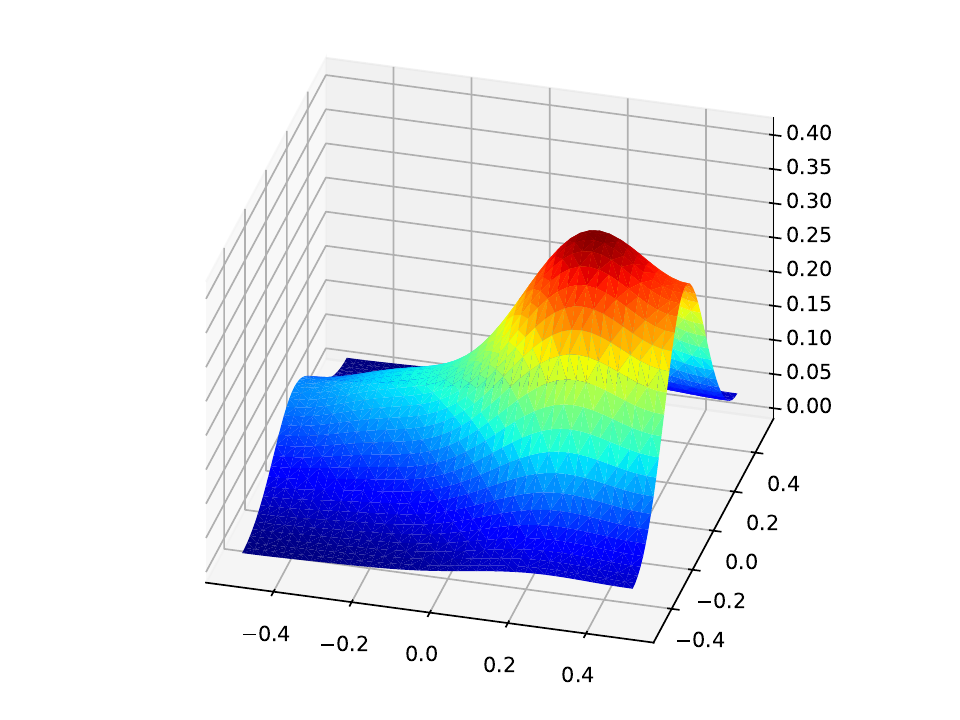}
    \end{subfigure}
    \begin{subfigure}[b]{0.22\textwidth}
        \centering
        \includegraphics[width=1.0\textwidth]{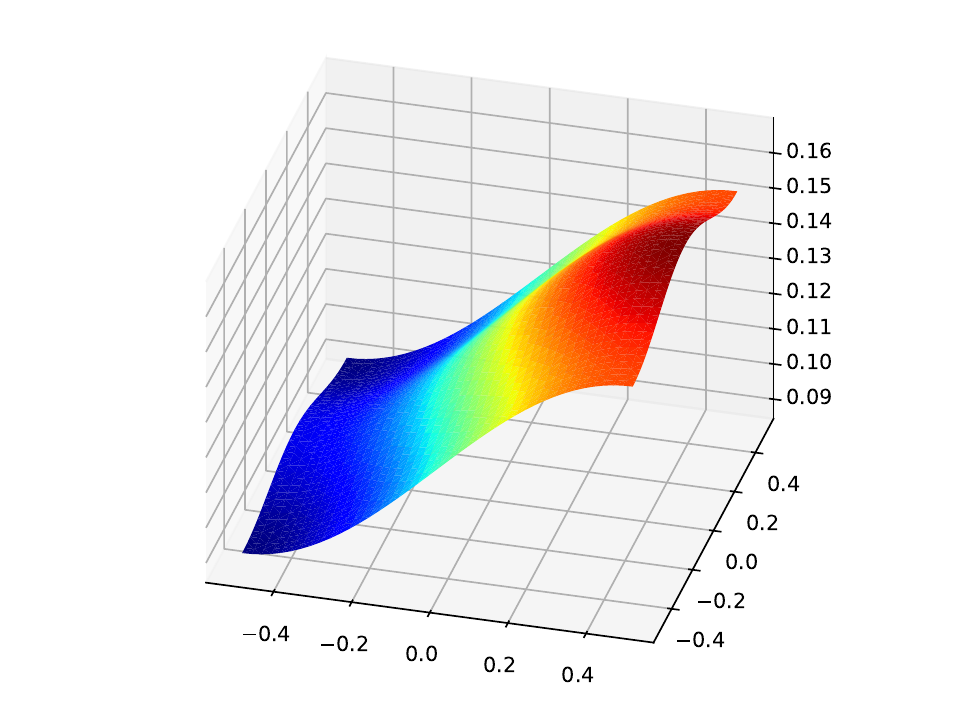}
    \end{subfigure}
    \begin{subfigure}[b]{0.22\textwidth}
        \centering
        \includegraphics[width=1.0\textwidth]{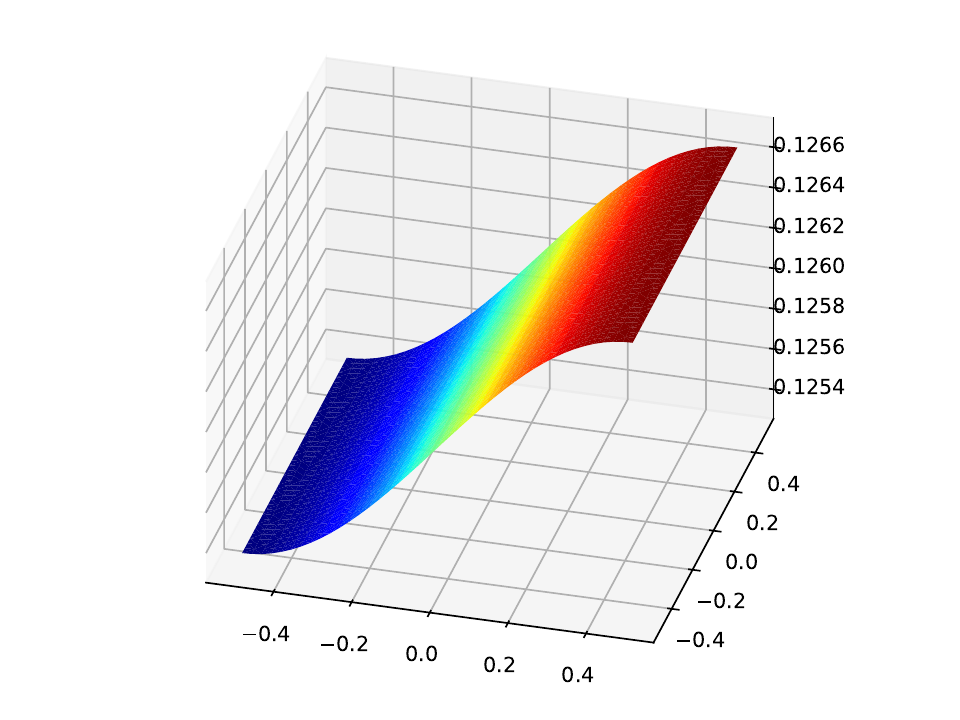}
    \end{subfigure}    
    \caption*{Snapshots of $p_h$ at times $t=0.01$, $0.02$, $0.1$ and $0.5$}
    \begin{subfigure}[b]{0.22\textwidth}
        \centering
        \includegraphics[width=1.0\textwidth]{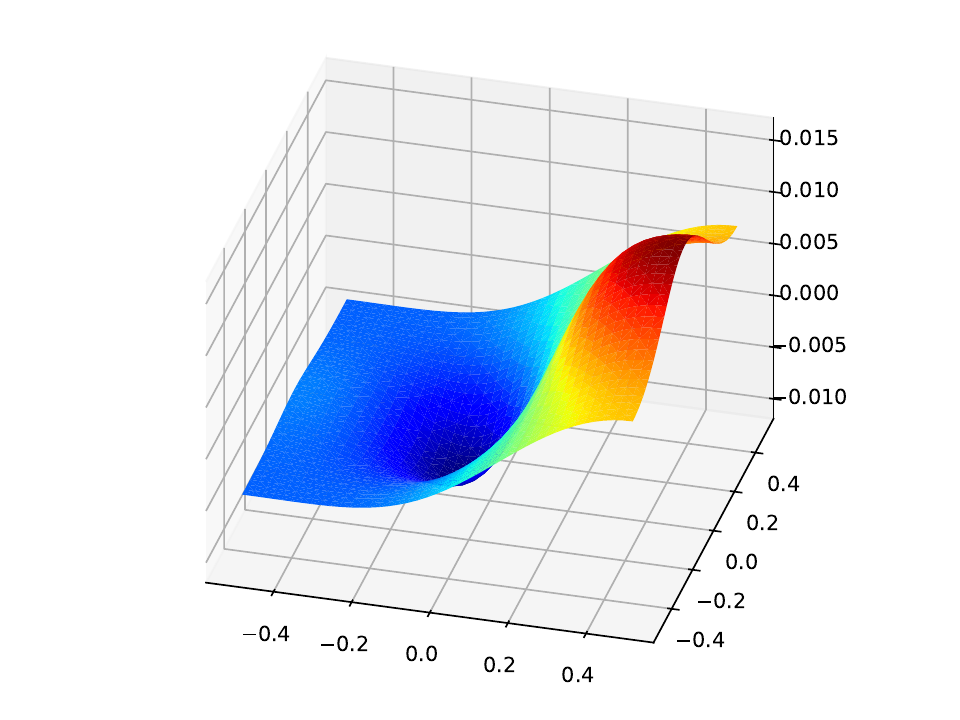}
    \end{subfigure}
    \begin{subfigure}[b]{0.22\textwidth}
        \centering
        \includegraphics[width=1.0\textwidth]{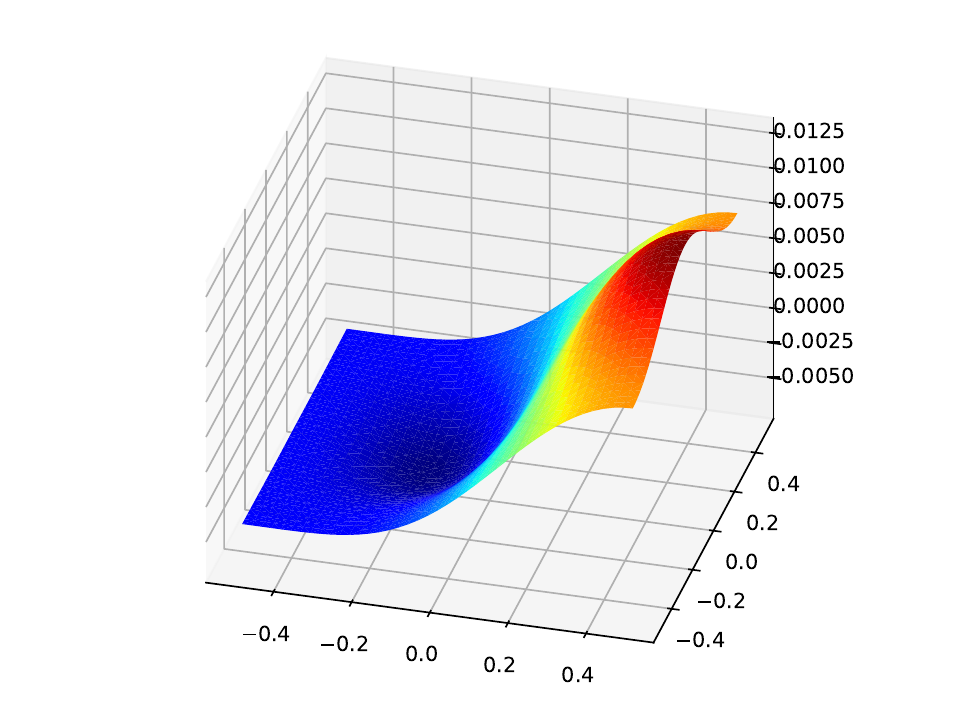}
    \end{subfigure}
    \begin{subfigure}[b]{0.22\textwidth}
        \centering
        \includegraphics[width=1.0\textwidth]{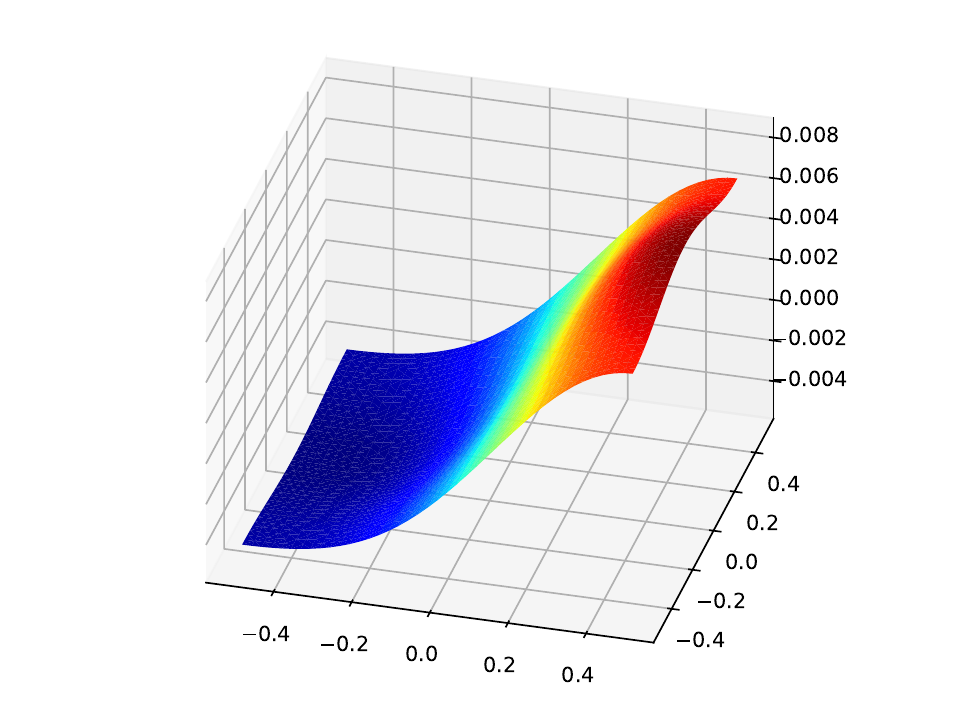}
    \end{subfigure}
    \begin{subfigure}[b]{0.22\textwidth}
        \centering
        \includegraphics[width=1.0\textwidth]{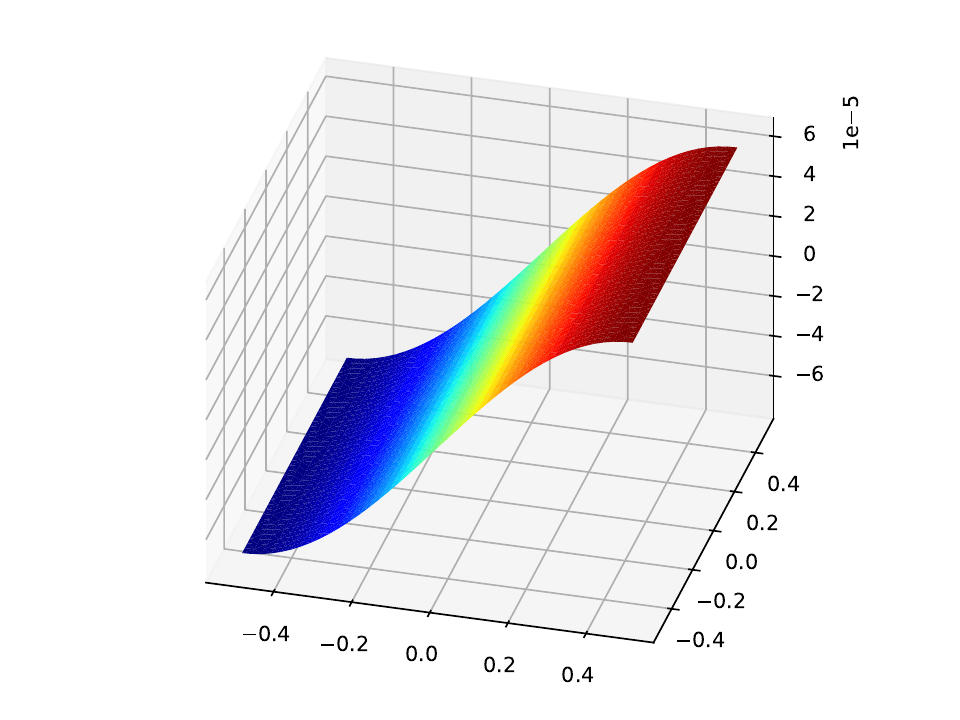}
    \end{subfigure}    
    \caption*{Snapshots of $\phi_h$ at times $t=0.01$, $0.02$, $0.04$ and $0.5$}
    \caption{Algorithm 1}
    \label{fig_smooth:snapshots_alg1}
\end{figure}

\begin{figure}    
\begin{subfigure}[b]{0.22\textwidth}
        \centering
        \includegraphics[width=1.0\textwidth]{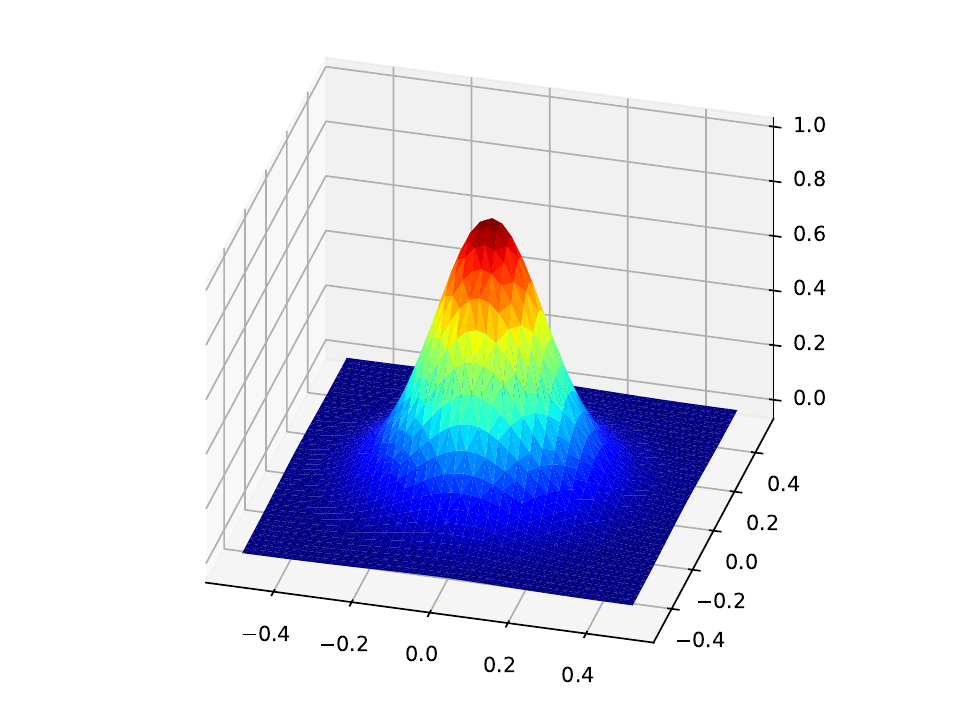}
    \end{subfigure}
    \begin{subfigure}[b]{0.22\textwidth}
        \centering
        \includegraphics[width=1.0\textwidth]{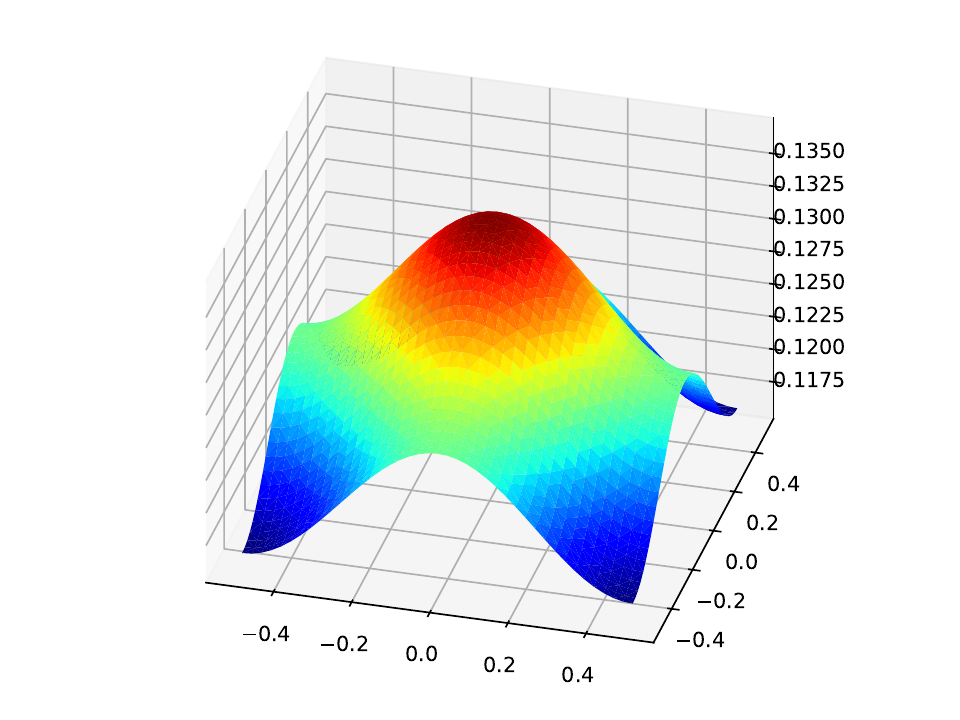}
    \end{subfigure}
    \begin{subfigure}[b]{0.22\textwidth}
        \centering
        \includegraphics[width=1.0\textwidth]{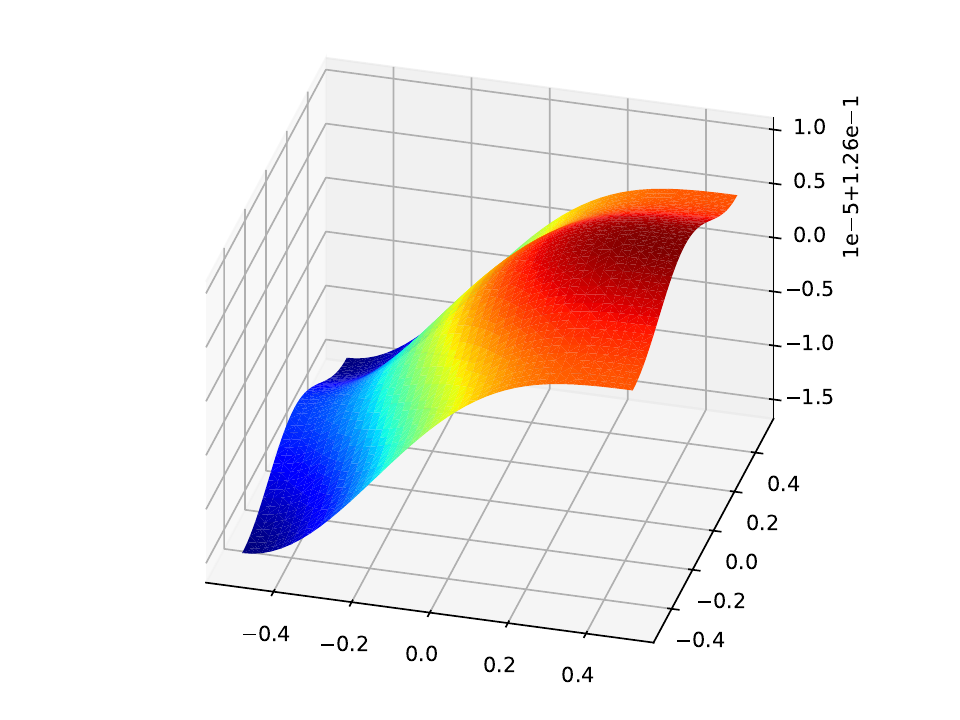}
    \end{subfigure}
    \begin{subfigure}[b]{0.22\textwidth}
        \centering
        \includegraphics[width=1.0\textwidth]{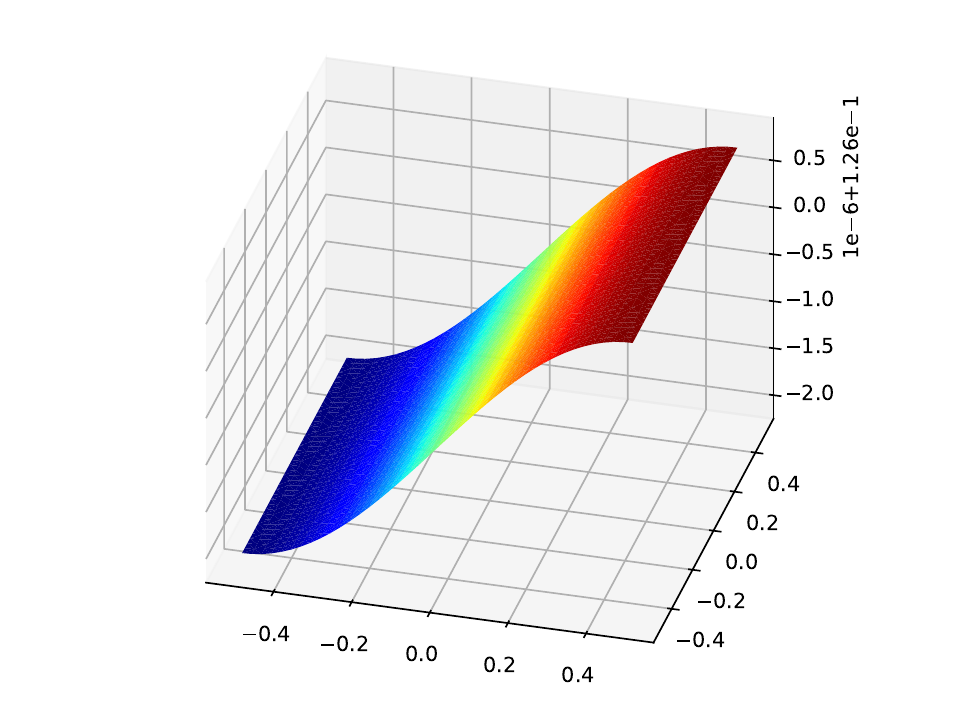}
    \end{subfigure}    
    \caption*{Snapshots of $n_h$ at times $t=0.01$, $0.1$, $0.3$ and $0.5$}
    \begin{subfigure}[b]{0.22\textwidth}
        \centering
        \includegraphics[width=1.0\textwidth]{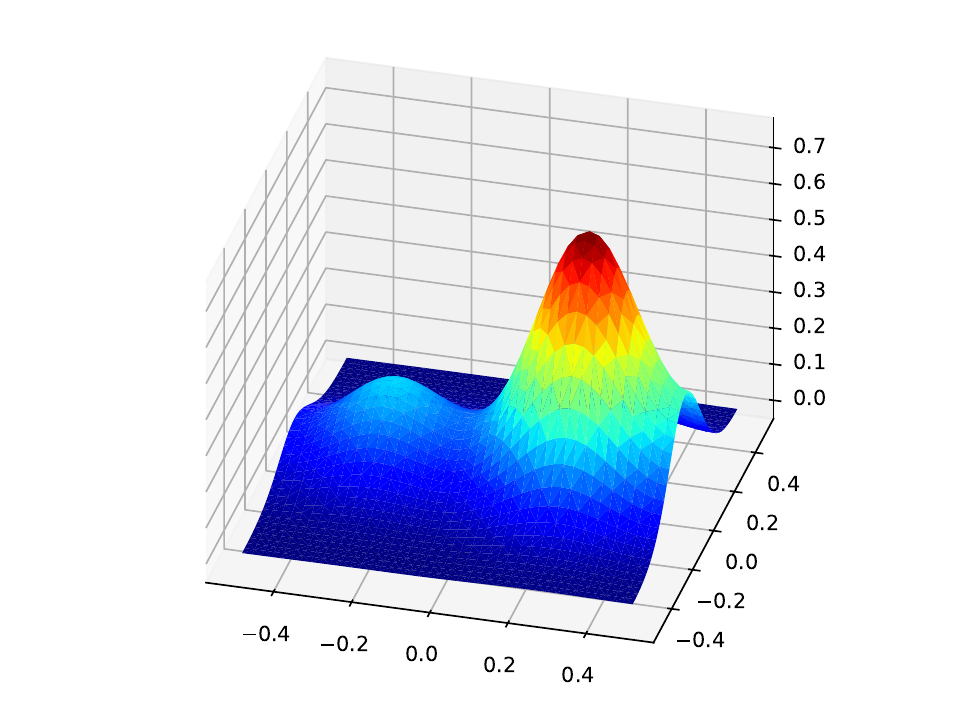}
    \end{subfigure}
    \begin{subfigure}[b]{0.22\textwidth}
        \centering
        \includegraphics[width=1.0\textwidth]{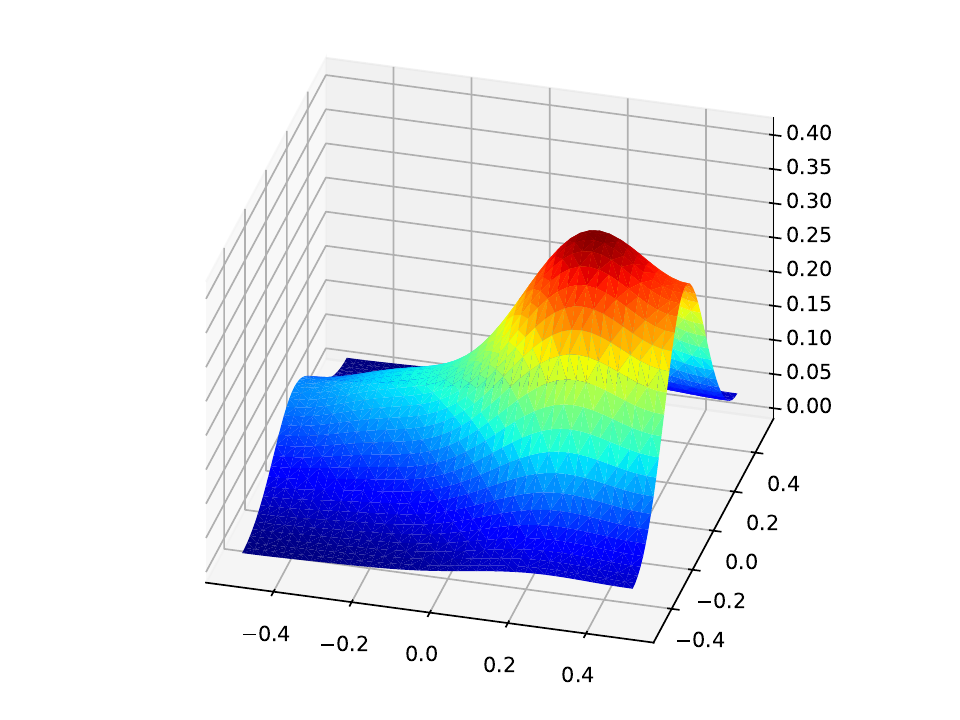}
    \end{subfigure}
    \begin{subfigure}[b]{0.22\textwidth}
        \centering
        \includegraphics[width=1.0\textwidth]{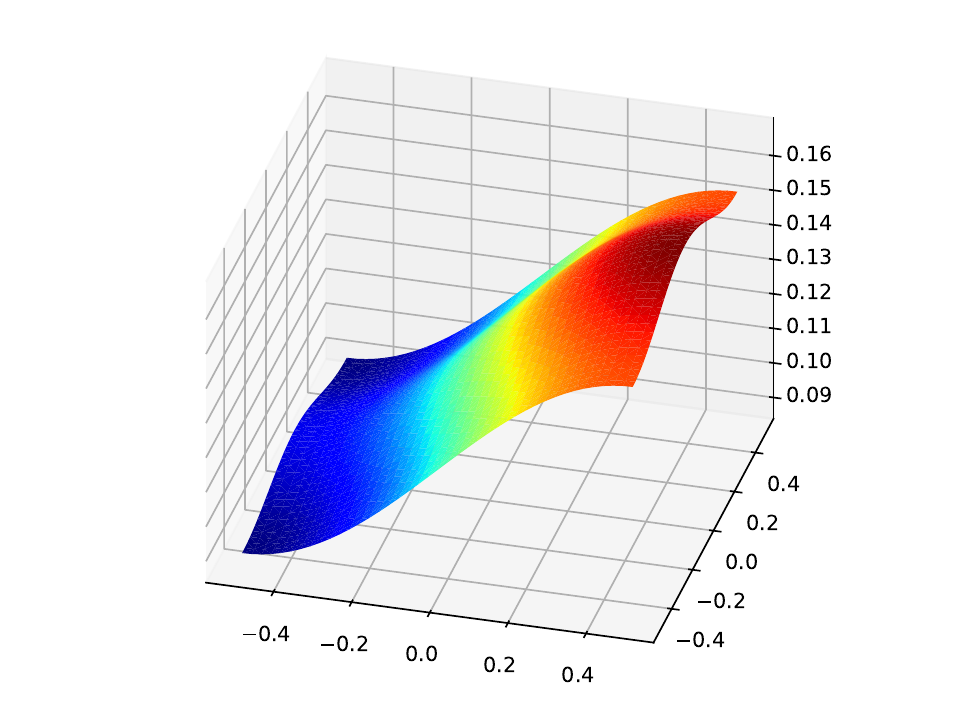}
    \end{subfigure}
    \begin{subfigure}[b]{0.22\textwidth}
        \centering
        \includegraphics[width=1.0\textwidth]{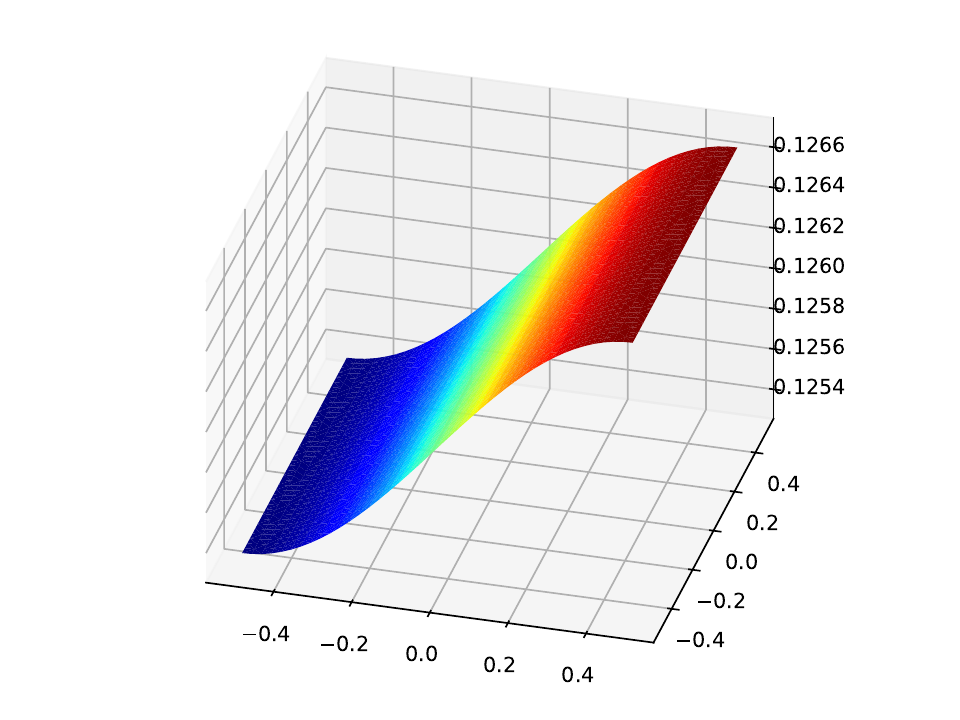}
    \end{subfigure}    
    \caption*{Snapshots of $p_h$ at times $t=0.01$, $0.02$, $0.1$ and $0.5$}
\begin{subfigure}[b]{0.22\textwidth}
        \centering
        \includegraphics[width=1.0\textwidth]{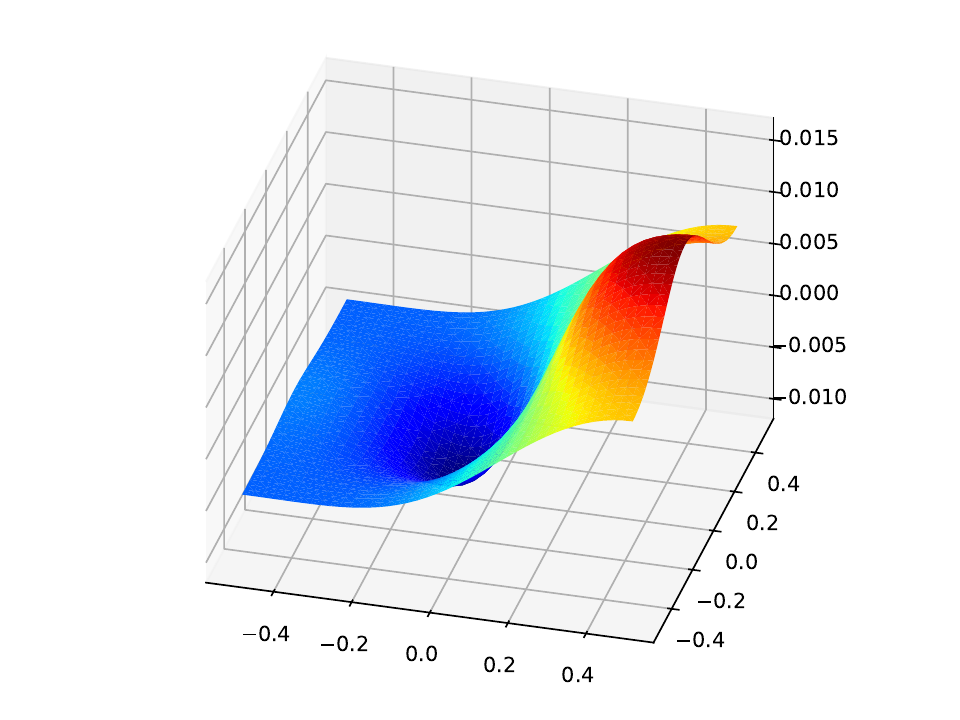}
    \end{subfigure}
    \begin{subfigure}[b]{0.22\textwidth}
        \centering
        \includegraphics[width=1.0\textwidth]{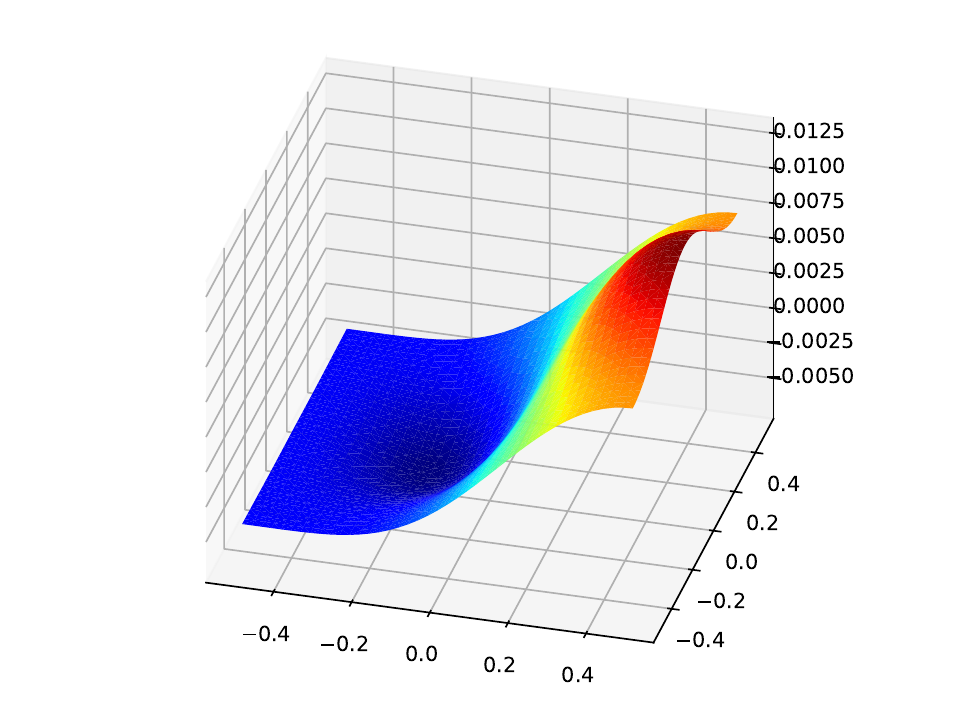}
    \end{subfigure}
    \begin{subfigure}[b]{0.22\textwidth}
        \centering
        \includegraphics[width=1.0\textwidth]{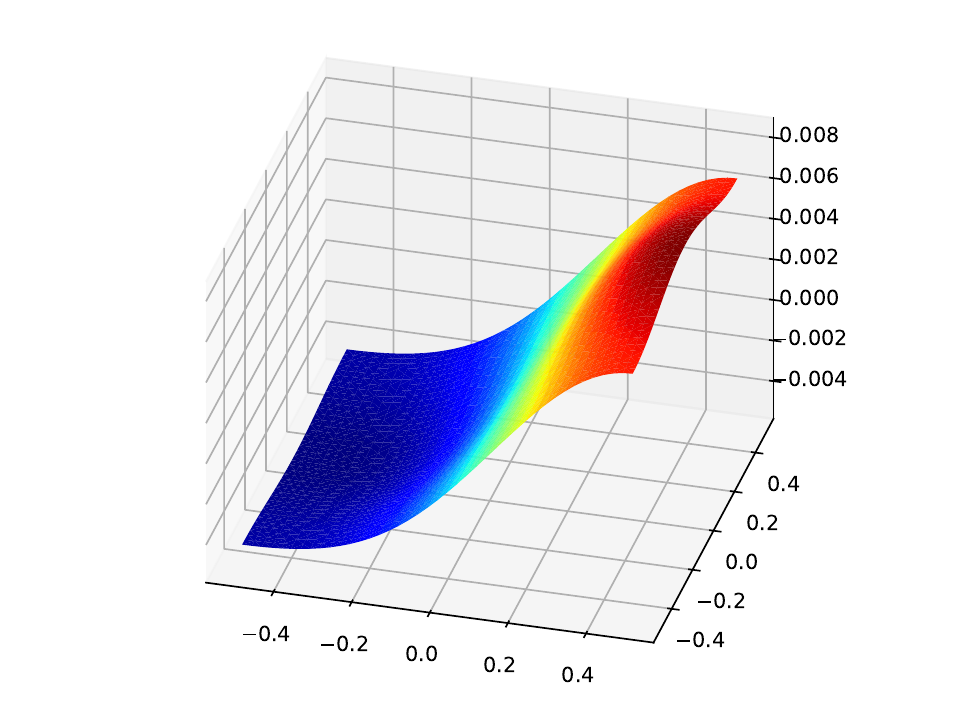}
    \end{subfigure}
    \begin{subfigure}[b]{0.22\textwidth}
        \centering
        \includegraphics[width=1.0\textwidth]{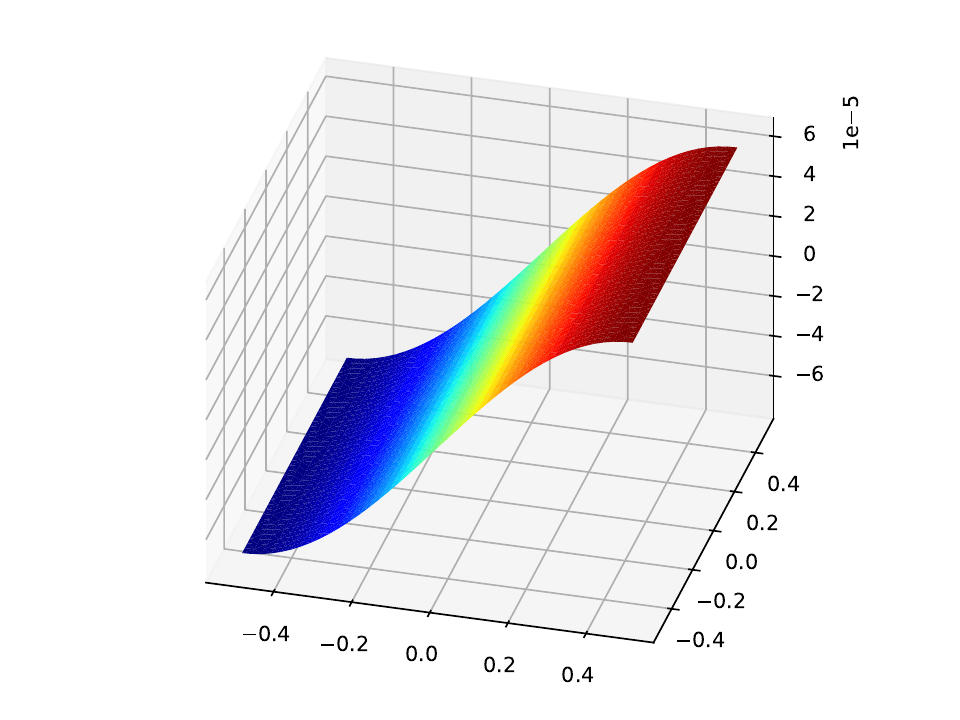}
    \end{subfigure}    
    \caption*{Snapshots of $\phi_h$ at times $t=0.01$, $0.02$, $0.04$ and $0.5$}
    \caption{Algorithm 2}
    \label{fig_smooth:snapshots_alg2}
\end{figure}

One must first observe that both algorithms preserve mass as shown in Figure \ref{fig_smooth:mass_and_energies} and demonstrate consistency with the energy and entropy dissipation. Figure  \ref{fig_smooth:max_and_min} illustrates how maxima and minima evolve over time, with maxima decreasing and minima increasing due to the effect of the diffusion.  Snapshots of $p_h$, $n_h$ and $\phi_h$ at different times are depicted in Figures \ref{fig_smooth:snapshots_alg1} and \ref{fig_smooth:snapshots_alg2} for Algorithms 1 and 2, respectively. No difference is appraised in the performance of both algorithms.    
 \subsection{Ion channel transport} 
\begin{figure}[!h]
\begin{subfigure}[b]{0.2\textwidth}
\centering
\begin{tikzpicture}[scale=0.71]
\draw [ultra thick] (-2,0) -- (2,0) -- (2, 1.5) -- (1 ,1.5) -- (1, 5.5) -- (2, 5.5) -- (2, 7) -- (-2,7) -- (-2, 5.5) -- (-1, 5.5) -- (-1, 1.5) -- (-2, 1.5) -- (-2, 0) ;

\draw  (-2,-0.5) -- (2, - 0.5);
\draw (-2, - 0.25) -- (-2, -0.7);
\draw (-2, - 0.25) -- (-2, -0.7);
\draw (2, - 0.25) -- (2, -0.7);
\draw (0, - 0.25) -- (0, -0.7);
\draw (-2.5, 0) -- (-2.5, 7);
\draw (-2.3, 0) -- (-2.7, 0);
\draw (-2.3, 1.5) -- (-2.7, 1.5);
\draw (-2.3, 5.5) -- (-2.7, 5.5);
\draw (-2.3, 7) -- (-2.7, 7);

\node at (-2,- 1) {$-2$};
\node at (-1,- 1) {$-1$};
\node at (2,- 1) {$2$};
\node at (1,- 1) {$1$};
\node at (0,- 1) {$0$};

\node at (-3, 0) {$0$};
\node at (-3.1, 1.5) {$1.5$};
\node at (-3.1, 5.5) {$5.5$};
\node at (-3, 7) {$7$};


\node at (0,0.5) {$\partial \Omega_b$};
\node at (0, 6.5) {$\partial \Omega_t$};
\node [rotate=90] at (1.5, 3.5) {$\partial\Omega_m$};
\node [rotate=90] at (-1.5, 3.5) {$\partial\Omega_m$}; 
\end{tikzpicture}

\end{subfigure}
\hspace {2cm}
\begin{subfigure}[b]{0.2\textwidth}
\centering
\includegraphics[width=0.8\textwidth]{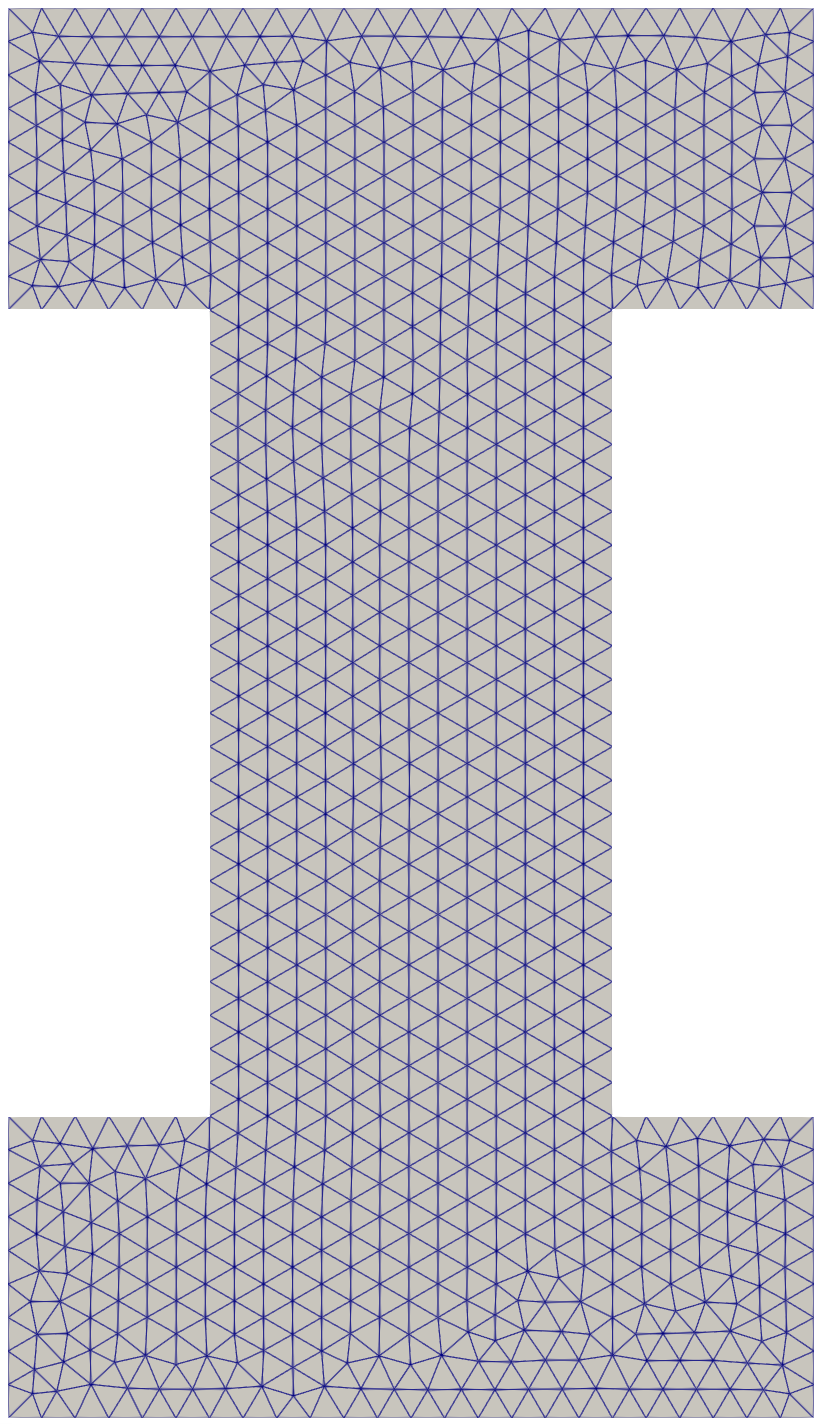}
\vspace{0.85cm}
\end{subfigure}
\caption{Channel and  mesh}
\label{fig:channel_mesh}
\end{figure}

 In this second example,  electrodiffusion of ions in a channel is considered. The phenomenon under consideration is the membrane transport process through a channel driven by an electric potential.   The computational domain is depicted in Figure \ref{fig:channel_mesh} and the electric potential drop for the ion transport is artificially attained by setting Dirichlet boundary conditions at the walls $\partial\Omega_t$ and $\partial\Omega_b$.  More precisely, 
$$
\phi = - 50 \quad\mbox{ on } \quad \partial\Omega_b\quad  \mbox{ and } \quad \phi= 50\quad \mbox{ on } \quad \partial\Omega_t
$$
and
$$
\partial_{\boldsymbol{n}}\phi=0\quad \mbox{ on }\quad \partial\Omega\backslash(\partial\Omega_b\cap \partial\Omega_t).
$$
Further Figure \ref{fig:channel_mesh} displays the mesh used for the forthcoming numerical experiments with mesh size $h=0.121854$. In all of these experiments, $k=10^{-2}$ will be the time step.  
\newpage

\subsubsection{Uniform ion distribution} 
\begin{figure}
    \begin{subfigure}[b]{0.25\textwidth}
        \centering
        \includegraphics[width=1.0\textwidth]{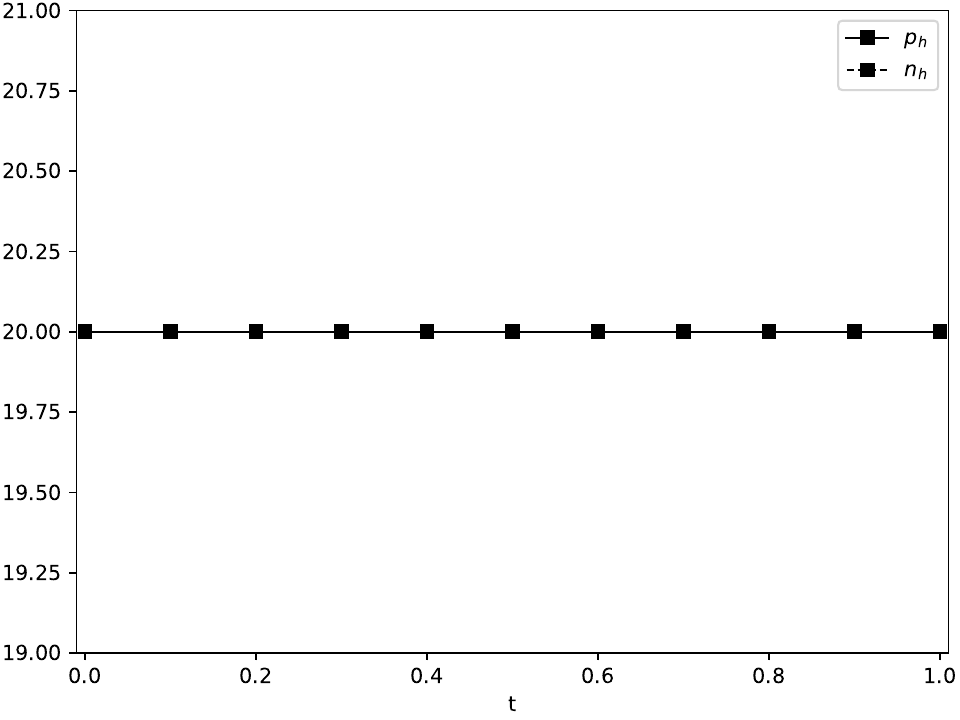}
    \end{subfigure}
    \begin{subfigure}[b]{0.25\textwidth}
        \centering
        \includegraphics[width=1.0\textwidth]{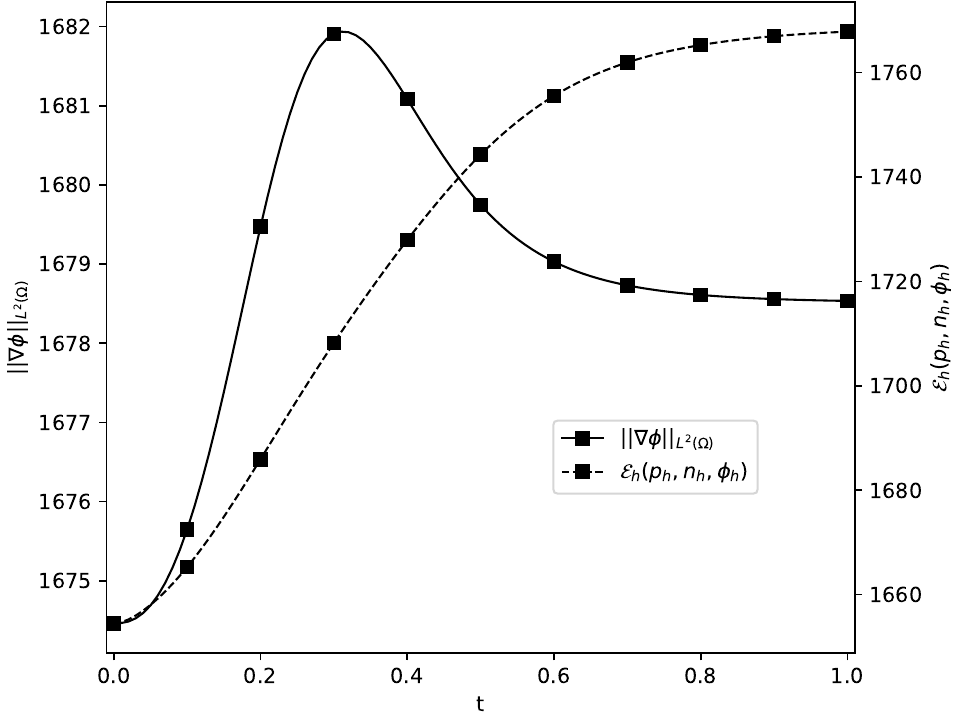}
    \end{subfigure}
    \caption*{Algorithm 1}
    \begin{subfigure}[b]{0.25\textwidth}
        \centering
        \includegraphics[width=1.0\textwidth]{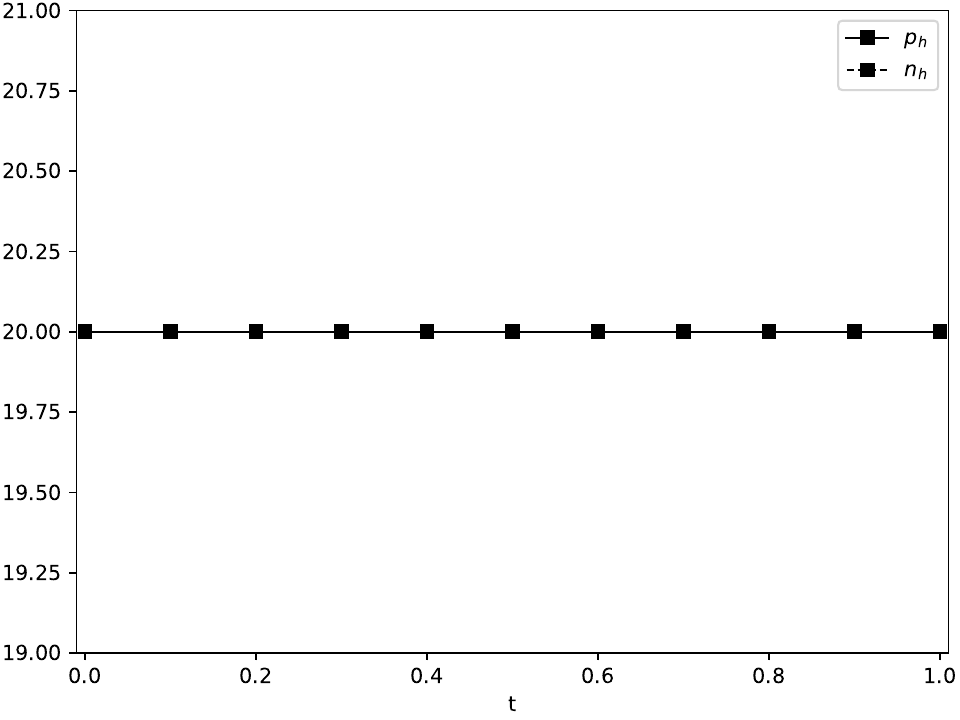}
    \end{subfigure}
    \begin{subfigure}[b]{0.25\textwidth}
        \centering
        \includegraphics[width=1.0\textwidth]{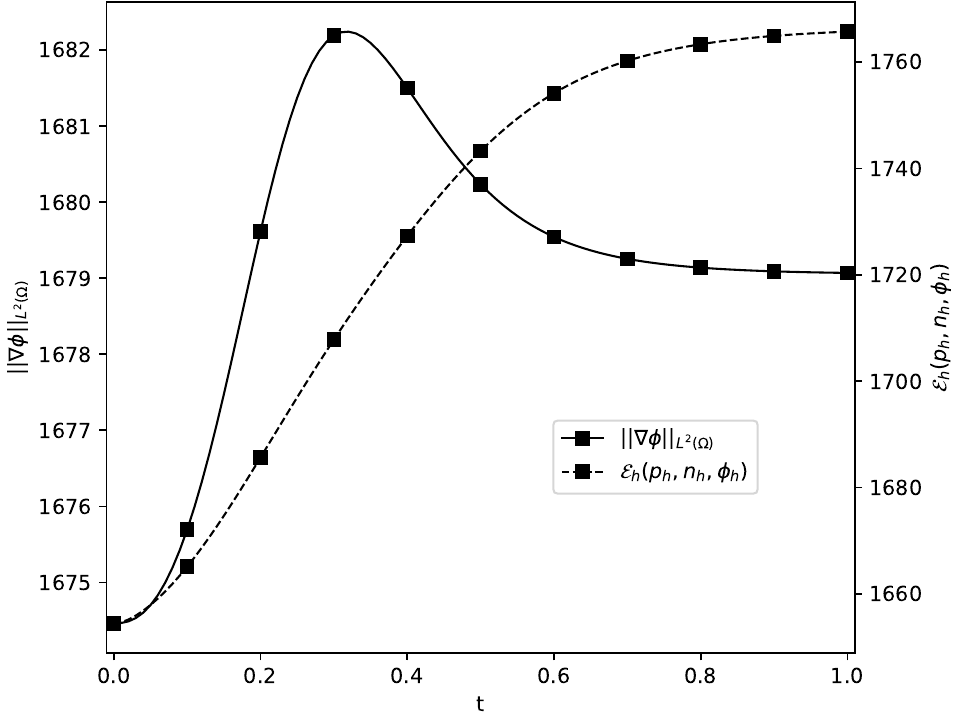}
    \end{subfigure}
    \caption*{Algorithm 2}
    \caption{Mass conservation (left), and energy and entropy evolutions (right)}
    \label{fig_channel:mass_and_energies}
\end{figure}

\begin{figure}
    \begin{subfigure}[b]{0.25\textwidth}
        \centering
        \includegraphics[width=1.0\textwidth]{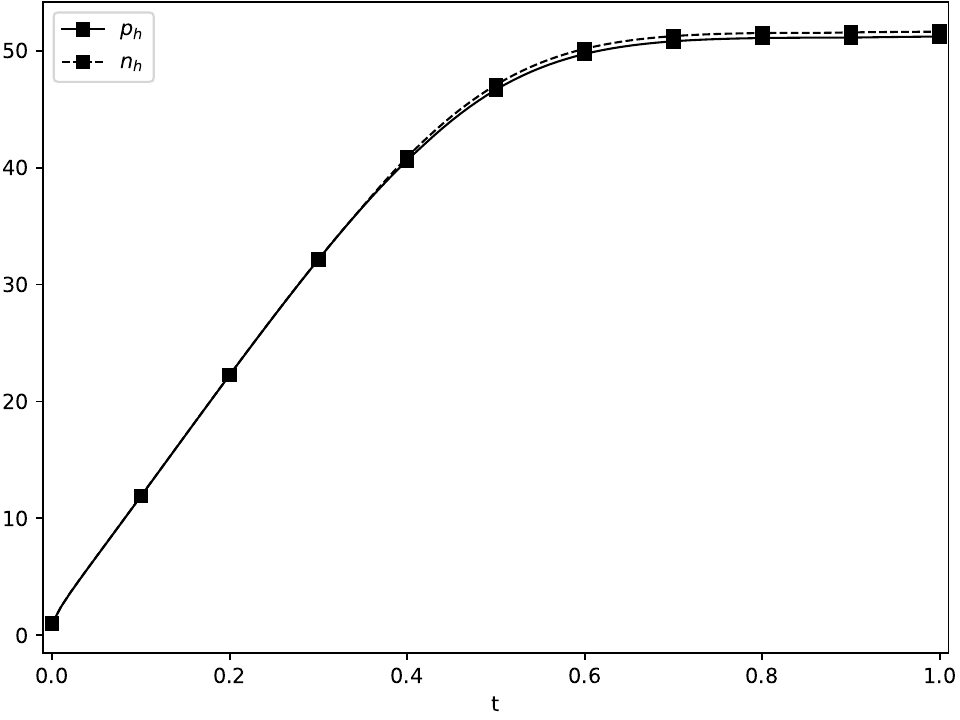}
     \end{subfigure}
    \begin{subfigure}[b]{0.25\textwidth}
        \centering
        \includegraphics[width=1.0\textwidth]{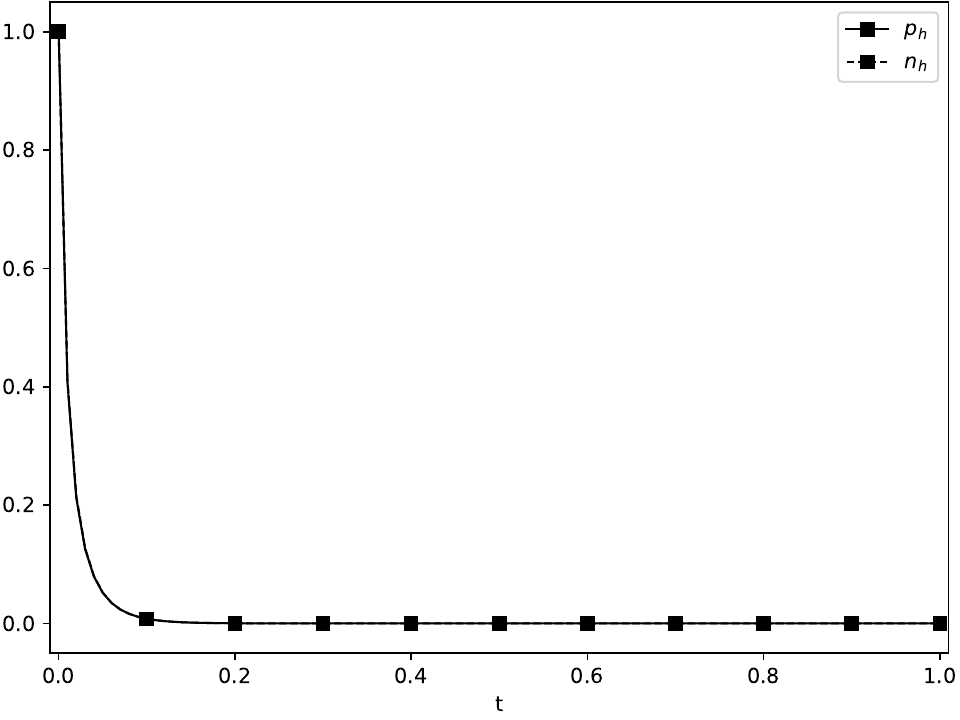}
     \end{subfigure}
  \caption*{Algorithm 1}\label{Alg1:smooth} 
       \begin{subfigure}[b]{0.25\textwidth}
        \centering
        \includegraphics[width=1.0\textwidth]{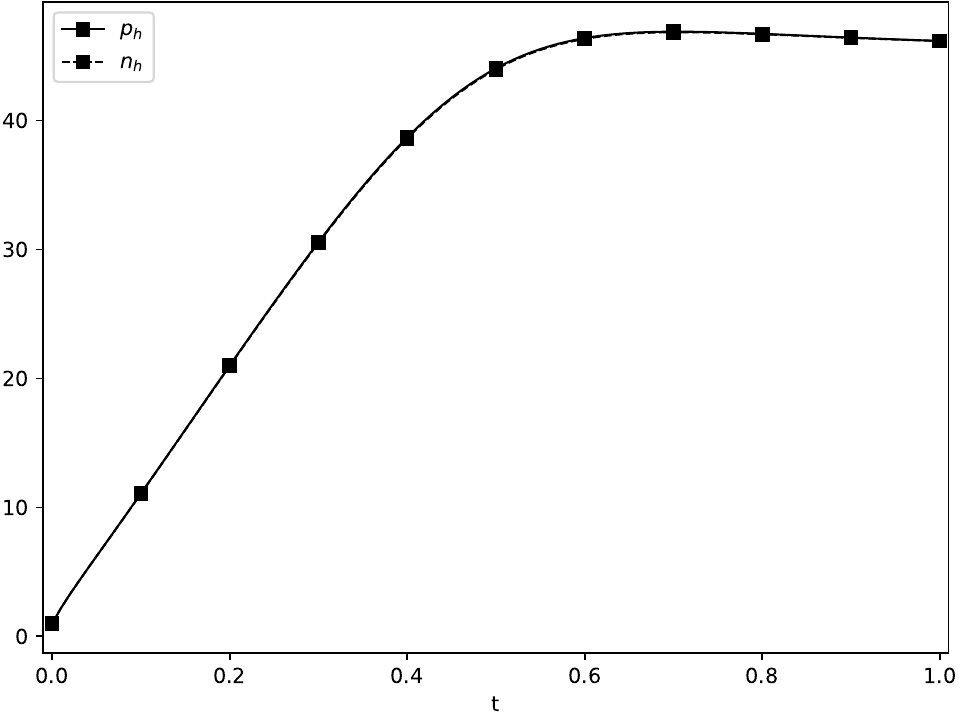}
    \end{subfigure}
    \begin{subfigure}[b]{0.25\textwidth}
        \centering
        \includegraphics[width=1.0\textwidth]{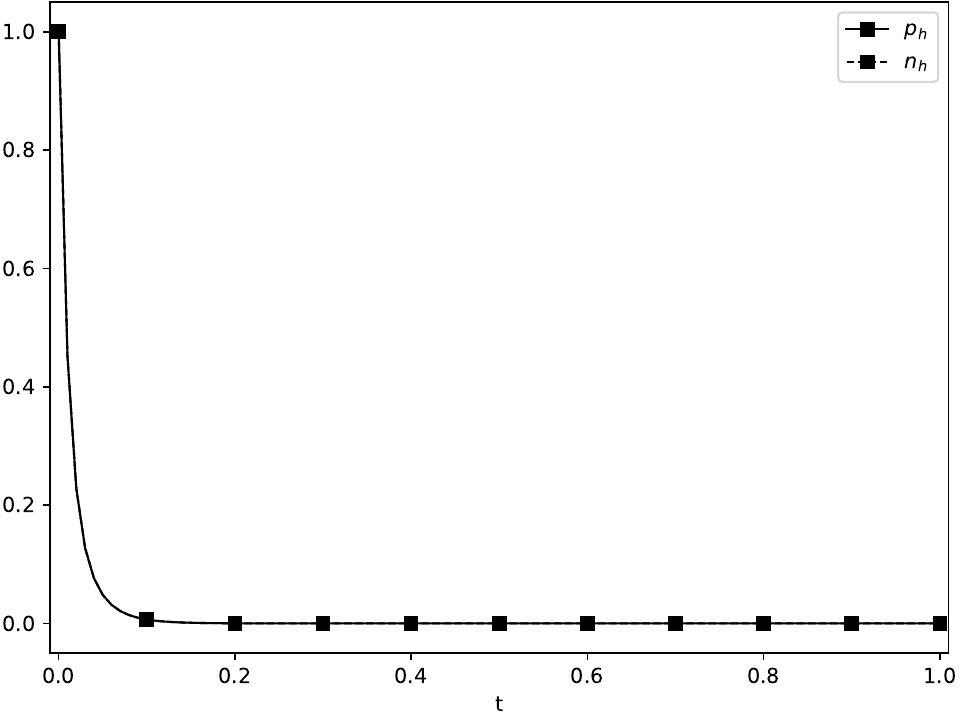}
    \end{subfigure}    
    \caption*{Algorithm 2}\label{Alg2:smooth}
    \caption{Maxima (left) and minima (right)}
    \label{fig_channel:max_and_min}
\end{figure}

\begin{figure} 
\begin{subfigure}[b]{0.22\textwidth}
        \centering
        \includegraphics[width=1.0\textwidth]{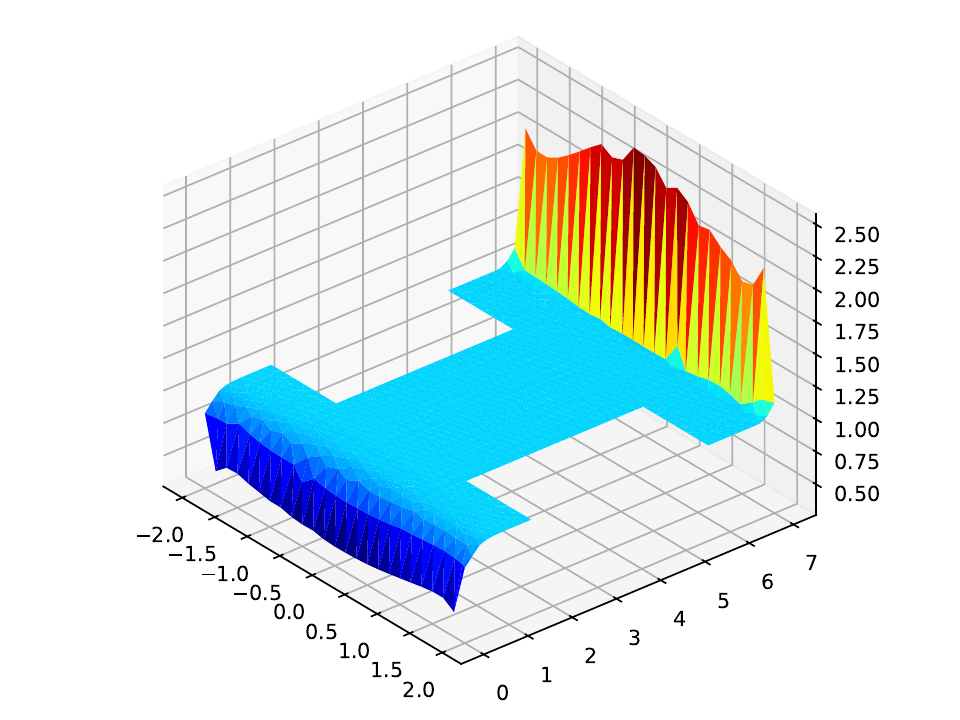}
    \end{subfigure}
    \begin{subfigure}[b]{0.22\textwidth}
        \centering
        \includegraphics[width=1.0\textwidth]{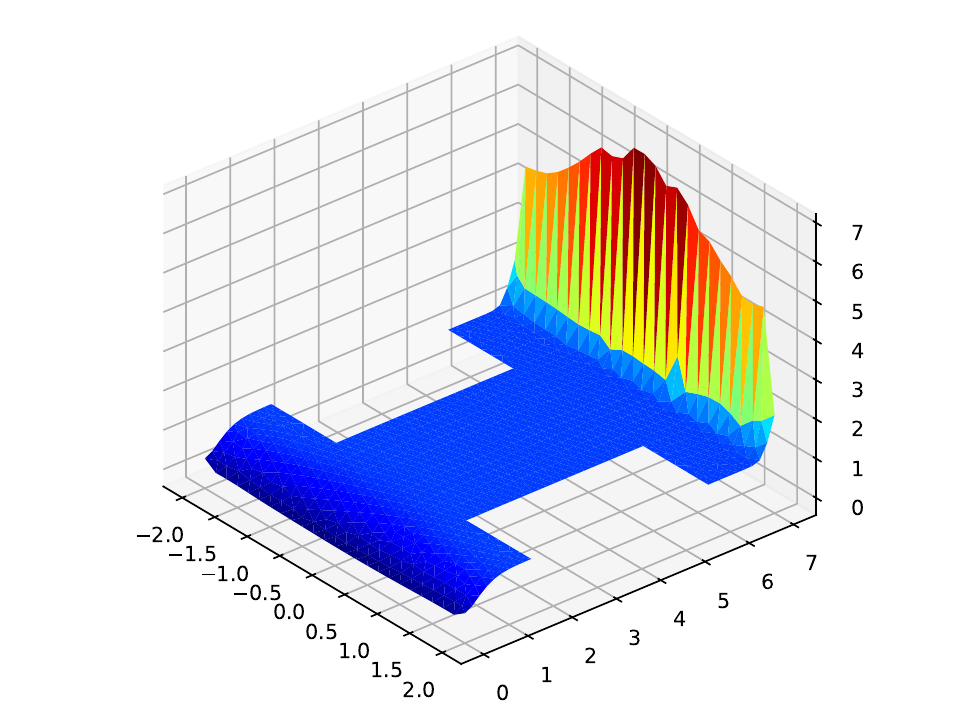}
    \end{subfigure}
    \begin{subfigure}[b]{0.22\textwidth}
        \centering
        \includegraphics[width=1.0\textwidth]{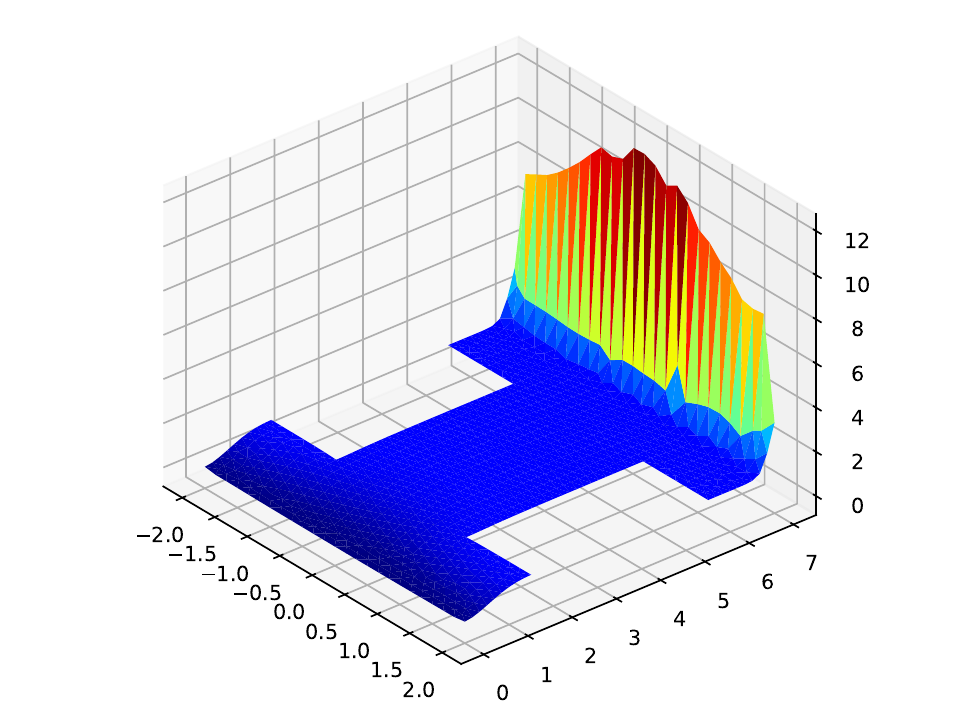}
    \end{subfigure}
    \begin{subfigure}[b]{0.22\textwidth}
        \centering
        \includegraphics[width=1.0\textwidth]{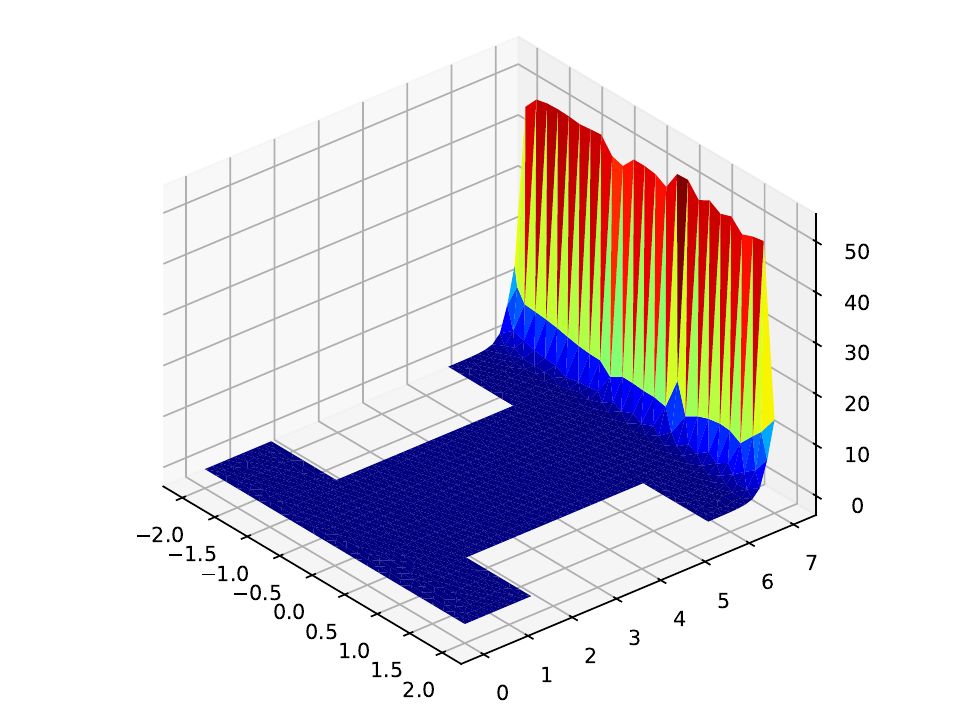}
    \end{subfigure}    
    \caption*{Snapshots of $n_h$ at times $t=0.01$, $0.05$, $0.1$ and $1.0$}
    \begin{subfigure}[b]{0.22\textwidth}
        \centering
        \includegraphics[width=1.0\textwidth]{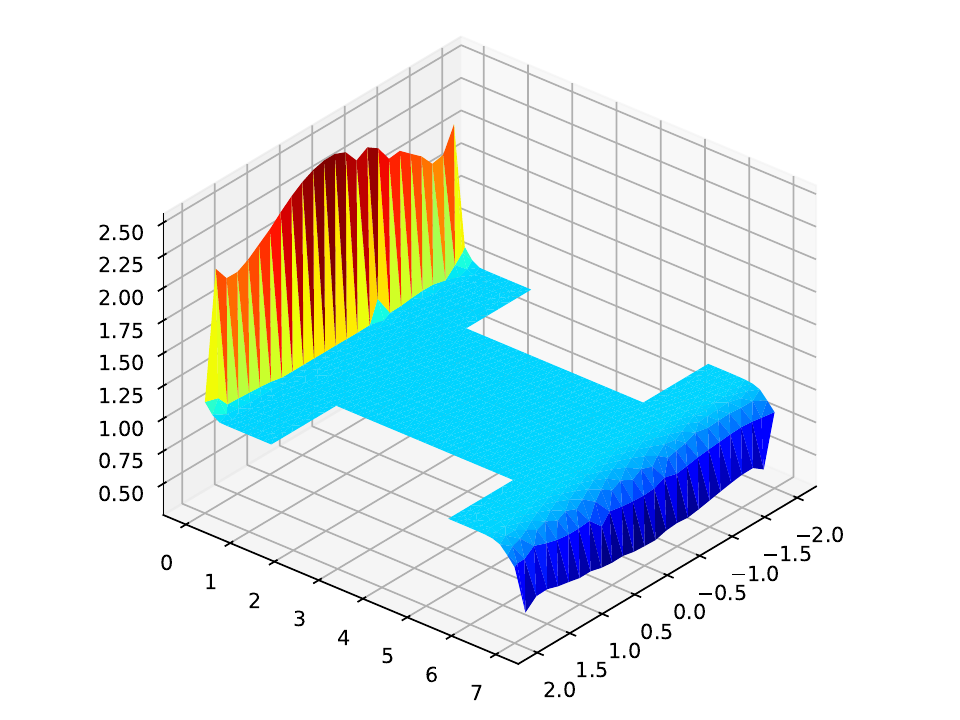}
    \end{subfigure}
    \begin{subfigure}[b]{0.22\textwidth}
        \centering
        \includegraphics[width=1.0\textwidth]{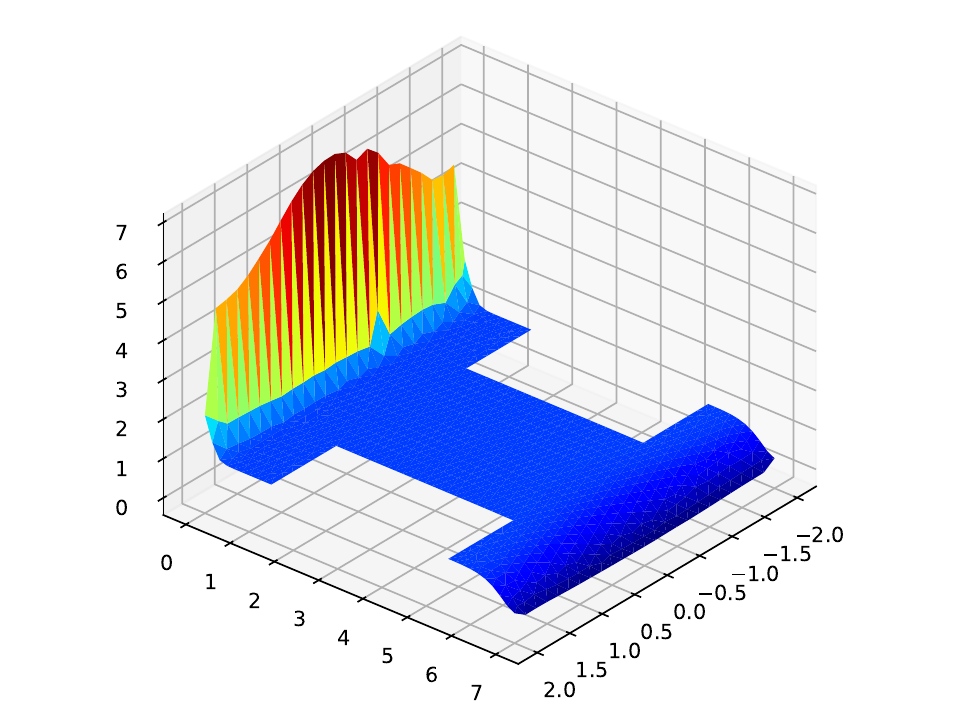}
    \end{subfigure}
    \begin{subfigure}[b]{0.22\textwidth}
        \centering
        \includegraphics[width=1.0\textwidth]{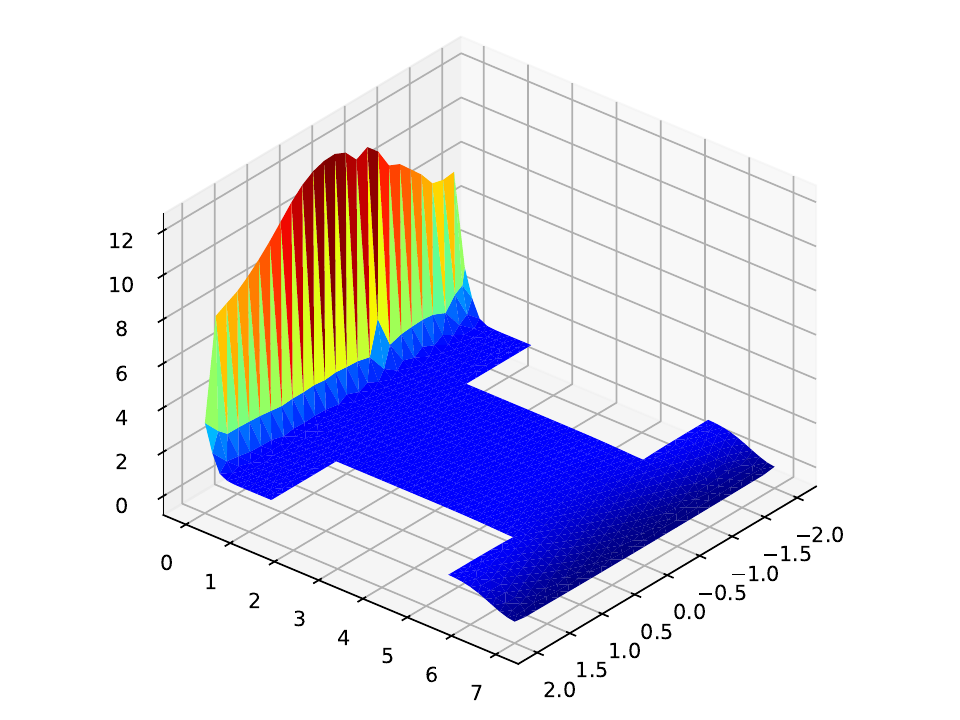}
    \end{subfigure}
    \begin{subfigure}[b]{0.22\textwidth}
        \centering
        \includegraphics[width=1.0\textwidth]{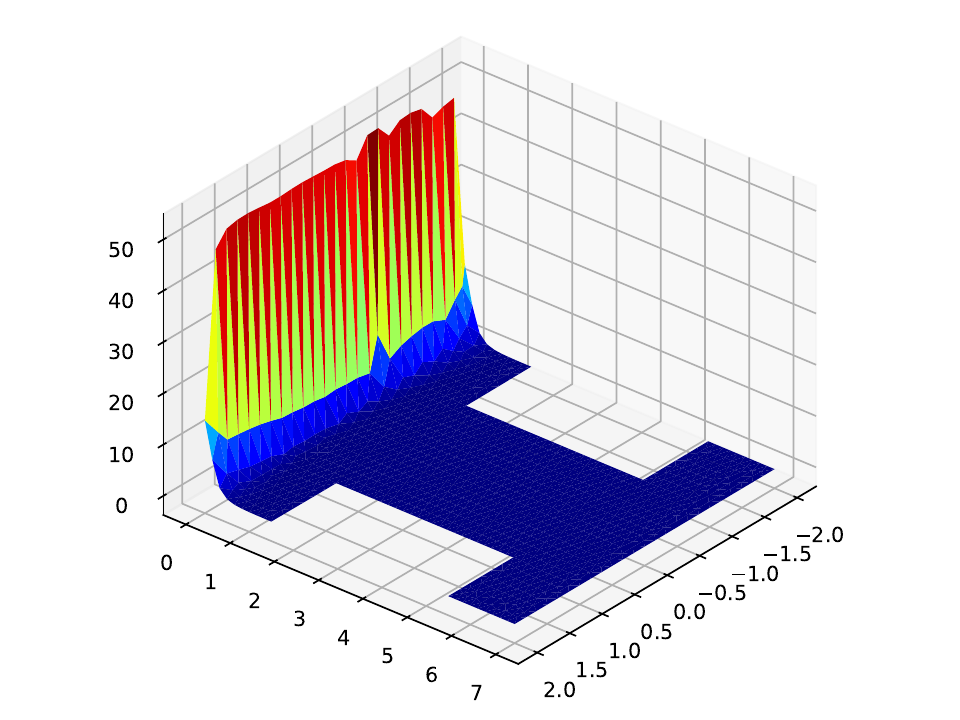}
    \end{subfigure}    
    \caption*{Snapshots of $p_h$ at times $t=0.01$, $0.05$, $0.1$ and $1.0$}
    \caption{Algorithm 1}
    \label{fig_channel:snapshots_alg1}
\end{figure}

\begin{figure}    
\begin{subfigure}[b]{0.22\textwidth}
        \centering
        \includegraphics[width=1.0\textwidth]{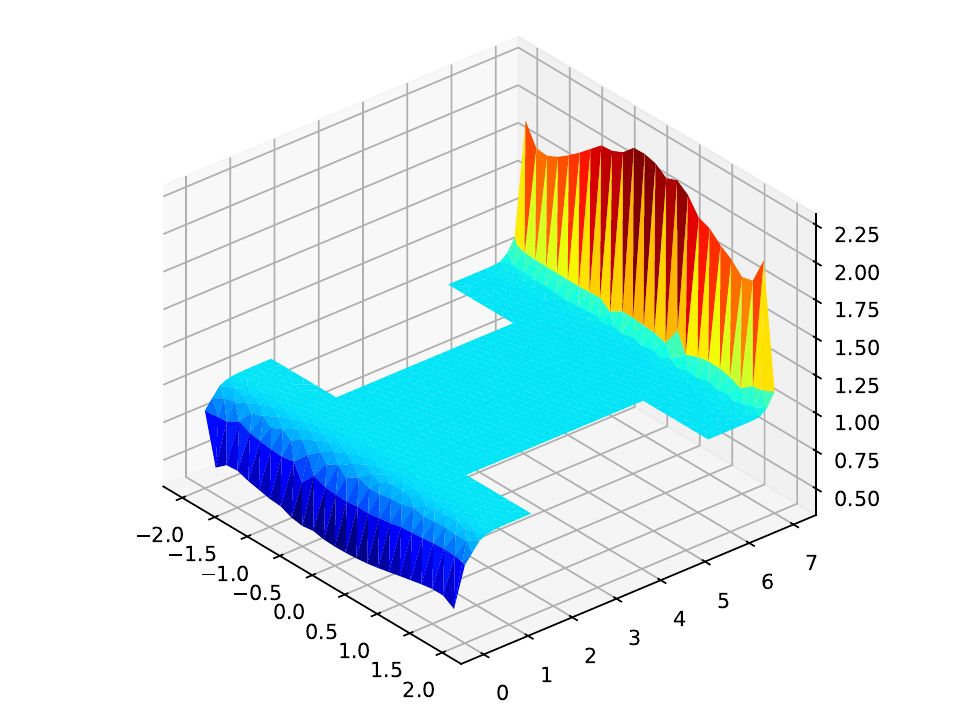}
    \end{subfigure}
    \begin{subfigure}[b]{0.22\textwidth}
        \centering
        \includegraphics[width=1.0\textwidth]{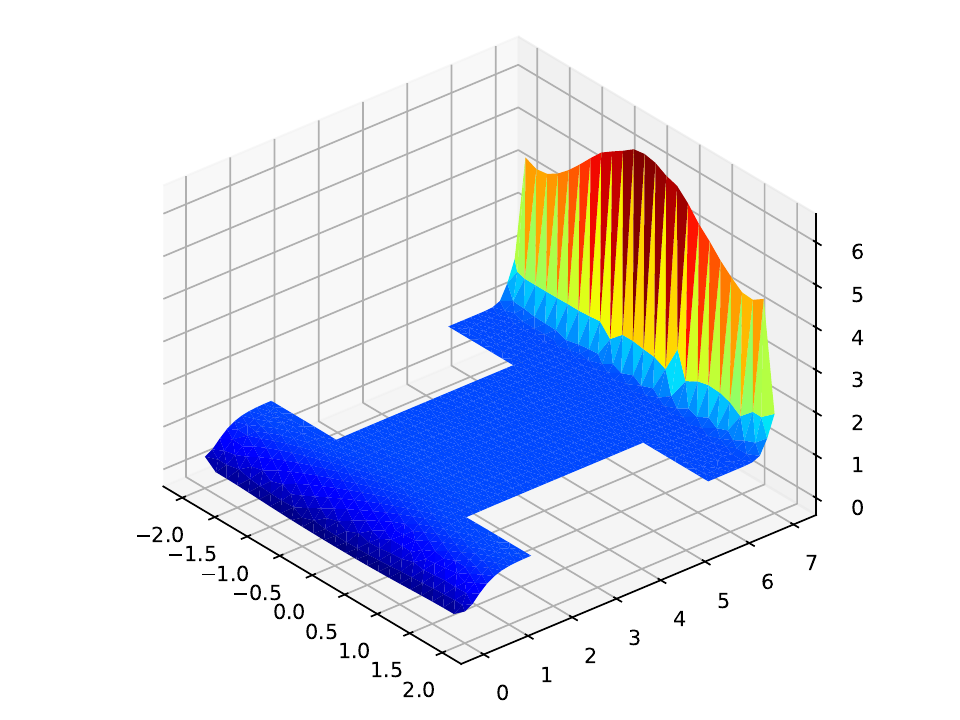}
    \end{subfigure}
    \begin{subfigure}[b]{0.22\textwidth}
        \centering
        \includegraphics[width=1.0\textwidth]{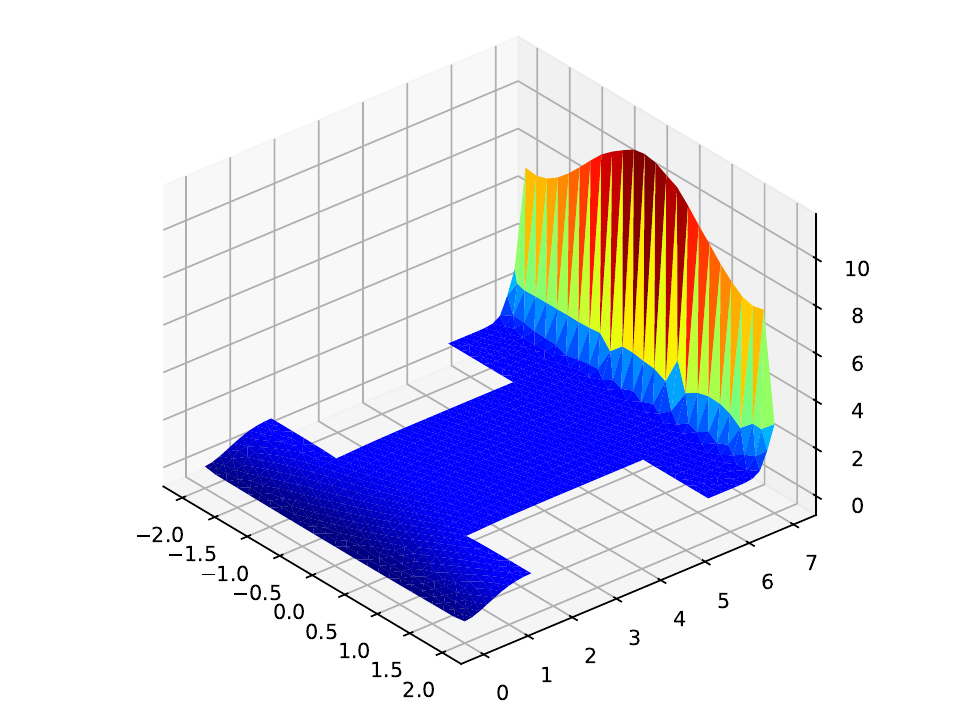}
    \end{subfigure}
    \begin{subfigure}[b]{0.22\textwidth}
        \centering
        \includegraphics[width=1.0\textwidth]{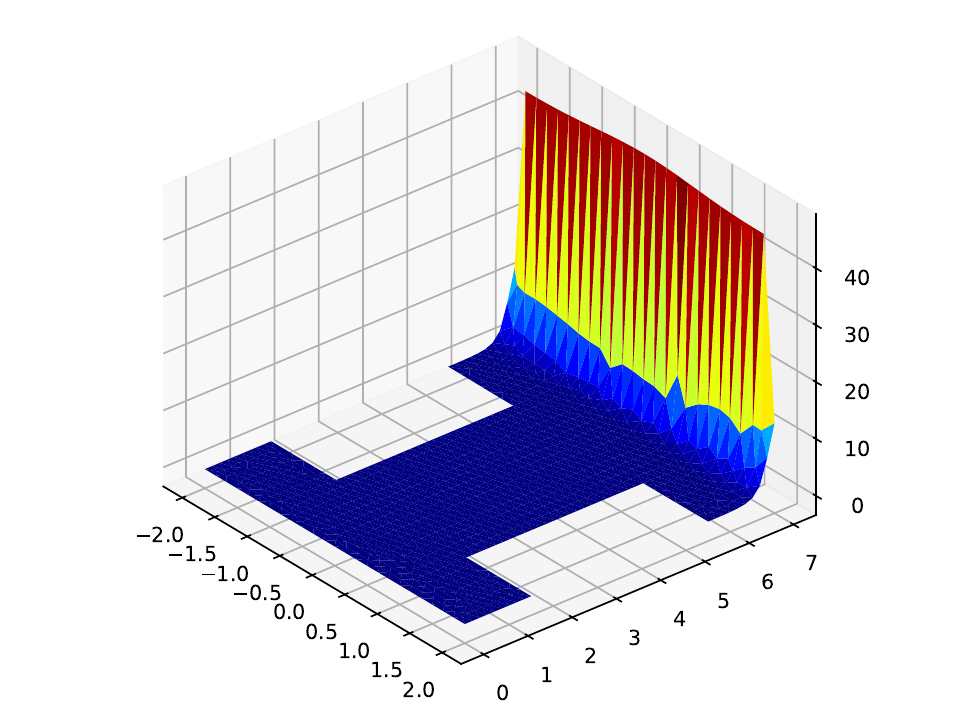}
    \end{subfigure}    
    \caption*{Snapshots of $n_h$ at times $t=0.01$, $0.05$, $0.1$ and $1.0$}
    \begin{subfigure}[b]{0.22\textwidth}
        \centering
        \includegraphics[width=1.0\textwidth]{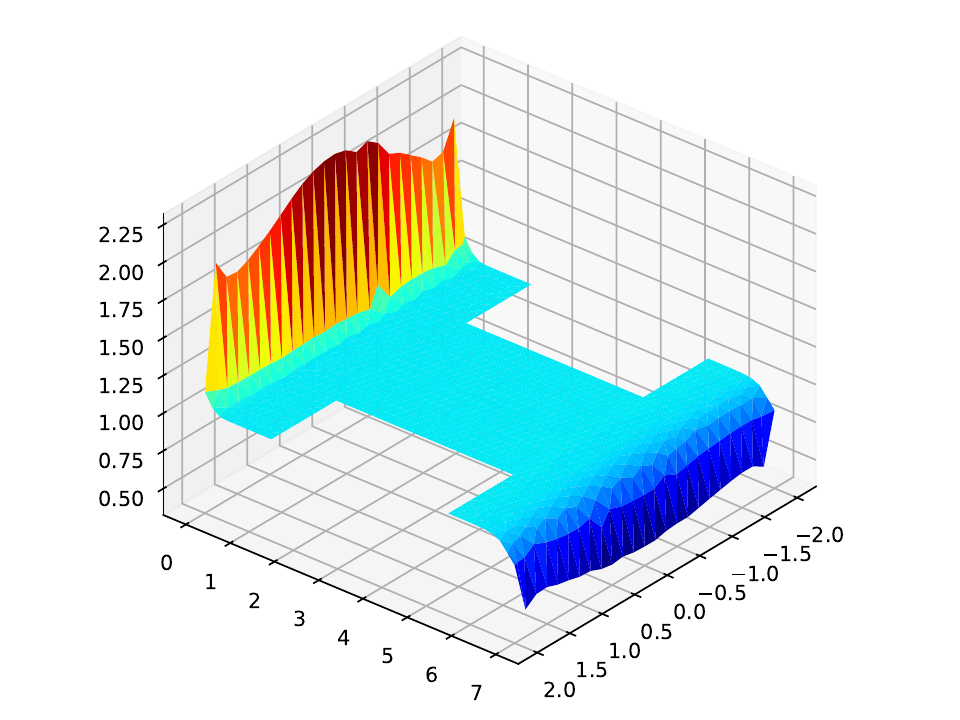}
    \end{subfigure}
    \begin{subfigure}[b]{0.22\textwidth}
        \centering
        \includegraphics[width=1.0\textwidth]{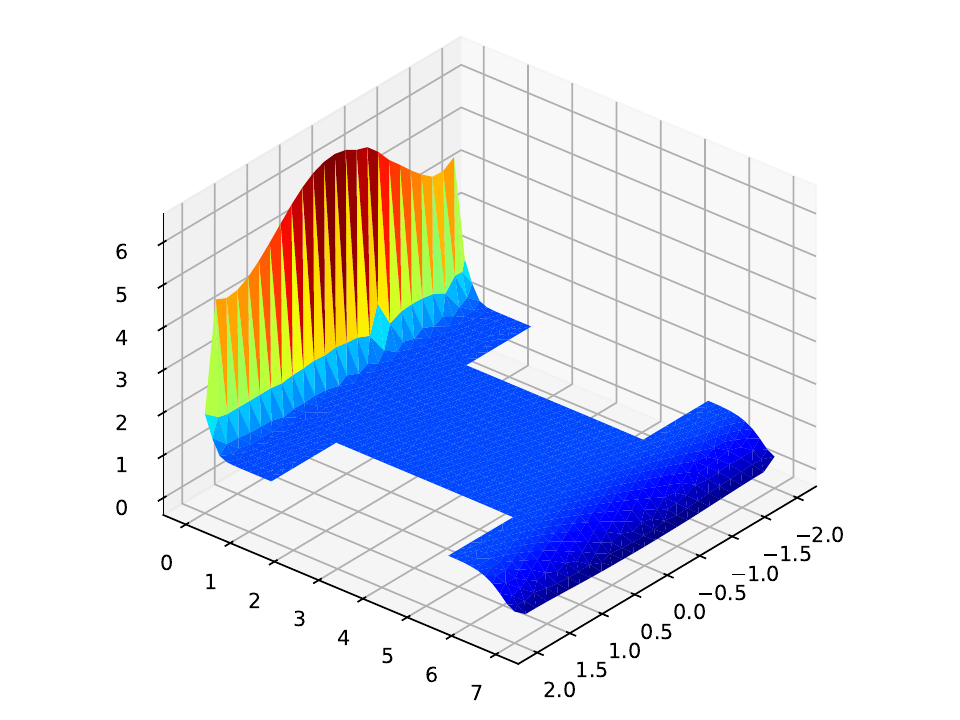}
    \end{subfigure}
    \begin{subfigure}[b]{0.22\textwidth}
        \centering
        \includegraphics[width=1.0\textwidth]{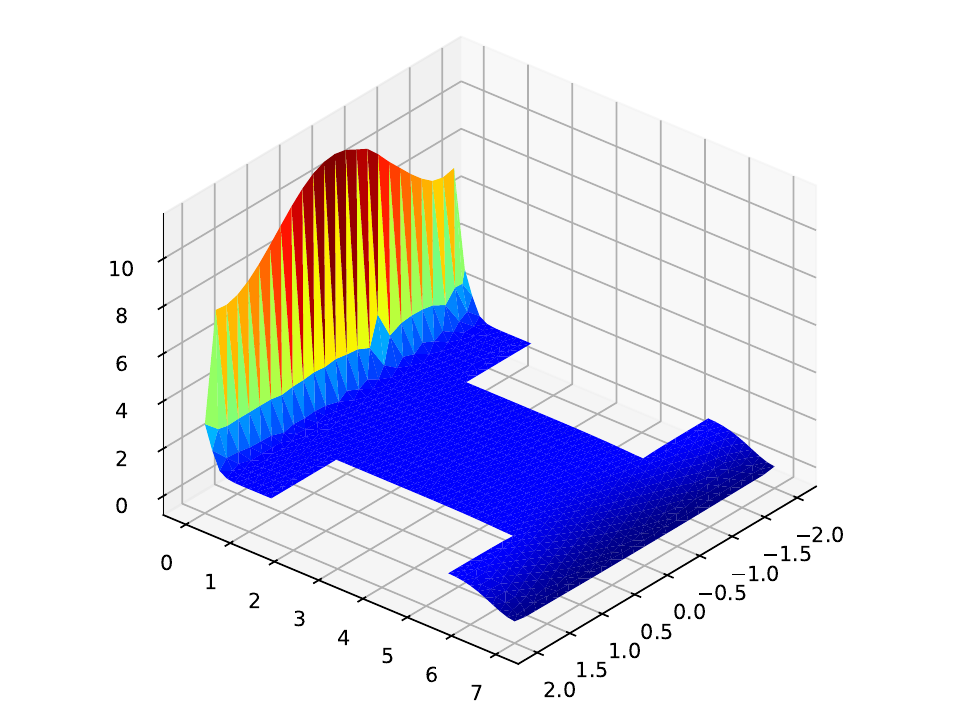}
    \end{subfigure}
    \begin{subfigure}[b]{0.22\textwidth}
        \centering
        \includegraphics[width=1.0\textwidth]{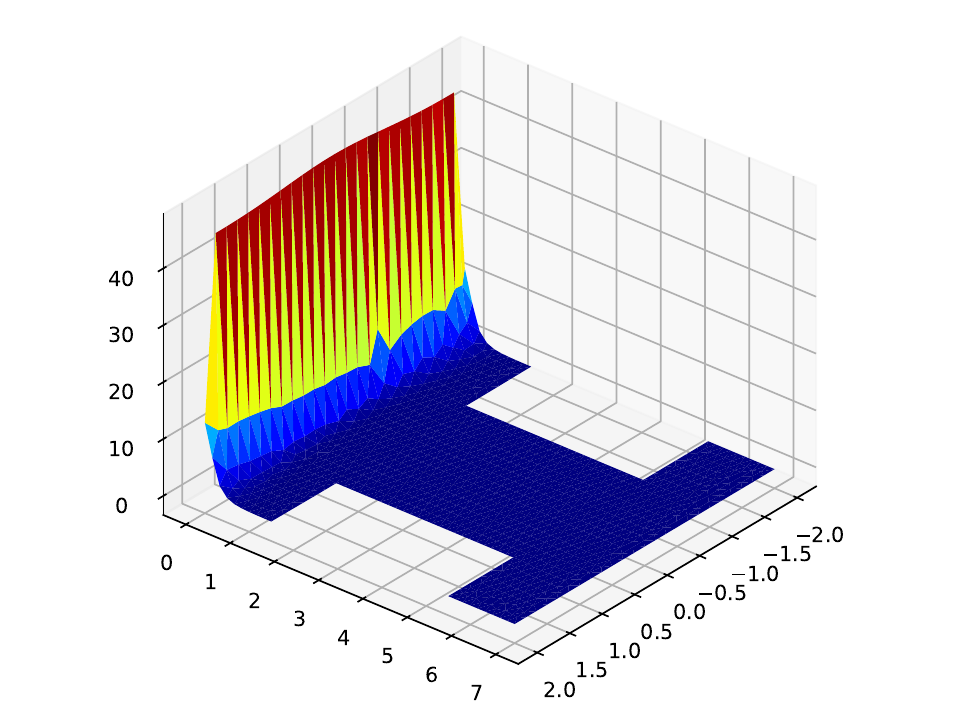}
    \end{subfigure}    
    \caption*{Snapshots of $p_h$ at times $t=0.01$, $0.05$, $0.1$ and $1.0$}
    \caption{Algorithm 2}
    \label{fig_channel:snapshots_alg2}
\end{figure}

As initial conditions for the positive and negative ions, it is assumed a homogeneous distribution, i. e.,  
$$
p_0=1\quad  \mbox{ and } \quad n_0=1\quad \mbox{ in } \quad \Omega.       
$$
It is obvious that the electroneutrality of the initial condition holds and is preserved in time throughout the  simulation as indicated in Figure \ref{fig_channel:mass_and_energies}. As for the energy and entropy of our system, they start becoming constant (Figure~\ref{fig_channel:mass_and_energies}) as both positive and negative ion maxima (Figure \ref{fig_channel:max_and_min}) are very close to its highest value; indicating that the system reaches equilibrium with the accumulation on $\partial\Omega_b$ and $\partial\Omega_t$ for the positive and negative ions, respectively, as indicated in Figures \ref{fig_channel:snapshots_alg1} and \ref{fig_channel:snapshots_alg2}. It is worthwhile noting that the electric potential barely changes. Additionally, the profile at the boundary where the ions accumulates is less regular for Algorithm 1 than for Algorithm 2, but Algorithm 1 reaches higher values. Minima evolve until practically zero; see Figure \ref{fig_channel:max_and_min}.  The evolution of the discrete solutions is given in Figures  \ref{fig_channel:snapshots_alg1} and \ref{fig_channel:snapshots_alg2}  at times at times $t=0.01$, $0.05$, $0.1$ and $1.0$. 

\subsubsection{Wave-like structures} 
\begin{figure}
    \begin{subfigure}[b]{0.25\textwidth}
        \centering
        \includegraphics[width=1.0\textwidth]{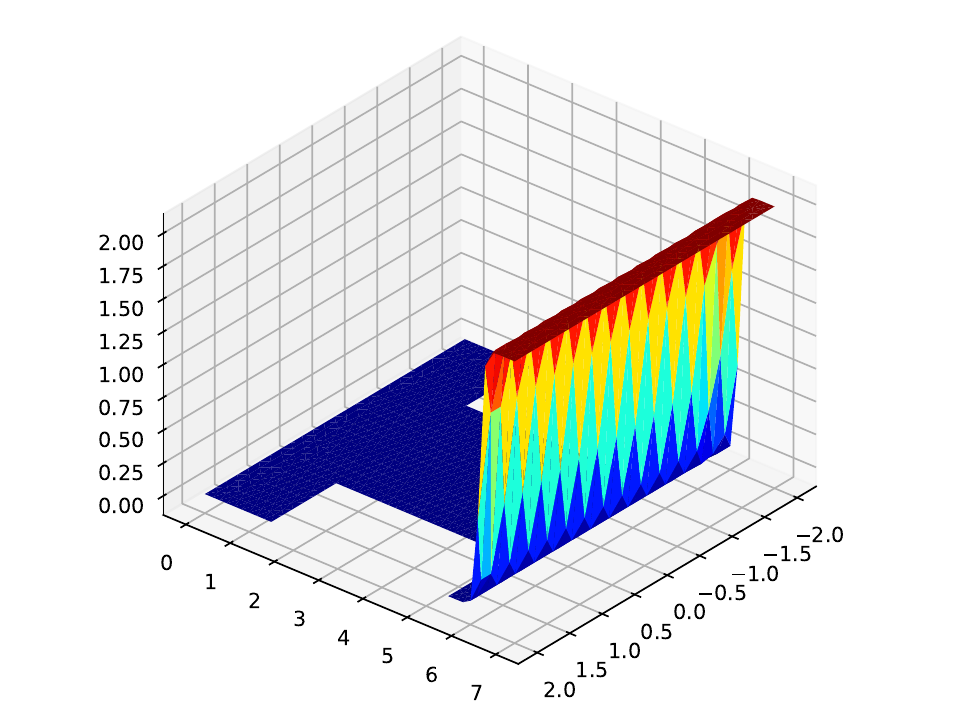}
    \end{subfigure}
    \begin{subfigure}[b]{0.25\textwidth}
        \centering
        \includegraphics[width=1.0\textwidth]{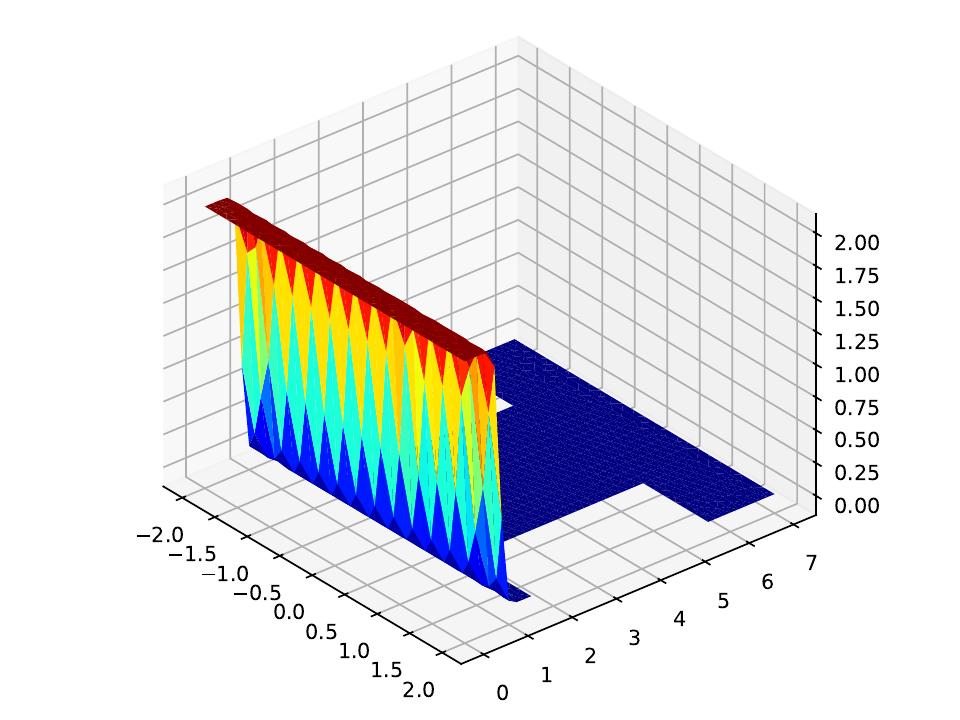}
    \end{subfigure}
    \begin{subfigure}[b]{0.25\textwidth}
\centering
\includegraphics[width=1.0\textwidth]{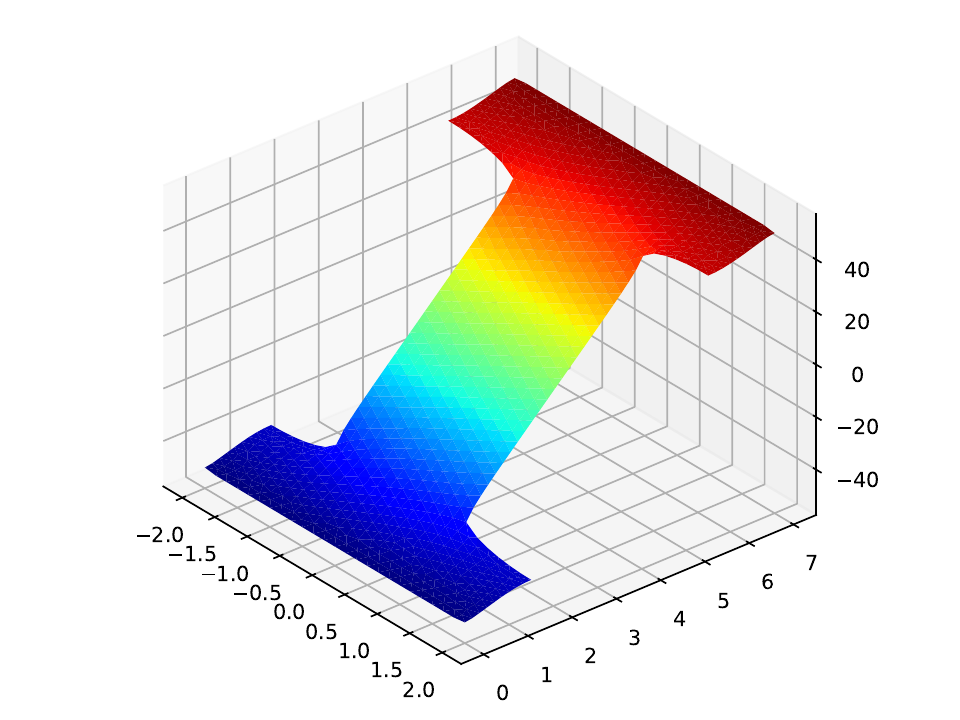}
\end{subfigure}
\caption{Initial conditions: $p_{0h}$, $n_{0h}$, and $\phi_{0h}$.}
\label{fig_channel_wave:initial_data}
\end{figure}

We are now interested in two initial conditions that, unlike the previous initial condition for ions, are concentrated on the boundaries $\Omega_t$ and $\Omega_b$ for both ions, respectively. With these initial conditions, we seek to simulate a diffusive travelling wave through the channel. To achieve this, we set: 
$$
p_0(x,y)=\tanh(10y -6.2)+1
$$
and
$$
n_0(x,y)=-\tanh(10y -0.8)+1.
$$
Figures  \ref{fig_channel_wave:initial_data} illustrates the nodal interpolations for $p_0$, $n_0$ and $\phi_0$.   

\begin{figure}
    \begin{subfigure}[b]{0.25\textwidth}
        \centering
        \includegraphics[width=1.0\textwidth]{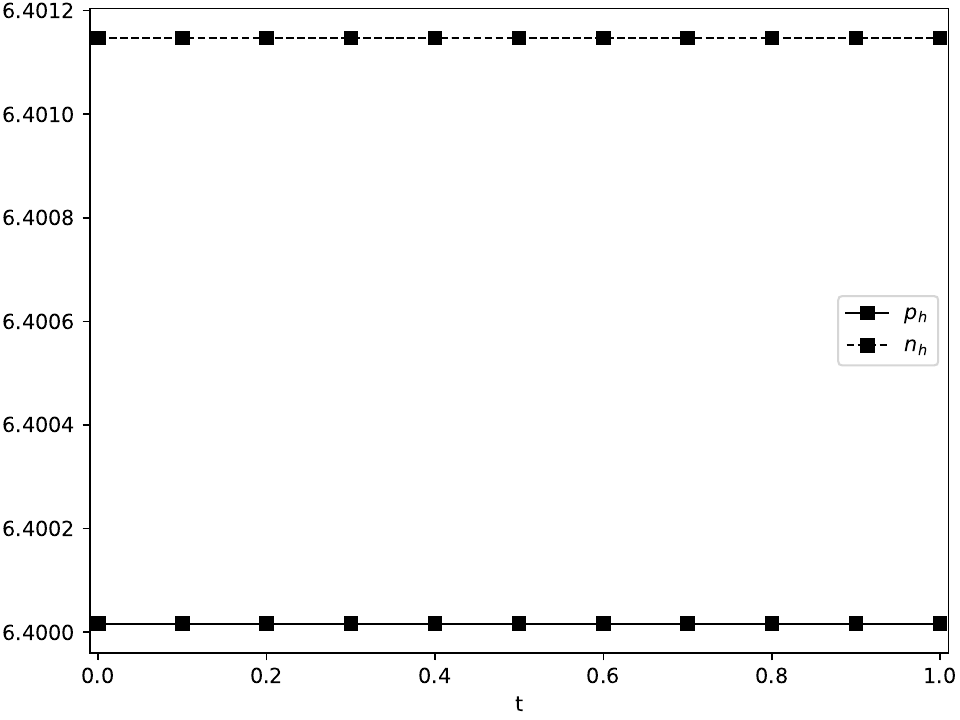}
    \end{subfigure}
    \begin{subfigure}[b]{0.25\textwidth}
        \centering
        \includegraphics[width=1.0\textwidth]{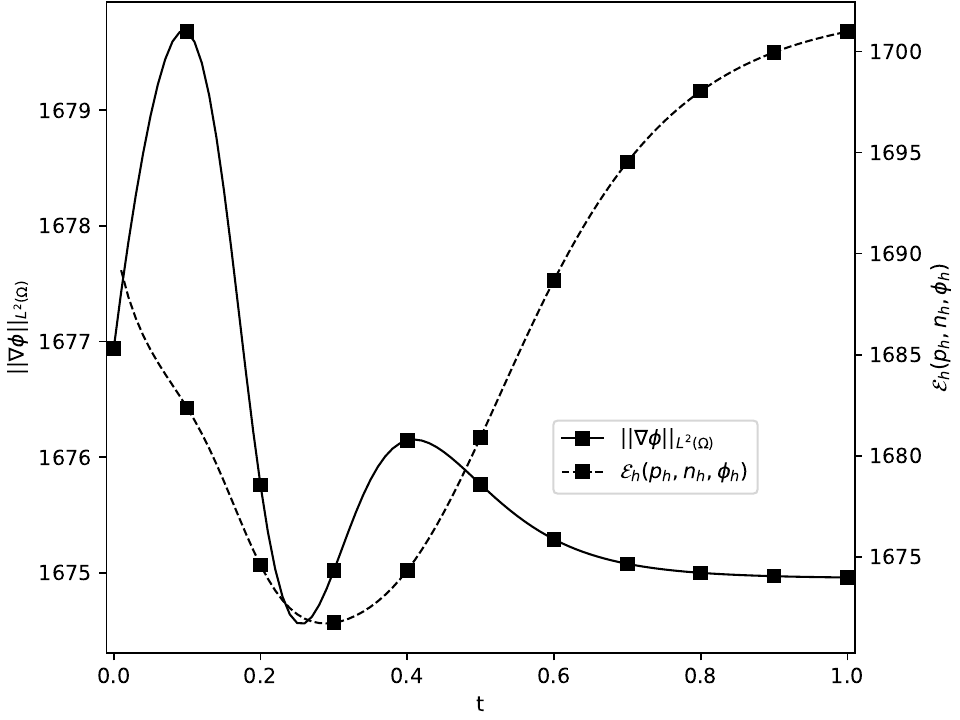}
    \end{subfigure}
    \caption*{Algorithm 1}
    \begin{subfigure}[b]{0.25\textwidth}
        \centering
        \includegraphics[width=1.0\textwidth]{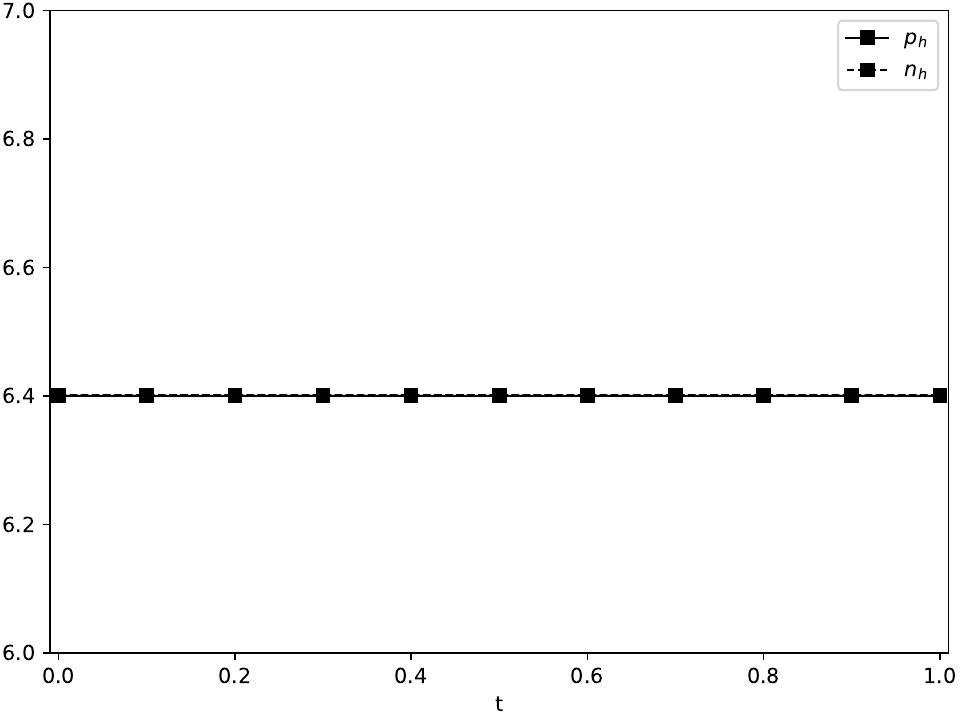}
    \end{subfigure}
    \begin{subfigure}[b]{0.25\textwidth}
        \centering
        \includegraphics[width=1.0\textwidth]{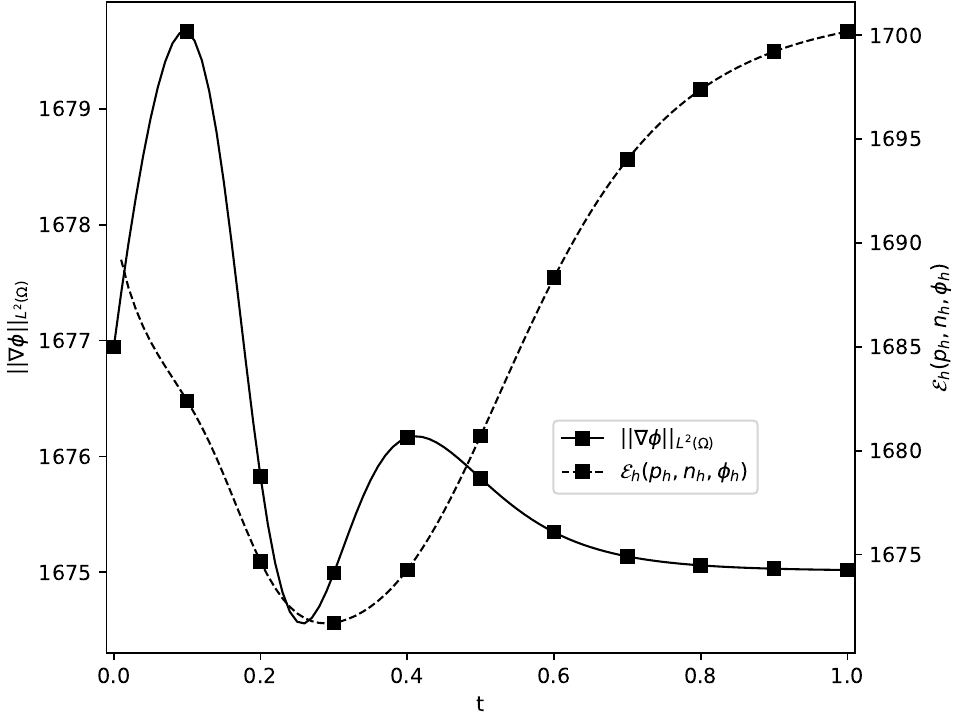}
    \end{subfigure}
    \caption*{Algorithm 2}
    \caption{Mass conservation (left), and energy and entropy evolutions (right)}
    \label{fig_channel_wave:mass_and_energies}
\end{figure}

\begin{figure}
    \begin{subfigure}[b]{0.25\textwidth}
        \centering
        \includegraphics[width=1.0\textwidth]{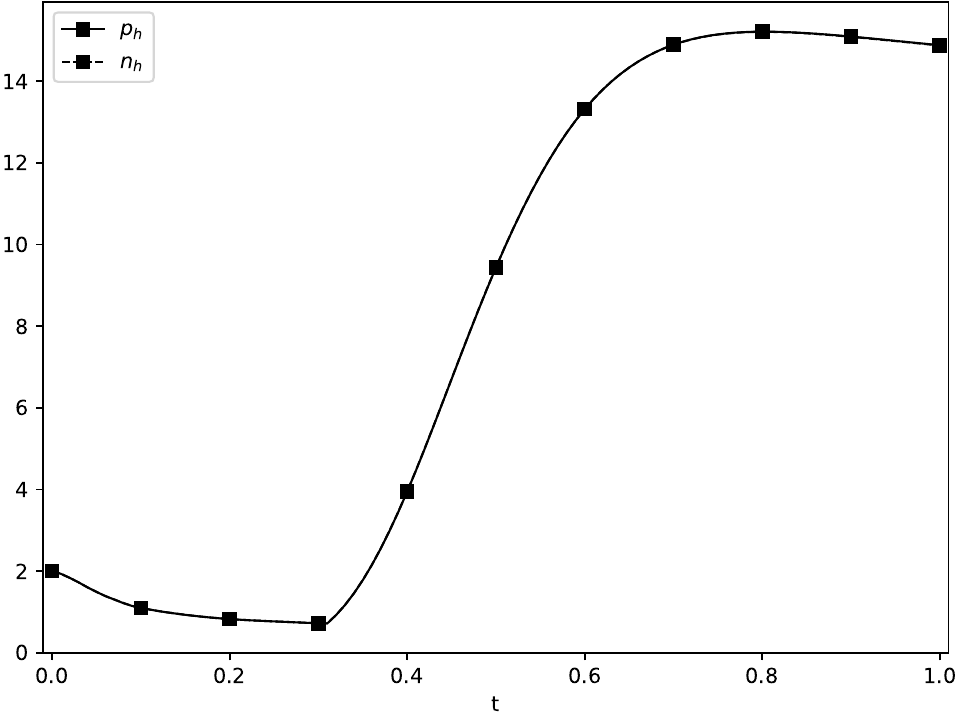}
     \end{subfigure}
    \begin{subfigure}[b]{0.25\textwidth}
        \centering
        \includegraphics[width=1.0\textwidth]{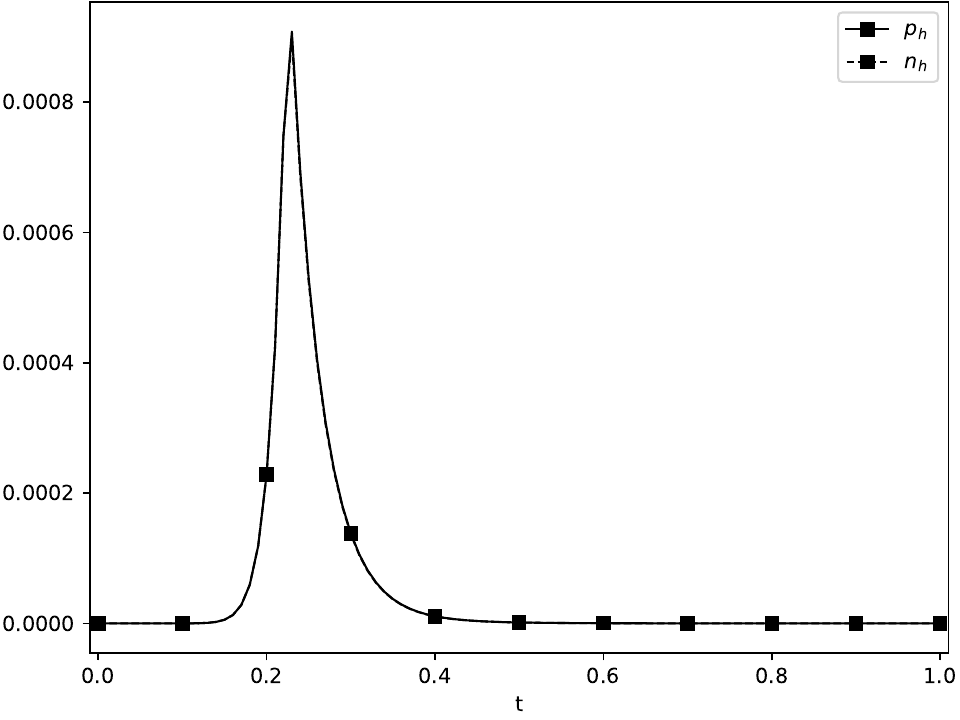}
     \end{subfigure}
  \caption*{Algorithm 1}
       \begin{subfigure}[b]{0.25\textwidth}
        \centering
        \includegraphics[width=1.0\textwidth]{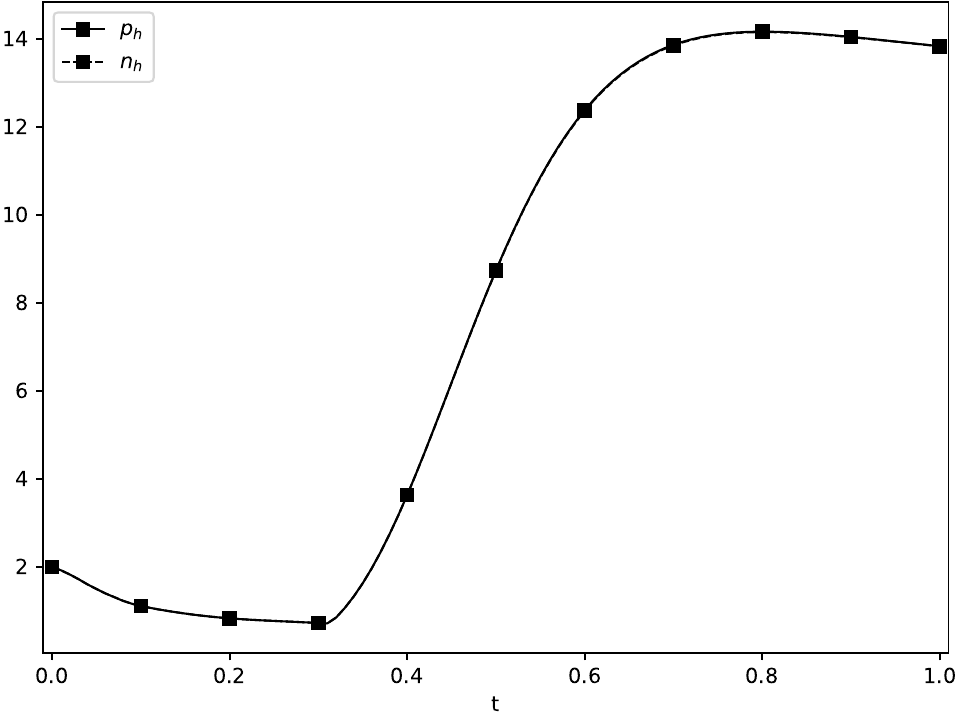}
    \end{subfigure}
    \begin{subfigure}[b]{0.25\textwidth}
        \centering
        \includegraphics[width=1.0\textwidth]{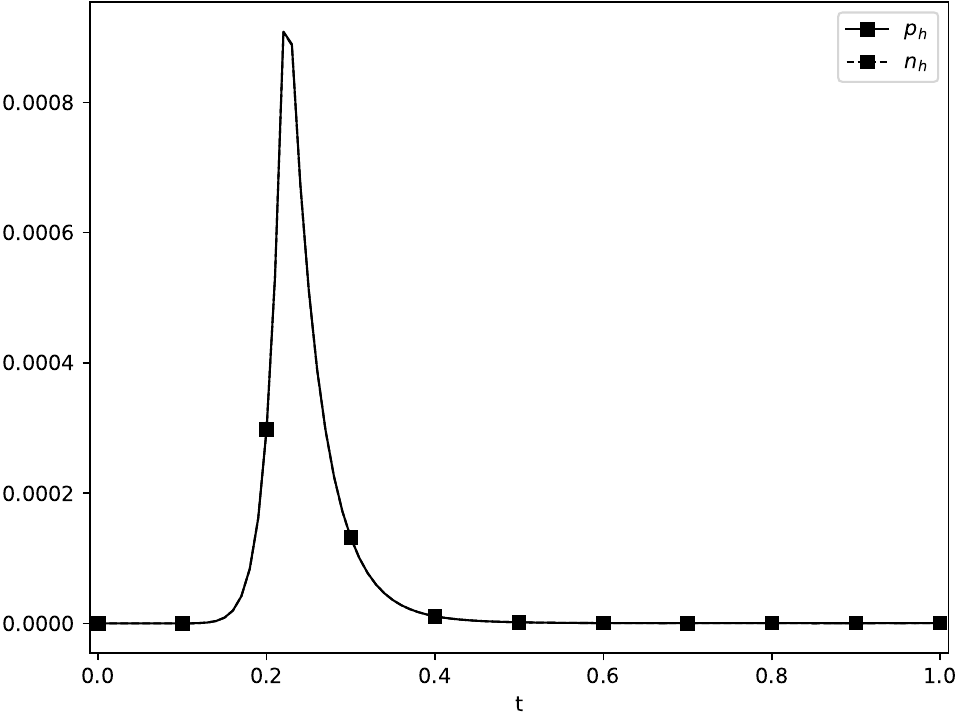}
    \end{subfigure}    
    \caption*{Algorithm 2}
    \caption{Maxima (left) and minima (right)}
    \label{fig_channel_wave:max_and_min}
\end{figure}

\begin{figure} 
\begin{subfigure}[b]{0.22\textwidth}
        \centering
        \includegraphics[width=1.0\textwidth]{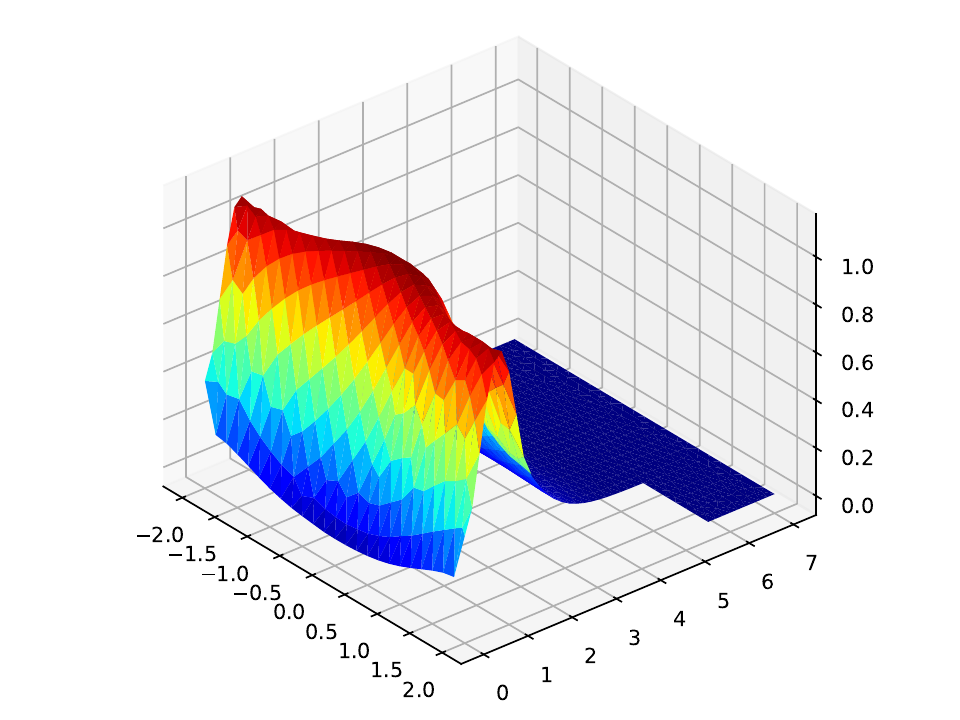}
    \end{subfigure}
    \begin{subfigure}[b]{0.22\textwidth}
        \centering
        \includegraphics[width=1.0\textwidth]{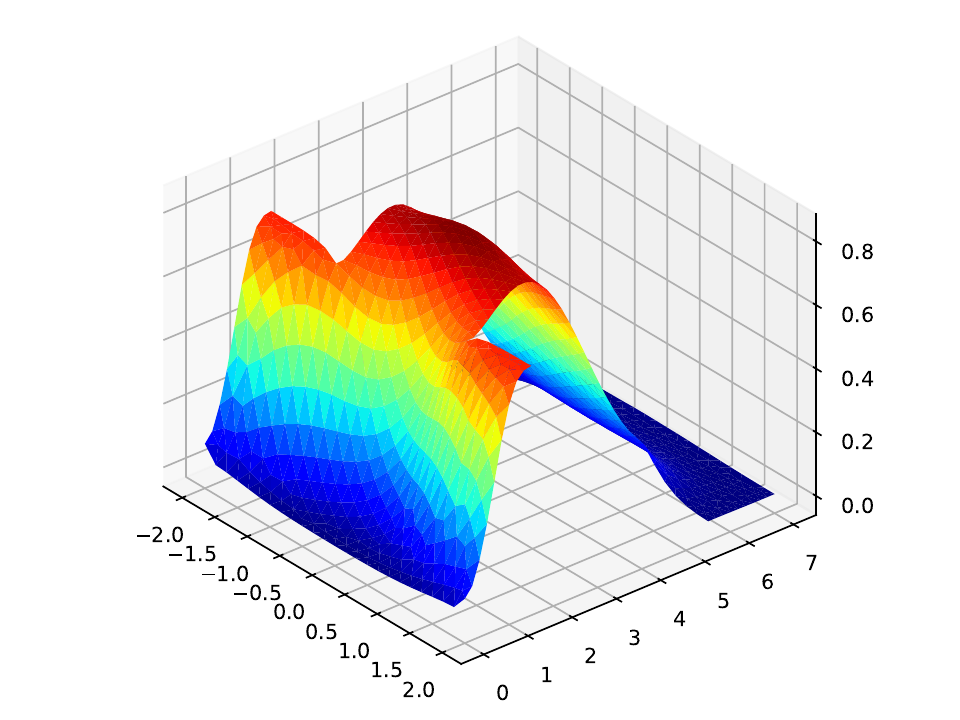}
    \end{subfigure}
    \begin{subfigure}[b]{0.22\textwidth}
        \centering
        \includegraphics[width=1.0\textwidth]{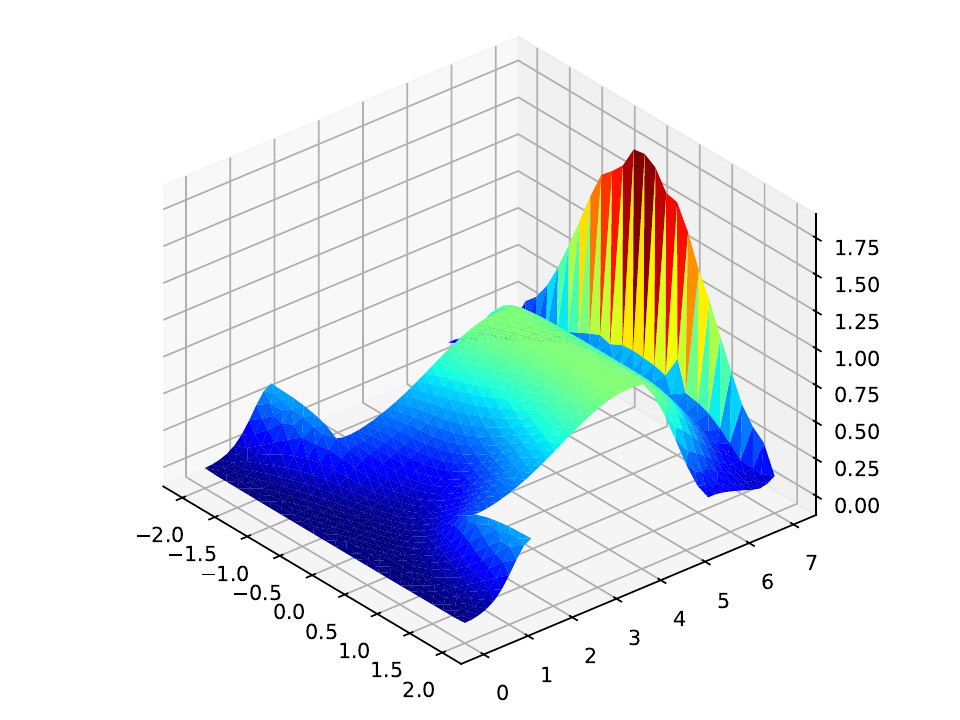}
    \end{subfigure}
    \begin{subfigure}[b]{0.22\textwidth}
        \centering
        \includegraphics[width=1.0\textwidth]{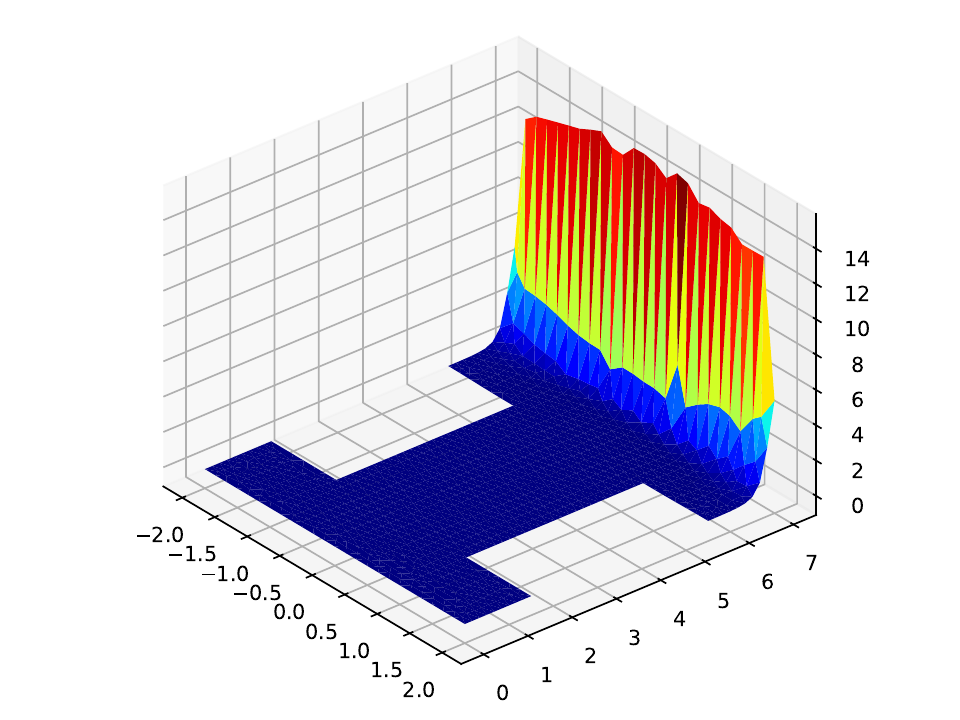}
    \end{subfigure}    
    \caption*{Snapshots of $n_h$ at times $t=0.1$, $0.20$, $0.35$ and $1.0$}
    \begin{subfigure}[b]{0.22\textwidth}
        \centering
        \includegraphics[width=1.0\textwidth]{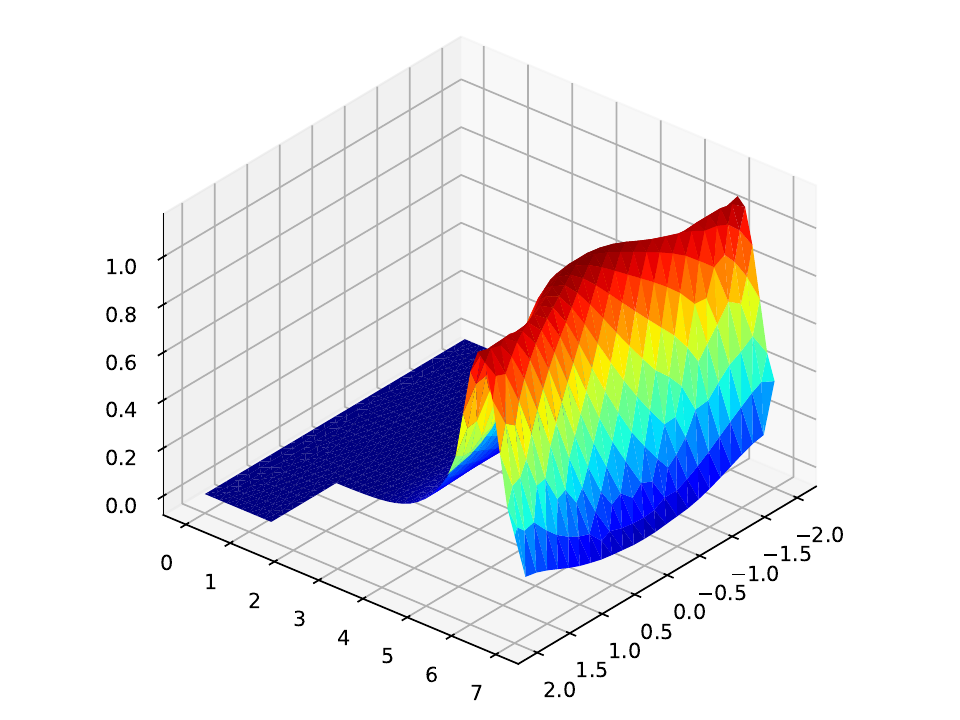}
    \end{subfigure}
    \begin{subfigure}[b]{0.22\textwidth}
        \centering
        \includegraphics[width=1.0\textwidth]{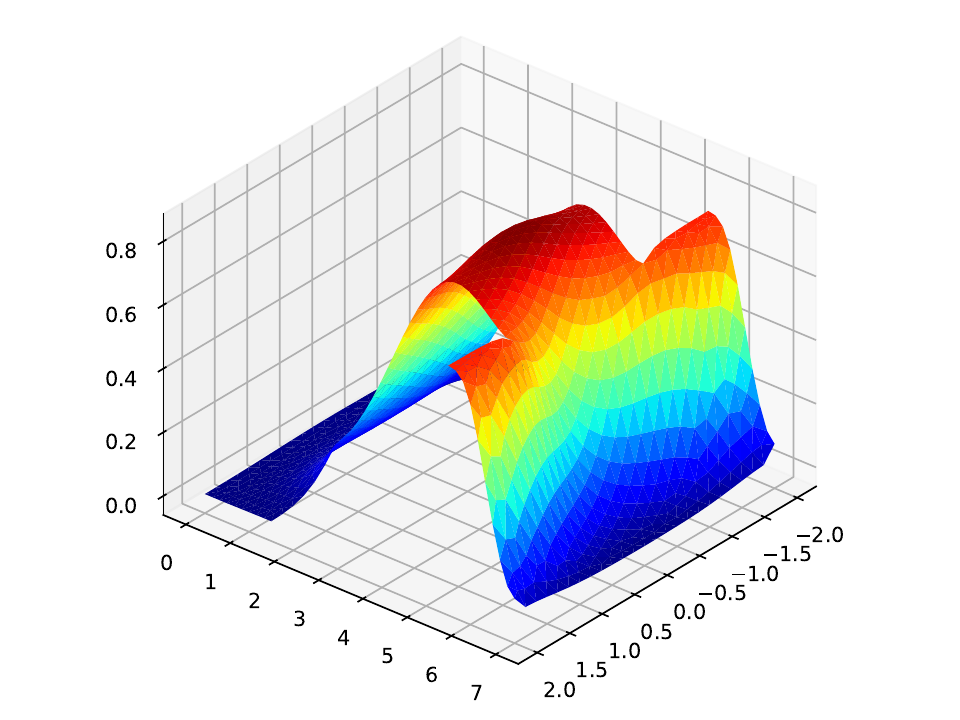}
    \end{subfigure}
    \begin{subfigure}[b]{0.22\textwidth}
        \centering
        \includegraphics[width=1.0\textwidth]{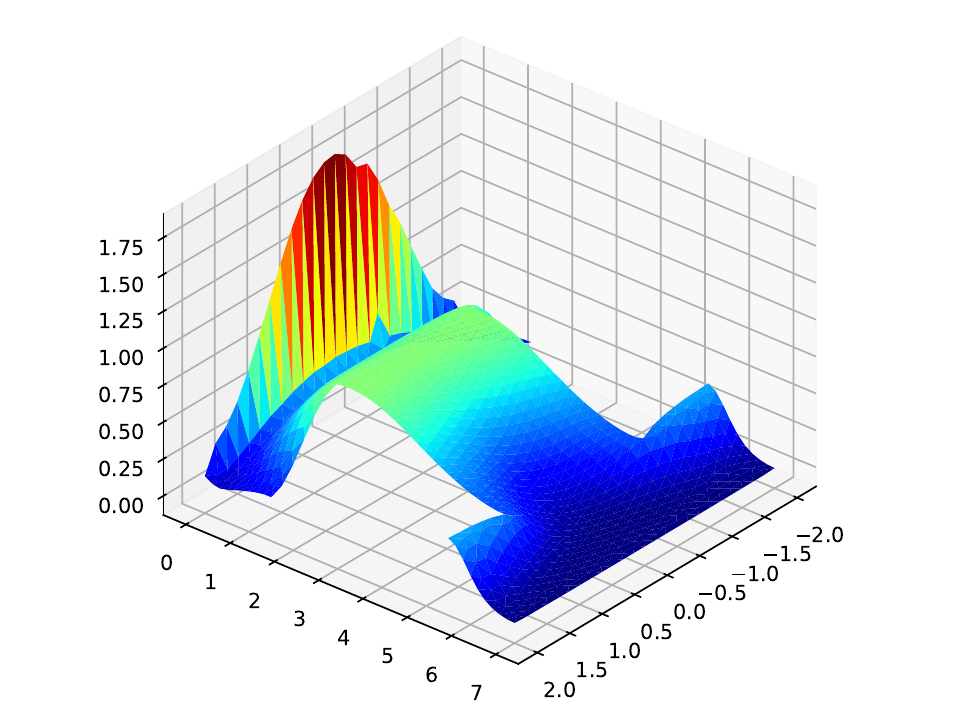}
    \end{subfigure}
    \begin{subfigure}[b]{0.22\textwidth}
        \centering
        \includegraphics[width=1.0\textwidth]{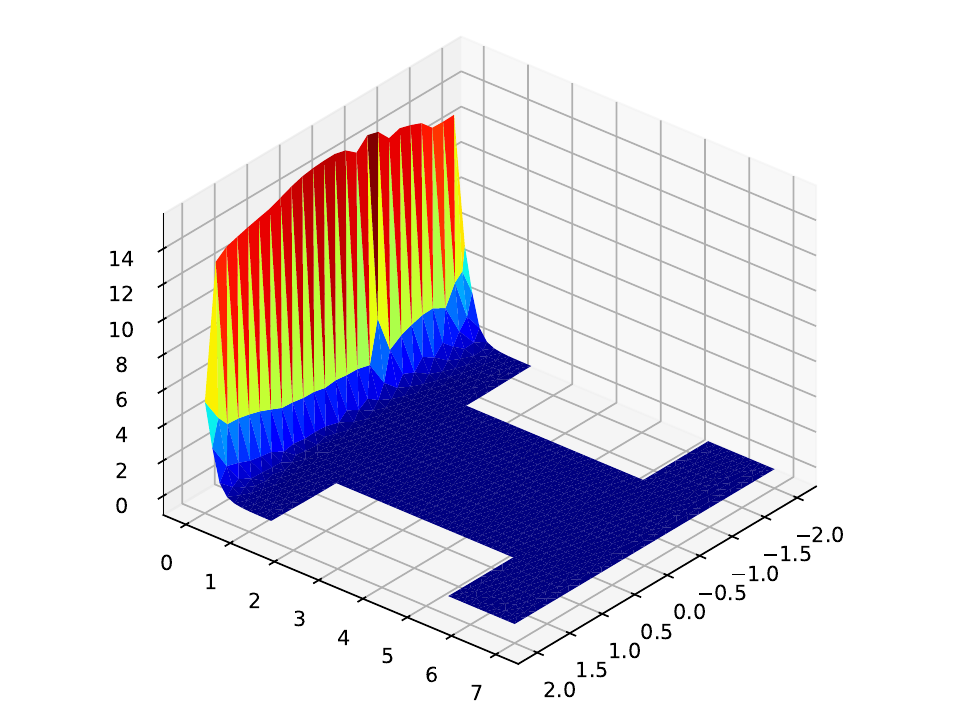}
    \end{subfigure}    
    \caption*{Snapshots of $p_h$ at times $t=0.1$, $0.2$, $0.35$ and $1.0$}
    \caption{Algorithm 1}
    \label{fig_channel_wave:snapshots_alg1}
\end{figure}

\begin{figure}
\begin{subfigure}[b]{0.22\textwidth}
        \centering
        \includegraphics[width=1.0\textwidth]{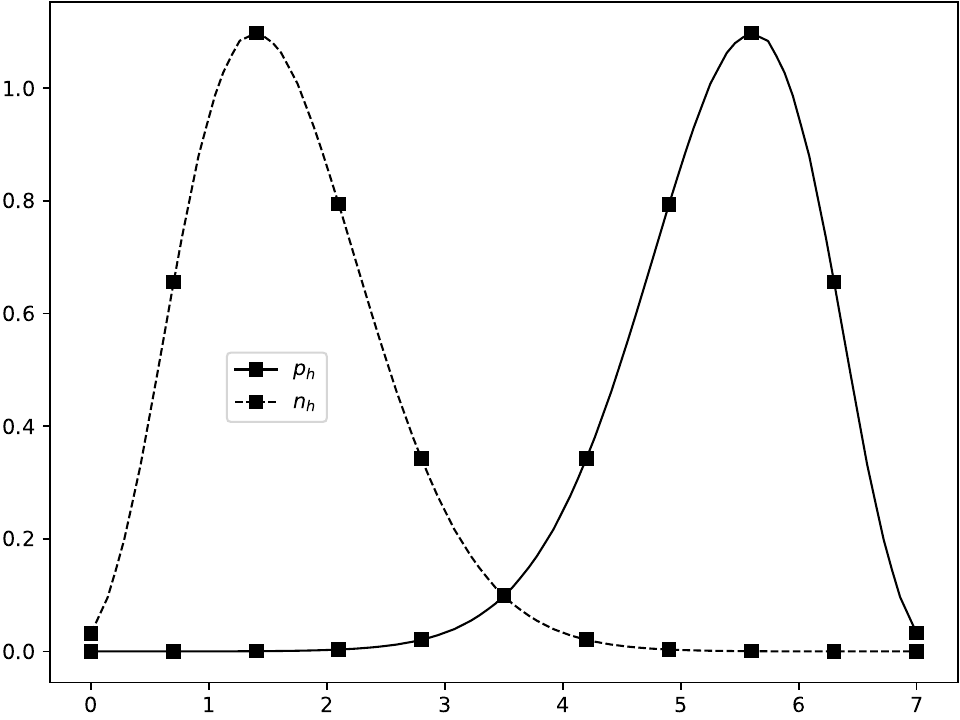}
    \end{subfigure}
    \begin{subfigure}[b]{0.22\textwidth}
        \centering
        \includegraphics[width=1.0\textwidth]{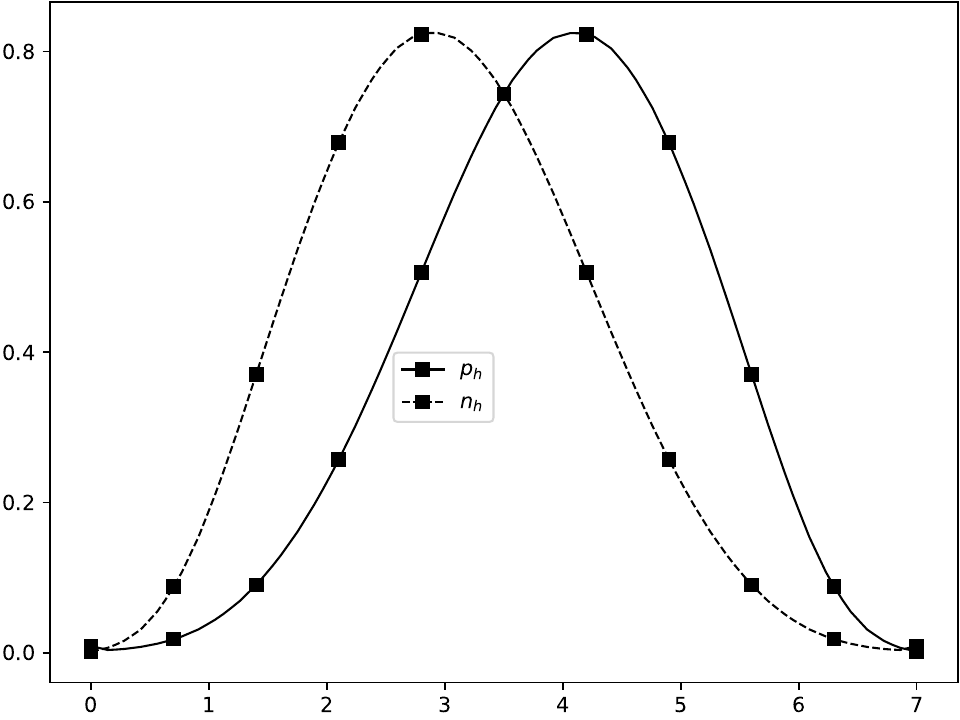}
    \end{subfigure}
    \begin{subfigure}[b]{0.22\textwidth}
        \centering
        \includegraphics[width=1.0\textwidth]{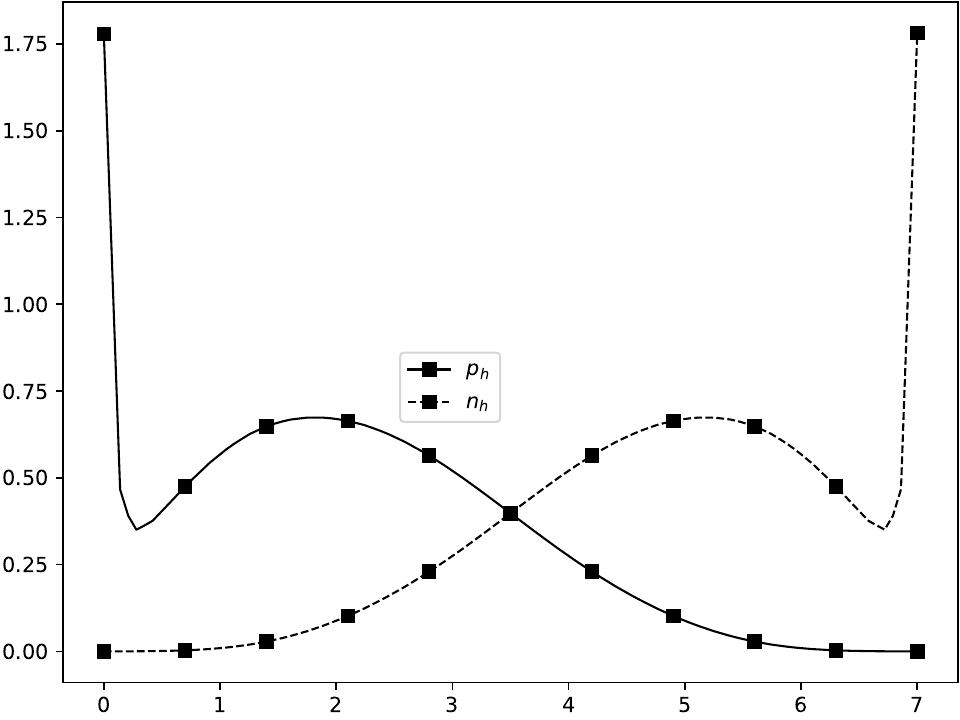}
    \end{subfigure}
    \begin{subfigure}[b]{0.22\textwidth}
        \centering
        \includegraphics[width=1.0\textwidth]{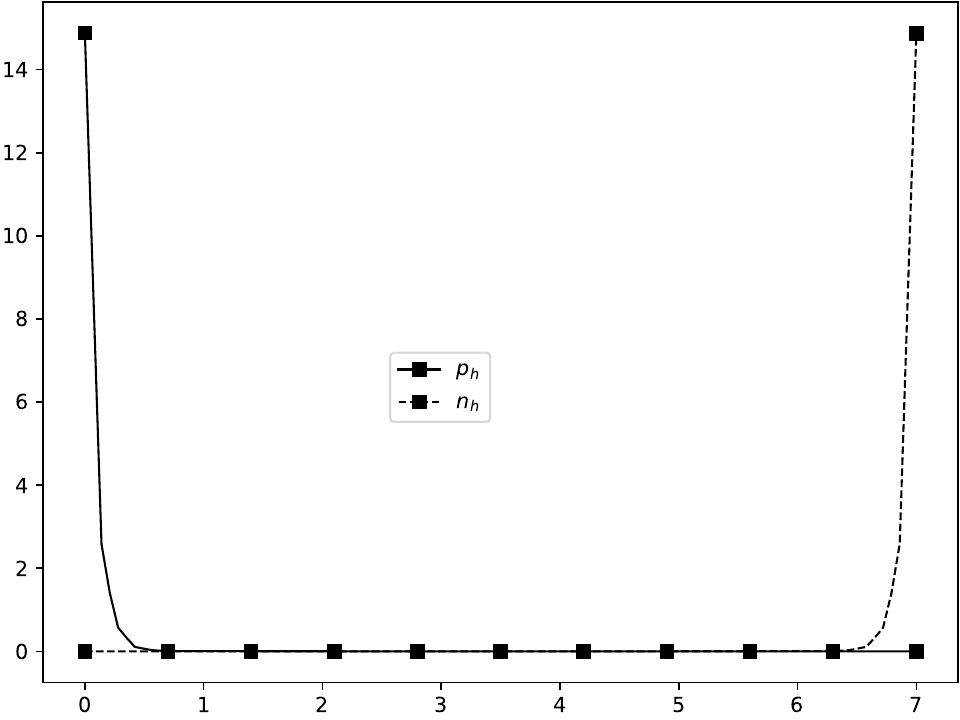}
    \end{subfigure}    
    \caption*{Profiles of $n_h$ and $p_h$ at times $t=0.1$, $0.20$, $0.35$ and $1.0$}
    \caption{Algorithm 1}
    \label{fig_channel_wave:profiles_alg1}
\end{figure}

\begin{figure}    
\begin{subfigure}[b]{0.22\textwidth}
        \centering
        \includegraphics[width=1.0\textwidth]{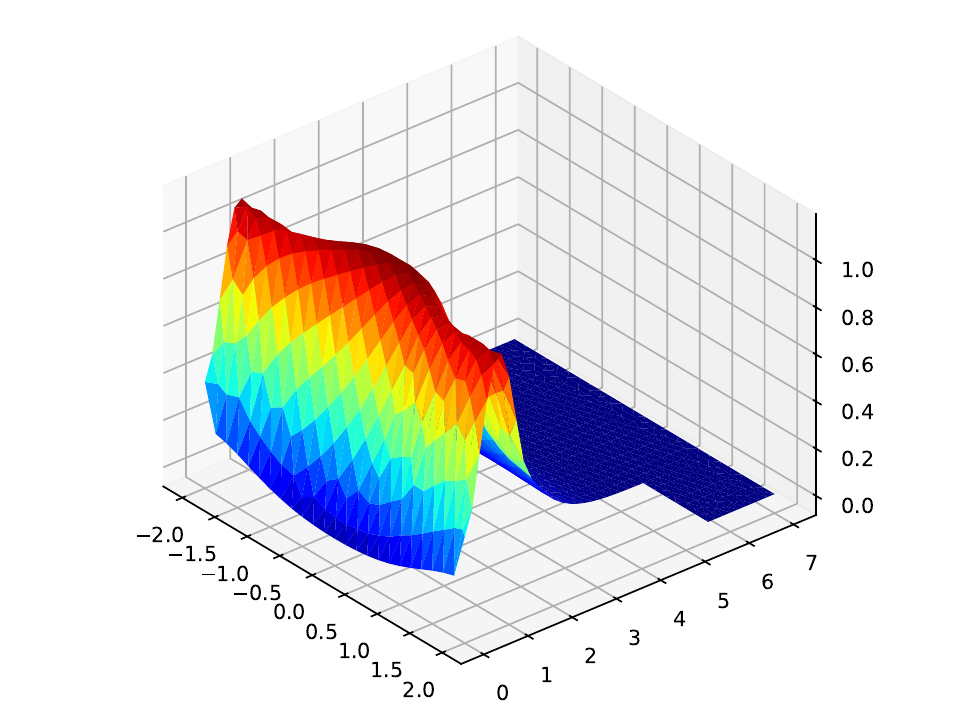}
    \end{subfigure}
    \begin{subfigure}[b]{0.22\textwidth}
        \centering
        \includegraphics[width=1.0\textwidth]{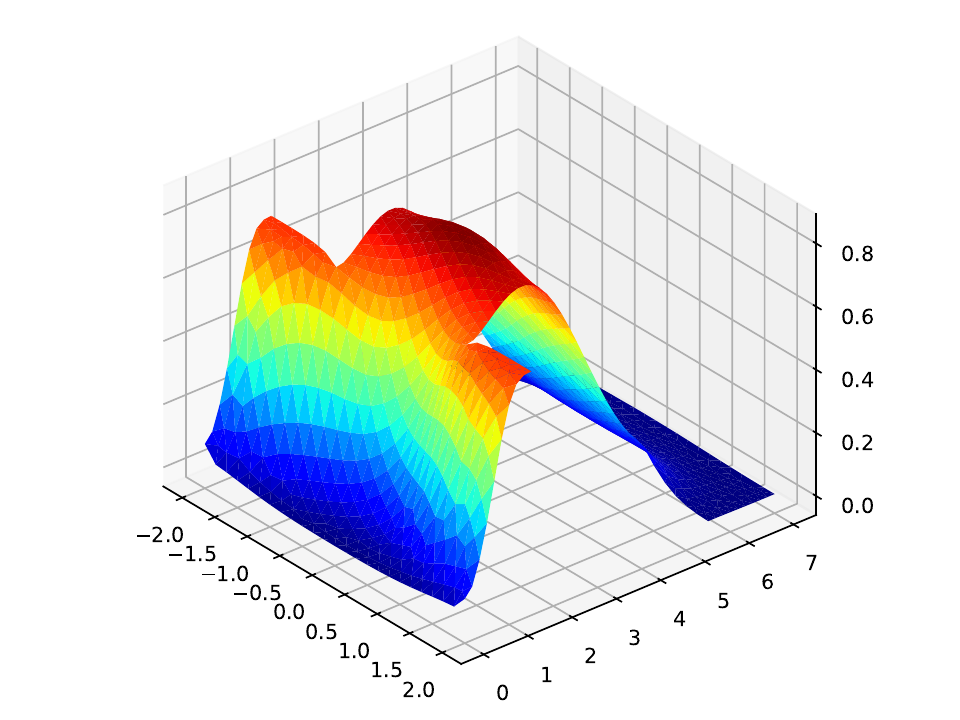}
    \end{subfigure}
    \begin{subfigure}[b]{0.22\textwidth}
        \centering
        \includegraphics[width=1.0\textwidth]{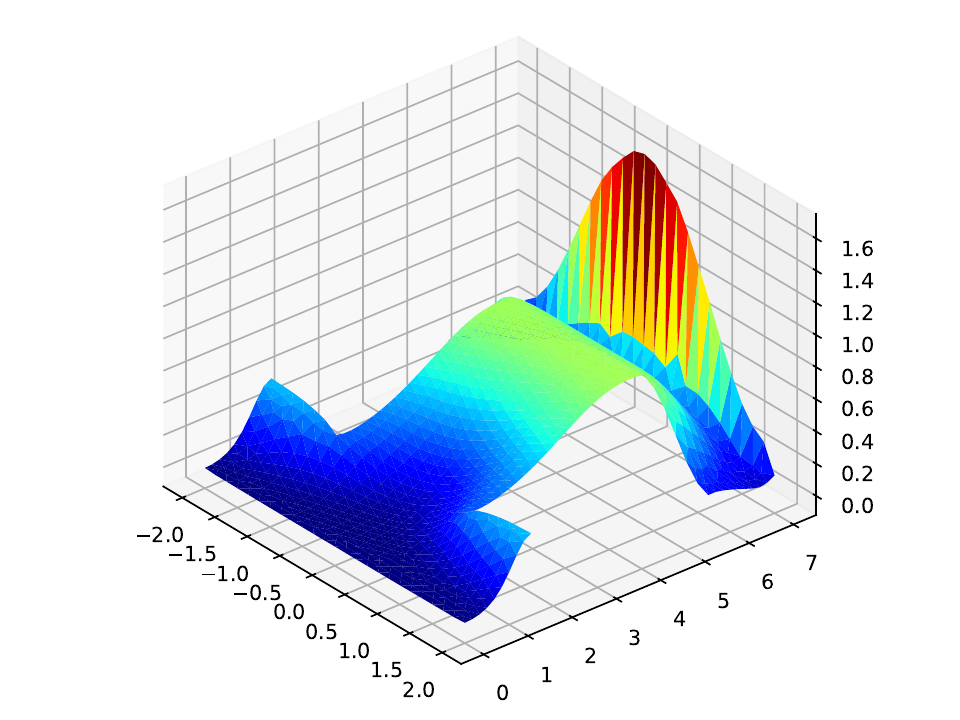}
    \end{subfigure}
    \begin{subfigure}[b]{0.22\textwidth}
        \centering
        \includegraphics[width=1.0\textwidth]{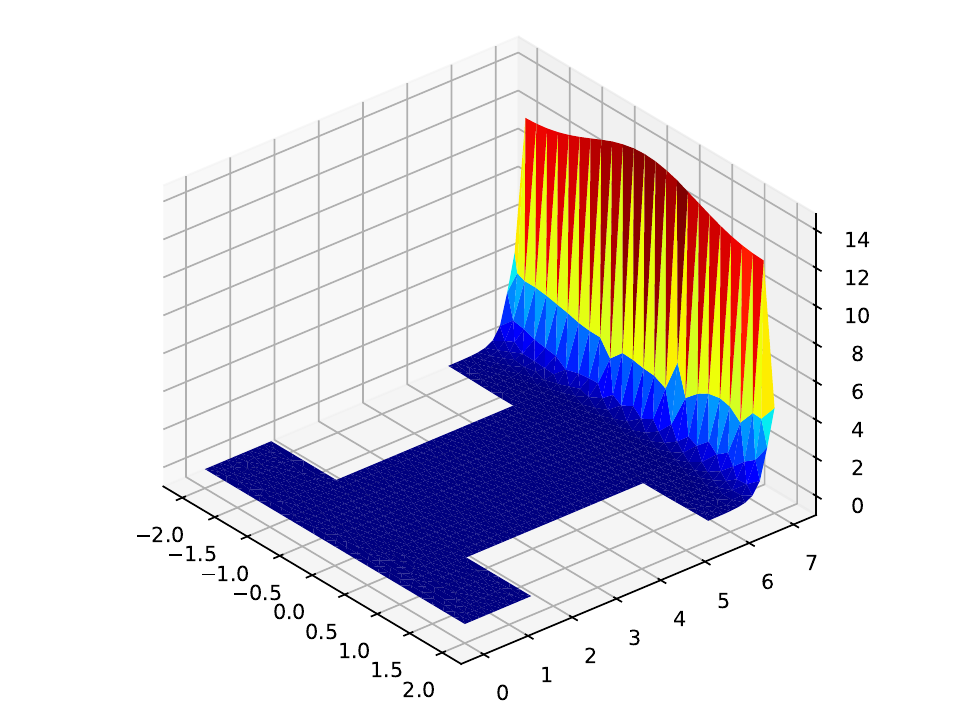}
    \end{subfigure}    
    \caption*{Algorithm 2: Snapshots of $n_h$ at times $t=0.1$, $0.2$, $0.35$ and $1.0$}
    \begin{subfigure}[b]{0.22\textwidth}
        \centering
        \includegraphics[width=1.0\textwidth]{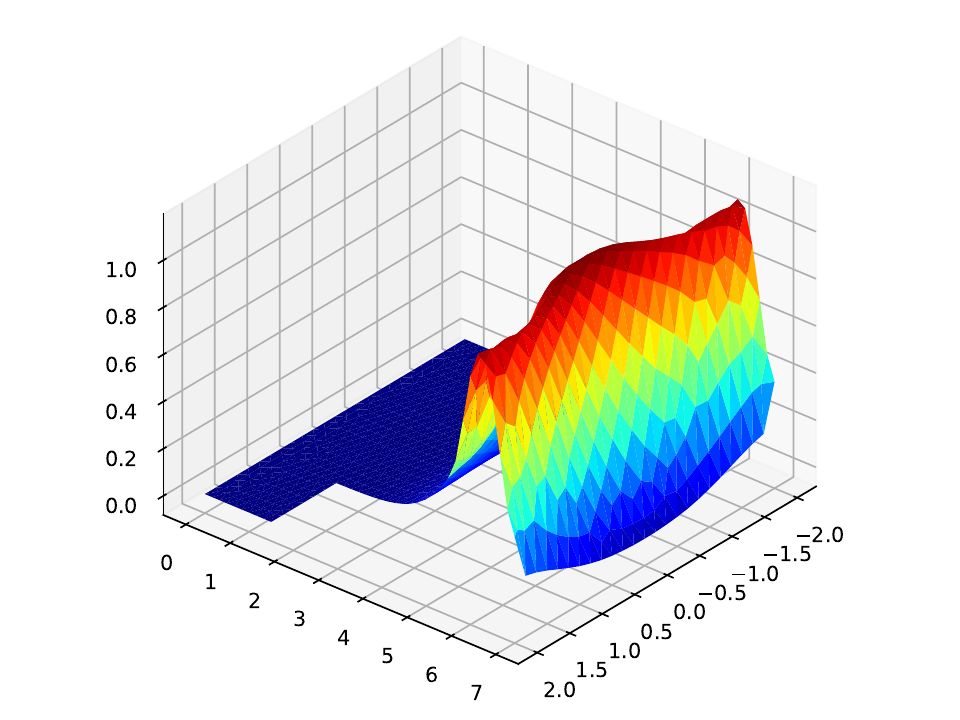}
    \end{subfigure}
    \begin{subfigure}[b]{0.22\textwidth}
        \centering
        \includegraphics[width=1.0\textwidth]{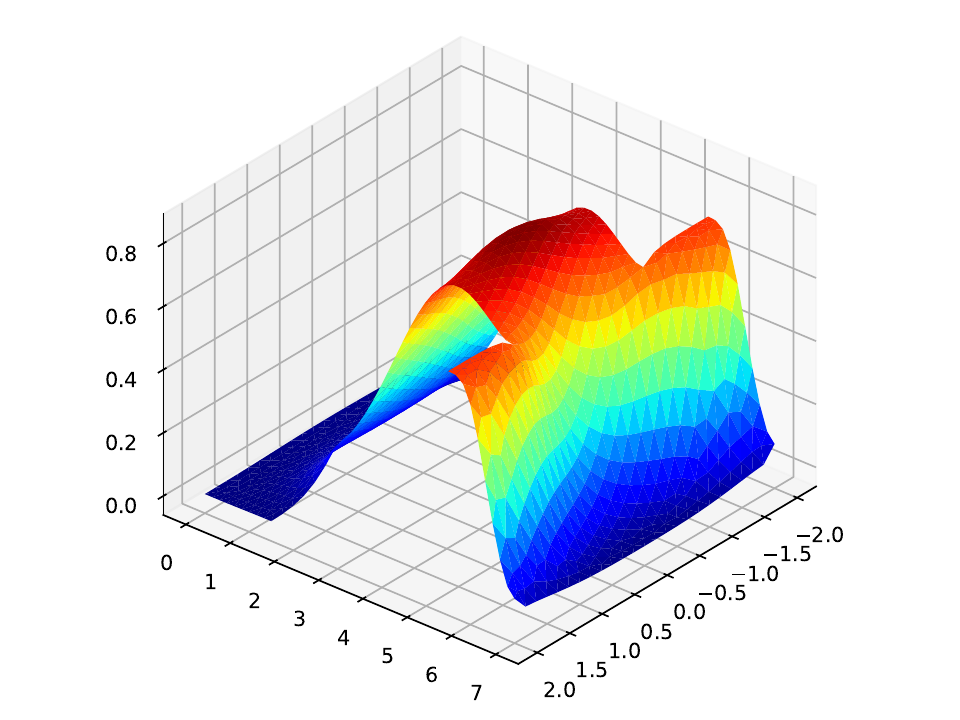}
    \end{subfigure}
    \begin{subfigure}[b]{0.22\textwidth}
        \centering
        \includegraphics[width=1.0\textwidth]{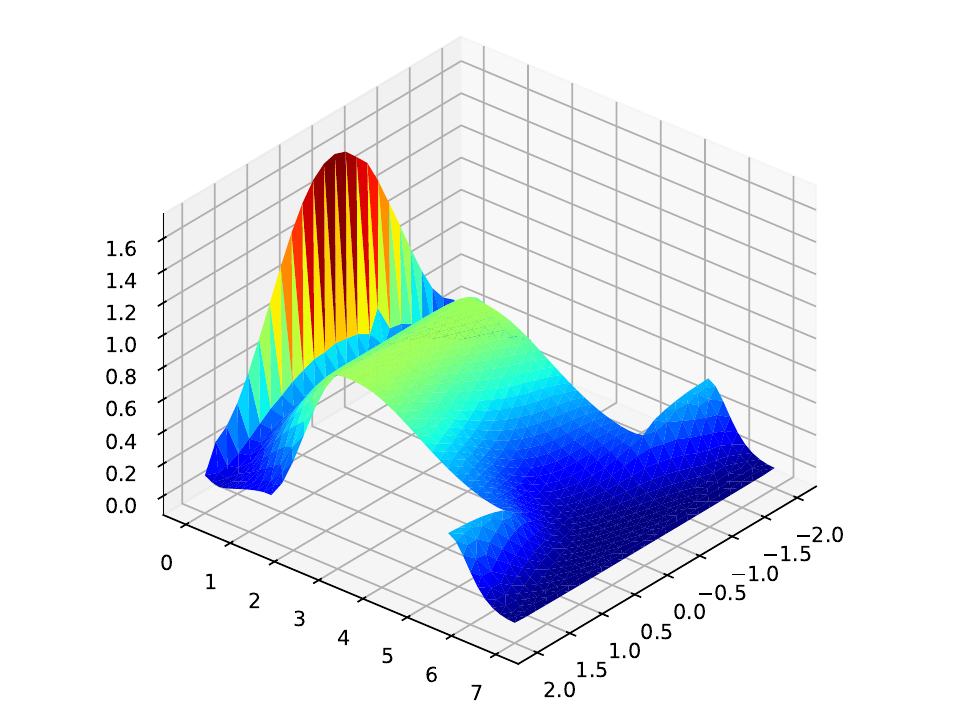}
    \end{subfigure}
    \begin{subfigure}[b]{0.22\textwidth}
        \centering
        \includegraphics[width=1.0\textwidth]{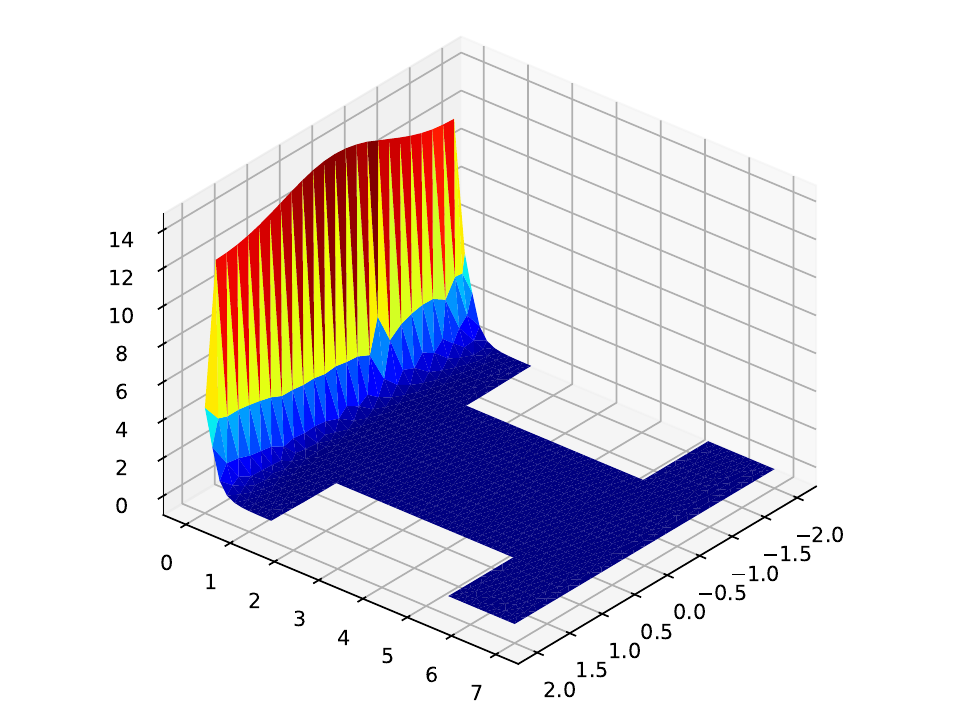}
    \end{subfigure}    
    \caption*{Algorithm 2: Snapshots of $p_h$ at times $t=0.1$, $0.2$, $0.35$ and $1.0$}
    \caption{Algorithm 2}
    \label{fig_channel_wave:snapshots_alg2}
\end{figure}

\begin{figure}
\begin{subfigure}[b]{0.22\textwidth}
        \centering
        \includegraphics[width=1.0\textwidth]{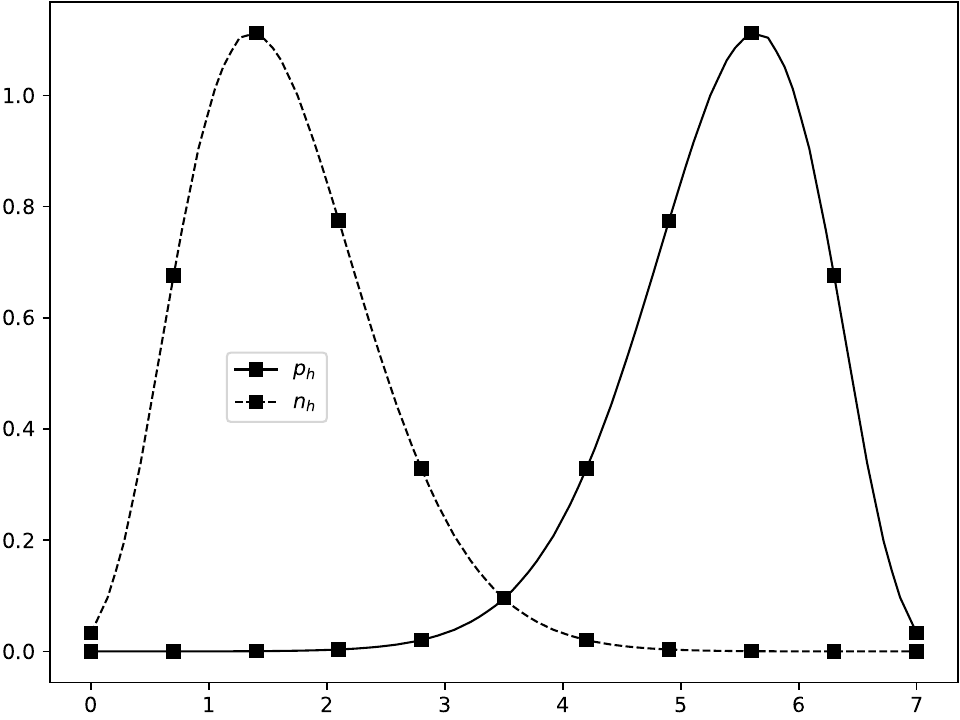}
    \end{subfigure}
    \begin{subfigure}[b]{0.22\textwidth}
        \centering
        \includegraphics[width=1.0\textwidth]{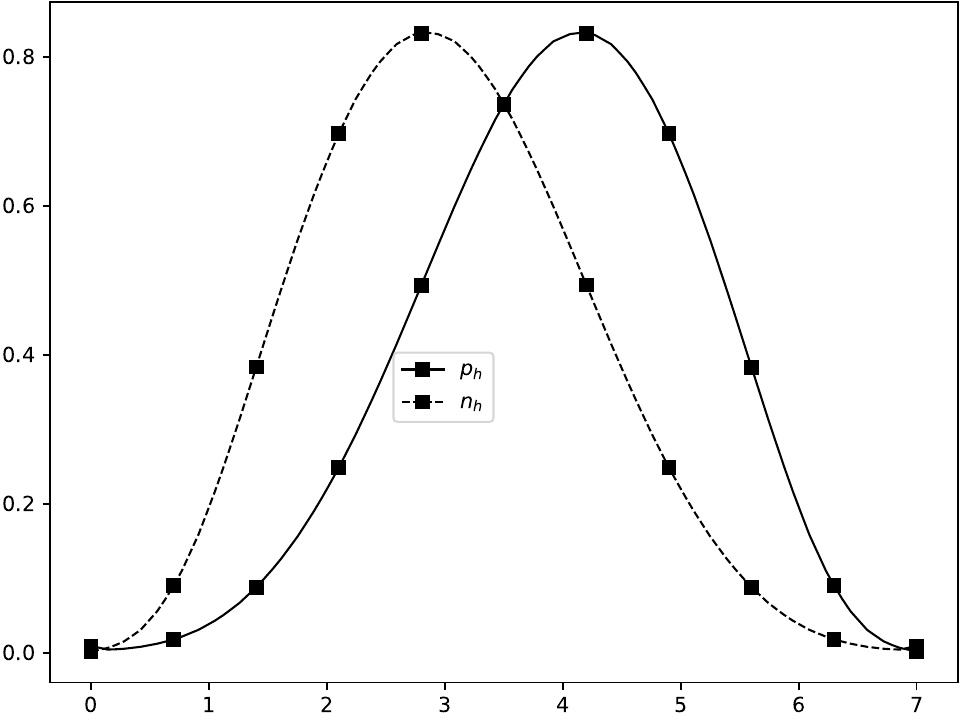}
    \end{subfigure}
    \begin{subfigure}[b]{0.22\textwidth}
        \centering
        \includegraphics[width=1.0\textwidth]{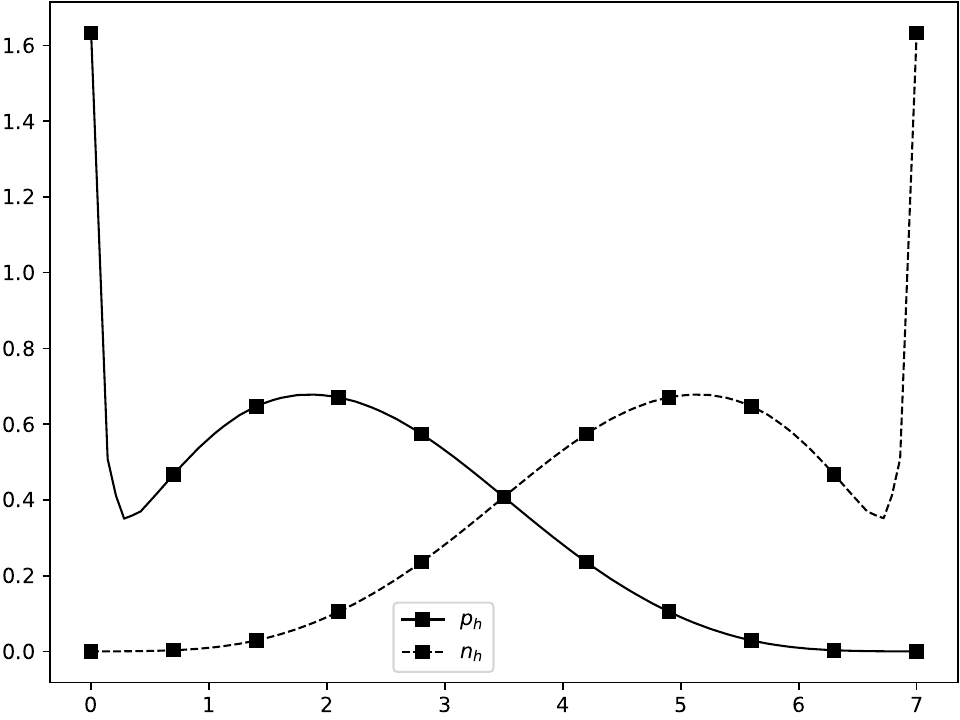}
    \end{subfigure}
    \begin{subfigure}[b]{0.22\textwidth}
        \centering
        \includegraphics[width=1.0\textwidth]{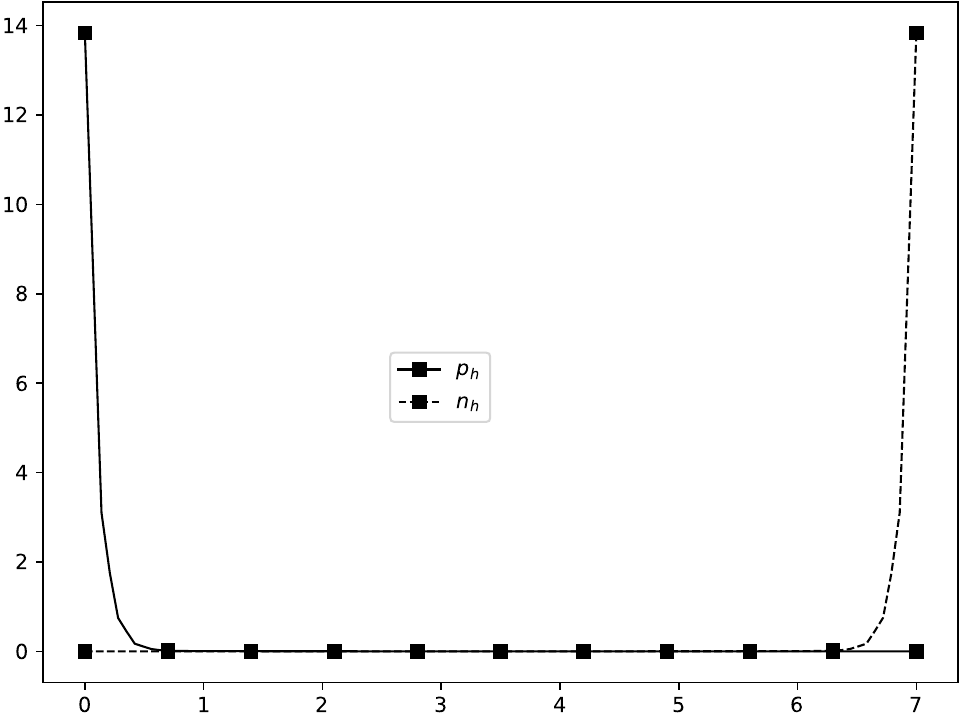}
    \end{subfigure}    
    \caption*{Profiles of $n_h$ and $p_h$ at times $t=0.1$, $0.2$, $0.35$ and $1.0$}
    \caption{Algorithm 2}
    \label{fig_channel_wave:profiles_alg2}
\end{figure}

Due to the Dirichlet boundary conditions fulfiled for the electric potential, it is not needed for the electriconeutrality to hold, but even so the initial total mass for both ions must be preserved as indicated in Figure \ref{fig_channel_wave:mass_and_energies}. Furthermore, Figure \ref{fig_channel_wave:mass_and_energies} shows that the evolution of the energy and entropy has a rich dynamics. In particular, the energy reaches its maximum value when the diffusive travelling wave enters the channel around $t=0.1$ and its minimum when the wave is in the middle of the channel  at approximately $t=0.25$. Algorithm $2$ differs from Algorithm $1$ in the fact that the dynamics is slightly delayed. Minima in Figures~\ref{fig_channel_wave:max_and_min} evolves toward a peak shape whose lowest value coincides with that of the energy and maxima reaches its minimum value prior to the wave to get out of the channel. The above-mentioned dynamics becomes evident in Figures \ref{fig_channel_wave:snapshots_alg1} and  \ref{fig_channel_wave:snapshots_alg2}. Also see Figures \ref{fig_channel_wave:profiles_alg1} and \ref{fig_channel_wave:profiles_alg2} for the profiles along the channel.                       

\subsubsection{Charge-selective transport} Finally, we investigate the impact of using Dirichlet boundary conditions for ions, allowing for the control of ion transport through charge-selective membranes. It is clear that the nature of the boundary conditions should be the source of significant changes in the transport of ions. Specifically, it is assumed that only one of both ions is subject to a selective boundary condition on $\partial\Omega_m$ and blocking boundary conditions for the rest $\Omega_n:=\partial\Omega\backslash\Omega_m$, i.e.,
$$
p= 1 \mbox{ on }\quad \partial\Omega_m\quad \mbox{ and }\quad \partial_{\boldsymbol{n}}p = 0\quad \mbox{ on }\quad\partial\Omega_n.   
$$

\begin{figure}
    \begin{subfigure}[b]{0.25\textwidth}
        \centering
        \includegraphics[width=1.0\textwidth]{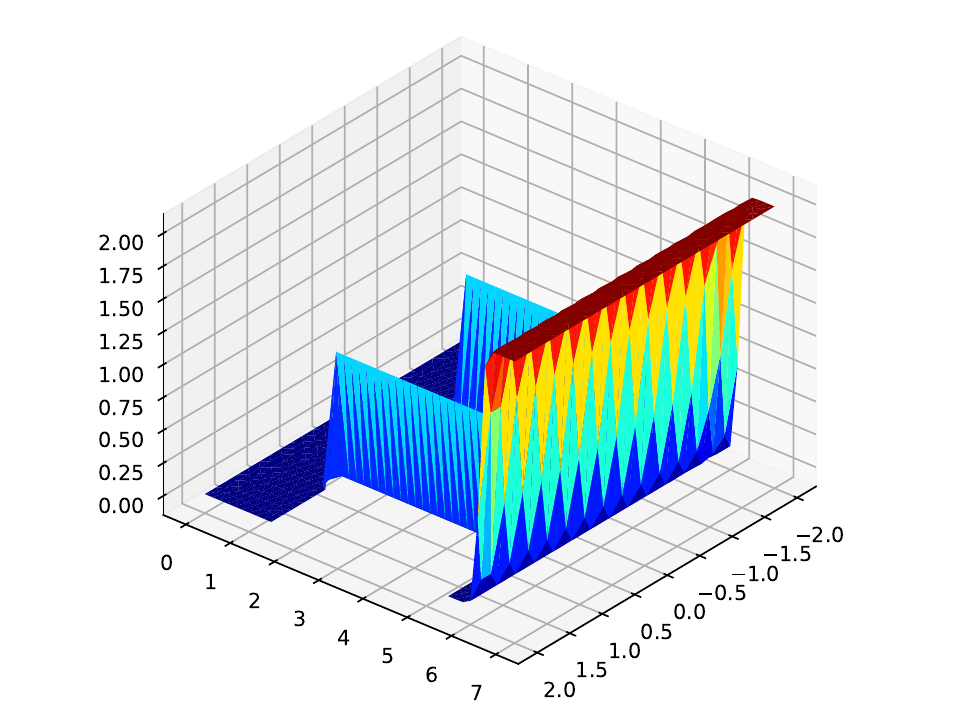}
    \end{subfigure}
    \begin{subfigure}[b]{0.25\textwidth}
        \centering
        \includegraphics[width=1.0\textwidth]{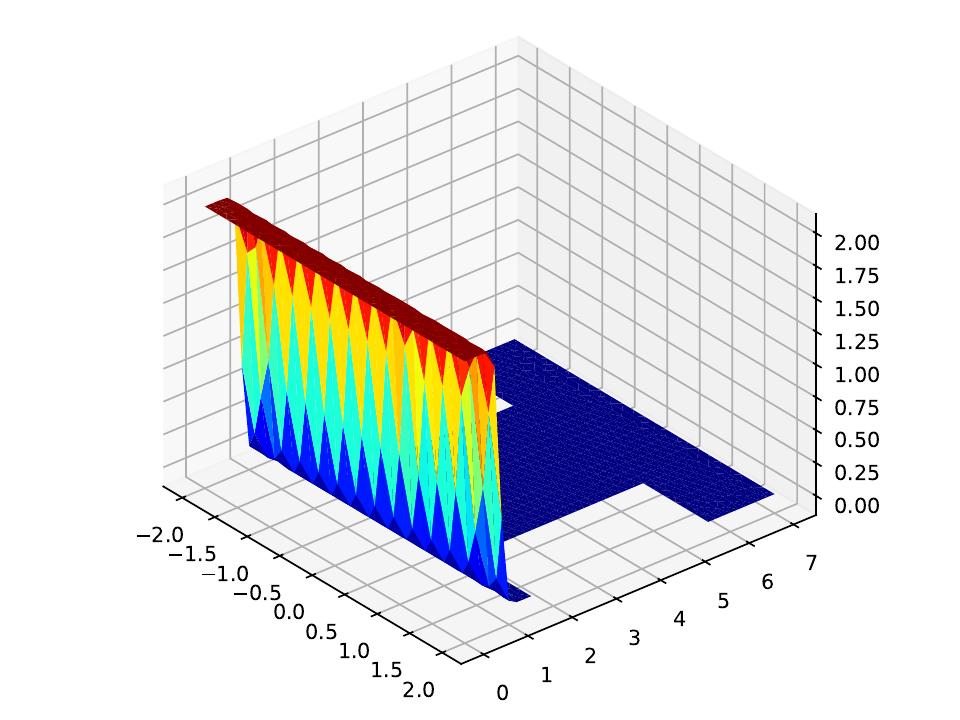}
    \end{subfigure}
    \begin{subfigure}[b]{0.25\textwidth}
\centering
\includegraphics[width=1.0\textwidth]{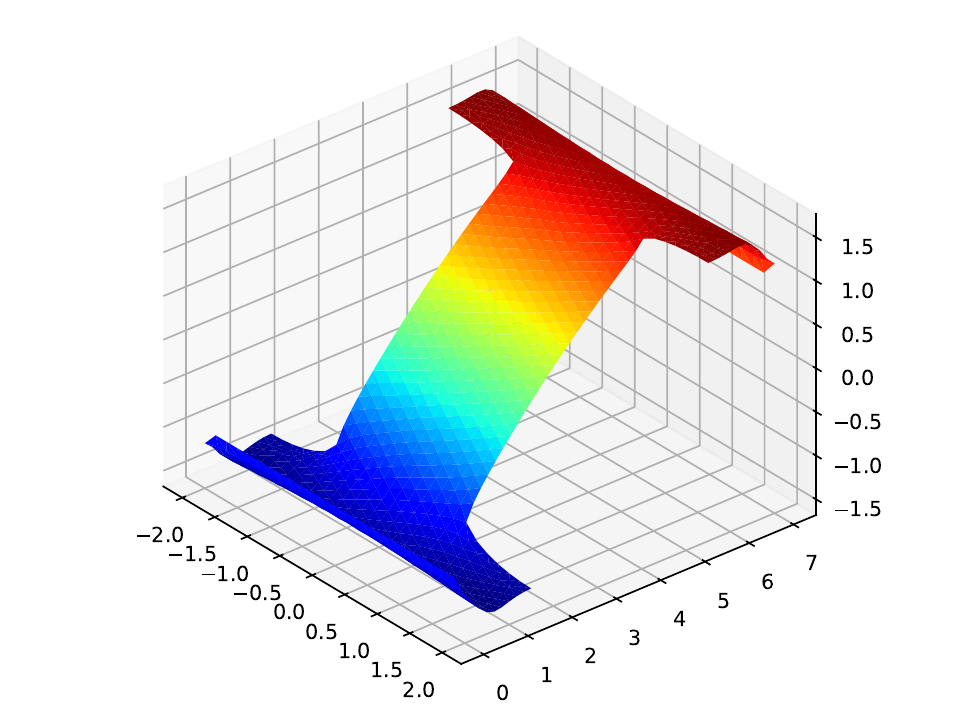}
\end{subfigure}
    \caption{Initial conditions: $p_{0h}$, $n_{0h}$ and $\phi_{0h}$}
    \label{fig_channel_wave_membrane:initial_conditions}
\end{figure}

For this example, we fix
$$
\phi = - 1 \quad\mbox{ on } \quad \partial\Omega_b\quad  \mbox{ and } \quad \phi= 1\quad \mbox{ on } \quad \partial\Omega_t
$$
and
$$
\partial_{\boldsymbol{n}}\phi=0\quad \mbox{ on }\quad \partial\Omega\backslash(\partial\Omega_b\cap \partial\Omega_t),
$$
see Figure \ref{fig_channel_wave_membrane:initial_conditions}.

\begin{figure}
    \begin{subfigure}[b]{0.25\textwidth}
        \centering
        \includegraphics[width=1.0\textwidth]{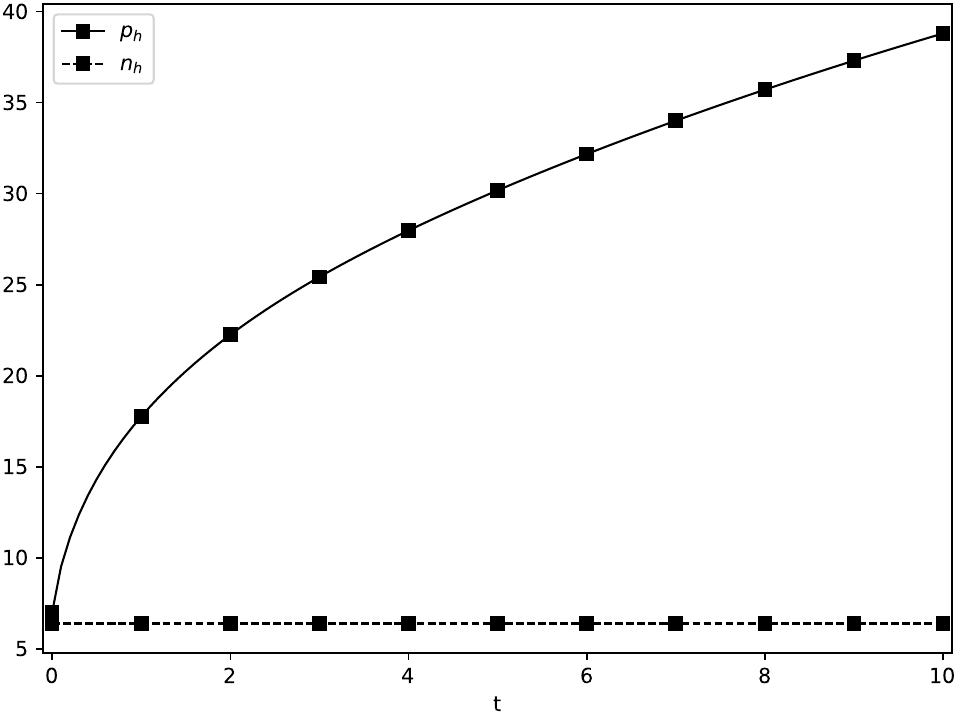}
    \end{subfigure}
    \begin{subfigure}[b]{0.25\textwidth}
        \centering
        \includegraphics[width=1.0\textwidth]{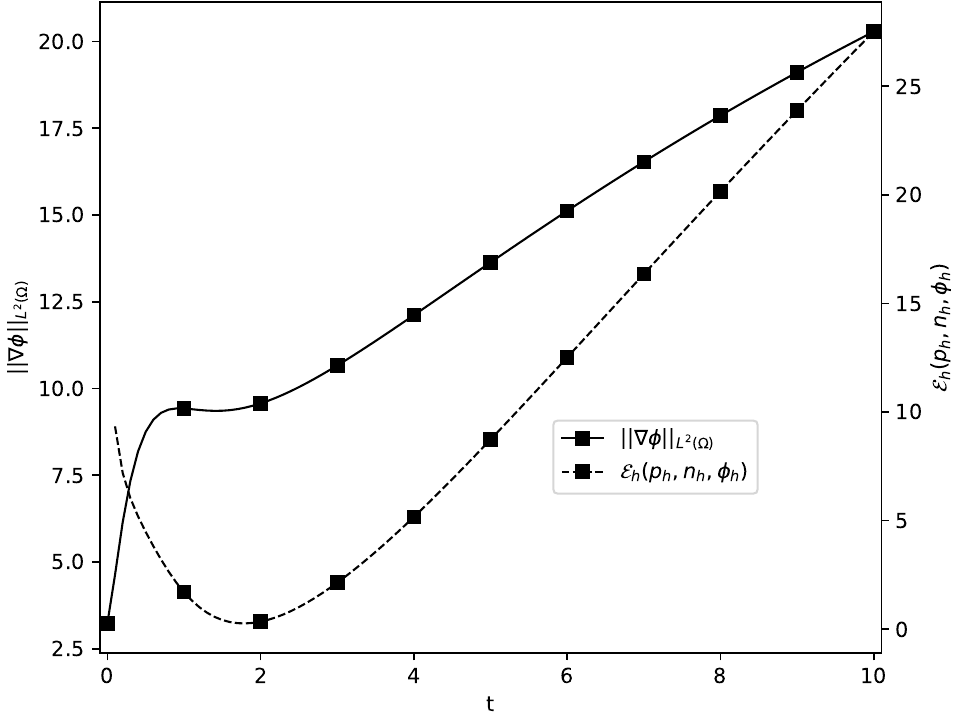}
    \end{subfigure}
    \caption*{Algorithm 1}
    \begin{subfigure}[b]{0.25\textwidth}
        \centering
        \includegraphics[width=1.0\textwidth]{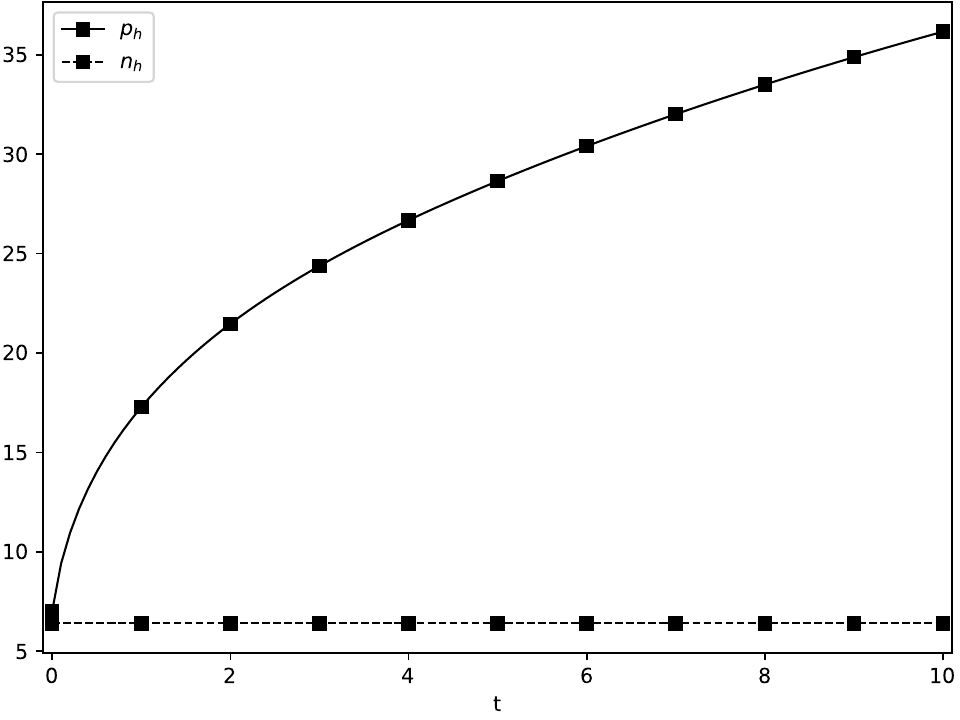}
    \end{subfigure}
    \begin{subfigure}[b]{0.25\textwidth}
        \centering
        \includegraphics[width=1.0\textwidth]{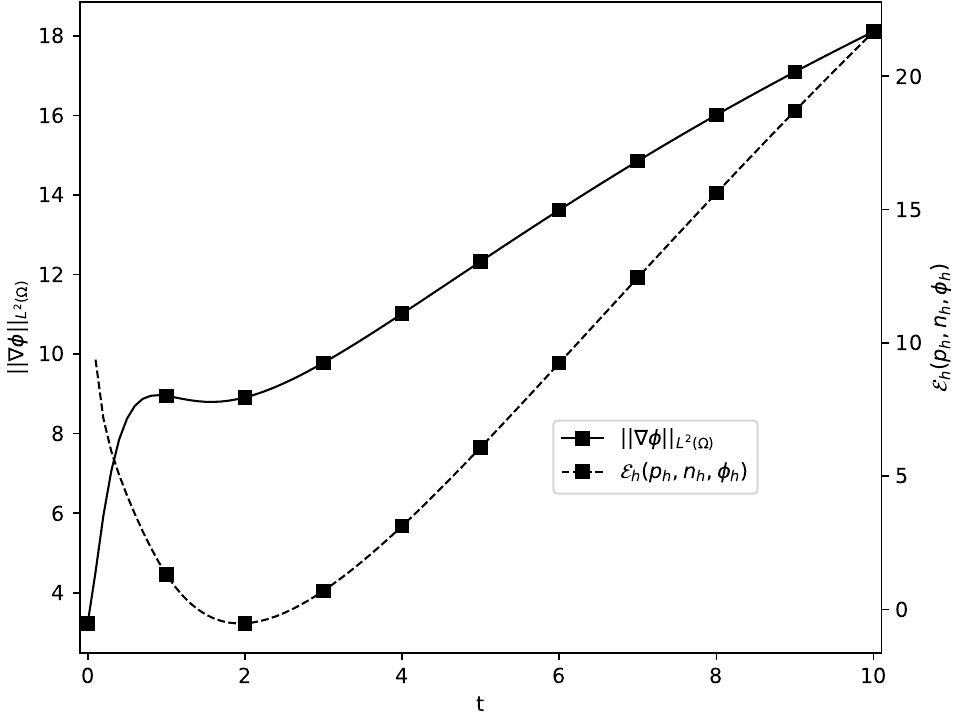}
    \end{subfigure}
    \caption*{Algorithm 2}
    \caption{Mass conservation (left), and energy and entropy evolutions (right)}
    \label{fig_channel_wave_membrane:mass_and_energies}
\end{figure}

\begin{figure}
    \begin{subfigure}[b]{0.25\textwidth}
        \centering
        \includegraphics[width=1.0\textwidth]{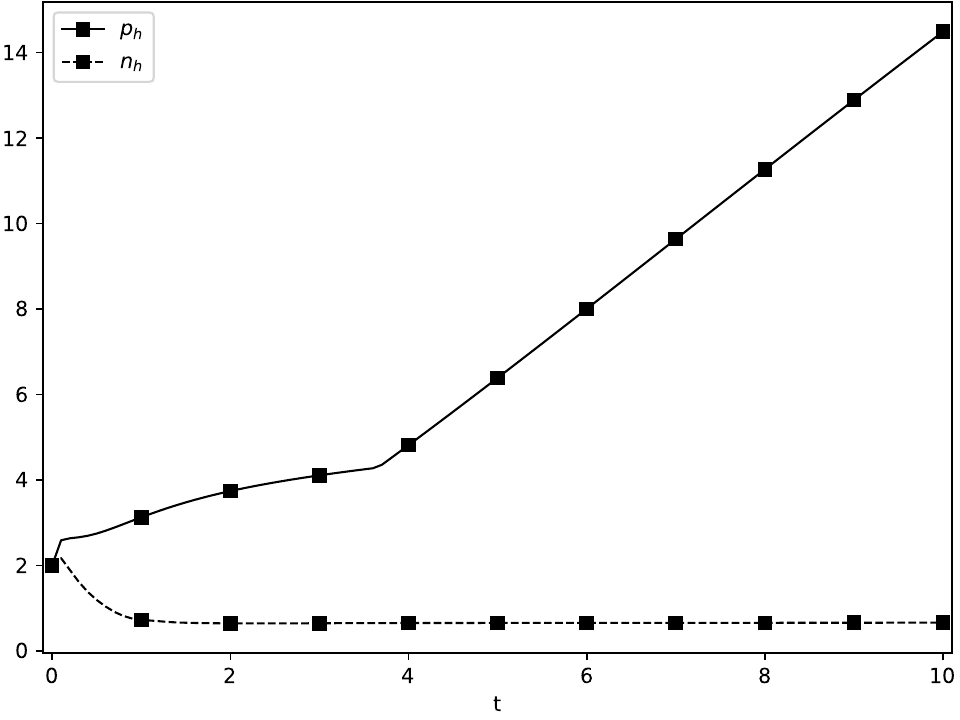}
     \end{subfigure}
    \begin{subfigure}[b]{0.25\textwidth}
        \centering
        \includegraphics[width=1.0\textwidth]{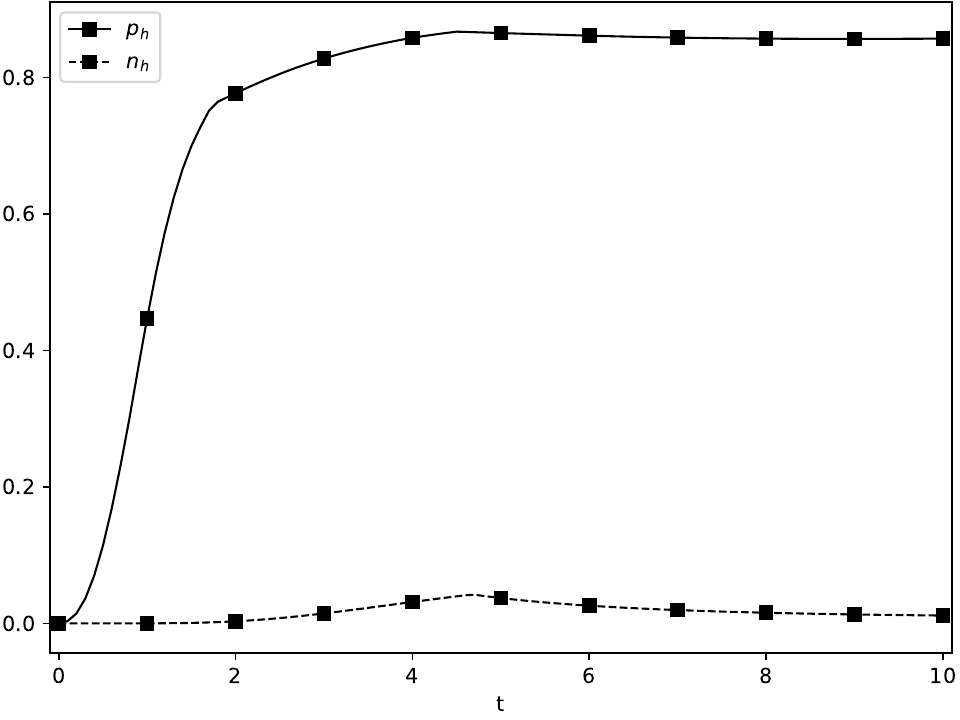}
     \end{subfigure}
  \caption*{Algorithm 1}
       \begin{subfigure}[b]{0.25\textwidth}
        \centering
        \includegraphics[width=1.0\textwidth]{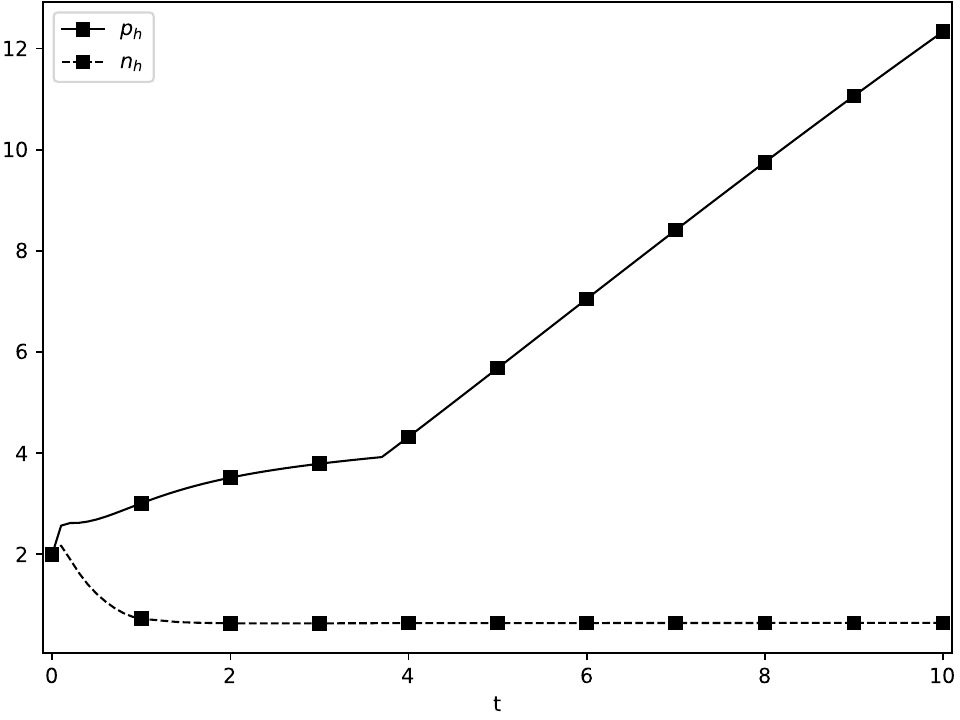}
    \end{subfigure}
    \begin{subfigure}[b]{0.25\textwidth}
        \centering
        \includegraphics[width=1.0\textwidth]{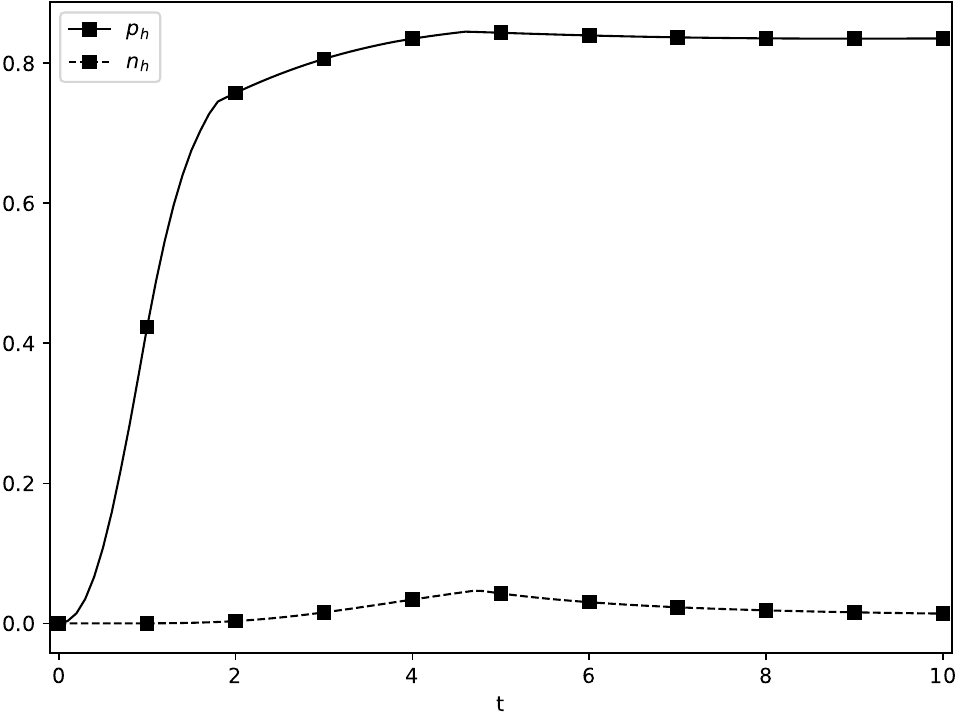}
    \end{subfigure}    
    \caption*{Algorithm 2}
    \caption{Maxima (left) and minima (right)}
    \label{fig_channel_wave_membrane:max_and_min}
\end{figure}

\begin{figure} 
\begin{subfigure}[b]{0.22\textwidth}
        \centering
        \includegraphics[width=1.0\textwidth]{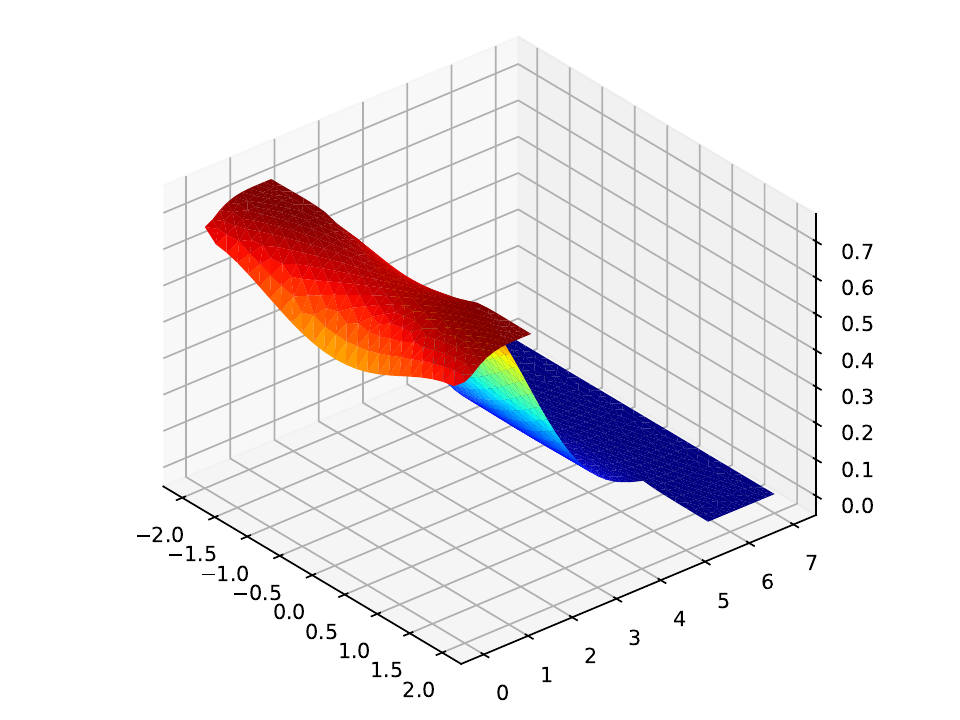}
    \end{subfigure}
    \begin{subfigure}[b]{0.22\textwidth}
        \centering
        \includegraphics[width=1.0\textwidth]{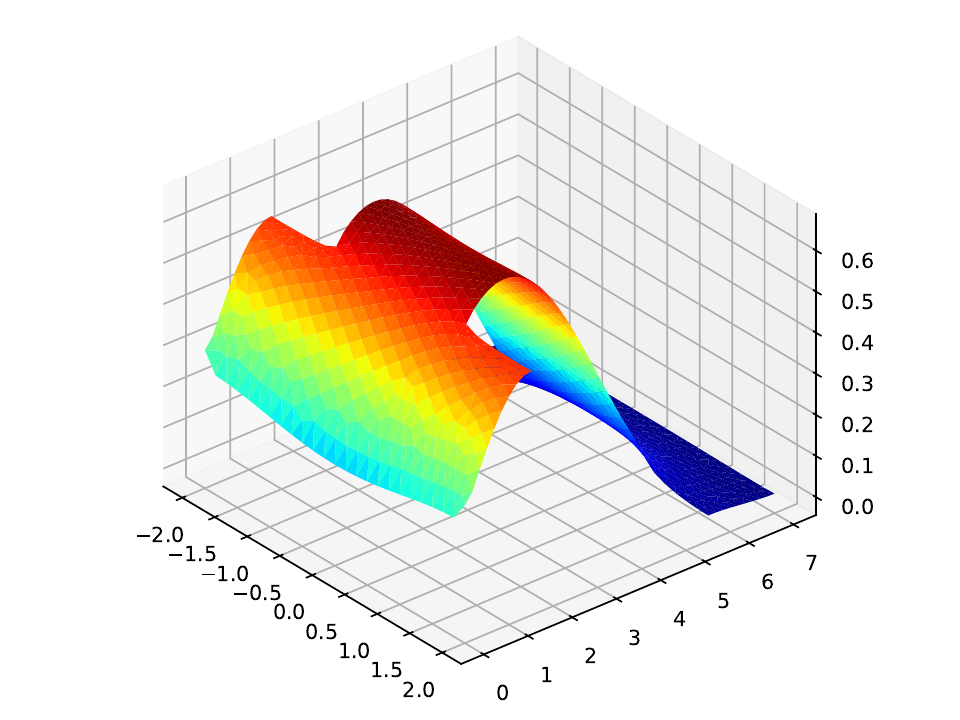}
    \end{subfigure}
    \begin{subfigure}[b]{0.22\textwidth}
        \centering
        \includegraphics[width=1.0\textwidth]{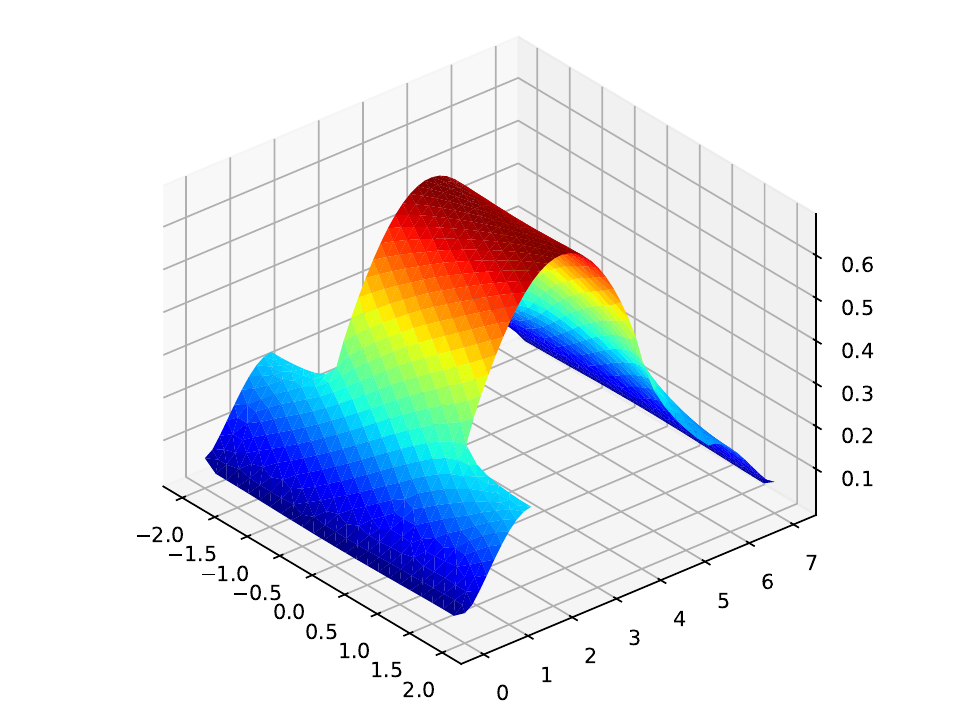}
    \end{subfigure}
    \begin{subfigure}[b]{0.22\textwidth}
        \centering
        \includegraphics[width=1.0\textwidth]{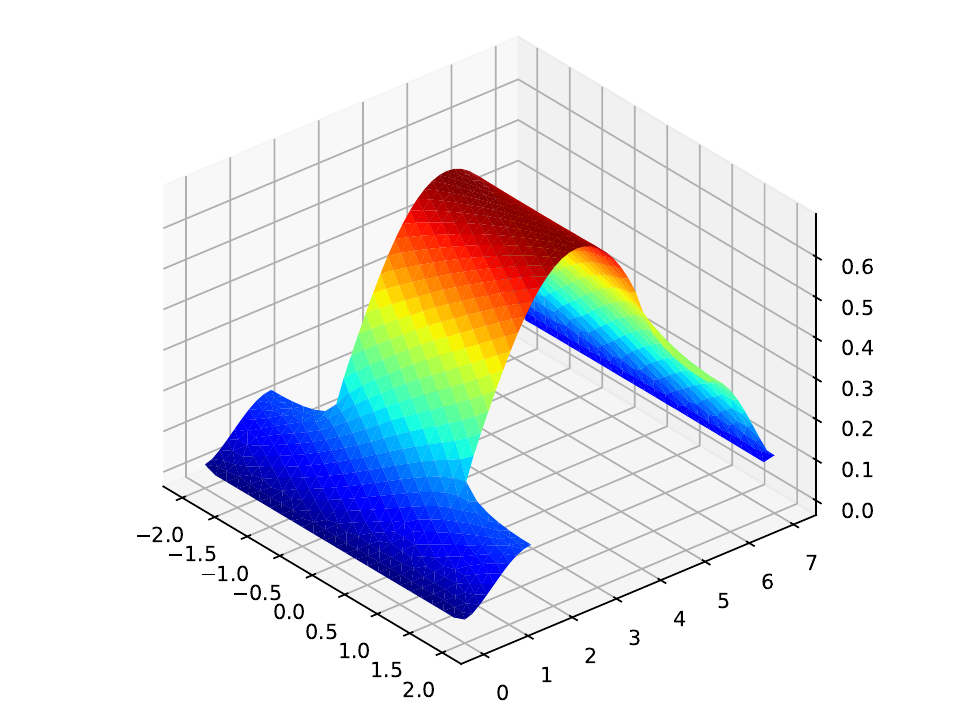}
    \end{subfigure}    
    \caption*{Snapshots of $n_h$ at times $t=1$, $2$, $5$ and $10$}
    \begin{subfigure}[b]{0.22\textwidth}
        \centering
        \includegraphics[width=1.0\textwidth]{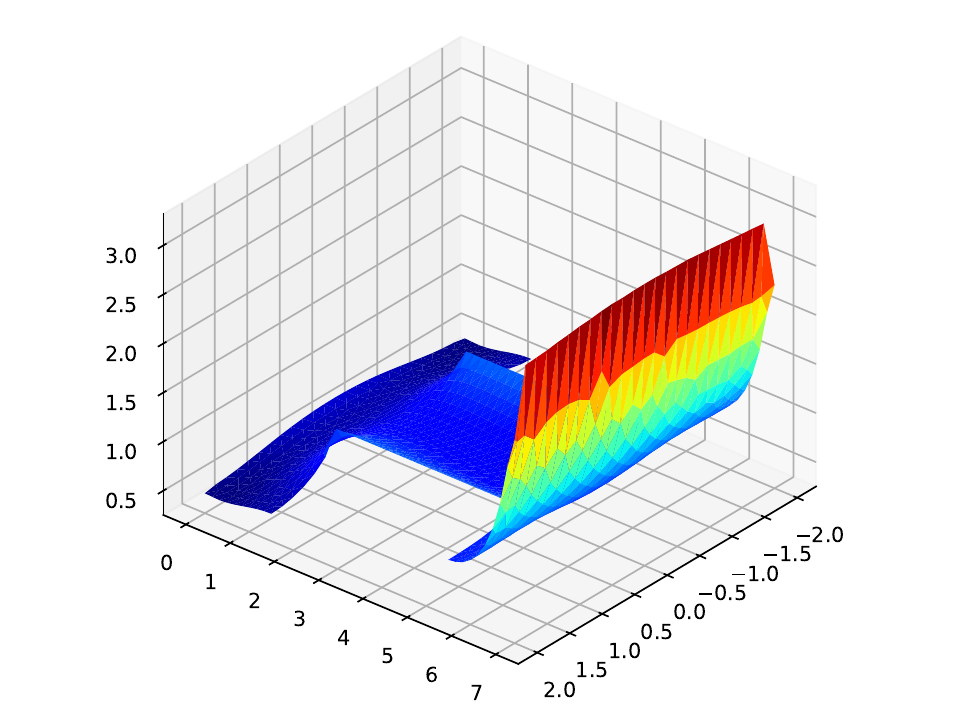}
    \end{subfigure}
    \begin{subfigure}[b]{0.22\textwidth}
        \centering
        \includegraphics[width=1.0\textwidth]{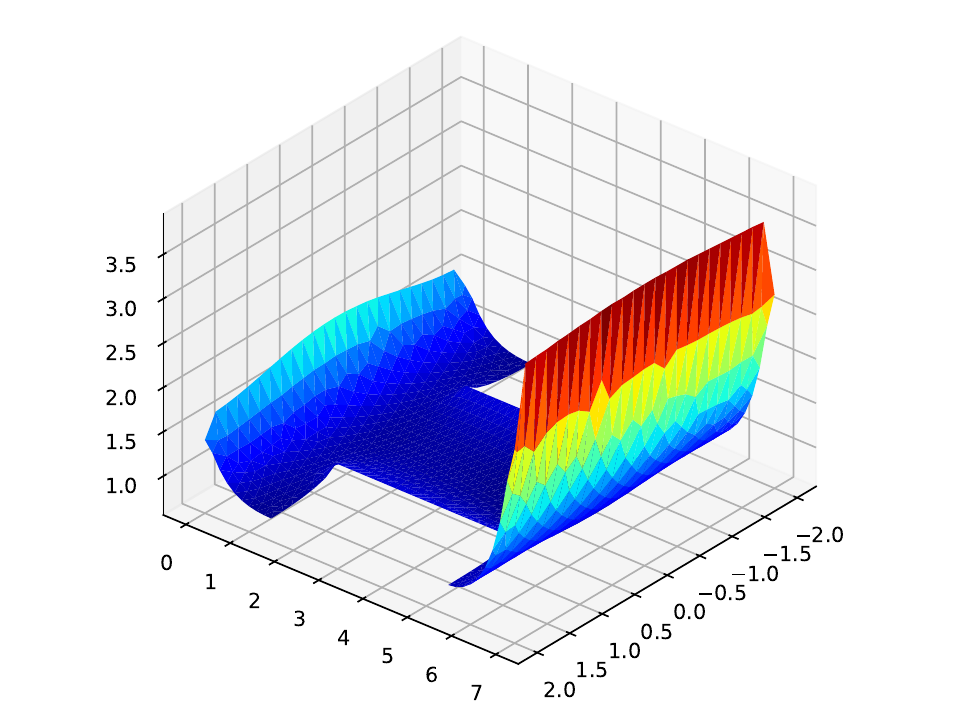}
    \end{subfigure}
    \begin{subfigure}[b]{0.22\textwidth}
        \centering
        \includegraphics[width=1.0\textwidth]{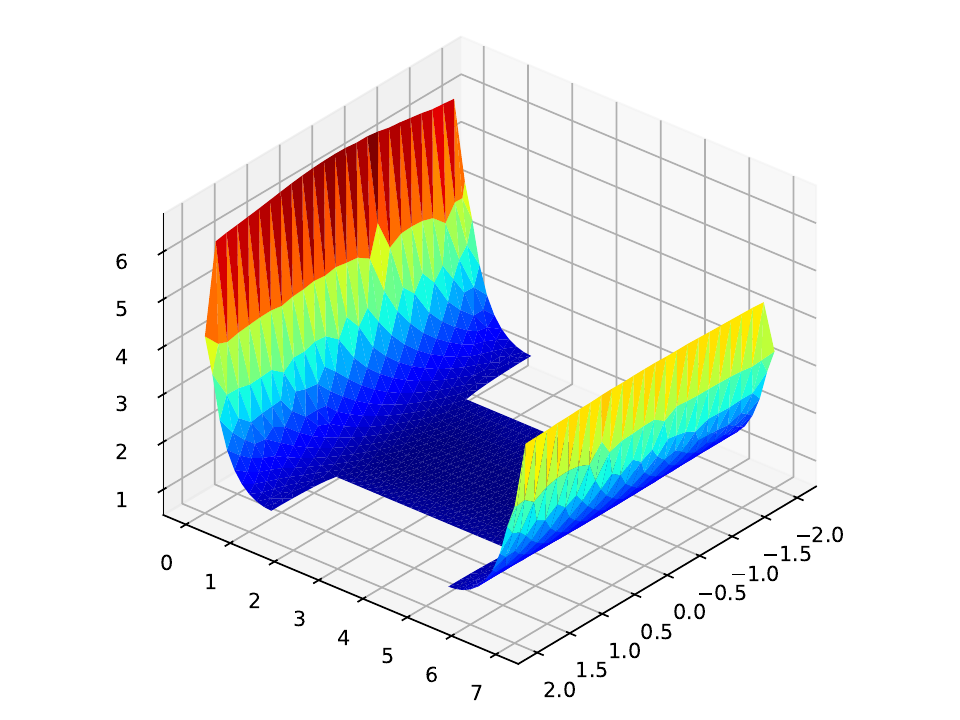}
    \end{subfigure}
    \begin{subfigure}[b]{0.22\textwidth}
        \centering
        \includegraphics[width=1.0\textwidth]{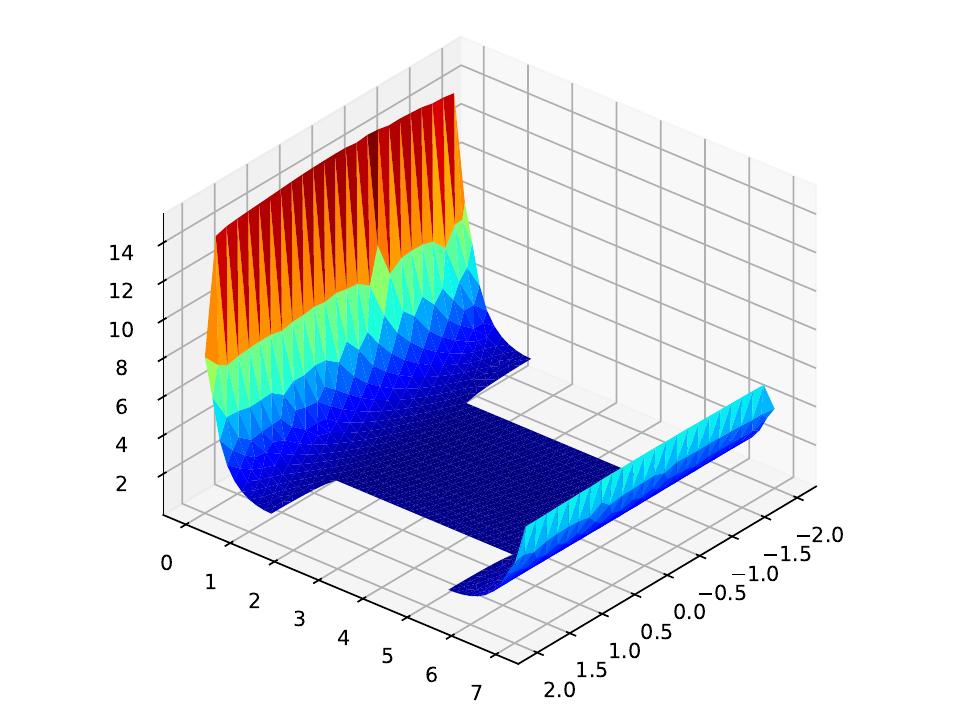}
    \end{subfigure}    
    \caption*{Snapshots of $p_h$ at times $t=1$, $2$, $5$ and $10$}
    \begin{subfigure}[b]{0.22\textwidth}
        \centering
        \includegraphics[width=1.0\textwidth]{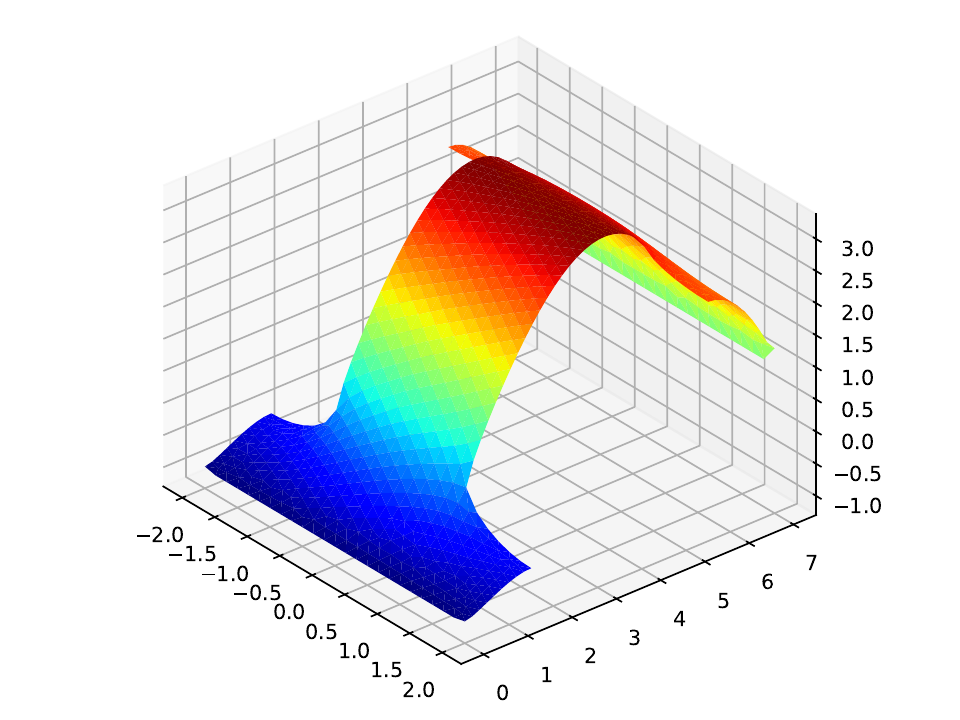}
    \end{subfigure}
    \begin{subfigure}[b]{0.22\textwidth}
        \centering
        \includegraphics[width=1.0\textwidth]{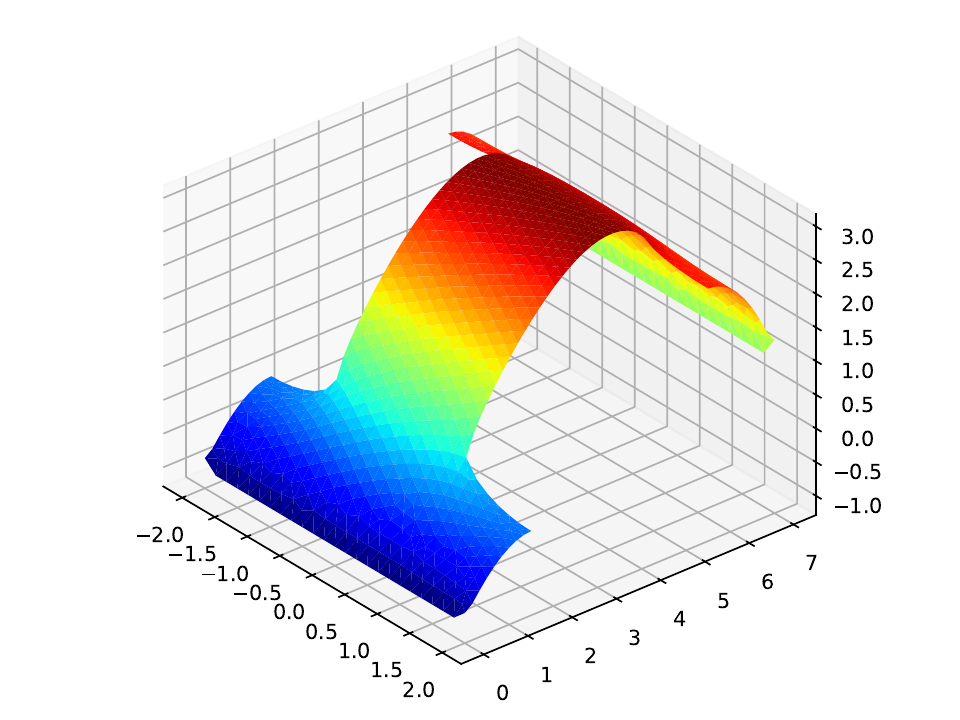}
    \end{subfigure}
    \begin{subfigure}[b]{0.22\textwidth}
        \centering
        \includegraphics[width=1.0\textwidth]{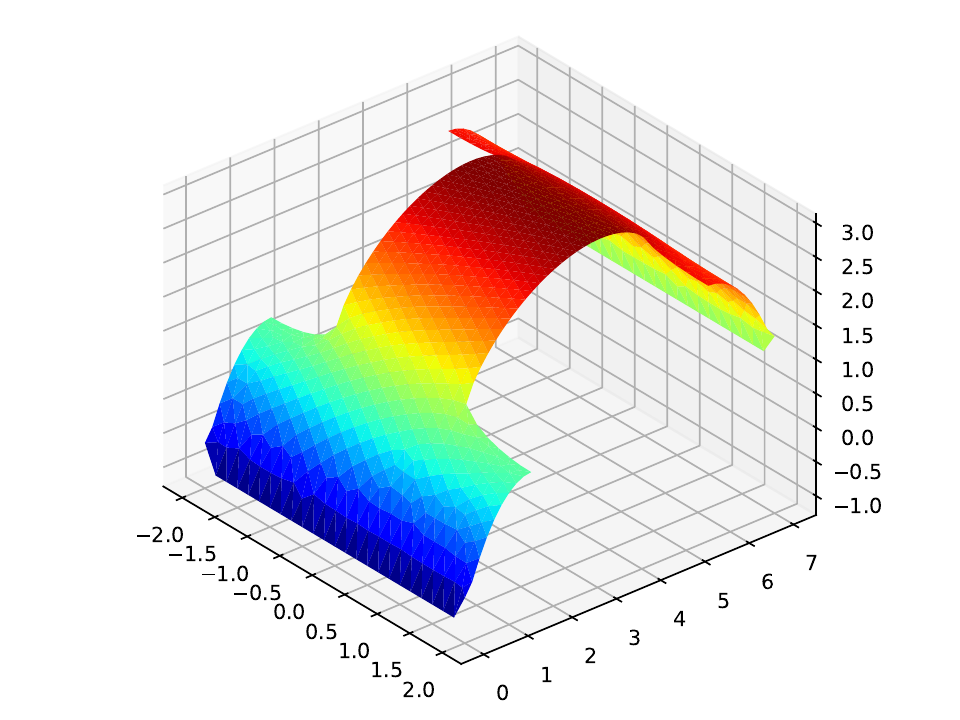}
    \end{subfigure}
    \begin{subfigure}[b]{0.22\textwidth}
        \centering
        \includegraphics[width=1.0\textwidth]{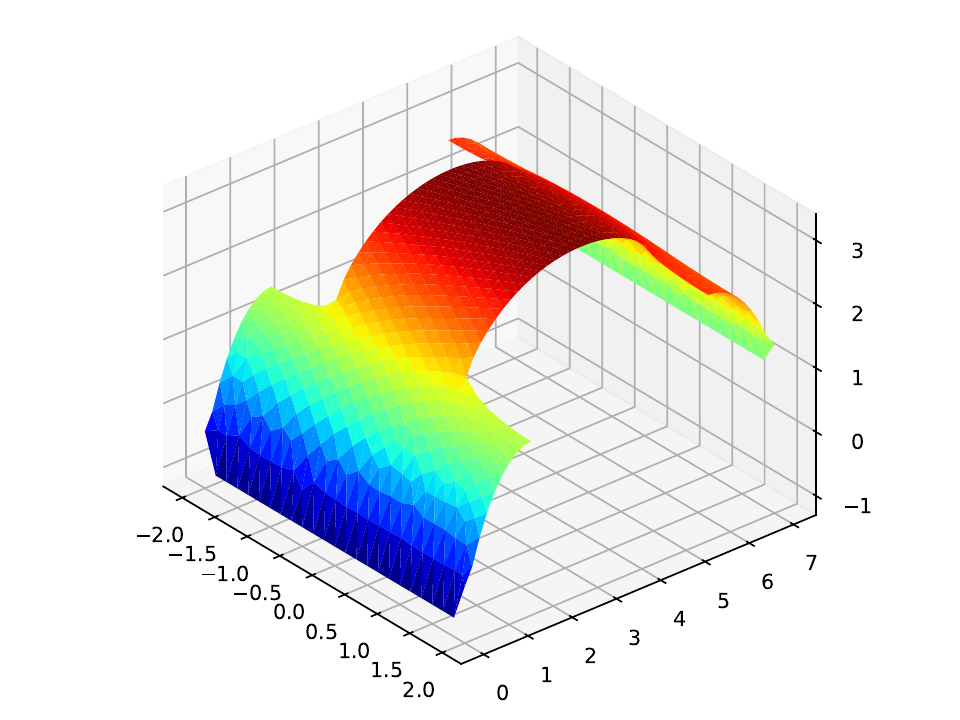}
    \end{subfigure}    
    \caption*{Snapshots of $\phi_h$ at times $t=1$, $2$, $5$ and $10$}
    \caption{Algorithm 1}
    \label{fig_channel_wave_membrane:snapshots_alg1}
\end{figure}

\begin{figure}
\begin{subfigure}[b]{0.22\textwidth}
        \centering
        \includegraphics[width=1.0\textwidth]{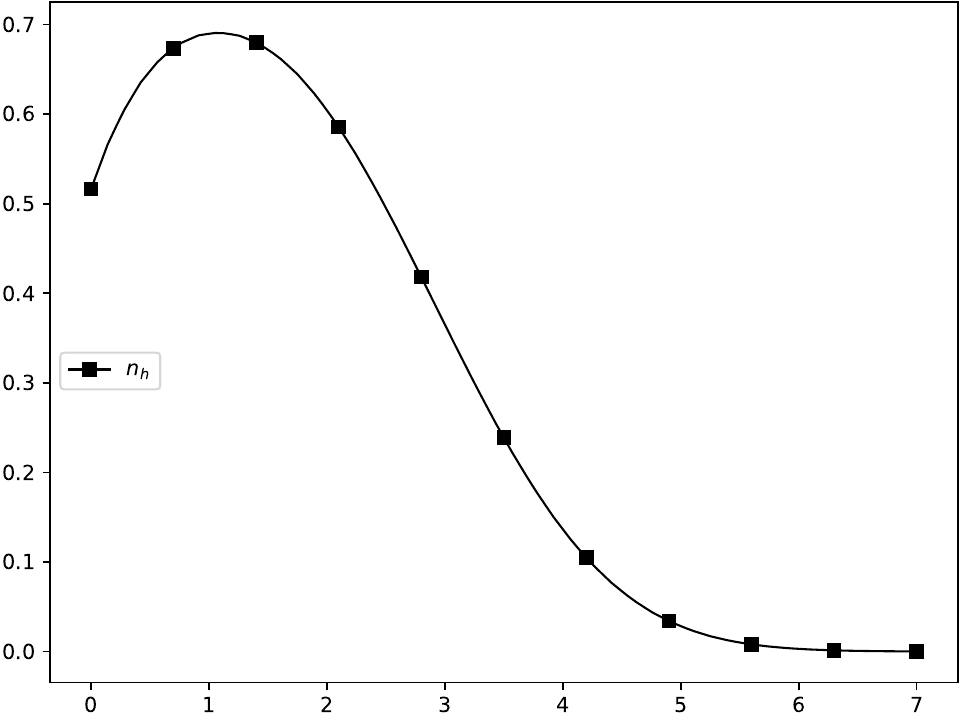}
    \end{subfigure}
    \begin{subfigure}[b]{0.22\textwidth}
        \centering
        \includegraphics[width=1.0\textwidth]{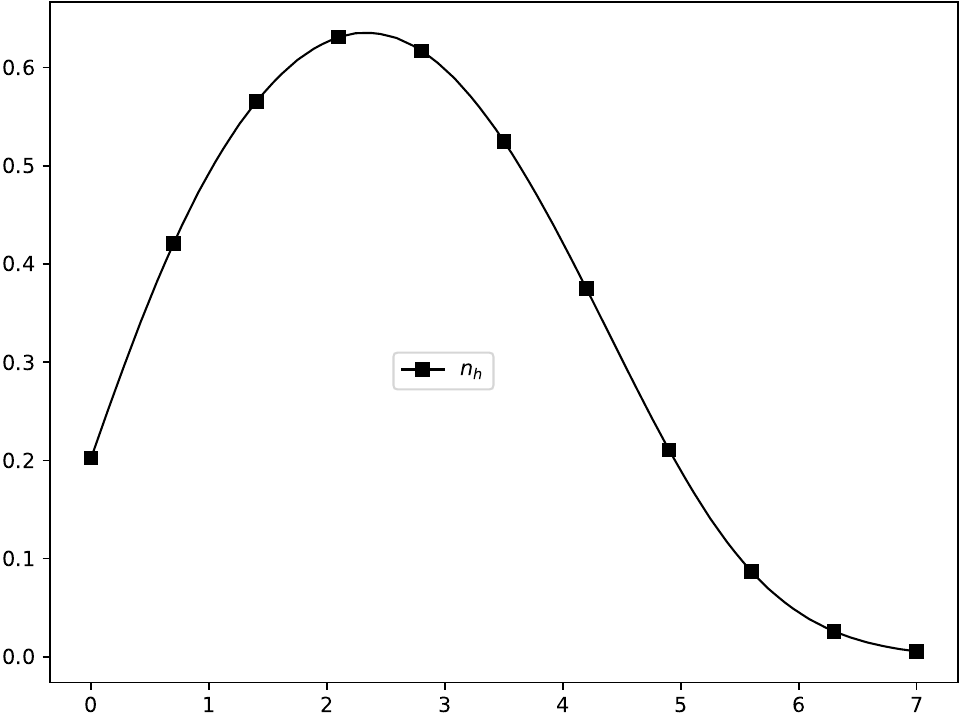}
    \end{subfigure}
    \begin{subfigure}[b]{0.22\textwidth}
        \centering
        \includegraphics[width=1.0\textwidth]{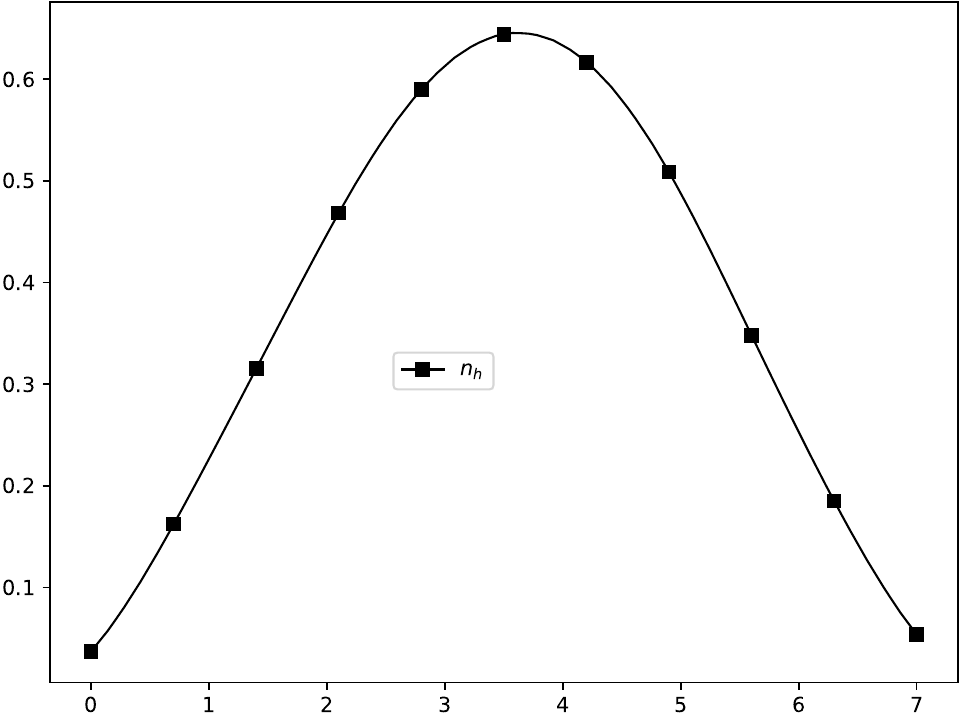}
    \end{subfigure}
    \begin{subfigure}[b]{0.22\textwidth}
        \centering
        \includegraphics[width=1.0\textwidth]{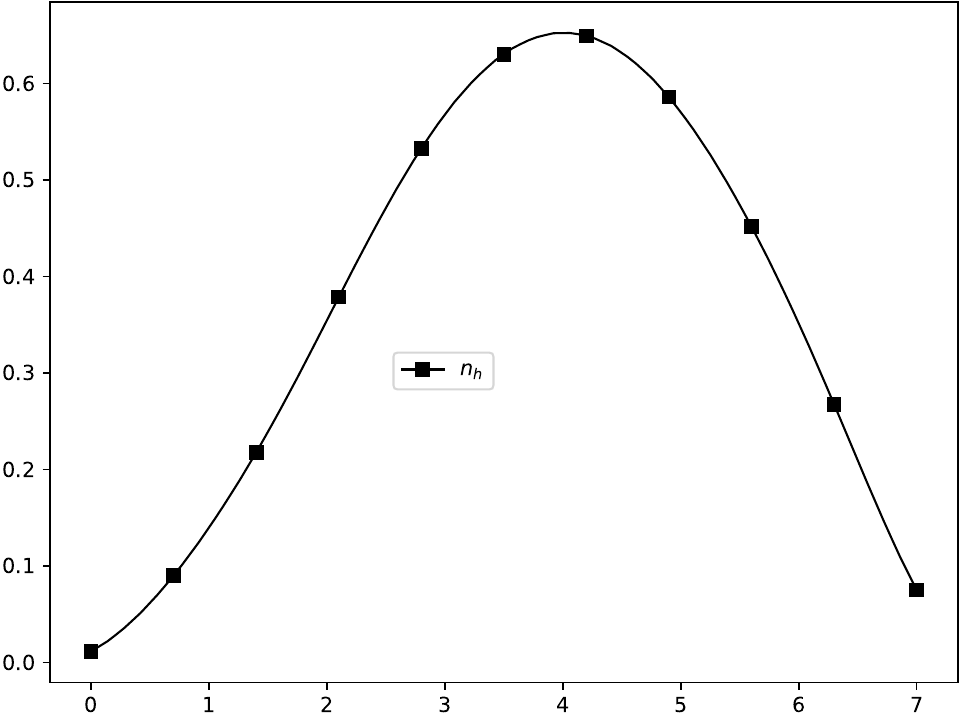}
    \end{subfigure}    
    \caption*{Profiles of $n_h$ at times $t=1$, $2$, $5$ and $10$}
\begin{subfigure}[b]{0.22\textwidth}
        \centering
        \includegraphics[width=1.0\textwidth]{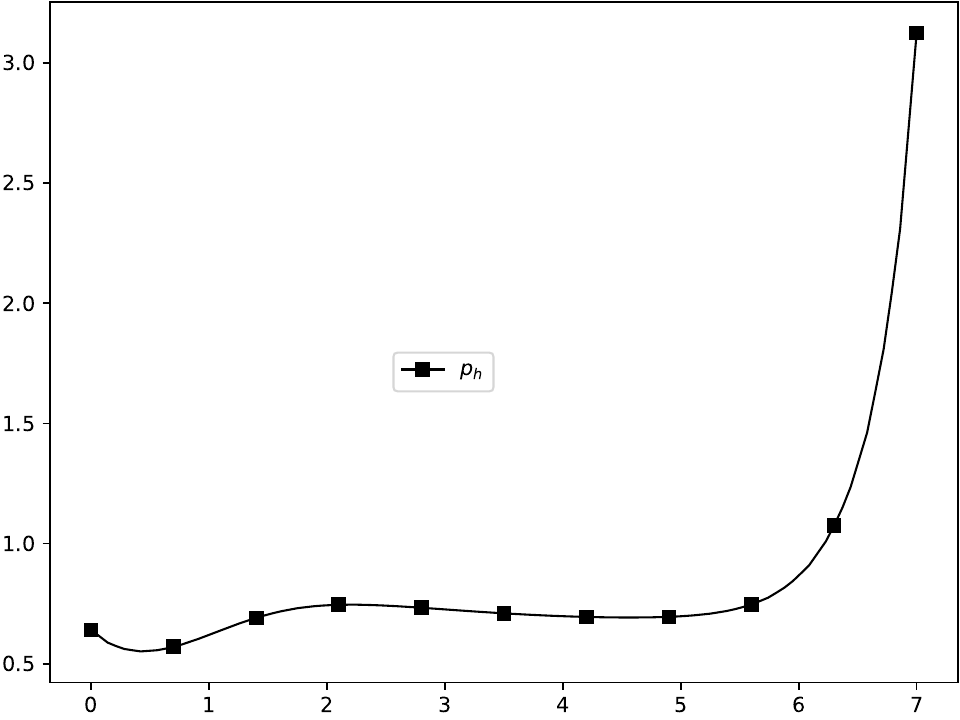}
    \end{subfigure}
    \begin{subfigure}[b]{0.22\textwidth}
        \centering
        \includegraphics[width=1.0\textwidth]{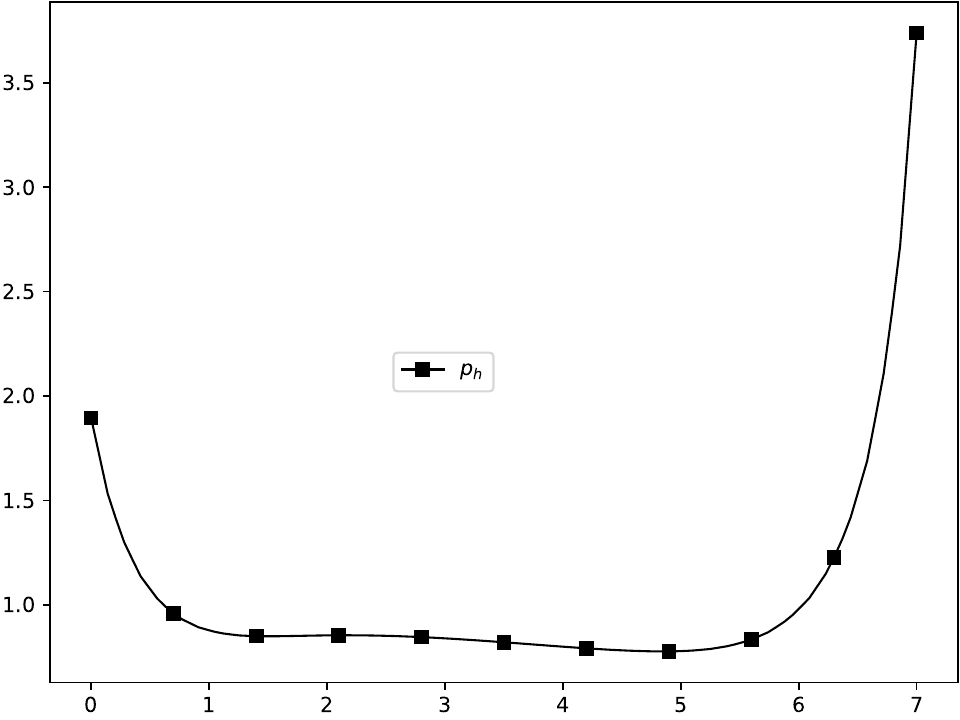}
    \end{subfigure}
    \begin{subfigure}[b]{0.22\textwidth}
        \centering
        \includegraphics[width=1.0\textwidth]{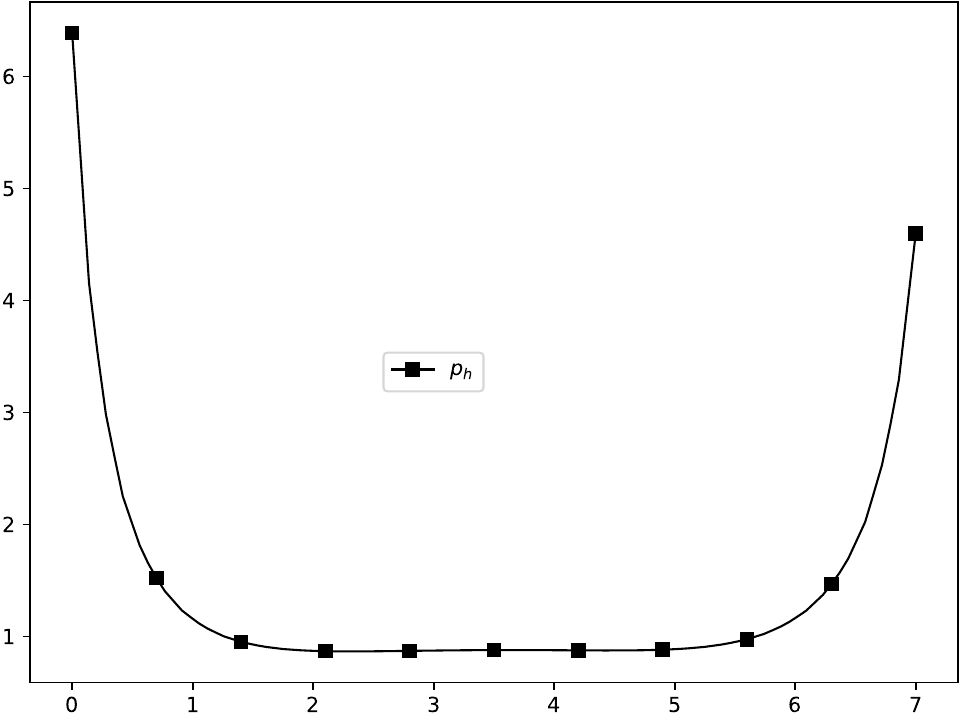}
    \end{subfigure}
    \begin{subfigure}[b]{0.22\textwidth}
        \centering
        \includegraphics[width=1.0\textwidth]{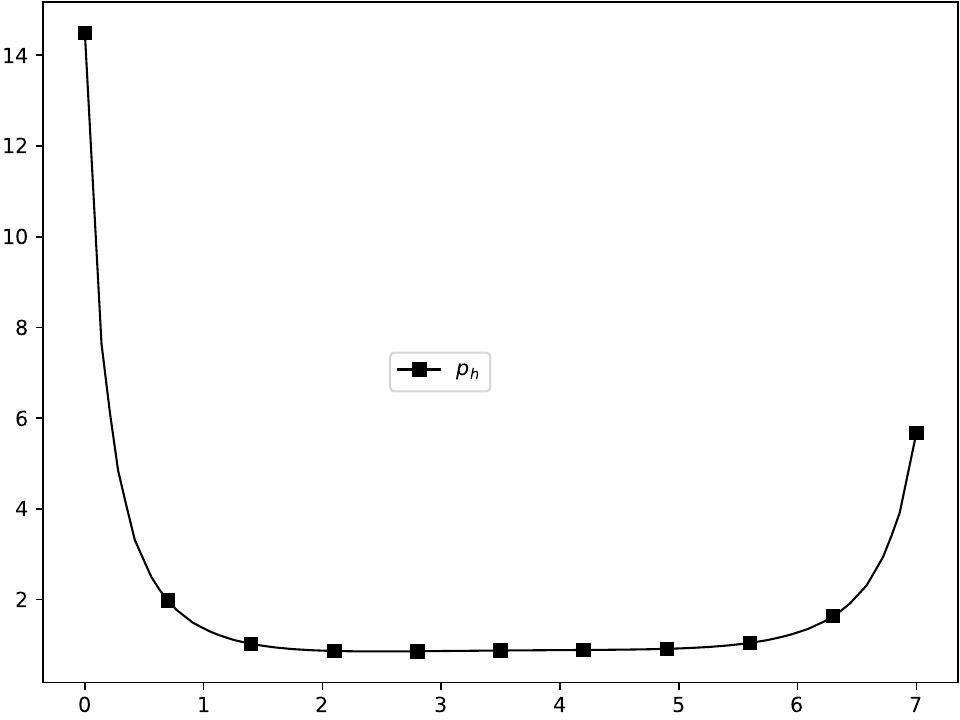}
    \end{subfigure}    
    \caption*{Profiles of $p_h$ at times $t=1$, $2$, $5$ and $10$}

    \begin{subfigure}[b]{0.22\textwidth}
        \centering
        \includegraphics[width=1.0\textwidth]{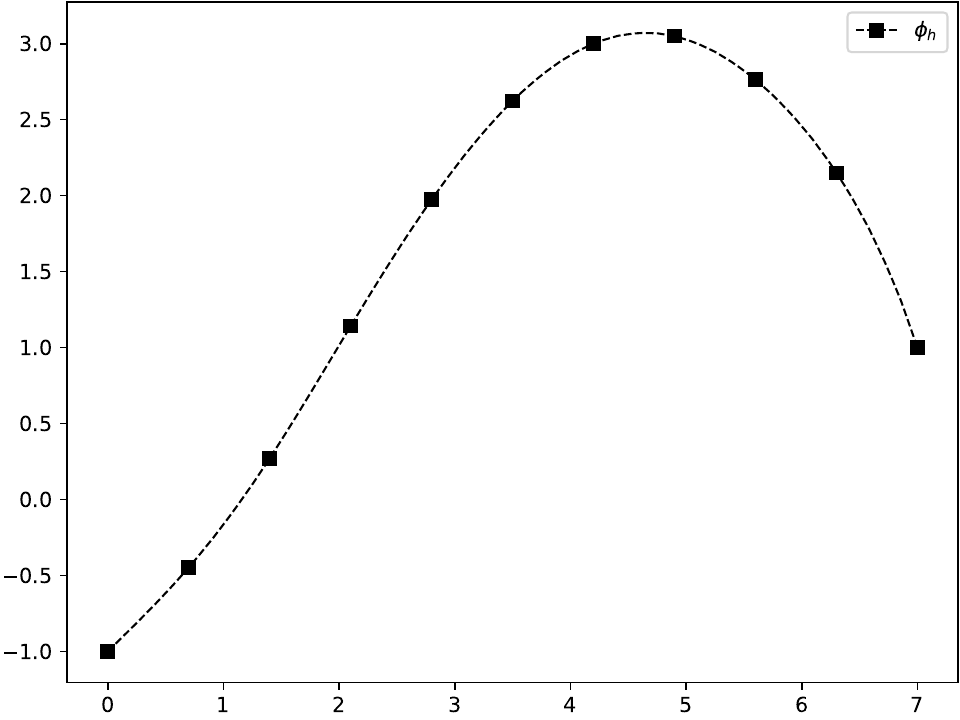}
    \end{subfigure}
    \begin{subfigure}[b]{0.22\textwidth}
        \centering
        \includegraphics[width=1.0\textwidth]{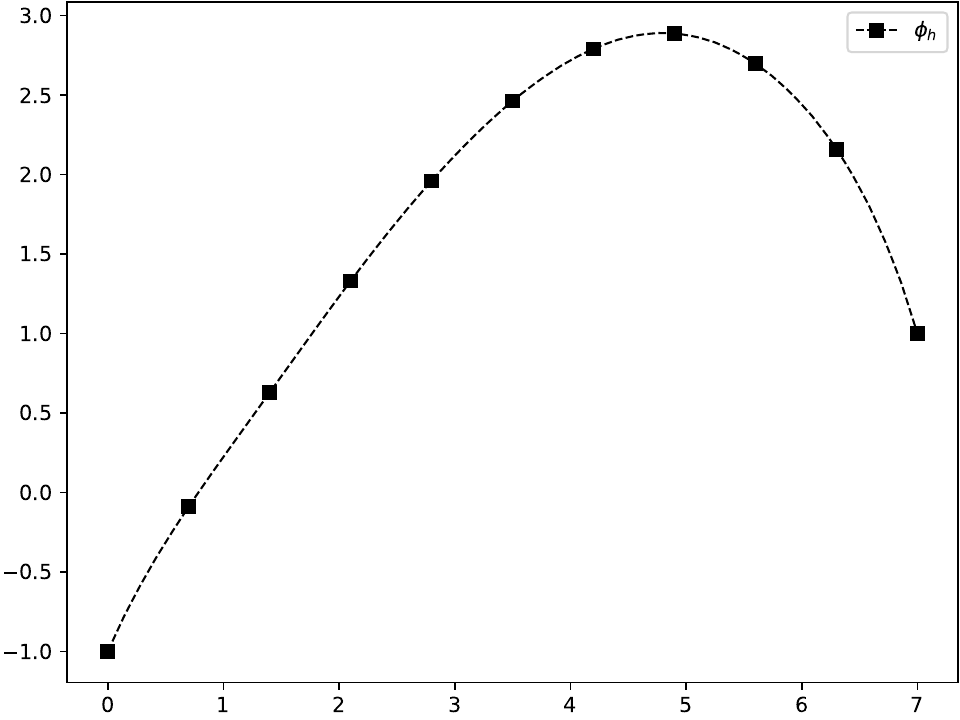}
    \end{subfigure}
    \begin{subfigure}[b]{0.22\textwidth}
        \centering
        \includegraphics[width=1.0\textwidth]{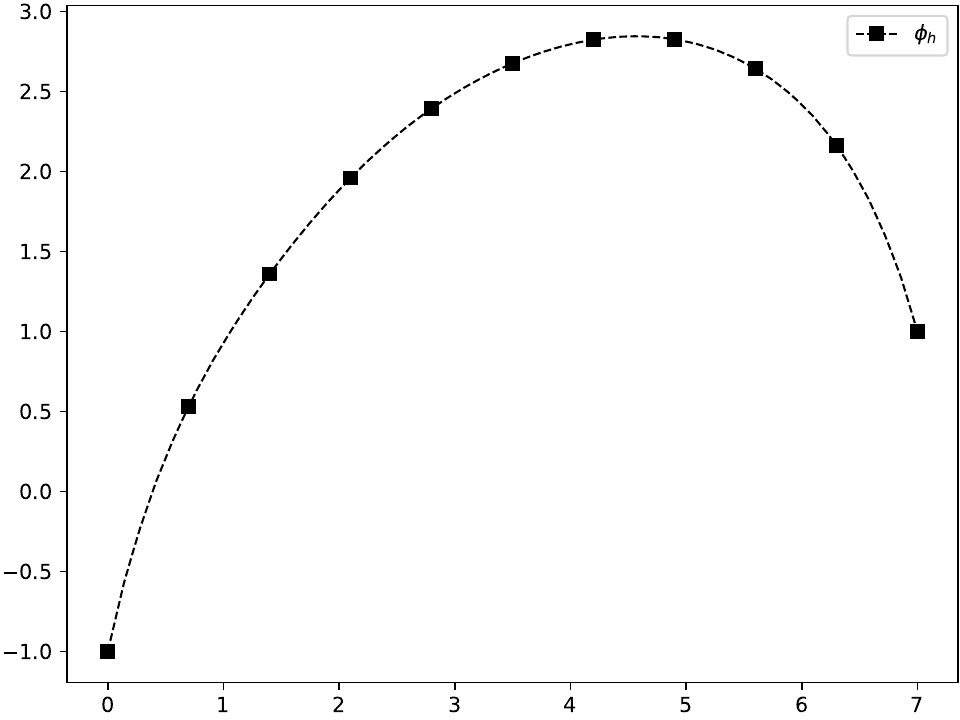}
    \end{subfigure}
    \begin{subfigure}[b]{0.22\textwidth}
        \centering
        \includegraphics[width=1.0\textwidth]{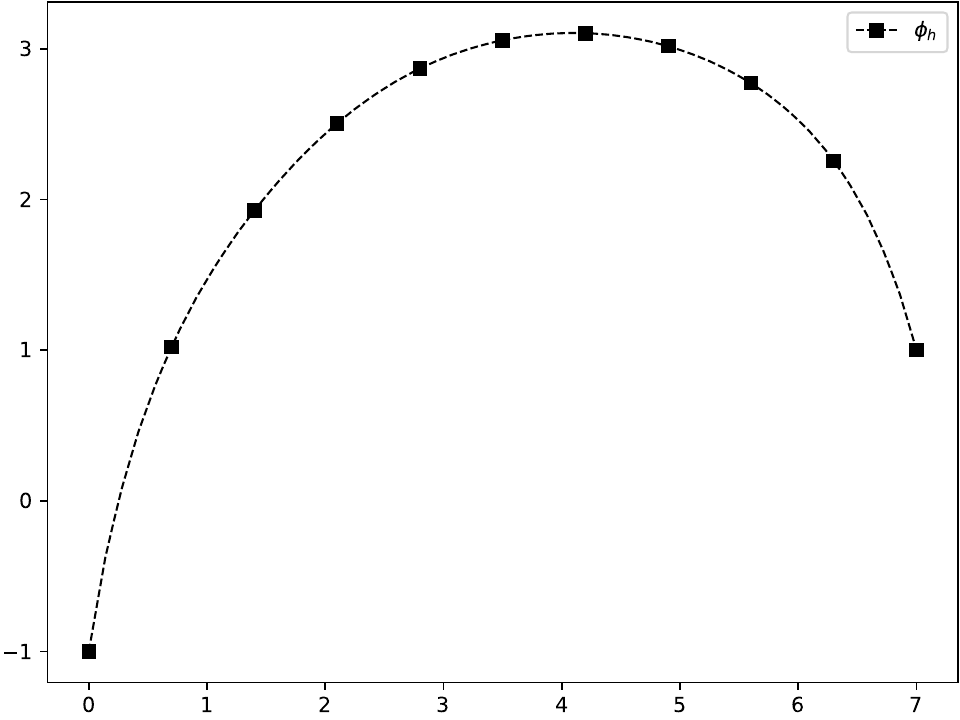}
    \end{subfigure}    
    \caption*{Snapshots of $\phi_h$ at times $t=1$, $2$, $5$ and $10$}
    \caption{Algorithm 1}
    \label{fig_channel_wave_membrane:profiles_alg1}
\end{figure}

\begin{figure}    
\begin{subfigure}[b]{0.22\textwidth}
        \centering
        \includegraphics[width=1.0\textwidth]{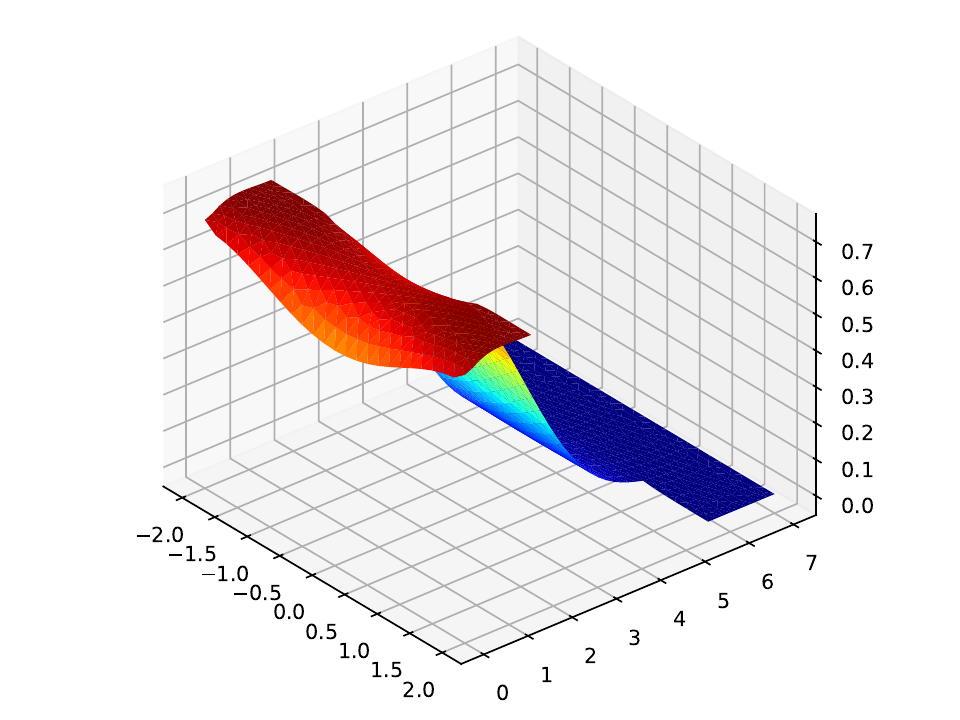}
    \end{subfigure}
    \begin{subfigure}[b]{0.22\textwidth}
        \centering
        \includegraphics[width=1.0\textwidth]{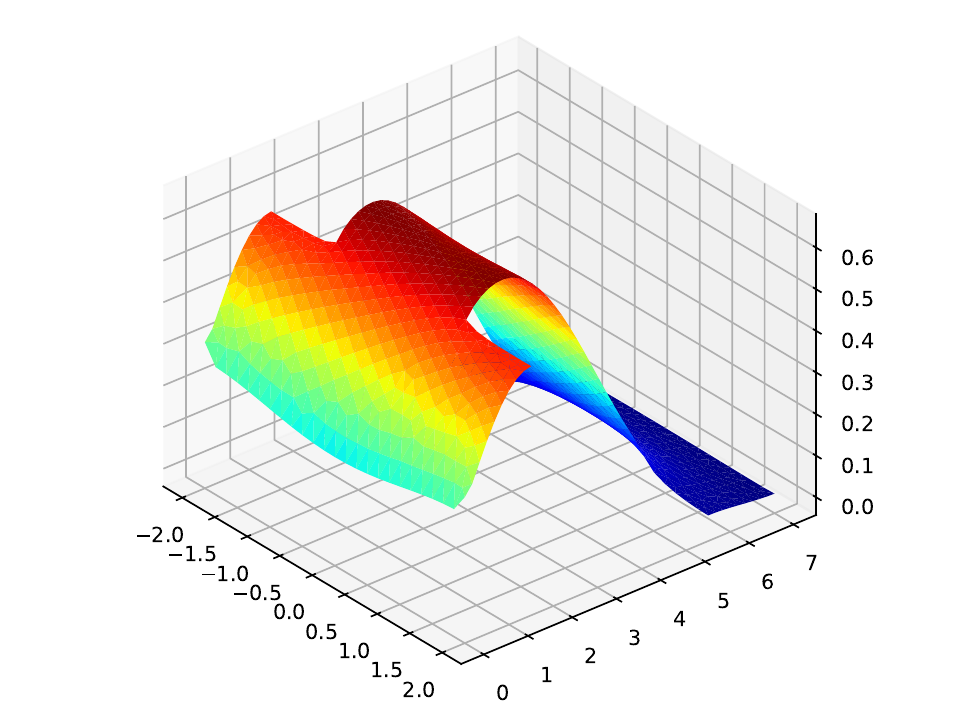}
    \end{subfigure}
    \begin{subfigure}[b]{0.22\textwidth}
        \centering
        \includegraphics[width=1.0\textwidth]{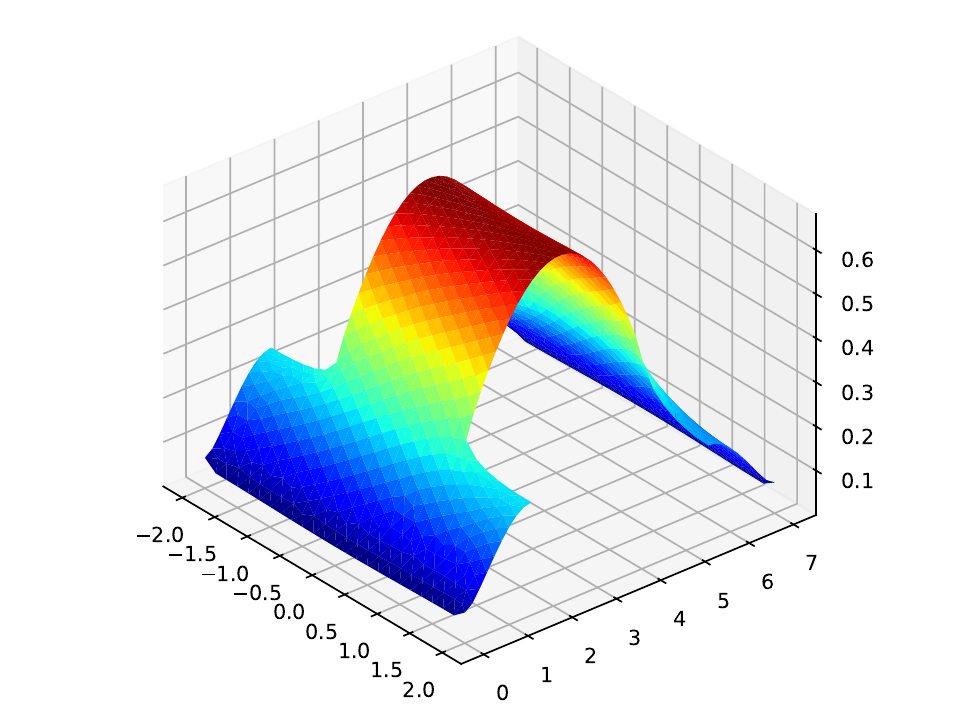}
    \end{subfigure}
    \begin{subfigure}[b]{0.22\textwidth}
        \centering
        \includegraphics[width=1.0\textwidth]{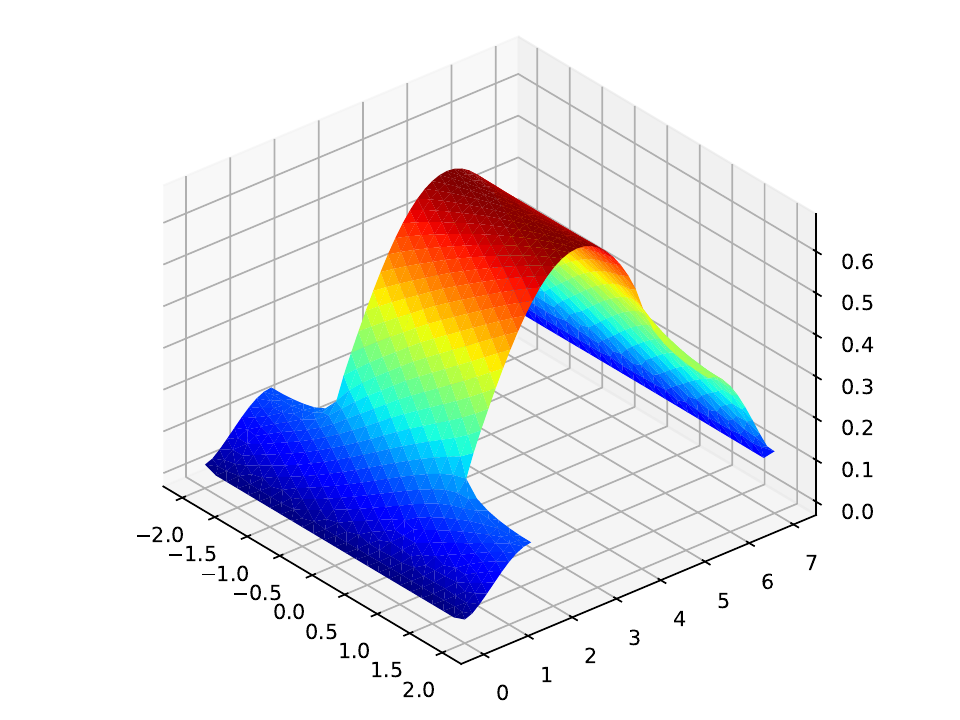}
    \end{subfigure}    
    \caption*{Snapshots of $n_h$ at times $t=1$, $2$, $5$ and $10$}
    \begin{subfigure}[b]{0.22\textwidth}
        \centering
        \includegraphics[width=1.0\textwidth]{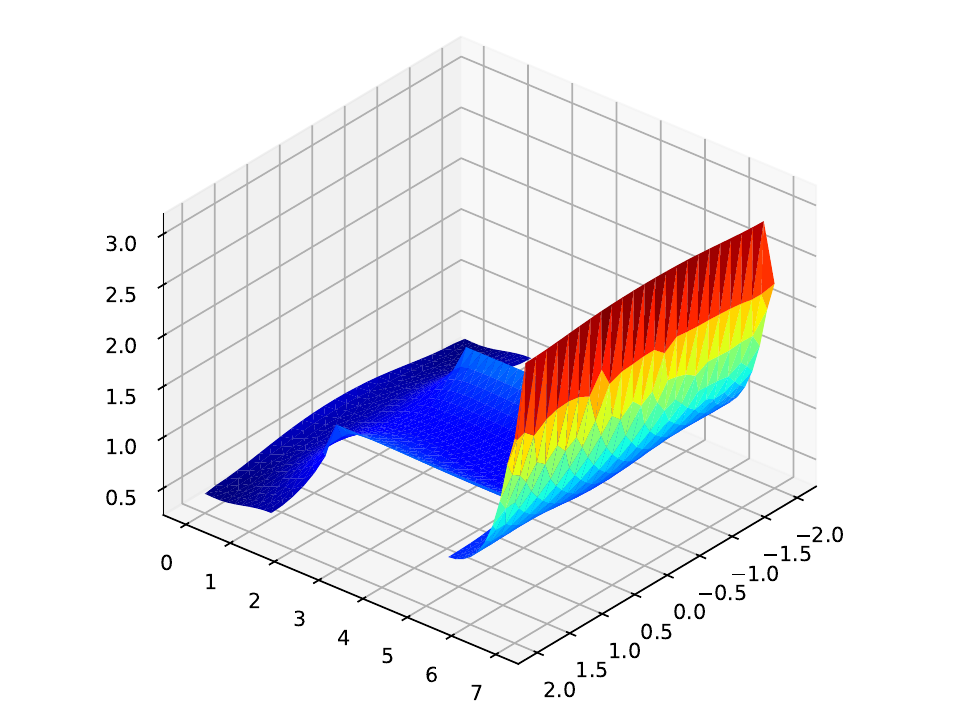}
    \end{subfigure}
    \begin{subfigure}[b]{0.22\textwidth}
        \centering
        \includegraphics[width=1.0\textwidth]{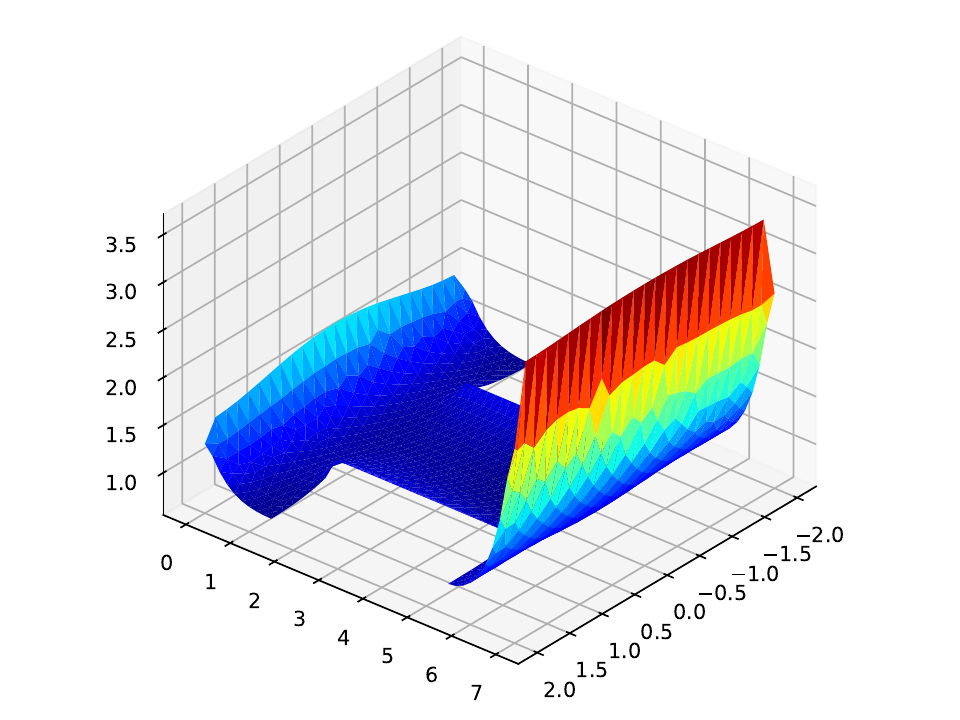}
    \end{subfigure}
    \begin{subfigure}[b]{0.22\textwidth}
        \centering
        \includegraphics[width=1.0\textwidth]{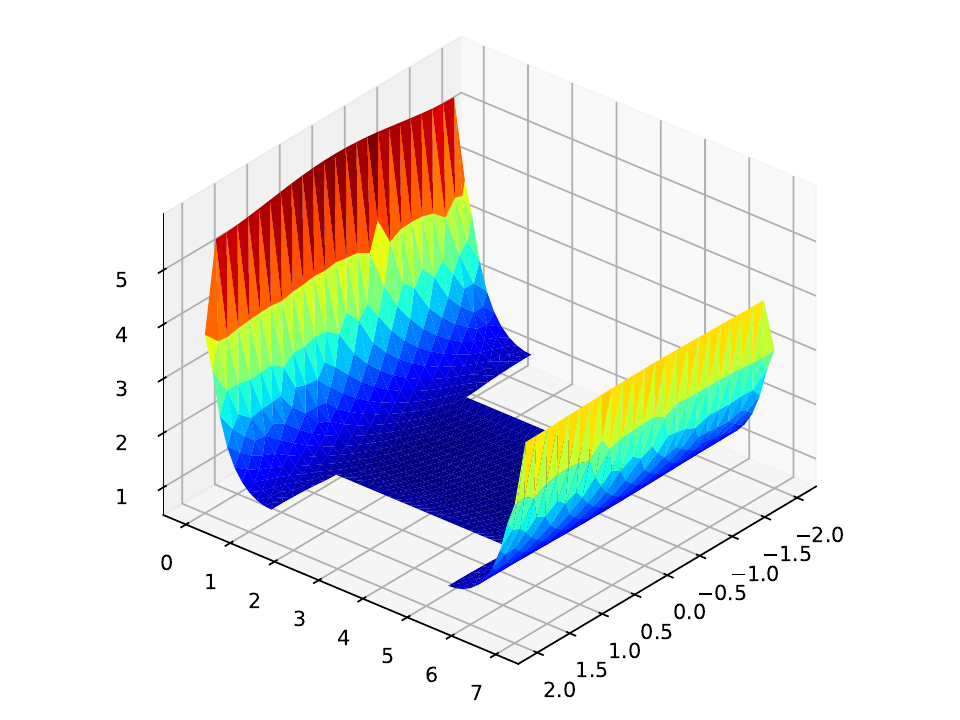}
    \end{subfigure}
    \begin{subfigure}[b]{0.22\textwidth}
        \centering
        \includegraphics[width=1.0\textwidth]{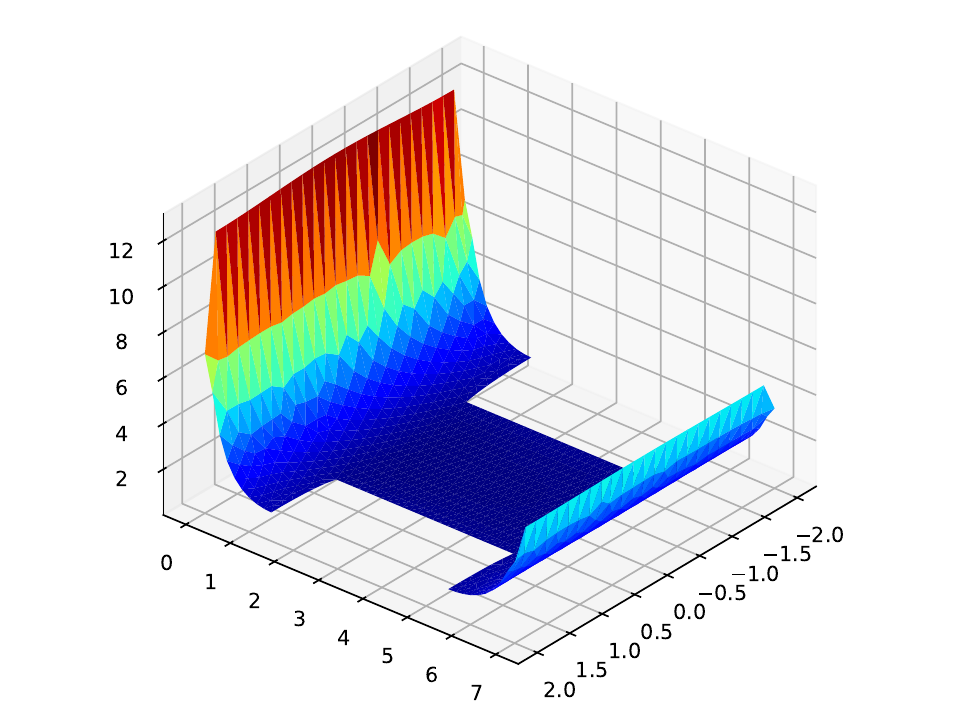}
    \end{subfigure}    
    \caption*{Snapshots of $p_h$ at times $t=1$, $2$, $5$ and $10$}
\begin{subfigure}[b]{0.22\textwidth}
        \centering
        \includegraphics[width=1.0\textwidth]{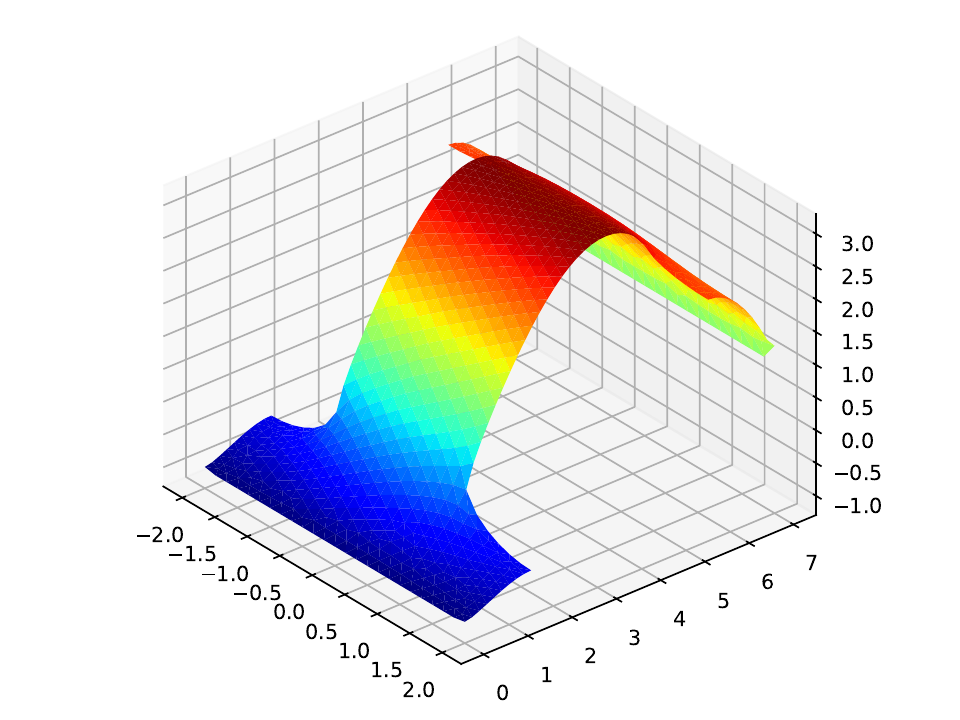}
    \end{subfigure}
    \begin{subfigure}[b]{0.22\textwidth}
        \centering
        \includegraphics[width=1.0\textwidth]{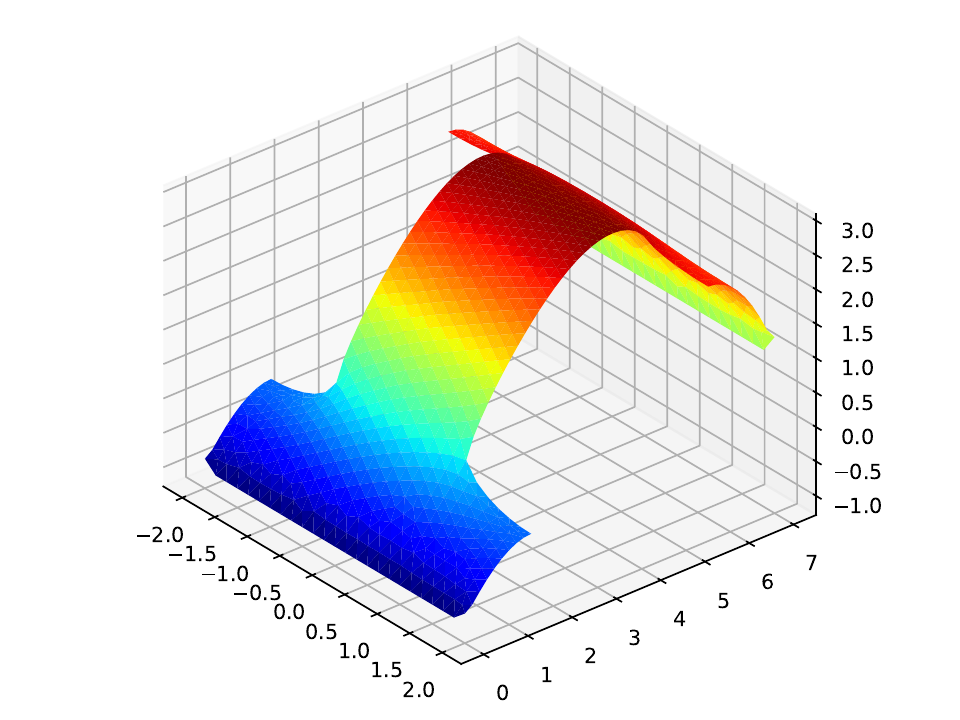}
    \end{subfigure}
    \begin{subfigure}[b]{0.22\textwidth}
        \centering
        \includegraphics[width=1.0\textwidth]{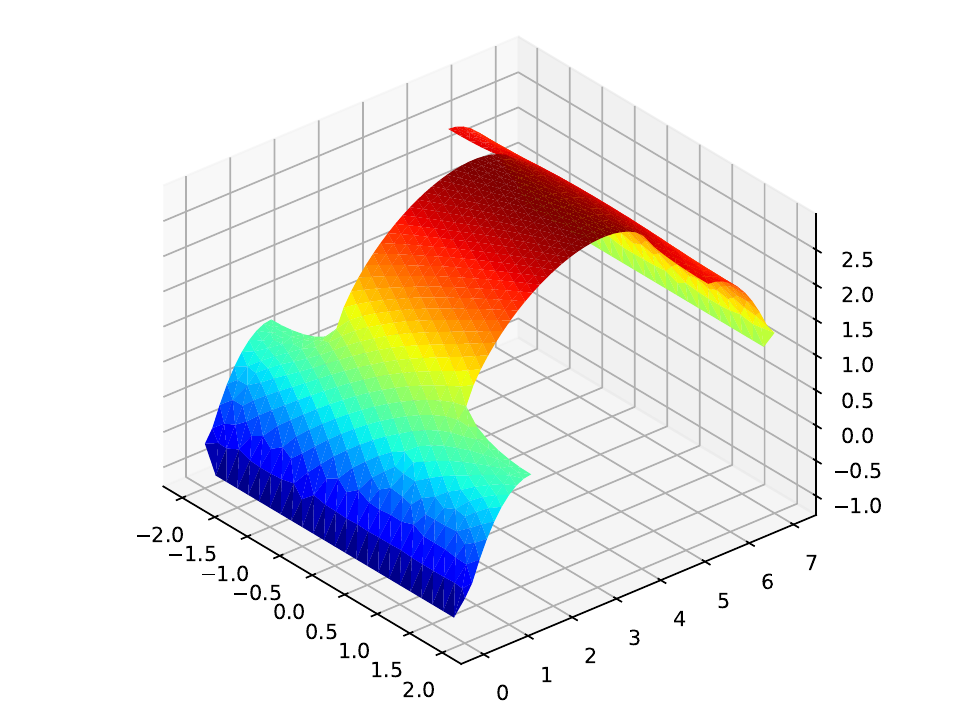}
    \end{subfigure}
    \begin{subfigure}[b]{0.22\textwidth}
        \centering
        \includegraphics[width=1.0\textwidth]{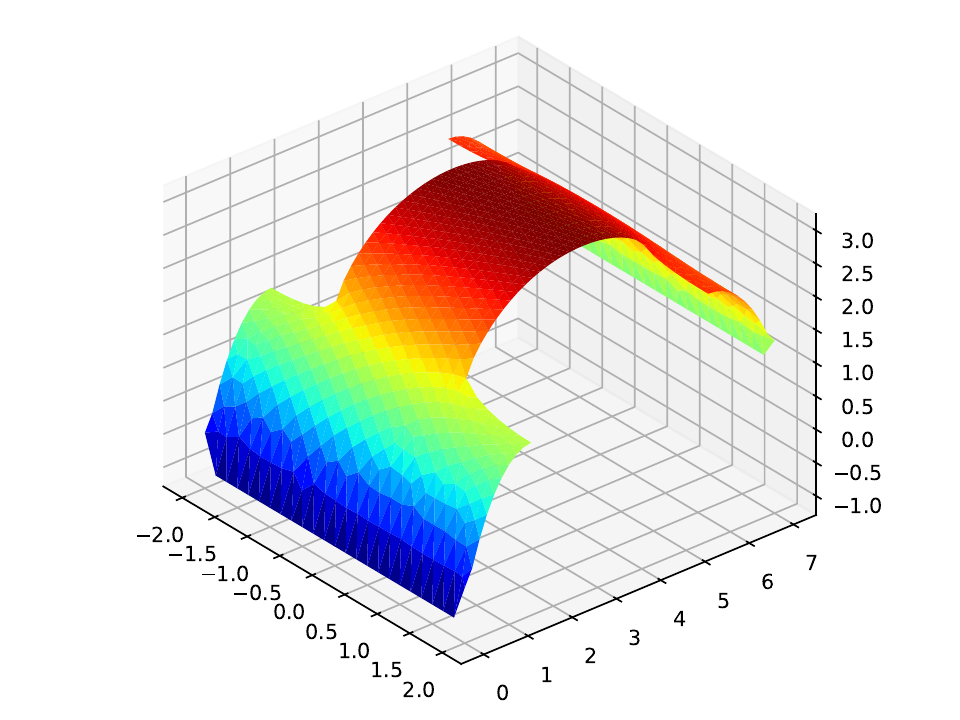}
    \end{subfigure}    
    \caption*{Algorithm 2: Snapshots of $\phi_h$ at times $t=1$, $2$, $5$ and $10$}
    \caption{Algorithm 2}
    \label{fig_channel_wave_membrane:snapshots_alg2}
\end{figure}

\begin{figure}
\begin{subfigure}[b]{0.22\textwidth}
        \centering
        \includegraphics[width=1.0\textwidth]{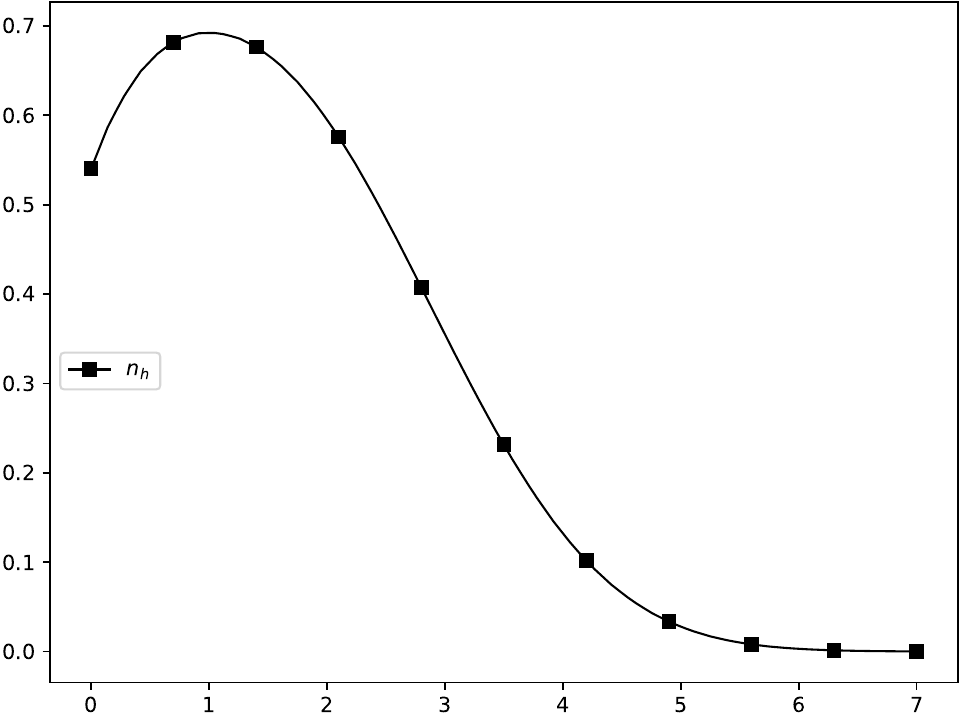}
    \end{subfigure}
    \begin{subfigure}[b]{0.22\textwidth}
        \centering
        \includegraphics[width=1.0\textwidth]{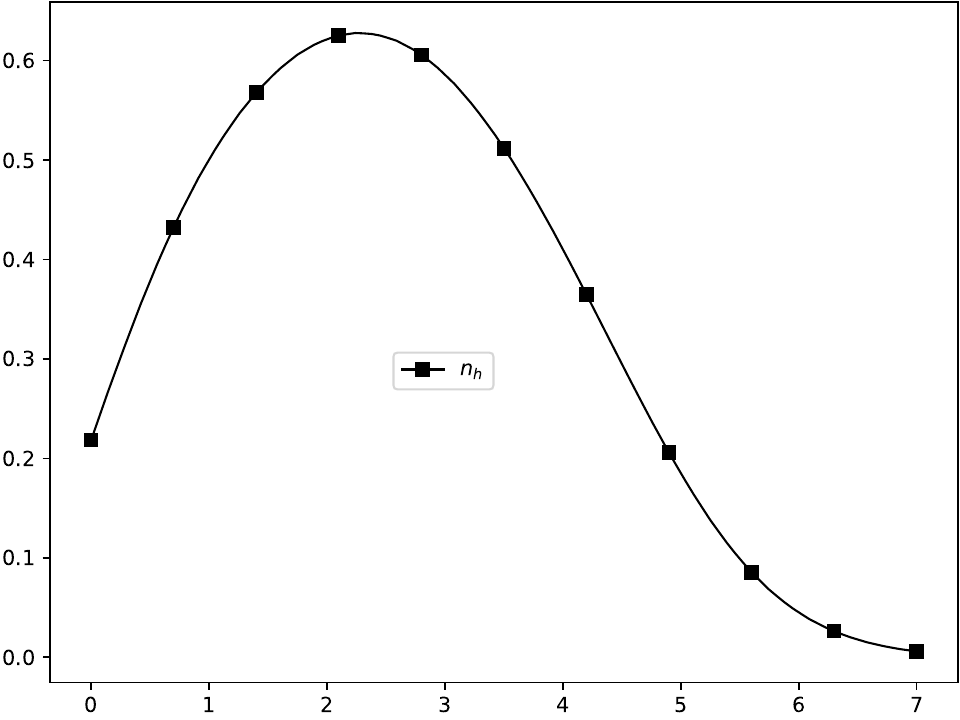}
    \end{subfigure}
    \begin{subfigure}[b]{0.22\textwidth}
        \centering
        \includegraphics[width=1.0\textwidth]{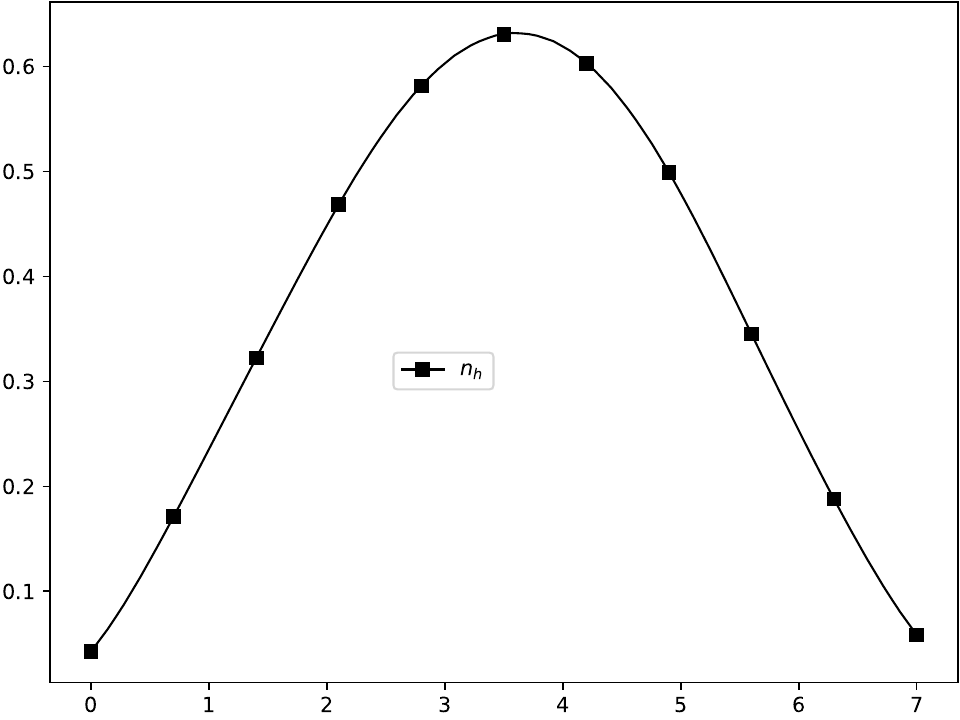}
    \end{subfigure}
    \begin{subfigure}[b]{0.22\textwidth}
        \centering
        \includegraphics[width=1.0\textwidth]{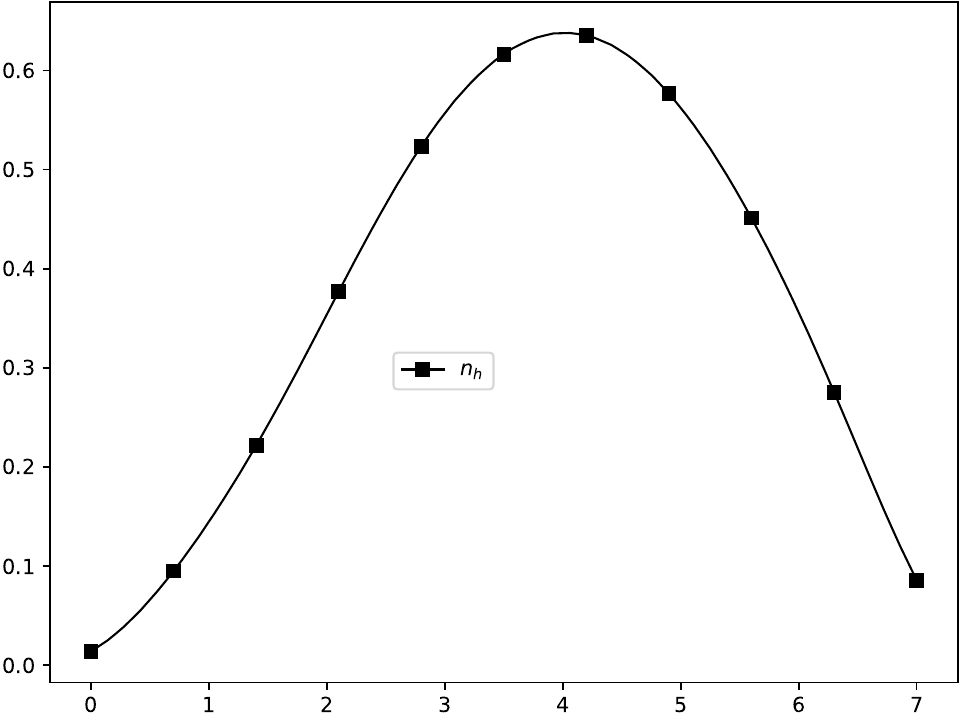}
    \end{subfigure}    
    \caption*{Profiles of $n_h$ at times $t=1$, $2$, $5$ and $10$}
\begin{subfigure}[b]{0.22\textwidth}
        \centering
        \includegraphics[width=1.0\textwidth]{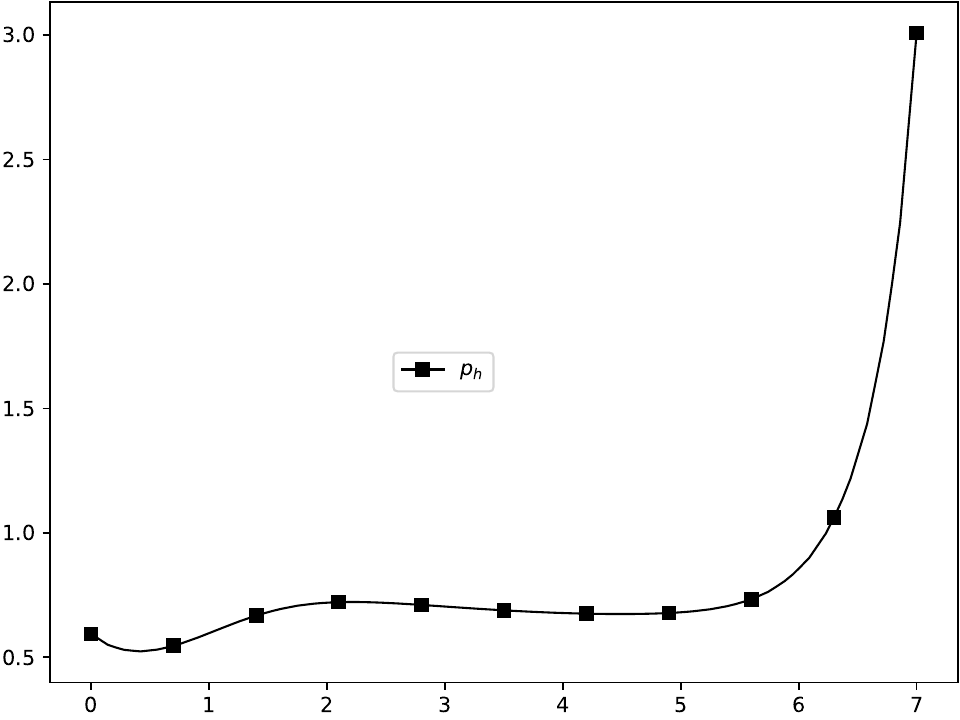}
    \end{subfigure}
    \begin{subfigure}[b]{0.22\textwidth}
        \centering
        \includegraphics[width=1.0\textwidth]{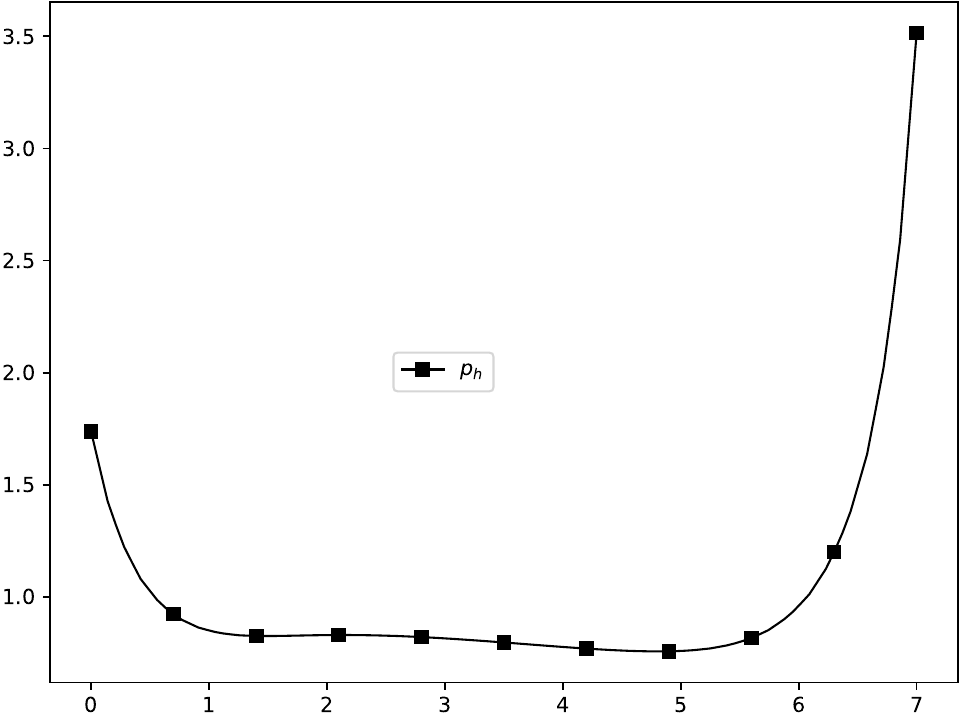}
    \end{subfigure}
    \begin{subfigure}[b]{0.22\textwidth}
        \centering
        \includegraphics[width=1.0\textwidth]{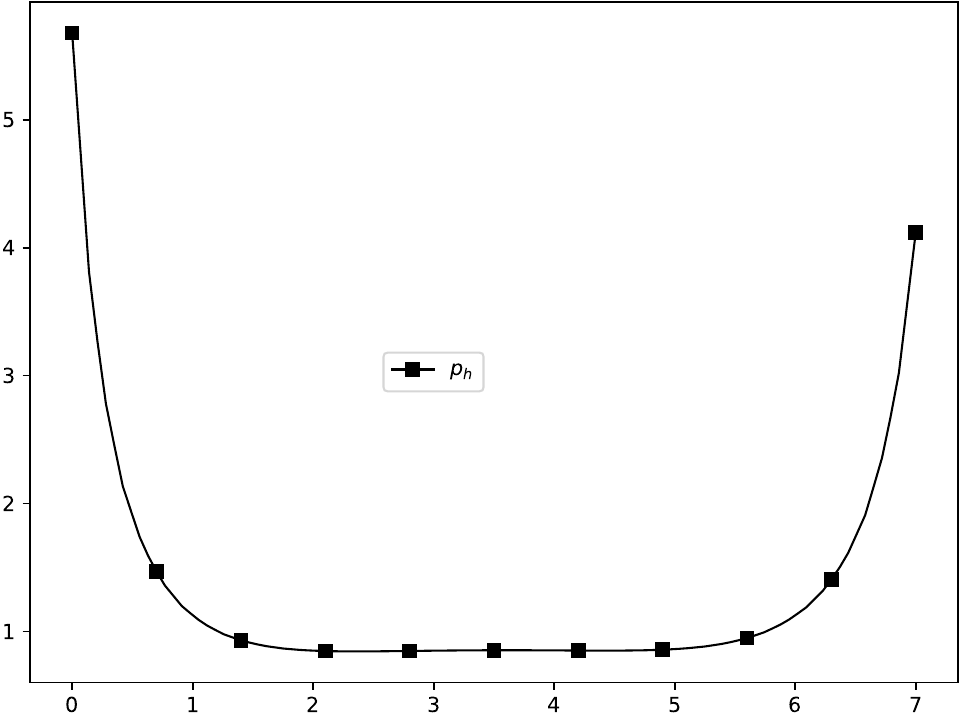}
    \end{subfigure}
    \begin{subfigure}[b]{0.22\textwidth}
        \centering
        \includegraphics[width=1.0\textwidth]{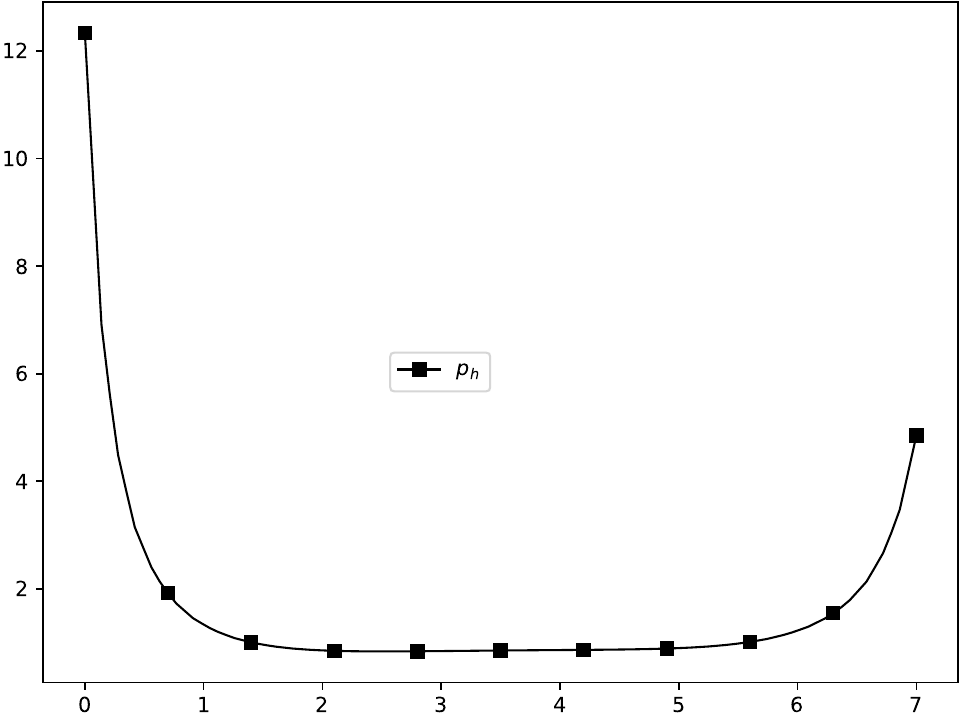}
    \end{subfigure}    
    \caption*{Profiles of $p_h$ at times $t=1$, $2$, $5$ and $10$}
    \begin{subfigure}[b]{0.22\textwidth}
        \centering
        \includegraphics[width=1.0\textwidth]{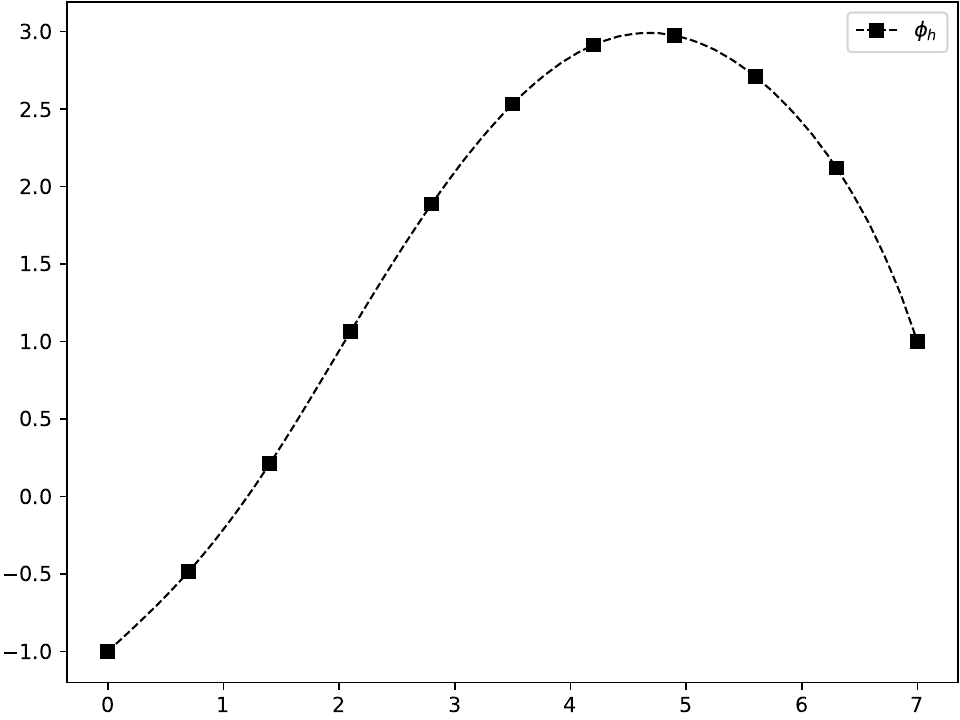}
    \end{subfigure}
    \begin{subfigure}[b]{0.22\textwidth}
        \centering
        \includegraphics[width=1.0\textwidth]{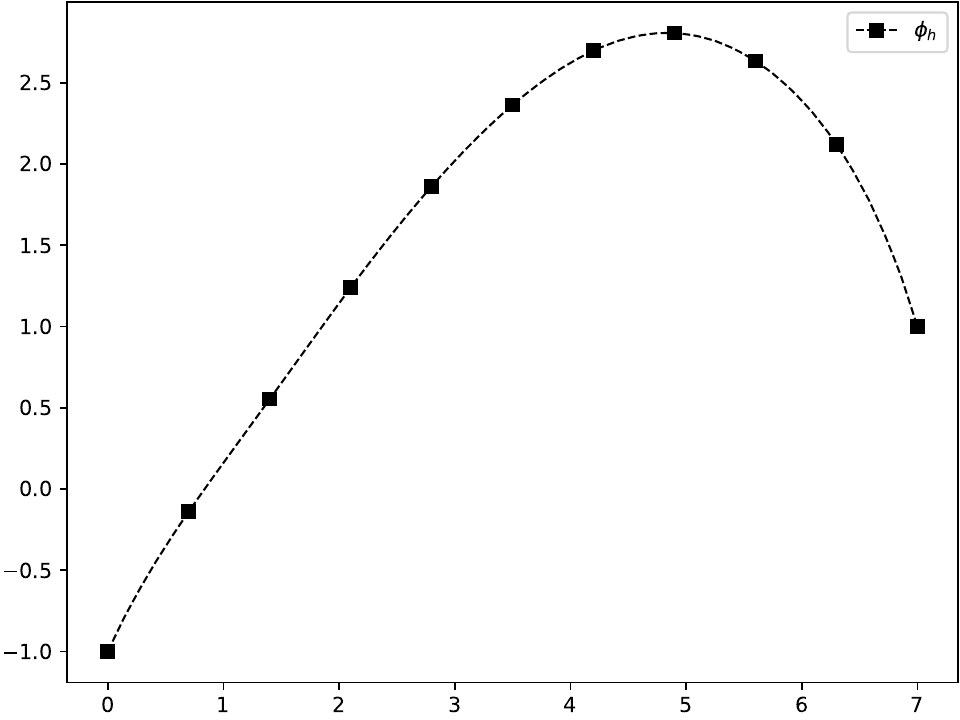}
    \end{subfigure}
    \begin{subfigure}[b]{0.22\textwidth}
        \centering
        \includegraphics[width=1.0\textwidth]{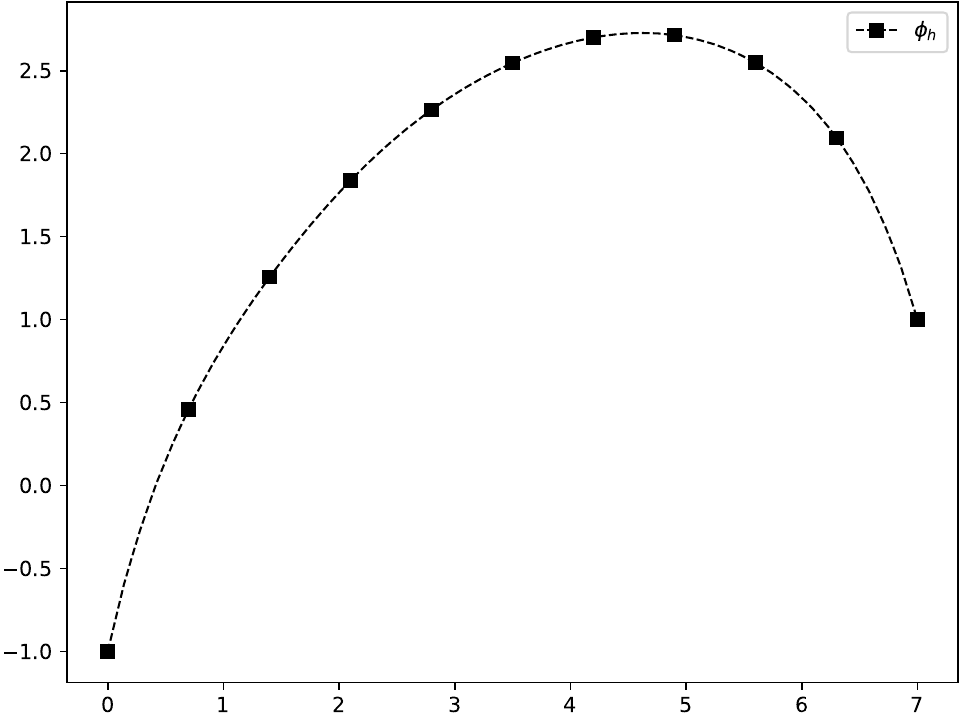}
    \end{subfigure}
    \begin{subfigure}[b]{0.22\textwidth}
        \centering
        \includegraphics[width=1.0\textwidth]{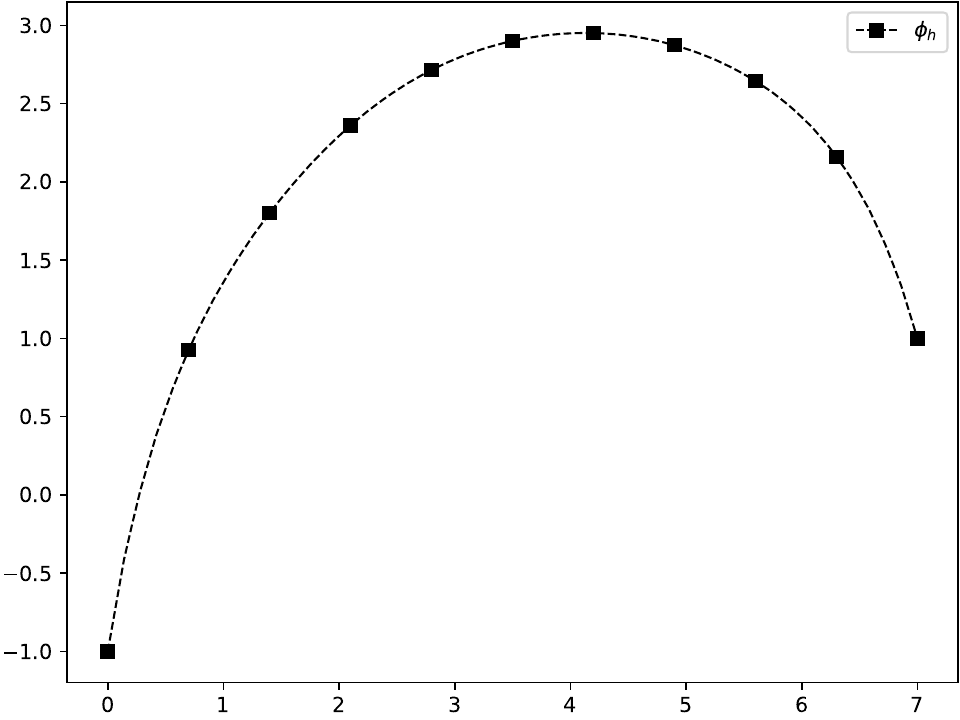}
    \end{subfigure}    
    \caption*{Snapshots of $\phi_h$ at times $t=1$, $2$, $5$ and $10$}
    \caption{Algorithm 2}
    \label{fig_channel_wave_membrane:profiles_alg2}
\end{figure}

In Figure \ref{fig_channel_wave_membrane:mass_and_energies} the total mass is conserved for $n_h$ and increases for $p_h$ due to the incoming of ions through the membrane $\partial\Omega_m$. One can see in Figures \ref{fig_channel_wave_membrane:profiles_alg1} and \ref{fig_channel_wave_membrane:profiles_alg2} how the wave created by $n_h$ stops in the middle of the channel; a consequence, negative ions do not go through the channel thoroughly,  and diffuses slowly toward the wall $\partial\Omega_t$, where $\phi=1$.  For $p_h$, the behavior seems to be that of  a diffusion effect from the channel directed toward both walls $\partial\Omega_t$ and $\partial\Omega_b$, rising the values of positive ions at the walls. It is obvious that maxima for $p_h$ will increase caused by the selective boundary conditions on $\partial\Omega_m$ and minima will increase as well to reach $1$ as shown in Figure \ref{fig_channel_wave_membrane:max_and_min}.  Figures \ref{fig_channel_wave_membrane:snapshots_alg1} and \ref{fig_channel_wave_membrane:snapshots_alg2} for Algorithms 1 and 2, respectively, depicts snapshots of $p_h$, $n_h$ and $\phi_h$ at times $t=1$, $2$, $5$ and $10$.

%

\section{Conclusion} 

In this paper we have constructed two physically consistent discretization of problem \eqref{PNP}-\eqref{IC}. The first algorithm has been devised to satisfy discrete counterparts of both a discrete maximum principle and a discrete minimum, for arbitrary meshes. In addition of the previous properties, the second algorithm enjoys a discrete entropy law for acute meshes, at the expense of being marginally more diffusive, i.e., slightly lower maximum values were observed in the numerical experiments. We have provided both theoretical results of the previous properties  as well as various numerical examples testing the performance of both algorithms. These tests show numerical evidence of the theoretical findings and that both algorithms are capable and well suited for dealing with different complex problem configurations.

\end{document}